\documentclass[12pt]{article}
\usepackage{amssymb}

\newtheorem{definition}{Definition}[section]
\newtheorem{theorem}[definition]{Theorem}
\newtheorem{lemma}[definition]{Lemma}
\newtheorem{corollary}[definition]{Corollary}
\newtheorem{remark}[definition]{Remark}
\newtheorem{notation}[definition]{Notation}
\newtheorem{example}[definition]{Example}
\newtheorem{conjecture}[definition]{Conjecture}
\newtheorem{problem}[definition]{Problem}
\newtheorem{note}[definition]{Note}

\def\R{\mathbb R}
\def\C{\mathbb C}
\def\K{\mathbb K}

\def\Z{\mathbb Z}
\def\fld{\mathbb K}
\def\Mdf{{\mbox{Mat}}_{d+1}(\K)}
\def\alg{\mathcal A}

\begin{document}
%%next line contains proof environment
%\newenvironment{proof}{\noindent{\it Proof\/}:}{\par\noindent $\Box$\par}

\title{ \bf Two linear transformations each tridiagonal
\\
with respect to an eigenbasis of the other; an algebraic
approach to the Askey scheme of orthogonal polynomials\footnote{
Lecture notes for the summer school on orthogonal polynomials
and special functions, Universidad Carlos III de Madrid, Leganes, Spain.
July 8--July 18, 2004.
\hfil\break
http://www.uc3m.es/uc3m/dpto/MATEM/summerschool/indice.html
%{\bf Keywords}. $q$-Racah polynomial,  Leonard pair,
%Tridiagonal pair,
% Askey scheme, 
% Askey-Wilson polynomials.
% \hfil\break
%\noindent {\bf 2000 Mathematics Subject Classification}. 
%05E30, 05E35, 33C45, 33D45. 
}}
\author{Paul Terwilliger  
}
\date{}
%to get date printout, comment out above line 
\maketitle
\begin{abstract} 
Let $\K$ denote a field, and let $V$ denote a  
vector space over $\K$ with finite positive dimension.
We consider a pair
of linear transformations
$A:V\rightarrow V$ and $A^*:V\rightarrow V$
that satisfy the following two conditions:
\begin{enumerate}
\item There exists a basis for $V$ with respect to which
the matrix representing $A$ is irreducible tridiagonal and the matrix
representing $A^*$ is diagonal.
\item There exists a basis for $V$ with respect to which
the matrix representing 
$A^*$ is irreducible tridiagonal and the  matrix representing
$A$ is diagonal.
\end{enumerate}

\medskip
\noindent
We call such a pair a {\it Leonard pair} on $V$. 
We give a correspondence between
Leonard pairs and a class of orthogonal polynomials.
This class
coincides with the 
terminating branch of the Askey scheme
and consists of the
$q$-Racah, $q$-Hahn, dual $q$-Hahn, 
$q$-Krawtchouk,
dual $q$-Krawtchouk,
quantum 
$q$-Krawtchouk,
affine 
$q$-Krawtchouk,
Racah, Hahn, dual Hahn, Krawtchouk,  Bannai/Ito, 
and 
orphan polynomials. 
We describe the above correspondence in detail. We  
show how, for the listed polynomials,
the
3-term recurrence, difference equation, Askey-Wilson
duality, and orthogonality can be expressed in a uniform 
and attractive manner using
the corresponding Leonard pair. 
We give some examples that indicate how
Leonard pairs arise in representation theory
and algebraic combinatorics. We discuss a mild generalization
of a Leonard pair called a tridiagonal pair.
At the end we list some open problems.
Throughout these notes our argument is
elementary and uses only linear algebra. No prior
exposure to the topic is assumed.
\end{abstract}

\newpage
%\tableofcontents
%- - - - - - - - - - - - - - - - - - - - - - - - - - - - - - - - - - -
\bigskip % \bigskip \bigskip
\def\leaderfill{\leaders\hbox to 1em{\hss.\hss}\hfill}
%{\narrower {\narrower {\narrower {\narrower \small
{\narrower {\narrower {\narrower  \small
		   \centerline{CONTENTS}

\medskip\noindent
 1. Leonard pairs                                      \leaderfill 2  \\
 2. An example                                 \leaderfill 3 \\
 3. Leonard systems                          \leaderfill  4 \\
 4. The $D_4$ action     \leaderfill 7 \\
5. The structure of a Leonard system                      \leaderfill 8 \\
 6. The antiautomorphism  $\dagger$                \leaderfill 10 \\
7. The scalars $a_i, x_i$                         \leaderfill 11 \\
8. The polynomials $p_i$                            \leaderfill 13 \\
9. The scalars $\nu, m_i$              \leaderfill 15 \\
10. The standard basis                                       \leaderfill 16 \\
11. The scalars $b_i, c_i$  \leaderfill 18 \\
12. The scalars $k_i$          \leaderfill 20 \\
13. The polynomials $v_i$                      \leaderfill 21 \\
14. The polynomials $u_i$           \leaderfill 22 \\
15. A bilinear form \leaderfill 23\\
16. Askey-Wilson duality \leaderfill 25\\ 
17. The three-term recurrence and the difference equation \leaderfill 26\\ 
18. The orthogonality relations \leaderfill 27\\ 
19. The matrix $P$ \leaderfill 28\\ 
20. The split decomposition \leaderfill 29\\ 
21. The split basis\leaderfill 32\\ 
22. The parameter array and the classifying space\leaderfill 33\\ 
23. Everything in terms of the parameter array\leaderfill 35\\ 
24. The terminating branch of the Askey scheme\leaderfill 38\\ 
25. A characterization of Leonard systems \leaderfill 41\\
26. Leonard pairs $A,A^*$ with $A$ lower bidiagonal and 
 $A^*$ upper bidiagonal
\leaderfill 43\\
27. Leonard pairs $A,A^*$ with $A$ tridiagonal and 
 $A^*$ diagonal
\leaderfill 44\\
28. A characterization of the parameter arrays I 
\leaderfill 44\\
29. A characterization of the parameter arrays II 
\leaderfill 45\\
30. The Askey-Wilson relations
\leaderfill 45\\
31. Leonard pairs and the Lie algebra $sl_2$
\leaderfill 46\\
32. Leonard pairs and the quantum algebra $U_q(sl_2)$
\leaderfill 47\\
33. Leonard pairs in combinatorics
\leaderfill 48\\
34. Tridiagonal pairs
\leaderfill 49\\
35. Appendix: List of parameter arrays
\leaderfill 53\\
36. Suggestions for further research
\leaderfill 63\\
References                       \leaderfill 68 \\

 \baselineskip=\normalbaselineskip \par}  \par} \par}
 %\baselineskip=\normalbaselineskip \par} \par} \par} \par}

\section{Leonard pairs}
These notes are based on the papers
\cite{TD00},
\cite{shape},
\cite{tdanduq},
\cite{LS99},
\cite{qSerre},
\cite{LS24},
\cite{conform},
\cite{lsint},
\cite{Terint},
\cite{TLT:split},
\cite{TLT:array},
\cite{qrac},
\cite{aw}.
We begin by recalling the notion of a {\it Leonard pair}. 
We will use the following terms. Let $X$ denote a square matrix.
Then $X$ is called {\it tridiagonal} whenever each nonzero entry
lies on either the diagonal, the subdiagonal, or the superdiagonal.
Assume $X$ is tridiagonal. Then $X$ is called {\it irreducible}\
whenever each entry on the subdiagonal is nonzero and each entry
on the superdiagonal is nonzero.

\medskip
\noindent We now define a Leonard pair. For the rest of this paper
$\K$ will denote a field.

\begin{definition} \cite{LS99}
\label{def:lprecall}
\rm
Let 
 $V$ denote a  
vector space over $\K$ with finite positive dimension.
By a {\it Leonard pair} on $V$,
we mean an ordered pair of linear transformations
$A:V\rightarrow V$ and $A^*:V\rightarrow V$ 
that
 satisfy both (i), (ii) below. 
\begin{enumerate}
\item There exists a basis for $V$ with respect to which
the matrix representing $A$ is irreducible tridiagonal and the matrix
representing $A^*$ is diagonal.
\item There exists a basis for $V$ with respect to which
the matrix representing 
$A^*$ is irreducible tridiagonal
and the matrix
representing 
$A$ is diagonal.
\end{enumerate}
\end{definition}

%%%%%%%%%%%%%%%
%%%%%%%%%
\begin{note}
\rm
According to a common notational convention
$A^*$ denotes the conjugate-transpose of $A$.
We are not using this convention. In a Leonard
pair $A,A^*$ the linear transformations $A$ and $A^*$
are arbitrary subject
to (i),  (ii) above.
\end{note}

\begin{note}
\rm 
Our use of the name ``Leonard pair'' is motivated
by a connection to a theorem of Doug Leonard
\cite{Leon},
\cite[p.~260]{BanIto}
that involves the $q$-Racah and related polynomials
of the Askey scheme.
\end{note}

%\noindent Leonard pairs play a role in representation theory.
%For instance, Leonard pairs arise naturally
%in the representation theory 
% of the Lie algebra $sl_2$ 
% \cite{TD00}, the quantum algebra
%$U_q(sl_2)$ 
%\cite{Koelink3},
%\cite{Koelink1},
%\cite{Koelink2}, 
%\cite{Koelink4},
%\cite{koo3},
%\cite[ch. 4]{Hjal},
%\cite{Terint}, 
%\cite{qSerre},
%the Askey-Wilson algebra 
%\cite{GYZnature},
% \cite{GYLZmut},
%\cite{GYZTwisted},
%\cite{GYZlinear},
%\cite{GYZspherical},
%\cite{Zhidd}, 
%\cite{ZheCart},
%\cite{Zhidden},
%and 
%the Tridiagonal algebra 
%\cite{TD00},
%\cite{qSerre},
%\cite{LS99}. 
%
%\medskip
%\noindent
%Leonard pairs play a role  in combinatorics. For instance,
%there is a combinatorial object  
%called a  $P$-and $Q$-polynomial
%association scheme \cite{BanIto}, \cite{bcn}, \cite{Leopandq},
%\cite{Tercharpq}, \cite{Ternew}.
%Leonard pairs have been used to describe
%certain irreducible modules for the subconstitutent algebra of
% these schemes \cite{TersubI}, 
%\cite{TersubII}, \cite{TersubIII}.  
%See \cite{Cau}, \cite{CurNom}, \cite{Curspin}, \cite{go}, \cite{HobIto},
%\cite{TD00},
%\cite{Tan}
%for more
%information on Leonard pairs and association  schemes.
%
%\medskip
%\noindent 
%Leonard pairs are closely related to the  work of 
%Grunbaum and Haine on the 
% ``bispectral problem''
%\cite{GH7},
%\cite{GH6}.
%See
%\cite{GH4},
%\cite{GH5},
%\cite{GH1}, 
%\cite{GH3},
%\cite{GH2} 
%for related work. 
%
%
\medskip

\section{An example}

Here is an example of a Leonard pair.
Set 
$V={\K}^4$ (column vectors), set 
\begin{eqnarray*}
A = 
\left(
\begin{array}{ c c c c }
0 & 3  &  0    & 0  \\
1 & 0  &  2   &  0    \\
0  & 2  & 0   & 1 \\
0  & 0  & 3  & 0 \\
\end{array}
\right), \qquad  
A^* = 
\left(
\begin{array}{ c c c c }
3 & 0  &  0    & 0  \\
0 & 1  &  0   &  0    \\
0  & 0  & -1   & 0 \\
0  & 0  & 0  & -3 \\
\end{array}
\right),
\end{eqnarray*}
%may and view $A$ and $A^*$  as linear transformations on $V$.
and view $A$ and $A^*$  as linear transformations from $V$ to $V$.
We assume 
the characteristic of $\K$ is not 2 or 3, to ensure
$A$ is irreducible.
Then $A, A^*$ is a Leonard
pair on $V$. 
Indeed, 
condition (i) in Definition
\ref{def:lprecall}
%\ref{def:leonardpairtalkphilS99}
is satisfied by the basis for $V$
consisting of the columns of the 4 by 4 identity matrix.
To verify condition (ii), we display an invertible  matrix  
$P$ such that 
$P^{-1}AP$ is 
diagonal and 
$P^{-1}A^*P$ is
irreducible tridiagonal.
Set 
\begin{eqnarray*}
P = 
\left(
\begin{array}{ c c c c}
1 & 3  &  3    &  1 \\
1 & 1  &  -1    &  -1\\
1  & -1  & -1  & 1  \\
1  & -3  & 3  & -1 \\
\end{array}
\right).
\end{eqnarray*}
%may By matrix multiplication $P^2=8I$,   
 By matrix multiplication $P^2=8I$, where $I$ denotes the identity,   
so $P^{-1}$ exists. Also by matrix multiplication,    
\begin{equation}
AP = PA^*.
%\qquad \qquad P^{-1}A^*P = A.
\label{eq:apeq}
\end{equation}
Apparently
$P^{-1}AP$ is equal to $A^*$ and is therefore diagonal.
By (\ref{eq:apeq}) and since $P^{-1}$ is
a scalar multiple of $P$, we find
$P^{-1}A^*P$ is equal to $A$ and is therefore irreducible tridiagonal.  Now 
condition (ii) of  Definition 
\ref{def:lprecall}
%\ref{def:leonardpairtalkphilS99}
is satisfied
by the basis for $V$ consisting of the columns of $P$. 

\medskip
\noindent The above example is a member of the following infinite
family of Leonard pairs.
For any nonnegative integer $d$  
the pair
\begin{equation}
A = 
\left(
\begin{array}{ c c c c c c}
0 & d  &      &      &   &{\bf 0} \\
1 & 0  &  d-1   &      &   &  \\
  & 2  &  \cdot    & \cdot  &   & \\
  &   & \cdot     & \cdot  & \cdot   & \\
  &   &           &  \cdot & \cdot & 1 \\
{\bf 0} &   &   &   & d & 0  
\end{array}
\right),
\qquad A^*= \hbox{diag}(d, d-2, d-4, \ldots, -d)
\label{eq:fam1}
\end{equation}
is a Leonard pair on the vector space $\K^{d+1} $,
provided the 
 characteristic of $\K$ is zero or an odd prime greater than $d$.
This can be  proved by modifying the 
 proof for $d=3$ given above. One shows  
$P^2=2^dI$  and $AP= PA^*$, where 
$P$ denotes the matrix with $ij$ entry
\begin{equation}
P_{ij} =  
\Biggl({{ d }\atop {j}}\Biggr) {{}_2}F_1\Biggl({{-i, -j}\atop {-d}}
\;\Bigg\vert \;2\Biggr)
\qquad \qquad (0 \leq i,j\leq d).
\label{eq:ex1}
\end{equation}
%\cite[Section 16]{LS24}.
%for $0 \;\leq i,j\leq d$.
We follow the standard notation for
hypergeometric series \cite[p.~3]{GR}. 
The details of the above calculations are given in Section 24
below.

\medskip

\section{Leonard systems}
\medskip
\noindent When working with a Leonard pair, 
it is often convenient to consider a closely related
and somewhat more abstract object called
a {\it Leonard system}.
In order to define this we first make an observation about
Leonard pairs.

%%%%%%%%%%%%%%%
\begin{lemma}
%\cite[Lemma 1.3]{LS99}
%In this section we review some results of \cite{TD00}, \cite{LS99}.
\label{lem:preeverythingtalkS99}
%With reference to Definition
%\ref{def:lprecall},
Let $V$ denote a vector space over $\K$ with finite positive
dimension and 
let $A, A^*$ denote a Leonard pair on $V$. Then
the eigenvalues of $A$ are mutually distinct
and contained in $\K$.
Moreover, the eigenvalues of
$A^*$ are mutually distinct and contained in $\fld$.
\end{lemma}
\noindent {\it Proof:}
 Concerning $A$, recall by
 Definition
\ref{def:lprecall}(ii)
that there exists
 a basis for $V$ consisting of eigenvectors for $A$.
 Consequently the eigenvalues of $A$ are all in $\fld$, and the minimal
 polynomial of $A$ has no repeated roots. To show the eigenvalues
 of $A$ are distinct, we show the minimal polynomial of $A$ has
 degree equal to $\hbox{dim}\,V$.
 By Definition
\ref{def:lprecall}(i),
 there exists
 a basis for $V$ with respect to which the matrix representing
 $A$ is
  irreducible tridiagonal. Denote this matrix by $B$. On one hand,
    $A$ and $B$ have the same minimal polynomial.
    On the other hand, 
    using the tridiagonal shape of $B$, we find
    $I, B, B^2, \ldots, B^d$ are linearly independent, where
    $d=\hbox{dim}\,V-1$,  so the minimal  polynomial of $B$
     has  degree $d+1=\hbox{dim}\,V$.
     We conclude the 
     mininimal polynomial of $A$  has degree equal to
       $\hbox{dim}\,V$,
       so the eigenvalues of  $A$ are distinct.
       We have now obtained our assertions about $A$,
       and the case of $A^*$ is similar.
\hfill $\Box $ \\

\noindent To prepare for our definition of  a Leonard system,
we recall a few concepts from  linear algebra. 
Let $d$  denote  a  nonnegative
integer 
and let 
$\Mdf$ denote the $\fld$-algebra consisting of all
$d+1$ by $d+1$ matrices that have entries in $\K$. We
index the rows and columns by $0,1,\ldots, d$.
We let ${\K}^{d+1}$ denote the $\K$-vector space consisting
of all $d+1$ by $1$ matrices that have entries in $\K$.
We index the rows by $0,1,\ldots, d$. We view 
 ${\K}^{d+1}$ as a left module for 
$\Mdf$.
We observe this module is irreducible.
For the rest of this paper we  
let  
$\mathcal A$
denote a $\fld$-algebra
isomorphic to 
$\hbox{Mat}_{d+1}(\K)$. 
When we refer to an $\mathcal A$-module we mean a left 
 $\mathcal A$-module.
%For the rest of this paper we let $V$ denote an irreducible
%$\mathcal A$-module.
Let $V$ denote an irreducible
$\mathcal A$-module.
We remark that $V$ is unique up to isomorphism of 
$\mathcal A$-modules, and 
that $V$ has dimension $d+1$.
Let $v_0, v_1, \ldots, v_d$ denote a basis for $V$.
For $X \in \mathcal A$ and
$Y \in 
\mbox{Mat}_{d+1}(\K)$,
we say 
$Y$ {\it represents $X$ with respect to $v_0, v_1, \ldots, v_d$}
whenever
$Xv_j = \sum_{i=0}^d Y_{ij}v_i$ for 
$0 \leq j\leq d$.
Let $A$ denote an element of $\mathcal A$.
We say $A$ is 
{\it multiplicity-free} whenever it has $d+1$ 
mutually distinct  eigenvalues in $\K$.
Let $A$ denote a  multiplicity-free element of $\alg$.
Let $\theta_0, \theta_1, \ldots, \theta_d$ denote an ordering of 
the eigenvalues
of $A$, and for $0 \leq i \leq d$   put 
\begin{equation}
E_i = \prod_{{0 \leq  j \leq d}\atop
{j\not=i}} {{A-\theta_j I}\over {\theta_i-\theta_j}},
\label{eq:primiddef}
\end{equation}
where $I$ denotes the identity of $\cal A$.
We observe
(i) 
$AE_i = \theta_iE_i \; (0 \leq i \leq d)$;
(ii) $E_iE_j = \delta_{ij}E_i \;(0 \leq i,j\leq d)$;
(iii) $\sum_{i=0}^d E_i = I $;
(iv) $A = \sum_{i=0}^d \theta_i E_i.$
Let $\mathcal D$ denote the subalgebra of $\mathcal A$ generated
by $A$. Using (i)--(iv) we find the sequence
 $E_0, E_1, \ldots, E_d$ is a  basis for the
$\K$-vector space $\cal D$.
We call $E_i$  the {\it primitive idempotent} of
$A$ associated with $\theta_i$.
It is helpful to think of these primitive idempotents as follows. 
Observe 
\begin{eqnarray}
V = E_0V + E_1V + \cdots + E_dV \qquad \qquad (\mbox{direct sum}).
\label{eq:VdecompS99}
\end{eqnarray}
For $0\leq i \leq d$, $E_iV$ is the (one dimensional) eigenspace of
$A$ in $V$ associated with the 
eigenvalue $\theta_i$, 
and $E_i$ acts  on $V$ as the projection onto this eigenspace. 
We remark that $\lbrace A^i |0 \leq i \leq d\rbrace $ 
is a basis for the 
$\K$-vector space $\mathcal D$ and that
$\prod_{i=0}^d (A-\theta_iI)=0$. By a {\it Leonard pair in ${\mathcal A}$}
we mean an ordered pair of elements taken from $\mathcal A$ that act
on $V$ as a Leonard pair in the sense of Definition
\ref{def:lprecall}. We call $\mathcal A$ the {\it ambient algebra} of
the pair and say the pair is {\it over $\K$}. 
We refer to $d$ as the {\it diameter} of the pair.
We now define a Leonard system.

\begin{definition} \cite[Definition 1.4]{LS99}
% \cite{TD00}, \cite{LS99}.
\label{def:defls}
\label{eq:ourstartingpt}
\rm
By a {\it Leonard system} in $\mathcal A$ we mean a 
sequence
$\Phi:=(A;A^*; \lbrace E_i\rbrace_{i=0}^d; $ $ 
\lbrace E^*_i\rbrace_{i=0}^d)$
that satisfies (i)--(v) below. 
\begin{enumerate}
\item Each of $A,A^*$ is a multiplicity-free element in $\mathcal A$.
\item $E_0,E_1,\ldots,E_d$ is an ordering of the primitive 
idempotents of $A$.
\item $E^*_0,E^*_1,\ldots,E^*_d$ is an ordering of the primitive 
idempotents of $A^*$.
\item ${\displaystyle{
E_iA^*E_j = \cases{0, &if $\;\vert i-j\vert > 1$;\cr
\not=0, &if $\;\vert i-j \vert = 1$\cr}
\qquad \qquad 
(0 \leq i,j\leq d)}}$.
\item ${\displaystyle{
 E^*_iAE^*_j = \cases{0, &if $\;\vert i-j\vert > 1$;\cr
\not=0, &if $\;\vert i-j \vert = 1$\cr}
\qquad \qquad 
(0 \leq i,j\leq d).}}$
\end{enumerate}
We refer to $d$ as the {\it diameter} of $\Phi$ and say 
$\Phi$ is {\it over } $\K$.  We call $\mathcal A$
the {\it ambient algebra} of $\Phi$. 
\end{definition}

\noindent We comment on how 
Leonard pairs and Leonard systems are related.
In the following discussion $V$ denotes
an irreducible $\mathcal A$-module. Let 
$(A;A^*; \lbrace E_i\rbrace_{i=0}^d;  
\lbrace E^*_i\rbrace_{i=0}^d) $
denote a Leonard system in $\mathcal A$. 
For $0 \leq i \leq d$ let $v_i$ denote a nonzero vector in 
$E_iV$. Then the sequence 
$v_0, v_1, \ldots, v_d$ is a basis for $V$ that satisfies
Definition 
\ref{def:lprecall}(ii).
For $0 \leq i \leq d$ let $v^*_i$ denote a nonzero vector in 
$E^*_iV$. Then the sequence 
$v^*_0, v^*_1, \ldots, v^*_d$ is a basis for $V$ that satisfies
Definition 
\ref{def:lprecall}(i). By these comments the pair $A,A^*$ is
a Leonard pair in $\mathcal A$. Conversely let $A,A^*$ denote
a Leonard pair in $\mathcal A$. Then each of $A, A^*$ is multiplicity-free
by Lemma
\ref{lem:preeverythingtalkS99}.
 Let 
$v_0, v_1, \ldots, v_d$ denote a basis for $V$ that satisfies
Definition 
\ref{def:lprecall}(ii). For $0 \leq i \leq d$ the vector
$v_i$ is an eigenvector for $A$;
 let
$E_i$ denote the corresponding primitive idempotent.
 Let 
$v^*_0, v^*_1, \ldots, v^*_d$ denote a basis for $V$ that satisfies
Definition 
\ref{def:lprecall}(i). For $0 \leq i \leq d$ the vector
$v^*_i$ is an eigenvector for $A^*$;
 let
$E^*_i$ denote the corresponding primitive idempotent.
Then 
$(A;A^*; \lbrace E_i\rbrace_{i=0}^d;  
\lbrace E^*_i\rbrace_{i=0}^d)$
is a Leonard system in $\mathcal A$.
In summary we have the following.

\begin{lemma} Let $A$ and $A^*$ denote elements of $\mathcal A$.
Then the pair $A,A^*$ is a Leonard pair in $\mathcal A$ if and only 
if the following (i), (ii) hold.
\begin{enumerate}
\item Each of $A,A^*$ is multiplicity-free.
\item There exists an ordering $E_0, E_1, \ldots, E_d$ of the primitive
idempotents of $A$ and there exists
an ordering $E^*_0, E^*_1, \ldots, E^*_d$ of the primitive
idempotents of $A^*$ such that
$(A;A^*; $ $\lbrace E_i\rbrace_{i=0}^d;  
\lbrace E^*_i\rbrace_{i=0}^d)$ is a Leonard system in $\mathcal A$.
\end{enumerate}
\end{lemma}

\medskip
\noindent We recall the notion of {\it isomorphism} for Leonard pairs and
Leonard systems.
%Let $A,A^*$ denote a Leonard pair in $\mathcal A$ and let
%$\sigma : {\mathcal A} \rightarrow {\mathcal A}'$ denote an isomorphism
%of $\K$-algebras. We observe the pair $A^\sigma, A^{*\sigma}$ is
%a Leonard pair in ${\mathcal A}'$.

\begin{definition} 
\rm
Let $A,A^*$ and $B,B^*$ denote Leonard pairs over
$\K$. By an {\it isomorphism of Leonard pairs from $A,A^*$ to $B,B^*$} 
we mean an isomorphism of $\K$-algebras from the ambient algebra of
$A,A^*$ to the ambient algebra of $B,B^*$ that sends $A$ to $B$
and $A^*$ to $B^*$.  The Leonard pairs $A,A^*$ and $B,B^*$ are
said to be {\it isomorphic} whenever there exists an isomorphism
of Leonard pairs from 
 $A,A^*$ to  $B,B^*$.
\end{definition} 
\noindent 
Let $\Phi$ 
denote the Leonard system from Definition
\ref{eq:ourstartingpt}
and let 
$\sigma :\alg \rightarrow {\cal A}'$ denote an isomorphism of
$\K$-algebras. We write 
$\Phi^{\sigma}:= 
(A^\sigma;A^{*\sigma};\lbrace E^\sigma_i\rbrace_{i=0}^d;  
\lbrace E^{*\sigma}_i\rbrace_{i=0}^d)$
and observe 
$\Phi^{\sigma}$ 
is a Leonard  system in ${\cal A }'$.

\begin{definition}
% \cite{LS99}
\label{def:isolsS99o}
\rm
Let $\Phi$ 
and  
 $\Phi'$ 
denote Leonard systems over $\K$.
 By an {\it isomorphism of Leonard  systems
 from $\Phi $ to $\Phi'$} we mean an isomorphism 
of $\K $-algebras
$\sigma $
from the ambient algebra of $\Phi$ to the ambient algebra of
$\Phi'$ such 
that  $\Phi^\sigma = \Phi'$. 
The Leonard systems $\Phi $, $\Phi'$
are said to be {\it isomorphic} whenever there exists
an isomorphism of Leonard  systems from $\Phi $ to $\Phi'$. 
\end{definition}

\noindent We have a remark.
Let  $\sigma :\mathcal A \rightarrow 
\mathcal A$ denote  any map.
By the Skolem-Noether theorem 
\cite[Corollary 9.122]{rotman},
$\sigma $ is an isomorphism of $\K$-algebras
if and only if there exists an invertible $S \in {\mathcal A}$ such that
$X^\sigma = S X S^{-1}$ for all  $X \in  {\mathcal A}$.

%\medskip 
%\noindent We finish this section by  recalling some parameters
%which will help us describe 
%a given Leonard system.
%
%\begin{definition}
% \cite{LS99}
%\label{def:evseq}
%Let $\Phi$ denote the Leonard system  from  
%Definition \ref{eq:ourstartingpt}.
%For $0 \leq i \leq d$, 
%we let $\theta_i $ (resp. $\theta^*_i$) denote the eigenvalue
%of $A$ (resp. $A^*$) associated with $E_i$ (resp. $E^*_i$).
%We refer to  $\theta_0, \theta_1, \ldots, \theta_d$ as the 
%eigenvalue sequence of $\Phi$.
%We refer to  $\theta^*_0, \theta^*_1, \ldots, \theta^*_d$ as the 
%dual eigenvalue sequence of $\Phi$. We observe 
% $\theta_0, \theta_1, \ldots, \theta_d$ are mutually distinct
% and contained in $\K$. Similarly
%  $\theta^*_0, \theta^*_1, \ldots, \theta^*_d$  
% are mutually distinct
% and contained in $\K$. 
%\end{definition}

\section{The $D_4$ action}

\medskip
\noindent A given Leonard system  can be modified in  several
ways to get a new Leonard system. For instance, 
let $\Phi$ 
 denote the Leonard system from 
Definition 
\ref{eq:ourstartingpt}, and let $\alpha, \alpha^*, \beta, \beta^*$
denote scalars in $\K$ such that $\alpha \not=0$, $\alpha^*\not=0$.
Then
the sequence
\begin{eqnarray*}
(\alpha A+\beta I;\alpha^*A^*+\beta^* I; \lbrace E_i\rbrace_{i=0}^d;  
\lbrace E^*_i\rbrace_{i=0}^d)
\end{eqnarray*}
 is a Leonard system
in $\mathcal A$. 
Also, each of the following three sequences is a Leonard system
in $\mathcal A$.
\begin{eqnarray*}
 \;\Phi^*&:=& (A^*; A; \lbrace E^*_i\rbrace_{i=0}^d; 
 \lbrace E_i\rbrace_{i=0}^d),
%\label{eq:lsdualS99}
\\
\Phi^{\downarrow}&:=& (A; A^*; \lbrace E_i\rbrace_{i=0}^d ;
\lbrace E^*_{d-i} \rbrace_{i=0}^d),
%\label{eq:lsinvertS99}
\\
\Phi^{\Downarrow} 
&:=& (A;A^*; \lbrace E_{d-i}\rbrace_{i=0}^d; \lbrace E^*_i\rbrace_{i=0}^d).
%\label{eq:lsdualinvertS99}
\end{eqnarray*}
Viewing $*, \downarrow, \Downarrow$
as permutations on the set of all Leonard systems,
\begin{eqnarray}
&&\qquad \qquad \qquad  *^2 \;=\;  
\downarrow^2\;= \;
\Downarrow^2 \;=\;1,
\qquad \quad 
\label{eq:deightrelationsAS99}
\\
&&\Downarrow *\; 
=\;
* \downarrow,\qquad \qquad   
\downarrow *\; 
=\;
* \Downarrow,\qquad \qquad   
\downarrow \Downarrow \; = \;
\Downarrow \downarrow.
\qquad \quad 
\label{eq:deightrelationsBS99}
\end{eqnarray}
The group generated by symbols 
$*, \downarrow, \Downarrow $ subject to the relations
(\ref{eq:deightrelationsAS99}),
(\ref{eq:deightrelationsBS99})
is the dihedral group $D_4$.  
We recall $D_4$ is the group of symmetries of a square,
and has 8 elements.
Apparently $*, \downarrow, \Downarrow $ induce an action of 
 $D_4$ on the set of all Leonard systems.
Two Leonard systems will be called {\it relatives} whenever they
are in the same orbit of this $D_4$ action.
The relatives of $\Phi$ are as follows:
\medskip

\centerline{
\begin{tabular}[t]{c|c}
        name &relative \\ \hline 
        $\Phi$ & $(A;A^*;\lbrace E_i\rbrace_{i=0}^d;\lbrace E^*_i\rbrace_{i=0}^d)$   \\ 
        $\Phi^\downarrow$ &
         $(A;A^*;\lbrace E_i\rbrace_{i=0}^d;\lbrace E^*_{d-i}\rbrace_{i=0}^d)$   \\ 
        $\Phi^\Downarrow$ &
         $(A;A^*;\lbrace E_{d-i}\rbrace_{i=0}^d;\lbrace E^*_i\rbrace_{i=0}^d)$   \\ 
        $\Phi^{\downarrow \Downarrow}$ &
         $(A;A^*;\lbrace E_{d-i}\rbrace_{i=0}^d;\lbrace E^*_{d-i}\rbrace_{i=0}^d)$   \\ 
	$\Phi^*$ &
       $(A^*;A;\lbrace E^*_i\rbrace_{i=0}^d;\lbrace E_i\rbrace_{i=0}^d)$   \\ 
        $\Phi^{\downarrow *}$ &
	 $(A^*;A;\lbrace E^*_{d-i}\rbrace_{i=0}^d;  
       \lbrace E_i\rbrace_{i=0}^d)$ \\
        $\Phi^{\Downarrow *}$ &
	 $(A^*;A;\lbrace E^*_i\rbrace_{i=0}^d;
         \lbrace E_{d-i}\rbrace_{i=0}^d)$ \\ 
	$\Phi^{\downarrow \Downarrow *}$ &
	 $(A^*;A;\lbrace E^*_{d-i} \rbrace_{i=0}^d;
         \lbrace E_{d-i}\rbrace_{i=0}^d)$
	\end{tabular}}
\medskip
\noindent 
There may be some isomorphisms among the above Leonard
systems.

\medskip
\noindent For the rest of this paper we will use the following
notational convention.

\begin{definition}
\label{def:notconv}
\rm
Let $\Phi$ denote a Leonard system. For any element $g$ in
the group $D_4$ and for any object $f$ that we associate with
$\Phi$, we let $f^g$ denote the corresponding object for the
Leonard system 
$\Phi^{g^{-1}}$. We have been using this
convention all along; an example is $E^*_i(\Phi)= 
E_i(\Phi^*)$.
\end{definition}

\section{The structure of a Leonard system}

\noindent In this section we establish a few basic 
facts concerning Leonard systems.
We begin with a definition and two routine lemmas.

\begin{definition}
%\cite{LS99}
\label{def:evseq}
\rm
Let $\Phi$ denote the Leonard system  from  
Definition \ref{eq:ourstartingpt}.
For $0 \leq i \leq d$, 
we let $\theta_i $ (resp. $\theta^*_i$) denote the eigenvalue
of $A$ (resp. $A^*$) associated with $E_i$ (resp. $E^*_i$).
We refer to  $\theta_0, \theta_1, \ldots, \theta_d$ as the 
{\it eigenvalue sequence} of $\Phi$.
We refer to  $\theta^*_0, \theta^*_1, \ldots, \theta^*_d$ as the 
{\it dual eigenvalue sequence} of $\Phi$. We observe 
 $\theta_0, \theta_1, \ldots, \theta_d$ are mutually distinct
 and contained in $\K$. Similarly
  $\theta^*_0, \theta^*_1, \ldots, \theta^*_d$  
 are mutually distinct
 and contained in $\K$. 
\end{definition}

\begin{lemma}
\label{lem:step0}
Let $\Phi$ denote the Leonard system
from Definition  
\ref{eq:ourstartingpt} and let $V$ denote an irreducible
$\mathcal A$-module. For $0 \leq i \leq d$ let $v_i$
denote a nonzero vector in $E^*_iV$ and observe
$v_0, v_1, \ldots, v_d$ is a basis for $V$.
Then  (i), (ii) hold below.
\begin{enumerate}
\item
For $0 \leq i \leq d$ the
matrix in 
$\mbox{Mat}_{d+1}(\K)$
that represents $E^*_i$ with respect to
$v_0, v_1, \ldots, v_d$ has $ii$ entry 1 and all other entries 0.
\item
The matrix in 
$\mbox{Mat}_{d+1}(\K)$
that represents $A^*$ with respect to
$v_0, v_1, \ldots, v_d$ is equal to
$\mbox{diag}(\theta^*_0, \theta^*_1, \ldots, \theta^*_d)$.
\end{enumerate}
\end{lemma}

\begin{lemma}
\label{lem:step1}
Let $A$ denote an irreducible tridiagonal matrix
in 
$\mbox{Mat}_{d+1}(\K)$.
Pick any integers $i,j$ $(0 \leq i,j\leq d)$.
Then (i)--(iii) hold below.
\begin{enumerate}
\item
The entry $(A^r)_{ij}= 0$ if  $r<\vert i-j \vert,
\qquad \qquad (0 \leq r \leq d).$
\item
Suppose $i\leq j$. Then 
the entry
$(A^{j-i})_{ij} = \prod_{h=i}^{j-1}A_{h,h+1}.$
%\begin{eqnarray*}
%(A^{j-i})_{ij} = \prod_{h=i}^{j-1}A_{h,h+1}.
%\end{eqnarray*}
Moreover 
$(A^{j-i})_{ij} \not=0$.
\item
Suppose $i\geq j$. Then 
the entry
$(A^{i-j})_{ij} = \prod_{h=j}^{i-1}A_{h+1,h}.$
%\begin{eqnarray*}
%(A^{i-j})_{ij} = \prod_{h=j}^{i-1}A_{h+1,h}.
%\end{eqnarray*}
Moreover 
$(A^{i-j})_{ij} \not=0$.
\end{enumerate}
\end{lemma}
%\noindent {\it Proof:}
%Routine using matrix multiplication.
%\hfill $\Box $ \\

\begin{theorem}
\label{eq:lsmatbasis}
Let $\Phi$ denote the Leonard system
from Definition  
\ref{eq:ourstartingpt}. Then the elements
\begin{equation}
A^rE^*_0A^s \qquad \qquad (0 \leq r,s\leq d)
\label{eq:lpbasis}
\end{equation}
form a basis for the $\K$-vector space $\cal A$.
\end{theorem}

\noindent {\it Proof:}
The number of elements in 
(\ref{eq:lpbasis}) is equal to $(d+1)^2$, and this number is the dimension
of $\cal A$. Therefore
it suffices
to show 
the elements
in (\ref{eq:lpbasis})
are linearly independent. To do this,
we represent 
the elements 
in (\ref{eq:lpbasis}) by matrices. 
Let $V$ denote an irreducible 
$\mathcal A$-module.
For $0 \leq i \leq d$ let $v_i$ denote a nonzero vector in 
$E^*_iV$, and observe $v_0, v_1,\ldots, v_d$ is a basis for
$V$. For the purpose of this proof, let us identify each element
of $\mathcal A$ with the matrix in 
$\hbox{Mat}_{d+1}(\K)$ that 
represents it with respect to the
basis $v_0, v_1, \ldots, v_d$. Adopting this point of view
we find 
$A$ is irreducible tridiagonal and $A^*$ is diagonal. 
For $0 \leq r,s\leq d$ we show  the
entries  of $A^rE^*_0A^s$ satisfy 
\begin{equation}
(A^rE^*_0A^s)_{ij} = 
 \cases{0, 
&$\qquad $if $\quad i>r \quad$ or $\quad j>s $ ;\cr
\not=0,  & $\qquad $if $\quad i=r\quad $ and $\quad j=s$\cr} 
\qquad \qquad (0 \leq i,j\leq d).
\label{eq:lowert}
\end{equation}
By Lemma
\ref{lem:step0}(i)
the matrix $E^*_0$ has $00$ entry 1 and all other entries $0$.
%Observe that for
%$0 \leq i,j\leq d$
%the $ij$ entry of $E^*_0$ is 
%one if both $i=0,j=0$, and zero otherwise. 
%From this we find
Therefore 
\begin{equation}
(A^rE^*_0A^s)_{ij} =  (A^r)_{i0} (A^s)_{0j} \qquad \qquad (0 \leq i,j\leq d).
\label{eq:redent}
\end{equation}
We mentioned $A$ is irreducible tridiagonal.
Applying Lemma \ref{lem:step1} we find that for 
$0 \leq i \leq d$ the entry $(A^r)_{i0}$ is zero if $i>r$, and  nonzero if
$i=r$. Similarly for $0 \leq j\leq d$ the entry
 $(A^s)_{0j}$ is zero 
if $j>s$, and nonzero if $j=s$. Combining these facts with 
(\ref{eq:redent}) we routinely obtain
(\ref{eq:lowert}) and it follows  the elements 
(\ref{eq:lpbasis}) are linearly independent.
Apparently the elements 
(\ref{eq:lpbasis}) form a basis  for $\cal A $, as desired.
\hfill $\Box $ \\

\begin{corollary} 
\label{cor:genset}
Let $\Phi$ denote the Leonard system
from Definition 
\ref{eq:ourstartingpt}. Then 
the elements 
$A, E^*_0$ together generate $\cal A$. Moreover
the elements $A,A^*$ together generate $\cal A$.
\end{corollary}

\noindent {\it Proof:}
The first assertion is immediate from Theorem 
\ref{eq:lsmatbasis}.  
The second assertion follows from the first assertion and
the observation that $E^*_0$ is a polynomial in $A^*$. 
\hfill $\Box $ \\

\noindent The following is immediate from Corollary
\ref{cor:genset}.

\begin{corollary} 
\label{cor:gensetlp}
Let $A,A^*$ denote a Leonard pair in 
$\mathcal A$. Then the 
elements $A,A^*$ together generate $\mathcal A$.
\end{corollary}

\noindent We mention a few implications of Theorem 
\ref{eq:lsmatbasis} that will be useful later in the paper.

\begin{lemma}
\label{lem:altbase}
Let $\Phi$ denote the Leonard system
from Definition 
\ref{eq:ourstartingpt}. Let $\cal D$ denote
the subalgebra of $\cal A$ generated by $A$.
 Let $X_0, X_1, \ldots, X_d$
denote a basis for the $\K$-vector space 
$\mathcal D$. Then the elements
\begin{equation}
X_r E^*_0 X_s \qquad \qquad (0 \leq r,s\leq d)
\label{altbase}
\end{equation}
form a basis for the $\K$-vector space 
$\mathcal A$.

\end{lemma}

\noindent {\it Proof:}
 The number of elements in    
(\ref{altbase}) is equal to  $(d+1)^2$, 
and this number is the dimension
of $\cal A$. Therefore it suffices to show the elements
(\ref{altbase}) span $\cal A$. 
But this is immediate from Theorem 
\ref{eq:lsmatbasis}, and since 
each  element 
in (\ref{eq:lpbasis}) is contained in the span of
the elements (\ref{altbase}).
\hfill $\Box $ \\

\begin{corollary}
\label{cor:eibasis}
Let $\Phi$ denote the Leonard system
from Definition
\ref{eq:ourstartingpt}.
Then the elements
\begin{equation}
E_r E^*_0 E_s \qquad \qquad (0 \leq r,s\leq d)
\label{eq:eibase}
\end{equation}
form a basis for the $\K$-vector space $\mathcal A$.
\end{corollary}

\noindent {\it Proof:}
Immediate from  
Lemma \ref{lem:altbase}, with  $X_i = E_i$ for $0 \leq i \leq d$.
\hfill $\Box $ \\

\begin{lemma}
\label{lem:oneiszero}
Let $\Phi$ denote the Leonard system
from Definition 
\ref{eq:ourstartingpt}. Let $\cal D$ denote
the subalgebra of $\cal A$ generated by $A$.
Let $X$ and $Y$ denote elements in $\mathcal D$
and assume $XE^*_0Y=0$. Then $X=0$ or $Y=0$.
\end{lemma}
\noindent {\it Proof:}
Let $X_0, X_1, \ldots, X_d$ denote a basis for 
the $\K$-vector space $\mathcal D$.
Since $X \in {\mathcal D}$ there exists
$\alpha_i \in \K$ $(0 \leq i \leq d)$
such that $X=\sum_{i=0}^d \alpha_i X_i$.
Similarly there exists 
$\beta_i \in \K$
$(0 \leq i \leq d)$
such that $Y=\sum_{i=0}^d \beta_i X_i$.
Evaluating 
$0=XE^*_0Y$ using these equations we get 
 $0= \sum_{i=0}^d \sum_{j=0}^d \alpha_i\beta_j X_iE^*_0X_j$.
From this and 
 Lemma
\ref{lem:altbase}
we find $\alpha_i \beta_j=0$ for $0 \leq i,j\leq d$.
We assume $X\not=0$ and show $Y=0$. 
Since $X \not=0$ there exists an integer $i$
$(0 \leq i \leq d)$ such that $\alpha_i\not=0$.
Now for $0 \leq j\leq d$ we have
$\alpha_i\beta_j=0$ so $\beta_j=0$. It follows
$Y=0$.
\hfill $\Box $ \\

\noindent We finish this section with a comment.

\begin{lemma}
\label{lem:absval}
Let $\Phi$ denote the Leonard system
from Definition 
\ref{eq:ourstartingpt}.
Pick any integers $i,j$ $(0 \leq i,j\leq d)$. Then (i)--(iv) hold
below.
\begin{enumerate}
\item
$
E^*_iA^rE^*_j = 0 \quad \hbox{if } \quad r< \vert i-j \vert,
\qquad (0 \leq r \leq d).
$
\item
Suppose $ i \leq j$. Then
\begin{eqnarray}
\label{eq:leftside}
E^*_i A^{j-i}E^*_j = E^*_iAE^*_{i+1}A \cdots E^*_{j-1}AE^*_j.
\end{eqnarray}
Moreover 
$E^*_i A^{j-i}E^*_j\not=0$.
\item Suppose 
$ i \geq j$. Then
\begin{eqnarray}
\label{eq:rightside}
E^*_i A^{i-j}E^*_j = E^*_iAE^*_{i-1}A \cdots E^*_{j+1}AE^*_j.
\end{eqnarray}
Moreover
$E^*_i A^{i-j}E^*_j\not=0$.
%\item Abbreviate $r=|i-j|$.
%Then $E^*_iA^rE^*_j$ is a  basis for
%the $\K$-vector space 
%$E^*_i{\mathcal A}E^*_j$.
\end{enumerate}
\end{lemma}
\noindent {\it Proof:}
Represent the elements of $\Phi$ by matrices 
as in the proof of Theorem
\ref{eq:lsmatbasis}, and use
Lemma
\ref{lem:step1}.
\hfill $\Box $ \\

%\noindent {\it Proof:}
%Let $V$ denote an irreducible $\mathcal A$-module.
%For $0 \leq h \leq d$ let $v_h$ denote a nonzero vector in
%$E^*_hV$ and observe $v_0, v_1, \ldots, v_d$ is a basis for $V$.
%For the purpose of this proof, let us identify each element
%of $\mathcal A$ with the matrix
%in
%$\hbox{Mat}_{d+1}(\K)$
%that represents it with respect to $v_0, v_1, \ldots, v_d$.
%Adopting this point of view we find $A$ is irreducible tridiagonal
%and $A^*$ is diagonal.
%Moroever for $0 \leq h \leq d$, $E^*_h$ 
%has $hh$ entry 1 and all other entries 0.
%\\
%\noindent (i) Let $r$ be given and assume
%$r<|i-j|$.
%We show $E^*_iA^rE^*_j=0$.
%In order to do this
%we show
%each entry of $E^*_iA^rE^*_j$ is zero. 
%The $ij$ entry
%of 
%$E^*_iA^rE^*_j$
%is equal to $A^r_{ij}$ and this is zero by
%Lemma \ref{lem:step1}(i).
%All other entries of 
%$E^*_iA^rE^*_j$ are zero by the construction.
%Therefore 
%$E^*_iA^rE^*_j=0$.
%\\
%\noindent (ii)
%For each side
%of
%(\ref{eq:leftside}),
%the $ij$ entry
%is equal to 
%$A_{i,i+1}A_{i+1,i+2}\cdots A_{j-1,j}$
%and this is not zero. All other entries are zero.
%\\
%\noindent (iii) 
%For each side of
%(\ref{eq:rightside}),
%the $ij$ entry
%is
%$A_{i,i-1}A_{i-1,i-2}\cdots A_{j+1,j}$
%and this is not zero. All other entries are zero.
%\\
%\noindent (iv) The space
%$E^*_i{\cal A}E^*_j$ consists
%of
%the matrices in 
%$\hbox{Mat}_{d+1}(\K)$ 
%that have $hk$ entry 0 for $(h,k)\not=(i,j)$.
%The dimension of 
%$E^*_i{\cal A}E^*_j$ is equal to 1.
%Observe $E^*_iA^rE^*_j$ is a nonzero 
%element of
%$E^*_i{\cal A}E^*_j$
%so 
%$E^*_iA^rE^*_j$ is a basis
%for 
%$E^*_i{\cal A}E^*_j$.
%\hfill $\Box $ \\

\section{The antiautomorphism $\dagger$}

We recall the notion of an {\it antiautomorphism} of $\mathcal A$.
Let $\gamma : {\mathcal A} \rightarrow {\mathcal A}$ denote any map.
We call $\gamma $
an {\it antiautomorphism} of $\mathcal A$
whenever $\gamma$ is an
isomorphism of $\K$-vector spaces  
and  $(XY)^\gamma = Y^\gamma X^\gamma $ for
all $X, Y \in {\mathcal A}$.
For example assume ${\mathcal A} = 
\hbox{Mat}_{d+1}(\K)$.
Then $\gamma $ is an antiautomorphism of $\mathcal A$ if and only 
if there exists an invertible element $R$ in ${\mathcal A}$ such that
$X^\gamma = R^{-1}X^t R$ for all $X \in {\mathcal A}$, where
$t$ denotes transpose.
This follows from the Skolem-Noether theorem \cite[Corollary 9.122]{rotman}.

\begin{theorem}
%\cite{TLT:split}
\label{thm:dagger}
Let $A,A^*$ denote a Leonard pair in $\mathcal A$.
Then there exists a unique antiautomorphism $\dagger$ of $\mathcal A$
such that $A^\dagger = A$ and $A^{*\dagger}=A^*$.
Moreover $X^{\dagger \dagger}=X$
for all $X \in {\mathcal A}$.
\end{theorem}
\noindent {\it Proof:}
Concerning existence, 
let $V$ denote an irreducible $\mathcal A$-module. 
By 
Definition
\ref{def:lprecall}(i) there exists a basis
for $V$ with respect to which the matrix representing
$A$ is irreducible tridiagonal and the matrix representing
$A^*$ is diagonal.
Let us denote this basis by  
$v_0, v_1, \ldots, v_d$.
For $X \in \mathcal A$ let $X^\sigma$ denote the matrix
in 
$\mbox{Mat}_{d+1}(\K)$ that represents $X$
with respect to
the basis
$v_0, v_1, \ldots, v_d$.
We observe
$\sigma : {\mathcal A}\rightarrow 
\mbox{Mat}_{d+1}(\K)$ 
is an isomorphism of $\K$-algebras.
We abbreviate $B=A^\sigma$ and observe
$B$ is irreducible tridiagonal.
We abbreviate $B^*=A^{*\sigma}$ and observe
$B^*$ is diagonal.
Let $D$ denote the diagonal matrix in  
$\mbox{Mat}_{d+1}(\K)$ 
that has $ii$ entry
\begin{eqnarray*}
D_{ii} = \frac{B_{01}B_{12} \cdots B_{i-1,i}}
{B_{10}B_{21} \cdots B_{i,i-1}}
\qquad \qquad (0 \leq i \leq d).
\end{eqnarray*}
It is routine to verify $D^{-1}B^tD=B$. Each of $D, B^*$ is
diagonal so $DB^*=B^*D$; also
$B^{*t}=B^*$ so $D^{-1}B^{*t}D=B^*$. Let
$
\gamma :
\mbox{Mat}_{d+1}(\K)
\rightarrow
\mbox{Mat}_{d+1}(\K)
$
denote the map that satisfies
$X^{\gamma} = D^{-1}X^tD$ for all $X \in 
\mbox{Mat}_{d+1}(\K)$.
We observe $\gamma $ is an antiautomorphism
of 
$\mbox{Mat}_{d+1}(\K)$ such that
$B^\gamma=B$ and $B^{*\gamma}=B^*$.
We define the map
$
\dagger :{ \mathcal A}
\rightarrow 
{ \mathcal A}
$
to be the  composition
$\dagger = \sigma \gamma \sigma^{-1}$.
We observe $\dagger $ is an antiautomorphism
of $\mathcal A$ such that $A^\dagger=A$
and $A^{*\dagger}=A^*$.
We have now shown there exists an
antiautomorphism
$\dagger $ of 
$\mathcal A$ such that
 $A^\dagger=A$
and $A^{*\dagger}=A^*$.
This antiautomorphism is unique since
$A,A^*$ together generate $\mathcal A$.
The map $X\rightarrow X^{\dagger \dagger}$ is
an isomorphism of $\K$-algebras from
$\mathcal A$ to itself.
This isomorphism is the identity
since
 $A^{\dagger \dagger}=A$,
 $A^{*\dagger \dagger}=A^*$,
 and since $A,A^*$ together generate 
$\mathcal A$.
\hfill $\Box $ \\

\begin{definition} 
\label{def:dag}
\rm
Let $A,A^*$ denote a Leonard pair in
 $\mathcal A$.
By the {\it antiautomorphism which corresponds to $A,A^*$}
 we mean the 
map $\dagger: {\mathcal A}\rightarrow {\mathcal A}$ from Theorem
\ref{thm:dagger}. Let 
$\Phi=(A;A^*;\lbrace E_i \rbrace_{i=0}^d; $ $  
\lbrace E^*_i\rbrace_{i=0}^d)$ denote a Leonard system in
 $\mathcal A$.
By the {\it antiautomorphism which corresponds to $\Phi$}
 we mean the
antiautomorphism which corresponds to the Leonard pair $A,A^*$. 
\end{definition}

\begin{lemma}
\label{cor:eistab}
Let $\Phi$ denote the Leonard system from Definition
\ref{eq:ourstartingpt} and let 
$\dagger $  denote the corresponding antiautomorphism.
Then the following (i), (ii) hold.
\begin{enumerate}
\item
Let $\mathcal D$ denote the subalgebra of $\cal A$ generated
by $A$. Then $X^\dagger=X$ for all $X \in \cal D$; in particular
$E_i^\dagger=E_i $ for $0 \leq i\leq d$. 
\item 
Let ${\mathcal D}^*$ denote the subalgebra of $\cal A$ generated
by $A^*$. Then $X^\dagger=X$ for all $X \in {\cal D}^*$; in particular
$E_i^{*\dagger}=E^*_i $ for $0 \leq i\leq d$. 
\end{enumerate}
\end{lemma}
\noindent {\it Proof:}
 (i) The sequence
$A^i$ $(0 \leq i \leq d)$ is a basis for
 the $\K$-vector space $\mathcal D$.
Observe $\dagger $ stabilizes $A^i$ for $0 \leq i \leq d$.
The result follows.
\\
\noindent (ii) Similar to the proof of (i) above. 
\hfill $\Box $ \\

\section{The scalars $a_i, x_i$}

\noindent In this section we introduce some scalars
that will help us describe Leonard systems.

\begin{definition}
\label{def:aixi}
\rm
Let $\Phi$ denote the Leonard system
from Definition 
\ref{eq:ourstartingpt}.
We define
\begin{eqnarray}
 a_i &=&\mbox{tr}(E^*_iA)\qquad    \quad   
  (0 \leq i \leq d),  \label{eq:aitr}
\\
x_i &=&\mbox{tr}(E^*_iAE^*_{i-1}A)     
\qquad  \quad   (1 \leq i \leq d),\label{eq:xitr} 
\end{eqnarray}
where $\mbox{tr}$ denotes trace.
For notational convenience we define $x_0=0$.
\end{definition}

\noindent 
We have a comment.

\begin{lemma}
\label{lem:bmat}
Let $\Phi$ denote the Leonard system
from Definition
\ref{eq:ourstartingpt} and let $V$ denote an irreducible
$\cal A$-module. For $0 \leq i \leq d$ let $v_i$ denote
a nonzero vector in $E^*_iV$ and observe $v_0, v_1, \ldots, v_d$
is a basis for $V$. Let $B$ denote the matrix in
$\hbox{Mat}_{d+1}(\K)$ that represents $A$ with respect to
$v_0, v_1, \ldots, v_d$. 
We observe $B$ is irreducible tridiagonal.
The following (i)--(iii) hold. 
\begin{enumerate}
\item
$B_{ii}=a_i$ $(0 \leq i \leq d)$.
\item
 $B_{i,i-1}B_{i-1,i}=x_i$  $(1 \leq i \leq d)$.
\item $x_i\not=0$ $(1 \leq i \leq d)$. 
\end{enumerate}
\end{lemma}

\noindent {\it Proof:}
(i), (ii) For $0 \leq i \leq d$ the matrix in
$\hbox{Mat}_{d+1}(\K)$ that represents $E^*_i$ with respect to
$v_0, v_1, \ldots, v_d$ has $ii$ entry
$1$ and all other entries 0. The result follows in
view of
Definition
\ref{def:aixi}.
\\
\noindent (iii) Immediate from (ii) and  since $B$ is irreducible.
\hfill $\Box $ \\

\begin{theorem}
\label{lem:monbasis}
Let $\Phi$ denote the Leonard system
from Definition
\ref{eq:ourstartingpt}. Let $V$ denote an irreducible
$\cal A$-module and  
let $v$ denote a nonzero vector in 
$E^*_0V$. Then
for $0 \leq i \leq d$ the vector
$E^*_iA^iv$ is nonzero and hence
a basis for $E^*_iV$. 
Moreover the sequence
\begin{eqnarray}
E^*_iA^iv \qquad (0 \leq i \leq d)
\label{eq:monbasis}
\end{eqnarray}
is a basis for $V$.
\end{theorem}
\noindent {\it Proof:}
We show $E^*_iA^iv\not=0$ for $0 \leq i \leq d$.
Let $i$ be given.
Setting $j=0$ in
Lemma
\ref{lem:absval}(iii)
we find
$E^*_iA^iE^*_0\not=0$.
Therefore
 $E^*_iA^iE^*_0V\not=0$.
The space $E^*_0V$  is spanned by
$v$
so $E^*_iA^iv\not=0$ as desired.
The remaining claims follow.
\hfill $\Box $ \\

\begin{theorem}
\label{lem:monbasis2}
Let $\Phi$ denote the Leonard system
from Definition
\ref{eq:ourstartingpt} and let
the scalars $a_i, x_i$ be as in
Definition
\ref{def:aixi}. Let
$V$ denote an irreducible
$\cal A$-module.
With respect to the basis for $V$ given
in 
(\ref{eq:monbasis})
 the matrix
that represents $A$ is equal to
\begin{equation}
\label{eq:monmat}
\left(
\begin{array}{ c c c c c c}
a_0 & x_1  &      &      &   &{\bf 0} \\
1 & a_1  &  x_2   &      &   &  \\
  & 1  &  \cdot    & \cdot  &   & \\
  &   & \cdot     & \cdot  & \cdot   & \\
  &   &           &  \cdot & \cdot & x_d \\
{\bf 0} &   &   &   & 1 & a_d  
\end{array}
\right) .
\end{equation} 
\end{theorem}

\noindent {\it Proof:}
With reference to
(\ref{eq:monbasis}) abbreviate
$v_i=E^*_iA^iv$ for $0 \leq i\leq d$.
Let $B$ denote the matrix in
$\hbox{Mat}_{d+1}(\K)$ that represents
$A$ with respect to
$v_0, v_1, \ldots, v_d$.
We show $B$ is equal to
(\ref{eq:monmat}). In view of
Lemma
\ref{lem:bmat} it suffices to show
$B_{i,i-1}=1$ for $1 \leq i \leq d$.
For $0 \leq i \leq d$ the 
matrix $B^i$ represents $A^i$ with
respect to 
$v_0, v_1, \ldots, v_d$; 
therefore 
$A^iv_0 = \sum_{j=0}^d(B^i)_{j0}v_j$.
Applying $E^*_i$ and using $v_0=v$
we find
$v_i = (B^i)_{i0}v_i$ so
$(B^i)_{i0}=1$, forcing
$
B_{i,i-1}
\cdots
B_{21}
B_{10}
=1$
by Lemma
\ref{lem:step1}.
We have shown
$
B_{i,i-1}
\cdots
B_{21}
B_{10}
=1$ for $1 \leq i \leq d$
so $B_{i,i-1}=1$ for $1 \leq i \leq d$. 
We now see $B$ is equal to
(\ref{eq:monmat}). 
\hfill $\Box $ \\

\noindent We finish this section with a few
comments.

\begin{lemma}
\label{lem:xiprod}
Let $\Phi$ denote the Leonard system
from Definition
\ref{eq:ourstartingpt} and
let the scalars $a_i, x_i$ be as in 
Definition
\ref{def:aixi}.
Then the following (i)--(iii) hold.
\begin{enumerate}
\item
$E^*_iAE^*_i = a_iE^*_i \qquad (0 \leq i \leq d)$.
\item 
$E^*_iAE^*_{i-1}AE^*_i = x_iE^*_i \qquad \qquad (1 \leq i \leq d)$.
\item
$E^*_{i-1}AE^*_{i}AE^*_{i-1} = x_iE^*_{i-1} \qquad \qquad (1 \leq i \leq d)$.
\end{enumerate}
\end{lemma}

\noindent {\it Proof:}
(i) Observe
$E^*_i$ is a basis for
$E^*_i{\cal A}E^*_i$. By this and since
$E^*_i AE^*_i$ is contained in
$E^*_i{\cal A}E^*_i$ we find
there exists
 $\alpha_i \in \K$ such that
$E^*_i AE^*_i = \alpha_i E^*_i.$
Taking the trace of both sides and using
$\mbox{tr}(XY)=\mbox{tr}(YX)$, $\mbox{tr}(E^*_i)=1$ 
we
find $a_i = \alpha_i$.
\\
\noindent (ii) We mentioned above that
 $E^*_i$ is a basis for
$E^*_i{\cal A}E^*_i$. By this
and since 
$E^*_i AE^*_{i-1}AE^*_i$ is contained in
$E^*_i{\cal A}E^*_i$ we find
there exists
 $\beta_i \in \K$ such that
$E^*_i AE^*_{i-1}AE^*_i=\beta_i E^*_i$.
Taking the trace of both sides we 
find 
$x_i=\beta_i$.
\\
\noindent (iii)  Similar to the proof of (ii) above.
\hfill $\Box $ \\

\begin{lemma}
\label{lem:xrec}
Let $\Phi$ denote the Leonard system
from Definition  
\ref{eq:ourstartingpt} and 
let the scalars $ x_i$ be as in 
Definition
\ref{def:aixi}. Then the following (i), (ii) hold.
\begin{enumerate}
\item 
$ 
E^*_jA^{j-i}E^*_iA^{j-i}E^*_j = x_{i+1} x_{i+2} \cdots x_j E^*_j
\qquad (0 \leq i \leq j\leq d)$.
\item
$ 
E^*_iA^{j-i}E^*_jA^{j-i}E^*_i = x_{i+1} x_{i+2} \cdots x_j E^*_i
\qquad (0 \leq i \leq j\leq d)$.
\end{enumerate}
\end{lemma}
\noindent {\it Proof:}
(i) Evaluate the expression on the left using
Lemma
\ref{lem:absval}(ii), (iii) 
and Lemma
\ref{lem:xiprod}(ii).
\\
(ii)
Evaluate the expression on the left using
Lemma
\ref{lem:absval}(ii), (iii) 
and Lemma
\ref{lem:xiprod}(iii). 
\hfill $\Box $ \\

\section{The polynomials $p_i$}

\medskip
\noindent  
In this section we begin our discussion of 
polynomials. 
We will use
the following notation. 
Let $\lambda $ denote
an indeterminate. We let 
 $\K \lbrack \lambda \rbrack $ 
denote the $\K$-algebra consisting of all polynomials
in $\lambda $ that have coefficients in $\K$.  
For the rest of this paper all polynomials that we discuss
are assumed to lie in 
 $\K \lbrack \lambda \rbrack $.

\begin{definition}
\label{def:pi}
\rm
Let $\Phi$ denote the Leonard system
from Definition 
\ref{eq:ourstartingpt} and let
 the scalars $a_i, x_i$ be as in 
Definition
\ref{def:aixi}. We define
a sequence of 
polynomials $p_0, p_1, \ldots, p_{d+1}$ by
\begin{eqnarray}
p_0&=&1,
\\
\lambda p_i &=& p_{i+1}+a_ip_i+x_ip_{i-1}
\qquad (0 \leq i \leq d),
\label{eq:prec}
\end{eqnarray}
where $p_{-1}=0$.
We observe $p_i$ is monic with degree exactly $i$ for $0 \leq i \leq d+1$.
\end{definition}

\begin{lemma}
\label{lem:pareq}
Let $\Phi$ denote the Leonard system
from Definition 
\ref{eq:ourstartingpt} and let the polynomials $p_i$ be as in
Definition 
\ref{def:pi}.
Let $V$ denote
an irreducible $\mathcal A$-module and 
let $v$ denote a nonzero vector in $E^*_0V$.
Then
$p_i(A)v=E^*_iA^iv$ for $0 \leq i  \leq d$ and
$p_{d+1}(A)v=0$.
\end{lemma}
\noindent {\it Proof:}
We
abbreviate
$v_i=p_i(A)v$ for $0 \leq i\leq d+1$.
We define
$v'_i=E^*_iA^iv$ for $0 \leq i \leq d$
and $v'_{d+1}=0$.
We show $v_i=v'_i$ for $0 \leq i \leq d+1$.
From the construction
$v_0=v$ and $v'_0=v$ so $v_0=v'_0$.
From
(\ref{eq:prec}) we obtain
%Lemma
%\ref{lem:monbasis2} we find
\begin{eqnarray}
Av_i = v_{i+1} + a_iv_i + x_iv_{i-1}
\qquad (0 \leq i \leq d)
\label{eq:rel}
\end{eqnarray}
where  $v_{-1}=0$.
From 
Theorem \ref{lem:monbasis2} we find
\begin{eqnarray}
Av'_i = v'_{i+1} + a_iv'_i + x_iv'_{i-1}
\qquad (0 \leq i \leq d)
\label{eq:rel2}
\end{eqnarray}
where  $v'_{-1}=0$.
Comparing
(\ref{eq:rel}), 
(\ref{eq:rel2}) and using $v_0=v'_0$
we find
$v_i=v'_i$ for $0 \leq i \leq d+1$. The result
follows.
\hfill $\Box $ \\

\noindent We mention a few consequences of Lemma
\ref{lem:pareq}.

\begin{theorem}
\label{thm:psend}
Let $\Phi$ denote the Leonard system
from Definition 
\ref{eq:ourstartingpt} and
let the 
polynomials $p_i$ be 
as in 
Definition
\ref{def:pi}. Let $V$ denote an irreducible 
$\mathcal A$-module.
 Then
$$
 p_i(A)E^*_0V = E^*_iV \quad (0 \leq i \leq d).
$$
\end{theorem}
\noindent {\it Proof:}
Let $v$ denote a nonzero vector in $E^*_0V$.
Then
$p_i(A)v=E^*_iA^iv$ by Lemma
\ref{lem:pareq}. 
Observe $v$ is a basis for $E^*_0V$.
By Theorem 
\ref{lem:monbasis} we find
$E^*_iA^iv$ is a basis for $E^*_iV$.
Combining these facts we find
 $p_i(A)E^*_0V = E^*_iV$.
\hfill $\Box $ \\

\begin{theorem}
\label{thm:pimon}
Let $\Phi$ denote the Leonard system
from Definition 
\ref{eq:ourstartingpt} and 
let the 
polynomials $p_i$ be 
as in 
Definition
\ref{def:pi}.
Then 
\begin{eqnarray} 
 p_i(A)E^*_0 = E^*_iA^iE^*_0 \qquad (0 \leq i \leq d).
\label{eq:pimon}
\end{eqnarray}
\end{theorem}
\noindent {\it Proof:}
Let the integer $i$ be given and abbreviate
$\Delta=
 p_i(A) -E^*_iA^i$.
We show
 $\Delta E^*_0=0$.
In order to do this we show
 $\Delta E^*_0V=0$,
 where $V$ denotes an irreducible $\mathcal A$-module.
Let $v$ denote a nonzero vector in $E^*_0V$ and recall
 $v$ is a basis for $E^*_0V$.
By Lemma
\ref{lem:pareq} we have
$\Delta v=0$ so
 $\Delta E^*_0V=0$.
Now 
 $\Delta E^*_0=0$ so
 $p_i(A)E^*_0 = E^*_iA^iE^*_0$.
\hfill $\Box $ \\

\begin{theorem}
\label{thm:pseq}
Let $\Phi$ denote the Leonard system
from Definition 
\ref{eq:ourstartingpt} and 
let the polynomial $p_{d+1}$ be as in
Definition
\ref{def:pi}. Then the following (i), (ii) hold.
\begin{enumerate}
\item 
 $p_{d+1}$ is both the minimal polynomial and the
 characteristic polynomial of
$A$.
\item
$
p_{d+1} = \prod_{i=0}^d (\lambda - \theta_i).
$
\end{enumerate}
\end{theorem}
\noindent {\it Proof:}
(i)
%(resp. $c$
%Recall the characteristic polynomial of $A$
%is $\mbox{det}(\lambda I-A)$.
We first show $p_{d+1}$ is equal to the
minimal polynomial of $A$.
Recall
 $I, A, \ldots, A^d$ are linearly
independent and that
$p_{d+1}$ is monic with degree $d+1$.
We show
$p_{d+1}(A)=0$.
Let $V$ denote an irreducible $\mathcal A$-module.
Let $v$ denote a nonzero vector in $E^*_0V$
and recall $v$ is a basis for $E^*_0V$.
From Lemma
\ref{lem:pareq}
we find
$p_{d+1}(A)v=0$.
It follows
$p_{d+1}(A)E^*_0V=0$ so $p_{d+1}(A)E^*_0=0$.
Applying Lemma
\ref{lem:oneiszero}
(with $X=p_{d+1}(A)$ and $Y=I$)
we find $p_{d+1}(A)=0$.
We have now shown $p_{d+1}$ is the minimal polynomial
of $A$.
By definition 
the characteristic polynomial of $A$
is equal to 
$\mbox{det}(\lambda I-A)$.
This polynomial is monic with degree $d+1$ and
has $p_{d+1}$ as a factor; therefore
it is  equal to $p_{d+1}$.
\\
\noindent (ii) 
For $0 \leq i \leq d$ the scalar $\theta_i$
is an eigenvalue
of $A$ and therefore a root of the characteristic polynomial
of $A$.
\hfill $\Box $ \\

\noindent The following result will be useful.

\begin{lemma}
\label{lem:eispoly}
Let $\Phi$ denote the Leonard system
from Definition 
\ref{eq:ourstartingpt} and
let the polynomials
 $p_i$ be 
as in 
Definition
\ref{def:pi}.
Let the 
scalars $x_i$ be as in
Definition
\ref{def:aixi}.
Then
\begin{eqnarray}
E^*_i = \frac{
p_i(A)E^*_0p_i(A)
}{
x_1x_2\cdots x_i
}
\qquad (0 \leq i \leq d).
\end{eqnarray}
\end{lemma}
\noindent {\it Proof:}
%Let $\dagger$ denote the antiautomorphism
%of 
%$\mathcal A$ associated with
%$\Phi$.
Let $\dagger:{\mathcal A}\rightarrow {\mathcal A}$ 
denote the antiautomorphism
which corresponds to 
$\Phi$.
From Theorem
\ref{thm:pimon}
we have 
$p_i(A)E^*_0
=
E^*_iA^iE^*_0
$.
Applying $\dagger$
we find
$E^*_0p_i(A)
=E^*_0A^iE^*_i
$.
From these comments  we find
\begin{eqnarray*}
p_i(A)E^*_0p_i(A) &=&
E^*_iA^iE^*_0
A^iE^*_i
\\
&=& x_1 x_2 \cdots x_iE^*_i
\end{eqnarray*}
in view of 
Lemma
\ref{lem:xrec}(i).
The result follows.
\hfill $\Box $ \\

%\noindent 
%We finish this section with a comment.
%
%\begin{lemma}
%\label{lem:aisum}
%Let $\Phi$ denote the Leonard system
%from Definition
%\ref{eq:ourstartingpt} and
%let the 
%polynomials $p_i$ be 
%as in 
%Definition
%\ref{def:pi}. Let the scalars $a_i$ be as in
%Definition \ref{def:aixi}.
%Then for $0 \leq i \leq d$ the coefficient of $\lambda^i$ in 
%$p_{i+1}$ is equal to
%$-\sum_{j=0}^i a_j $.
%\end{lemma}
%
%\noindent {\it Proof:}
%Let $\alpha_i$ denote the coefficient of $\lambda^i$ in $p_{i+1}$.
%Computing the coefficient of $\lambda^i$ in
%(\ref{eq:prec}) we find $\alpha_{i-1}=\alpha_i+a_i$
%for $0 \leq i \leq d$, where $\alpha_{-1}=0$. It follows
%$\alpha_i = -\sum_{j=0}^i a_j$
%for $0 \leq i \leq d$.
%\hfill $\Box $ \\

\section{The scalars $\nu, m_i$}

\medskip
\noindent
In this section 
we introduce some more scalars that 
will help us describe Leonard systems.

\begin{definition} 
\label{def:mid}
\rm
Let $\Phi$ denote the Leonard system
from Definition  
\ref{eq:ourstartingpt}.
We define 
\begin{equation}
\label{eq:mid}
m_i =  \mbox{tr}(E_iE^*_0) \qquad \qquad (0 \leq i \leq d).
\end{equation}
\end{definition}

\begin{lemma}
\label{lem:mid}
Let $\Phi$ denote the Leonard system
from Definition
\ref{eq:ourstartingpt}.
Then (i)--(v) hold below.
\begin{enumerate}
\item $E_iE^*_0E_i = m_i E_i
\qquad (0 \leq i\leq d)$.
\item $E^*_0E_iE^*_0 = m_i E^*_0
\qquad (0 \leq i\leq d)$.
\item $m_i\not=0
\qquad (0 \leq i\leq d)$.
\item
$\sum_{i=0}^d m_i = 1$.
\item $m_0 = m^*_0$.
\end{enumerate}
\end{lemma}
\noindent {\it Proof:}
(i) Observe 
$E_i$ is a basis for $E_i{\mathcal A}E_i$.
By this and since
$E_iE^*_0E_i$ is contained in
 $E_i{\mathcal A}E_i$, there exists
 $\alpha_i \in \K$ 
such that 
 $E_iE^*_0E_i=\alpha_i E_i$. 
Taking the trace of both
 sides in this equation and using
$\hbox{tr}(XY)= \hbox{tr}(YX)$, $\mbox{tr}(E_i)=1$
  we find $\alpha_i=m_i$.
\\
\noindent (ii) Similar to the proof of (i).
\\
\noindent (iii) 
Observe $m_iE_i$ is equal to $E_iE^*_0E_i$ by part (i) above
and 
$E_iE^*_0E_i$ is nonzero by 
Corollary
\ref{cor:eibasis}. 
It follows $m_iE_i\not=0$ so $m_i\not=0$.
\\
\noindent (iv) 
Multiply each term
in the equation
$\sum_{i=0}^d E_i = I$
 on the right by $E^*_0$, and then take the trace.
Evaluate the result using
Definition
\ref{def:mid}.
\\
\noindent (v) 
The elements $E_0E^*_0$ and $E^*_0E_0$ have the same trace. 
\hfill $\Box $ \\

\begin{definition}
\label{def:n} 
\rm
 Let $\Phi$ denote the Leonard system
from Definition 
\ref{eq:ourstartingpt}. 
Recall $m_0=m^*_0$ by Lemma
\ref{lem:mid}(v); 
 we let $\nu$ denote the multiplicative inverse
of this common value. We observe
$\nu= \nu^*$. We 
emphasize
\begin{equation}
\mbox{tr}(E_0E^*_0) = \nu^{-1}.
\label{eq:sizedef}
\end{equation}
\end{definition}

\begin{lemma}
\label{lem:threeone}
Let $\Phi$ denote the Leonard 
system from Definition
\ref{eq:ourstartingpt}  and let 
the scalar $\nu$ be as in Definition
\ref{def:n}. Then the following (i), (ii) hold.
\begin{enumerate}
\item
$\nu E_0E^*_0E_0 = E_0$.
\item 
$\nu E^*_0E_0E^*_0 = E^*_0.$
\end{enumerate}
\end{lemma}
\noindent {\it Proof:}
(i) Set $i=0$ in Lemma
\ref{lem:mid}(i) and recall $m_0={\nu}^{-1}$.
\\
(ii) Set $i=0$ in Lemma
\ref{lem:mid}(ii) and recall $m_0={\nu}^{-1}$.
\hfill $\Box $ \\

%\begin{theorem}
%Let $\Phi$ denote the Leonard system
%from Definition 
%\ref{eq:ourstartingpt} and let the
%polynomials $p_i$ be as in Definition \ref{def:pi}.
%Let the scalars $\theta_i$ be as in
%Definition \ref{def:evseq} and let the 
%scalars $m_i$
%be as in
%Definition \ref{def:mid}.
%Then 
%\begin{eqnarray}
%p_i(\theta_j) = m^{-1}_j \mbox{tr}(E_jE^*_iA^iE^*_0)
%\qquad (0 \leq i,j\leq d).
%\end{eqnarray}
%\end{theorem}
%\noindent {\it Proof:}
%Using Theorem 
%\ref{thm:pimon}
%we find
%\begin{eqnarray*}
%\mbox{tr}(E_jE^*_iA^iE^*_0) &=& 
%\mbox{tr}(E_jp_i(A)E^*_0)
%\\
%&=& p_i(\theta_j)\mbox{tr}(E_jE^*_0)
%\\
%&=& p_i(\theta_j)m_j.
%\end{eqnarray*}
%The result follows.
%\hfill $\Box $ \\
%

\section{The standard basis}

In this section we 
 discuss the notion of a  {\it standard basis}.
We begin with a comment.

\begin{lemma}
\label{lem:stbasisp}
Let $\Phi$ denote the Leonard system
from Definition 
\ref{eq:ourstartingpt} and let $V$ denote an irreducible
$\mathcal A$-module.
Then 
\begin{eqnarray}
E^*_iV = E^*_iE_0V \qquad (0 \leq i \leq d).
\end{eqnarray}
\end{lemma}
\noindent {\it Proof:}
The space $E^*_iV$ has dimension 1 and contains
$E^*_iE_0V$. We show
$E^*_iE_0V\not=0$.
Applying
Corollary
\ref{cor:eibasis}
to $\Phi^*$ we find $E^*_iE_0\not=0$. It follows
 $E^*_iE_0V\not=0$.
We conclude
 $E^*_iV=E^*_iE_0V$.
\hfill $\Box $ \\

\begin{lemma}
\label{lem:stbasis}
Let $\Phi$ denote the Leonard system
from Definition 
\ref{eq:ourstartingpt} and
let
$V$ denote an irreducible $\mathcal A$-module.
Let 
$u$ denote a nonzero vector in  $E_0V$. Then
for $0 \leq i \leq d$ the vector $E^*_iu$ is nonzero
and hence a basis
for $E^*_iV$. Moreover the sequence
\begin{equation}
E^*_0u, E^*_1u, \ldots, E^*_du
\label{eq:stbasisint}
\end{equation}
is a basis for $V$.
\end{lemma}
\noindent {\it Proof:}
Let the integer $i$ be given.
We show $E^*_iu\not=0$.
Recall $E_0V$ has dimension 1 and
$u$ is a nonzero vector in 
 $E_0V$  
so $u$ spans  $E_0V$. Applying $E^*_i$ we find
$E^*_iu$ spans $E^*_iE_0V$. The space $E^*_iE_0V$ is nonzero by
Lemma
\ref{lem:stbasisp}
so 
$E^*_iu$ is nonzero.
The remaining assertions are clear.
\hfill $\Box $ \\

\begin{definition}
\label{def:stbasis}
\rm
Let $\Phi$ denote the Leonard system
from Definition 
\ref{eq:ourstartingpt} and let $V$ denote an irreducible
$\mathcal A$-module.
By a {\it $\Phi$-standard basis} for $V$,
we mean
a sequence 
$$
E^*_0u, E^*_1u, \ldots, E^*_du,
$$
where $u$ is a nonzero vector in $E_0V$.
\end{definition}

\noindent We give a few characterizations of the standard
basis.

\begin{lemma}
\label{lem:eggechar}
Let $\Phi$ denote the Leonard system from
Definition
\ref{eq:ourstartingpt}  and let $V$ denote an irreducible
$\mathcal A$-module.
Let $v_0, v_1, \ldots, v_d$ denote
a sequence of vectors in $V$, not all 0. Then this sequence is a 
$\Phi$-standard basis for $V$ if and only if both 
(i), (ii)  hold below.
\begin{enumerate}
\item $v_i \in E^*_iV$ for $0 \leq i \leq d$.
\item $\sum_{i=0}^d v_i \in E_0V$.
\end{enumerate}
\end{lemma} 
\noindent {\it Proof:}
To prove the lemma in one direction,
 assume 
$v_0, v_1, \ldots, v_d$ is a
$\Phi$-standard basis for $V$.
By Definition
\ref{def:stbasis}
there exists a nonzero $u \in E_0V$ such
that $v_i = E^*_iu$ for $0 \leq i \leq d$.
Apparently $v_i \in E^*_iV$ for $0 \leq i \leq d$ so (i) holds.
Let $I$ denote the identity element of $\mathcal A$ and
recall $I=\sum_{i=0}^d E^*_i$.
Applying this to $u$ we find
$u=\sum_{i=0}^d v_i $ and (ii) follows.
We have now proved the lemma in one direction. To prove the lemma
in the other direction, assume $v_0, v_1, \ldots, v_d$ satisfy
(i), (ii) above.
We define
$u=\sum_{i=0}^d v_i$ and observe $u\in E_0V$.
Using (i) 
we find $E^*_iv_j=\delta_{ij}v_j$ for
$0 \leq i,j\leq d$; it follows
$v_i = E^*_iu$ for $0 \leq i \leq d$.
Observe $u \not=0$  since at least one
of $v_0, v_1, \ldots, v_d$ is nonzero.
Now 
$v_0, v_1, \ldots, v_d$ is 
a 
$\Phi$-standard  basis for $V$ by Definition
\ref{def:stbasis}.
\hfill $\Box $ \\

\noindent  We recall some notation.
Let $d$ denote a nonnegative
integer and let $B$ denote a matrix in 
$\hbox{Mat}_{d+1}(\K)$. Let $\alpha $ denote a scalar in $\K$.
Then $B$ is said to have {\it constant row  sum $\alpha$} whenever
$B_{i0}+B_{i1}+\cdots + B_{id}=\alpha$ for $0 \leq i \leq d$.

\begin{lemma}
\label{lem:rowsum}
Let $\Phi$ denote the Leonard system from  Definition
\ref{eq:ourstartingpt} and let the scalars $\theta_i, \theta^*_i$
be as in Definition \ref{def:evseq}.
Let $V$ denote an irreducible $\mathcal A$-module and
let $v_0, v_1, \ldots, v_d$ denote a basis
for $V$.
Let $B$ (resp. $B^*$) denote the matrix in 
$\hbox{Mat}_{d+1}(\K)$
that represents
$A$ (resp. $A^*)$ with respect to 
this basis.
Then 
 $v_0, v_1, \ldots, v_d$ is a $\Phi$-standard basis
 for $V$ if and only if both (i), (ii) hold below.
\begin{enumerate}
\item $B$ has constant row sum $\theta_0$.
\item $B^*=\hbox{diag}(\theta^*_0, \theta^*_1, \ldots, \theta^*_d)$.
\end{enumerate}
\end{lemma}
\noindent {\it Proof:}
Observe $A \sum_{j=0}^d v_j = \sum_{i=0}^d v_i(B_{i0}+B_{i1}+\cdots B_{id})$.
Recall $E_0V$ is the eigenspace for $A$ and eigenvalue $\theta_0$.
Apparently 
$B$ has constant row sum $\theta_0$ if and only if
$\sum_{i=0}^d v_i \in E_0V$.
Recall that for $0 \leq i \leq d$,
$E^*_iV$ is the eigenspace for $A^*$ and eigenvalue $\theta^*_i$.
Apparently 
$B^*=\hbox{diag}(\theta^*_0, \theta^*_1, \ldots, \theta^*_d)$
if and only if $v_i \in E^*_iV$ for $0 \leq i \leq d$. 
The result follows in view of Lemma
\ref{lem:eggechar}.
\hfill $\Box $ \\

\begin{definition}
\label{def:flatcon}
\rm
Let $\Phi$ denote the Leonard system from 
Definition
\ref{eq:ourstartingpt}. We define
a map $\flat : {\mathcal A}\rightarrow 
\mbox{Mat}_{d+1}(\K)$ as follows.
Let 
$V$ denote an irreducible $\mathcal A$-module.
For all $X \in {\mathcal  A}$ we let $X^\flat $ denote
the matrix 
in 
$\mbox{Mat}_{d+1}(\K)$
that represents $X$ with respect to 
a $\Phi$-standard basis for $V$.
We observe $\flat : {\mathcal A} \rightarrow 
\mbox{Mat}_{d+1}(\K)$ is an isomorphism of
$\K$-algebras.
\end{definition}

\begin{lemma}
\label{lem:firstcom}
Let $\Phi$ denote the Leonard system from 
Definition
\ref{eq:ourstartingpt}  and let the scalars
$\theta_i, \theta^*_i$
be as in
Definition
\ref{def:evseq}.
Let the map $\flat : {\mathcal A}\rightarrow \mbox{Mat}_{d+1}(\K)$
be as in Definition
\ref{def:flatcon}.
Then (i)--(iii) hold below.
\begin{enumerate}
\item $A^\flat $ has constant row sum $\theta_0$.
\item $A^{*\flat}=\hbox{diag}(\theta^*_0,\theta^*_1, \ldots, \theta^*_d)$.
\item For $0 \leq i \leq d$ the matrix
$E_i^{*\flat}$ has $ii$ entry 1
and all other entries 0.
\end{enumerate}
\end{lemma}
\noindent {\it Proof:}
(i), (ii) Combine 
Lemma 
\ref{lem:rowsum} and
Definition
\ref{def:flatcon}.
\\
(iii) Immediate from
Lemma
\ref{lem:step0}(i).
%Definition
%\ref{def:stbasis}.
\hfill $\Box $ \\

\section{The scalars $b_i, c_i$}

\noindent In this section we consider some scalars
that arise naturally in the context of the standard basis.

\begin{definition}
\label{def:sharpmp}
\rm
Let $\Phi$ denote the Leonard system
from Definition 
\ref{eq:ourstartingpt}
and let the map
$\flat :{\mathcal A}\rightarrow
\mbox{Mat}_{d+1}(\K)
$ be as in
Definition
\ref{def:flatcon}.
For $0 \leq i \leq d-1$ we let $b_i$ denote the
$i,i+1$ entry of $A^\flat$.
For
 $1 \leq i \leq d$ we let $c_i$ denote the
$i,i-1$ entry of $A^\flat$.
We observe 
\begin{equation}
\label{eq:matrepls}
A^\flat =  \left(
\begin{array}{ c c c c c c}
a_0 & b_0  &      &      &   &{\bf 0} \\
c_1 & a_1  &  b_1   &      &   &  \\
  & c_2  &  \cdot    & \cdot  &   & \\
  &   & \cdot     & \cdot  & \cdot   & \\
  &   &           &  \cdot & \cdot & b_{d-1} \\
{\bf 0} &   &   &   & c_d & a_d  
\end{array}
\right), 
\end{equation} 
where the $a_i$ are from 
Definition
\ref{def:aixi}.
For notational convenience we define $b_d=0$ and $c_0=0$.
\end{definition}

\begin{lemma}
\label{def:bici}
Let $\Phi$ denote the Leonard system
from Definition 
\ref{eq:ourstartingpt}
and let the scalars $b_i, c_i$ be as in
Definition
\ref{def:sharpmp}.
 Then with reference
to Definition
\ref{def:evseq} and Definition
\ref{def:aixi} the following (i), (ii) hold.
\begin{enumerate}
\item $b_{i-1}c_i=x_i \qquad (1 \leq i \leq d)$.
\item $c_i + a_i + b_i=\theta_0 \qquad (0 \leq i \leq d)$.
\end{enumerate}
\end{lemma}
\noindent {\it Proof:}
(i) Apply Lemma
\ref{lem:bmat}(ii) with $B=A^\flat$.
\\
(ii) Combine
(\ref{eq:matrepls})
and
Lemma
\ref{lem:firstcom}(i).
\hfill $\Box $ \\

\begin{lemma}
\label{def:bici2}
Let $\Phi$ denote the Leonard system
from Definition
\ref{eq:ourstartingpt} and let the scalars $b_i, c_i$
be as in
Definition
\ref{def:sharpmp}.
Let the polynomials $p_i$ be as in
Definition
\ref{def:pi}
and let the scalar $\theta_0$ be as in
Definition \ref{def:evseq}.
Then the following (i)--(iii) hold.
\begin{enumerate}
\item $b_i \not=0 \qquad (0 \leq i \leq d-1)$.
\item $c_i \not=0 \qquad (1 \leq i \leq d)$.
\item $b_0b_1 \cdots b_{i-1} = p_i(\theta_0) \qquad (0 \leq i \leq d+1)$.
%\item $x_1x_2 \cdots x_{i} = c_1c_2\cdots c_ip_i(\theta_0) 
%\qquad (0 \leq i \leq d)$.
\end{enumerate}
\end{lemma}
\noindent {\it Proof:}
(i), (ii) Immediate from 
Lemma
\ref{def:bici}(i) and
since each of $x_1, x_2, \ldots, x_d$ is nonzero.
\\
\noindent (iii) 
Assume $0 \leq i \leq d$; otherwise each side
is zero.
Let $\dagger:{\mathcal A}\rightarrow {\mathcal A}$
denote the antiautomorphism which
corresponds to $\Phi$.
Applying 
$\dagger $ to
both sides of
(\ref{eq:pimon})
we get
$E^*_0p_i(A)=
E^*_0A^{i}E^*_i$.
Let $u$ denote a nonzero vector in $E_0V$ and observe
$Au=\theta_0u$.
Recall $E^*_0u, E^*_1u, \ldots, E^*_du$ is a $\Phi$-standard basis
for $V$, and that $A^{\flat}$ represents $A$ with respect to this
basis. From
(\ref{eq:matrepls}) we find $b_0b_1\cdots b_{i-1}$ is
the $0i$ entry of $A^{i\flat}$. Now
\begin{eqnarray*}
b_0b_1\cdots b_{i-1}E^*_0u
&=&
E^*_0A^iE^*_iu 
\\
&=&E^*_0p_i(A)u
\\
&=& p_i(\theta_0)E^*_0u
\end{eqnarray*}
and it follows
$b_0b_1\cdots b_{i-1}
=
p_i(\theta_0)$.
\hfill $\Box $ \\

\begin{theorem}
\label{def:bici3}
Let $\Phi$ denote the Leonard system
from Definition 
\ref{eq:ourstartingpt} and let the polynomials
$p_i$ be as in
Definition
\ref{def:pi}. Let the scalar $\theta_0$ be as in
Definition \ref{def:evseq}.
Then 
$p_i(\theta_0)\not=0$ for 
$0 \leq i \leq d$.
Let the scalars $b_i, c_i$ be
as in
Definition
\ref{def:sharpmp}.
Then 
\begin{eqnarray}
\label{eqbi}
b_i&=& \frac{p_{i+1}(\theta_0)}
{p_i(\theta_0)} 
\qquad  \qquad 
(0 \leq i \leq d)
\end{eqnarray}
and 
\begin{eqnarray}
\label{eqci}
c_i&=& \frac{x_ip_{i-1}(\theta_0)}
{p_i(\theta_0)} 
\qquad  \qquad 
(1 \leq i \leq d).
\end{eqnarray}
\end{theorem}
\noindent {\it Proof:}
Observe $p_i( \theta_0)\not=0$ for $0 \leq i \leq d$
by Lemma
\ref{def:bici2}(i), (iii).
Line 
(\ref{eqbi}) is immediate from
Lemma
\ref{def:bici2}(iii).
To get 
(\ref{eqci}) combine
(\ref{eqbi}) and
Lemma \ref{def:bici}(i).
\hfill $\Box $ \\

%\begin{lemma}
%\label{def:bici4a}
%Let $\Phi$ denote the Leonard system
%from Definition 
%\ref{eq:ourstartingpt} and let the scalars $b_i, c_i$
%be as in Definition
%\ref{def:sharpmp}.
%Then the following (i), (ii) hold.
%\begin{enumerate}
%\item 
%$
%E^*_{i}AE^*_{i+1}E_0 = b_iE^*_iE_0 \qquad 
%(0 \leq i \leq d-1).
%$
%\item
%$
%E^*_{i}AE^*_{i-1}E_0
%=c_iE^*_iE_0
%\qquad 
%(1 \leq i \leq d).
%$
%\end{enumerate}
%\end{lemma}
%\noindent {\it Proof:}
%(i) 
%This is (\ref{eq:xentrypre})
%with
%$X=A$ and $j=i+1$.
%\\
%\noindent (ii) 
%This is
%(\ref{eq:xentrypre})
%with
%$X=A$ and $j=i-1$.
%\hfill $\Box $ \\
%
%
%\begin{theorem}
%\label{def:bici4}
%Let $\Phi$ denote the Leonard system
%from Definition 
%\ref{eq:ourstartingpt} and let the scalars $b_i, c_i$
%be as in Definition
%\ref{def:sharpmp}. Let the scalars $m^*_i$ be as in
%Definition
%\ref{def:mid}.
%Then the following (i), (ii) hold.
%\begin{enumerate}
%\item 
%$
%b_i = m^{*-1}_i tr(E^*_{i}AE^*_{i+1}E_0) \qquad 
%(0 \leq i \leq d-1).
%$
%\item
%$
%c_i = m^{*-1}_i tr(E^*_{i}AE^*_{i-1}E_0) \qquad 
%(1 \leq i \leq d).
%$
%\end{enumerate}
%\end{theorem}
%\noindent {\it Proof:}
%(i) 
%This is
%(\ref{eq:xentry})
%with $X=A$ and $j=i+1$.
%\\
%\noindent (ii)
%This is
%(\ref{eq:xentry})
%with $X=A$ and $j=i-1$.
%\hfill $\Box $ \\

\begin{theorem}
\label{lem:longeq}
Let $\Phi$ denote the Leonard system
from Definition 
\ref{eq:ourstartingpt} and let the scalars $ c_i$ be as in
Definition
\ref{def:sharpmp}.
Let the scalars $\theta_i$
be as in Definition \ref{def:evseq}
and let the scalar $\nu$ be as in
Definition \ref{def:n}.
 Then
\begin{equation}
\label{frame}
(\theta_0-\theta_1)
(\theta_0-\theta_2) \cdots
(\theta_0-\theta_d)
= 
\nu c_1c_2\cdots c_d.
\end{equation}
\end{theorem}
\noindent {\it Proof:}
Let $\delta$ denote the expression on the left-hand side
of 
(\ref{frame}). Setting $i=0$ in
(\ref{eq:primiddef}) we find
$\delta E_0=\prod_{j=1}^d(A-\theta_jI)$.
We multiply both sides of this equation on the left by
$E^*_d$ and on the right by $E^*_0$. We evaluate
the resulting equation using
Lemma
\ref{lem:absval}(i)
to obtain
$\delta E^*_dE_0E^*_0=E^*_dA^dE^*_0$.
We multiply both sides of this equation on the right by $E_0$
and use
Lemma
\ref{lem:threeone}(i) to obtain
\begin{eqnarray}
\label{eq:fourterm}
\delta \nu^{-1}E^*_dE_0=E^*_dA^dE^*_0E_0.
\end{eqnarray}
Let $u$ denote a nonzero vector in $E_0V$
and
observe $E_0u=u$.
Recall $E^*_0u, E^*_1u, \ldots, E^*_du$ is
a $\Phi$-standard basis for $V$, and that
$A^\flat$ represents $A$ with respect to this basis.
From
(\ref{eq:matrepls}) we find
$c_1c_2\cdots c_d$ is the $d0$ entry of $A^{d\flat}$. Now
\begin{eqnarray*}
c_1c_2 \cdots c_d E^*_du
&=&
E^*_dA^dE^*_0u
\\
&=&
E^*_dA^dE^*_0E_0u
\\
&=&
\delta \nu^{-1}E^*_du
\end{eqnarray*}
so $c_1c_2\cdots c_d
=\delta \nu^{-1}$.
The result follows.
\hfill $\Box $ \\

\section{The scalars $k_i$}

\noindent In this section we consider some scalars that are closely
related to the scalars from Definition
\ref{def:mid}.

\begin{definition}
\label{def:ki1}
\rm
Let $\Phi$ denote the Leonard system
from Definition
\ref{eq:ourstartingpt}.
We define
\begin{eqnarray}
\label{eq:ki1}
k_i = m^*_i \nu \qquad   (0 \leq i \leq d),
\end{eqnarray}
where the $m^*_i$ are from
Definition
\ref{def:mid} and $\nu$ is from
Definition
\ref{def:n}. 
\end{definition}

\begin{lemma}
\label{def:ki2}
Let $\Phi$ denote the Leonard system
from
Definition \ref{eq:ourstartingpt} and let the scalars $k_i$
be as in
Definition
\ref{def:ki1}. Then
%the following (i)--(iii) hold.
%\begin{enumerate} 
%\item $k_0=1$.
%\item $k_i\not=0\quad  (0 \leq i \leq d)$.
%\item $\sum_{i=0}^dk_i  = \nu$.
%\end{enumerate}
(i) $k_0=1$;
(ii) $k_i\not=0$ for  $0 \leq i \leq d$;
(iii) $\sum_{i=0}^d k_i  = \nu$.
%\end{enumerate}

\end{lemma}
\noindent {\it Proof:}
(i) Set $i=0$ in 
(\ref{eq:ki1}) and recall $m^*_0=\nu^{-1}$.
\\
(ii)  Applying Lemma
\ref{lem:mid}(iii) to $\Phi^*$ we find
$m^*_i\not=0$ for $0 \leq i \leq d$.
We have $\nu\not=0$ by Definition
\ref{def:n}. 
The result follows in view of
(\ref{eq:ki1}).
\\
\noindent (iii)
  Applying Lemma
\ref{lem:mid}(iv) to $\Phi^*$ we find
$\sum_{i=0}^d m^*_i=1$. The result follows
in view of 
(\ref{eq:ki1}).
\hfill $\Box $ \\

\begin{lemma}
\label{def:ki3}
Let $\Phi$ denote the Leonard system
from
Definition \ref{eq:ourstartingpt} and let the scalars $k_i$
be as in
Definition
\ref{def:ki1}.
Then with reference to Definition
\ref{def:evseq}, Definition 
\ref{def:aixi}, and Definition
\ref{def:pi},
\begin{eqnarray}
\label{eq:kipi}
k_i = \frac
{
p_i(\theta_0)^2
}
{
x_1 x_2 \cdots x_i
}
\qquad (0 \leq i \leq d).
\end{eqnarray}
\end{lemma}
\noindent {\it Proof:}
We show that each side of 
(\ref{eq:kipi}) is equal to
$\nu {\rm tr}(E^*_iE_0)$.
Using 
(\ref{eq:mid}) and
(\ref{eq:ki1}) we find
$\nu {\rm tr}(E^*_iE_0)$ is equal
to the left-hand side of
(\ref{eq:kipi}).
Using Lemma \ref{lem:eispoly} we 
find 
$\nu {\rm tr}(E^*_iE_0)$ is equal
to the right-hand side of
(\ref{eq:kipi}).
%Observe
%\begin{eqnarray*}
%m^*_i &=& \mbox{tr}(E^*_iE_0)
%\\
%&=&  \mbox{tr}\; \frac{p_i(A)E^*_0p_i(A)E_0}{x_1x_2\cdots x_i}
%\qquad \qquad (\mbox{by Theorem 
%\ref{lem:eispoly}})
%\\
%&=&  \mbox{tr}\; \frac{E^*_0p_i(A)E_0p_i(A)}{x_1x_2\cdots x_i}
%\\
%&=&   \frac{p_i(\theta_0)^2}{x_1x_2\cdots x_i} \mbox{tr}(E^*_0E_0)
%\\
%&=&   \frac{p_i(\theta_0)^2 m^*_0}{x_1x_2\cdots x_i}.
%\end{eqnarray*}
%The result follows in view of
%(\ref{eq:ki1}) and since $m^*_0=\nu^{-1}$.
\hfill $\Box $ \\

\begin{theorem}
\label{def:ki4}
Let $\Phi$ denote the Leonard system
from
Definition \ref{eq:ourstartingpt} and let the scalars $k_i$
be as in
Definition
\ref{def:ki1}. Let the scalars $b_i, c_i$ be as in
Definition 
\ref{def:sharpmp}.
Then
\begin{eqnarray}
k_i = \frac
{
b_0b_1\cdots b_{i-1}
}
{
c_1 c_2 \cdots c_i
}
\qquad (0 \leq i \leq d).
\label{eq:kibc}
\end{eqnarray}
\end{theorem}
\noindent {\it Proof:}
Evaluate the expression on the right in
(\ref{eq:kipi}) using
Lemma \ref{def:bici}(i) and 
Lemma \ref{def:bici2}(iii).
\hfill $\Box $ \\

\section{The polynomials  $v_i$}

\noindent Let $\Phi$ denote the Leonard system from Definition
\ref{eq:ourstartingpt} and let the polynomials $p_i$ be
as in
Definition \ref{def:pi}.
The  $p_i$ have two normalizations
 of interest; we call these the
$u_i$ and the $v_i$. In this section we discuss the
$v_i$.
In the next section we will discuss the $u_i$.

\begin{definition}
\label{def:vi1}
\rm
Let $\Phi$ denote the Leonard system
from Definition 
\ref{eq:ourstartingpt} and let the polynomials $p_i$ be as in
 Definition
\ref{def:pi}.
For $0 \leq i \leq d$ we define the polynomial $v_i $ 
by 
\begin{equation}
\label{vicl}
v_i = \frac{p_i}{c_1c_2\cdots c_i},
\end{equation}
where the 
 $c_j$ are
from Definition
\ref{def:sharpmp}.
We observe  $v_0=1$.
\end{definition}

\begin{lemma}
\label{def:vi1a}
Let $\Phi$ denote the Leonard system
from Definition 
\ref{eq:ourstartingpt} and let the polynomials $v_i$ be as in
 Definition 
\ref{def:vi1}. Let the 
scalar $\theta_0$ be as in Definition \ref{def:evseq} and let the
scalars $k_i$ be as in
Definition
\ref{def:ki1}.
Then
\begin{eqnarray}
\label{eq:vnorm}
 v_i(\theta_0)= k_i
\qquad ( 0 \leq i \leq d).
\end{eqnarray}
\end{lemma}
\noindent {\it Proof:}
Use Lemma
\ref{def:bici2}(iii),
Theorem
\ref{def:ki4}, and
(\ref{vicl}).
\hfill $\Box $ \\

\begin{lemma}
\label{def:vi2}
Let $\Phi$ denote the Leonard system
from Definition 
\ref{eq:ourstartingpt} 
and let the polynomials $v_i$ be as in
 Definition 
\ref{def:vi1}. Let the scalars $a_i, b_i, c_i$ be as in
Definition
\ref{def:aixi} and Definition
\ref{def:sharpmp}.
Then
\begin{eqnarray}
\lambda v_i = c_{i+1}v_{i+1}
+ a_i v_i 
+ b_{i-1} v_{i-1}
\qquad (0 \leq i \leq d-1),
\end{eqnarray}
where $b_{-1}=0$ and $v_{-1}=0$. Moreover
\begin{eqnarray}
\lambda v_d 
- a_d v_d 
-b_{d-1}v_{d-1} = (c_1 c_2 \cdots c_d)^{-1}p_{d+1}.
\end{eqnarray}
\end{lemma}
\noindent {\it Proof:}
In 
(\ref{eq:prec}), 
divide both sides  by $c_1c_2 \cdots c_i$.
Evaluate the result using
Lemma \ref{def:bici}(i)
and 
(\ref{vicl}).
\hfill $\Box $ \\

\begin{theorem}
\label{def:vi3}
Let $\Phi$ denote the Leonard system
from Definition
\ref{eq:ourstartingpt}
and let the polynomials $v_i$ be as in 
 Definition 
\ref{def:vi1}.
Let $V$ denote an irreducible $\mathcal A$-module and let
$u$ denote a nonzero vector in $E_0V$.
Then
\begin{equation}
\label{eq:viau}
v_i(A)E^*_0u = E^*_iu \qquad \qquad (0 \leq i \leq d).
\end{equation}
\end{theorem}
\noindent {\it Proof:}
For $0 \leq i \leq d$ we define
$w_i=v_i(A)E^*_0u$ 
and $w'_i=E^*_iu$.
We show $w_i=w'_i$.
Each of $w_0$,  $w'_0$ is equal to
$E^*_0u$ so $w_0=w'_0$.
Using
Lemma \ref{def:vi2}
we obtain
\begin{eqnarray}
Aw_i = c_{i+1}w_{i+1} + a_iw_i + b_{i-1}w_{i-1}
\qquad (0 \leq i \leq d-1)
\label{eq:releta}
\end{eqnarray}
where  $w_{-1}=0$ and $b_{-1}=0$.
By
Definition 
\ref{def:flatcon},
Definition
\ref{def:sharpmp}, and since
$w'_0, w'_1, \ldots, w'_d$
is a $\Phi$-standard basis,
\begin{eqnarray}
Aw'_i = c_{i+1}w'_{i+1} + a_iw'_i + b_{i-1}w'_{i-1}
\qquad (0 \leq i \leq d-1)
\label{eq:rel2eta}
\end{eqnarray}
where  $w'_{-1}=0$.
Comparing
(\ref{eq:releta}), 
(\ref{eq:rel2eta}) and using $w_0=w'_0$
we find
$w_i=w'_i$ for $0 \leq i \leq d$.
The result follows.
\hfill $\Box $ \\

%\noindent We finish this section with a comment.
%
%\begin{lemma}
%\label{def:vi4}
%Let $\Phi$ denote the Leonard system
%from Definition
%\ref{eq:ourstartingpt}
%and let the polynomials $v_i$ be as in 
% Definition 
%\ref{def:vi1}. Let the scalar $\nu$ be as in 
%Definition \ref{def:n}.
%Then the following (i), (ii) hold.
%\begin{enumerate}
%\item
%$v_i(A)E^*_0E_0 = E^*_iE_0 \qquad (0 \leq i \leq d).$
%\item
%$v_i(A)E^*_0 = \nu E^*_iE_0E^*_0 \qquad (0 \leq i \leq d).$
%\end{enumerate}
%\end{lemma}
%\noindent {\it Proof:}
%(i) Let the integer $i$ be given and abbreviate
%$
%\Delta=
%v_i(A)E^*_0-E^*_i
%$.
%We show $\Delta E_0=0$.
%In order to to do this we show  
%$\Delta E_0V=0$, where $V$ denotes an irreducible 
%$\mathcal A$-module.
%Let $u$ denote a nonzero vector in $E_0V$ and 
%recall
%$u$ spans $E_0V$.
%Observe
%$\Delta u=0$ by 
%Theorem \ref{def:vi3} so
%$\Delta E_0V=0$. Now $\Delta E_0=0$ so
%$v_i(A)E^*_0E_0 = E^*_iE_0$.
%\\
%\noindent (ii) In the equation of (i) above, multiply  both sides
%on the right by $E^*_0$ and simplify the result using
%Lemma \ref{lem:threeone}(ii).
%\hfill $\Box $ \\

\section{The polynomials $u_i$}

\noindent Let $\Phi$ denote the Leonard system from Definition
\ref{eq:ourstartingpt} and let the polynomials $p_i$ be
as in
Definition \ref{def:pi}. In the previous section we gave
a 
 normalization of the $p_i$ that we called the $v_i$.
 In this section we give a normalization for
 the $p_i$ that we call the $u_i$.

\begin{definition}
\label{def:ui1}
\rm
Let $\Phi$ denote the Leonard system
from 
Definition 
\ref{eq:ourstartingpt} and let the polynomials $p_i$ be
as in Definition
\ref{def:pi}.
For $0 \leq i \leq d$ we define 
the polynomial $u_i$ by
\begin{equation}
\label{uicl}
u_i = \frac{p_i}{p_i(\theta_0)},
\end{equation}
where $\theta_0$ is from
Definition
\ref{def:evseq}.
We observe $u_0=1$. Moreover 
\begin{eqnarray}
\label{eq:unorm}
u_i(\theta_0)= 1 
\qquad ( 0 \leq i \leq d).
\end{eqnarray}
\end{definition}

\begin{lemma}
\label{def:ui2}
Let $\Phi$ denote the Leonard system
from 
Definition 
\ref{eq:ourstartingpt}
and let the polynomials
$u_i$ be as in
Definition 
\ref{def:ui1}.
Let the scalars $a_i, b_i, c_i$ be as in
Definition
\ref{def:aixi} and Definition
\ref{def:sharpmp}.
Then  
\begin{eqnarray}
\lambda u_i = b_iu_{i+1}
+ a_i u_i 
+ c_i u_{i-1}
\qquad (0 \leq i \leq d-1),
\label{eq:uirec1}
\end{eqnarray}
where $u_{-1}=0$.
Moreover
%Let the scalar $\theta_0$ be as in
%Definition \ref{def:evseq}.
%Then
\begin{eqnarray}
\lambda u_d 
- c_d u_{d-1} 
-a_du_d = p_d(\theta_0)^{-1}p_{d+1},
\label{eq:uirec2}
\end{eqnarray}
where $\theta_0$ is from 
Definition \ref{def:evseq}.
%Then
\end{lemma}
\noindent {\it Proof:}
In 
(\ref{eq:prec}), 
divide both sides  by $p_i(\theta_0)$
and evaluate the result using
Lemma \ref{def:bici}(i),
Theorem
\ref{def:bici3},
%(\ref{eqbi}),
and (\ref{uicl}).
\hfill $\Box $ \\

\noindent 
The $u_i$ and $v_i$ are related as follows.

\begin{lemma}
\label{def:ui4}
Let $\Phi$ denote the Leonard system
from 
Definition 
\ref{eq:ourstartingpt}.
Let the polynomials $u_i, v_i$ be as in
Definition
\ref{def:ui1} and 
Definition
\ref{def:vi1} respectively.
Then 
\begin{equation}
\label{eq:uvk}
v_i = k_i u_i  \qquad (0 \leq i \leq d),
\end{equation}
where the $k_i$ are from
Definition \ref{def:ki1}.
\end{lemma}
\noindent {\it Proof:}
Compare 
(\ref{vicl}) and
(\ref{uicl}) in light of
Lemma \ref{def:bici2}(iii) and
Theorem
\ref{def:ki4}.
\hfill $\Box $ \\

\section{A bilinear form}

\medskip
\noindent
In this section we associate with each Leonard pair a certain
bilinear form. To prepare for this   
we recall a few concepts from linear
algebra.

\medskip
\noindent 
Let $V$ denote a finite dimensional
vector space over $\K$.
By a {\it bilinear form on $V$} we mean
a  map $\langle \,,\,\rangle : V \times V \rightarrow \K$
that satisfies the following four conditions for all
$u,v,w \in V$ and for all $\alpha \in \K$:
(i) $\langle u+v,w \rangle = 
 \langle u,w \rangle 
+
\langle v,w \rangle$;
(ii)
$
\langle \alpha u,v \rangle 
=
\alpha \langle u,v \rangle 
$;
(iii)
 $\langle u,v+w \rangle = 
 \langle u,v \rangle 
+
\langle u,w \rangle$;
(iv)
$
\langle u, \alpha v \rangle 
=
\alpha \langle u,v \rangle 
$.
We observe 
that a scalar multiple of 
a bilinear form on $V$ is a bilinear form on $V$.
Let
$
\langle\, ,\, \rangle 
$
denote a bilinear form on $V$.
This form is said to be {\it symmetric}
whenever
$
\langle u,v \rangle 
= 
\langle v,u \rangle 
$
for all $u, v \in V$.
Let 
$
\langle\, ,\, \rangle 
$
denote a bilinear form on $V$.
Then the following
are equivalent: (i) there exists a nonzero $u \in V$ such
that $
\langle u,v \rangle 
= 0 
$
for all $v \in V$;
 (ii) there exists a nonzero $v \in V$ such
that $
\langle u,v \rangle 
= 0 
$
for all $u \in V$.
The form 
$\langle \,,\,\rangle $
is 
said to be {\it degenerate }
whenever (i), (ii) hold and {\it nondegenerate}
otherwise.
Let $\gamma :{\mathcal A} \rightarrow {\mathcal A}$
denote an antiautomorphism
and let $V$ denote an irreducible $\mathcal A$-module.
Then there exists a nonzero bilinear form 
$\langle \,,\,\rangle $ on $V$ such that
$\langle Xu,v\rangle 
=
\langle u,X^\gamma v\rangle 
$
for all $u,v \in V$ and for all $X \in \mathcal A$.
The form is unique
up to multiplication by a nonzero scalar in $\K$.
The form in nondegenerate.
We refer to this form as the {\it bilinear form on $V$ associated
with $\gamma $}. This form is not symmetric in general.

\medskip
\noindent We now return our attention to Leonard pairs.

\begin{definition}
\label{def:ip}
\rm
Let
$\Phi=(A;A^*; \lbrace E_i\rbrace_{i=0}^d; $ $ 
\lbrace E^*_i\rbrace_{i=0}^d)$
denote a Leonard system in
$\mathcal A$.
Let
$\dagger :{\mathcal A}\rightarrow {\mathcal A} $ 
denote the corresponding antiautomorphism 
from Definition
\ref{def:dag}.
Let $V$ denote an 
irreducible $\cal A$-module.
For the rest of this paper  we let  
$\langle \,,\,\rangle $ denote the bilinear
form on $V$ 
associated with $\dagger$. 
We abbreviate $\Vert u \Vert^2 = \langle u,u\rangle $
for all $u \in V$.
By the construction,
 for $X \in {\mathcal A}$  we have
\begin{eqnarray}
\langle Xu, v \rangle = \langle u, X^\dagger v \rangle
\qquad (\forall u \in V, \forall v \in V).
\label{eq:daggerstart}
\end{eqnarray}
\end{definition}

\noindent We make an observation.

\begin{lemma}
With reference to Definition \ref{def:ip},
let $\mathcal D$ (resp. $\mathcal D^*$) denote
the subalgebra of $\mathcal A$ generated  by $A$ (resp. $A^*$).
Then for $X \in {\mathcal D} \cup {\mathcal D}^*$ we have
\begin{eqnarray}
\label{eq:bilinv}
\langle Xu, v \rangle = \langle u, Xv \rangle
\qquad (\forall u \in V, \forall v \in V).
\end{eqnarray}
\end{lemma}
\noindent {\it Proof:}
Combine  (\ref{eq:daggerstart})
and
Lemma \ref{cor:eistab}.
\hfill $\Box $ \\

%given $u,v \in V$ we say $u,v$ are {\it orthogonal}
%with respect to 
%$\langle \,,\,\rangle $ whenever 
%$\langle u,v\rangle = 0 $.
%We abbreviate $\Vert u \Vert^2 = \langle u,u\rangle $
%for all $u \in V$.

\begin{theorem}
\label{thm:orthog}
With reference to Definition
\ref{def:ip},
let $u$ denote a nonzero vector in $E_0V$ and
recall 
%\begin{eqnarray*}
%\label{eq:orthogset}
$E^*_0u, E^*_1u,\ldots, E^*_du $
%$\end{eqnarray*}
is a $\Phi$-standard basis for $V$.
%With respect to 
%$\langle \,,\,\rangle $  the vectors 
% (\ref{eq:orthogset})
%are mutually orthogonal. Moreover 
We have
\begin{eqnarray}
\langle E^*_iu, E^*_ju \rangle = \delta_{ij}k_i \nu^{-1}   
\Vert u \Vert^2  
\qquad \qquad 
(0 \leq i,j \leq d),
\label{eq:vert}
\end{eqnarray}
where the $k_i$ are from
Definition \ref{def:ki1}
and $\nu$ is from
Definition \ref{def:n}.
\end{theorem}
\noindent {\it Proof:}
By 
(\ref{eq:bilinv}) and since $E_0u=u$ we find
$\langle E^*_iu, E^*_ju\rangle
=
\langle u, E_0E^*_iE^*_jE_0u\rangle$.
Using
Lemma \ref{lem:mid}(ii) and
(\ref{eq:ki1}) we find
$\langle u, E_0E^*_iE^*_jE_0u\rangle=
\delta_{ij}k_i \nu^{-1}   
\Vert u \Vert^2 $. 
%For $0 \leq i,j\leq d$,
%\begin{eqnarray*}
%\langle E^*_iu, E^*_ju\rangle
%&=&
%\langle E^*_iE_0u, E^*_jE_0u\rangle \\
%&=&
%\langle u, (E^*_iE_0)^\dagger E^*_jE_0u\rangle \\
%&=&
%\langle u, E_0E^*_iE^*_jE_0u\rangle \\
%&=&
%\delta_{ij}\langle u, E_0E^*_iE_0u\rangle \\
%&=&
%\delta_{ij}m^*_i\langle u, E_0u\rangle  \\
%&=&
%\delta_{ij}m^*_i\langle u, u\rangle.
%\end{eqnarray*}
%The result follows from this and 
%Definition \ref{def:ki1}.
\hfill $\Box $ \\

\begin{corollary}
\label{cor:sym}
With reference to
Definition
\ref{def:ip}, the bilinear form 
 $\langle \,, \,\rangle $ is symmetric.
\end{corollary}
\noindent {\it Proof:}
Let $u$ denote a nonzero vector in $E_0V$ and abbreviate
$v_i=E^*_iu$ for $0 \leq i \leq d$.
From Theorem
\ref{thm:orthog} we find
$\langle v_i,v_j\rangle= \langle v_j,v_i\rangle$
for $0 \leq i,j\leq d$. The result follows since
$v_0, v_1,\ldots, v_d$ is a basis for $V$.
\hfill $\Box $ \\

\noindent We have a comment.

\begin{lemma}
\label{lem:anglebasic}
With reference to Definition
\ref{def:ip}, let $u$ denote a nonzero vector in 
$E_0V$ and let  $v$ denote a nonzero vector in $E^*_0V$. Then 
the following (i)--(iv) hold.
\begin{enumerate}
\item Each of
 $\Vert u \Vert^2, 
\Vert v \Vert^2,
\langle u,v\rangle $
is nonzero.
\item
$E^*_0u = 
\langle u,v\rangle 
\Vert v\Vert^{-2} 
 v.$
\item
$
E_0v = 
\langle u,v\rangle 
\Vert u\Vert^{-2} u.$ 
\item $\nu 
\langle u,v\rangle^2 =
\Vert u\Vert^2
\Vert v\Vert^2
$.
\end{enumerate}
\end{lemma}
\noindent {\it Proof:}
(i)
Observe $\Vert u \Vert^2\not=0$ by
Theorem \ref{thm:orthog}
and since $\langle \,,\,\rangle $ is not 0.
Similarly
$\Vert v \Vert^2\not=0$.
To see that $\langle u,v\rangle \not=0$, observe that
$v$ is a basis for $E^*_0V$ so there exists $\alpha \in \K$ 
such that
$E^*_0u = \alpha v$.
Recall $E^*_0u\not=0$ by
Lemma
\ref{lem:stbasis} so $\alpha\not=0$.
Using (\ref{eq:bilinv}) and $E^*_0v=v$ we routinely find
$\langle u,v \rangle =
\alpha \Vert v\Vert^2$ and it follows
$\langle u,v \rangle \not=0$. 
%We now see
%\begin{eqnarray*}
%\langle u,v \rangle &=&
%\langle u,E^*_0v \rangle
%\\
%&=&
%\langle E^*_0u,v \rangle
%\\
%&=&
%\alpha \Vert v\Vert^2
%\\
%&\not=& 0.
%\end{eqnarray*}
\\
(ii) In the proof of part (i) we 
found $E^*_0u=\alpha v$ where $\langle u,v\rangle =\alpha \Vert v\Vert^2$.
The result follows.
\\
(iii) Similar to the proof of (ii) above.
\\
(iv)
Using $u=E_0u$ and $\nu E_0E^*_0E_0=E_0$ we find  
$\nu^{-1} u = E_0E^*_0u$.
To finish the proof,
evaluate $E_0E^*_0u$ using 
(ii) above and then (iii) above.
%Using $u=E_0u$ and $\nu E_0E^*_0E_0=E_0$ we find  
%\begin{eqnarray*}
%\nu^{-1} u &=& E_0E^*_0E_0u
%\\
%&=& E_0E^*_0u
%\\
%&=& \frac{\langle u,v\rangle}{\Vert v \Vert^2}E_0v
%\\
%&=& \frac{\langle u,v\rangle^2}{\Vert u \Vert^2 \Vert v \Vert^2}u.
%\end{eqnarray*}
%The result follows.
\hfill $\Box $ \\

\section{Askey-Wilson duality}

In this section we show the polynomials $u_i, v_i, p_i$
satisfy a relation known as {\it Askey-Wilson duality}.
We begin with a lemma.

\begin{lemma}
\label{lem:predu}
With reference to Definition
\ref{def:ip}, let $u$ denote a nonzero vector in 
$E_0V$ and let  $v$ denote a nonzero vector in $E^*_0V$. Then 
\begin{equation}
\langle E^*_iu, E_jv\rangle = 
\nu^{-1}k_i k^*_j u_i(\theta_j)  \langle u,v \rangle 
\qquad  \qquad 
 (0 \leq i,j\leq d).
\label{eq:predu}
\end{equation}
\end{lemma}
\noindent {\it Proof:}
Using Theorem
\ref{def:vi3} we find
\begin{eqnarray}
\langle E^*_iu, E_jv\rangle &=&
\langle v_i(A)E^*_0u, E_jv\rangle 
\nonumber
\\
&=&
\langle E^*_0u,v_i(A)E_jv\rangle 
\nonumber
\\
&=&
v_i(\theta_j)\langle E^*_0u,E_jv\rangle 
\nonumber
\\
&=&
v_i(\theta_j)\langle E^*_0u,v^*_j(A^*)E_0v\rangle 
\nonumber
\\
&=&
v_i(\theta_j)\langle v^*_j(A^*)E^*_0u,E_0v\rangle 
\nonumber
\\
&=&
v_i(\theta_j)v^*_j(\theta^*_0)\langle E^*_0u,E_0v\rangle.
\label{eq:almostdone}
\end{eqnarray}
Using Lemma
\ref{lem:anglebasic}(ii)--(iv) we find
$\langle E^*_0u,E_0v\rangle =\nu^{-1}\langle u,v\rangle$.
Observe
$v_i(\theta_j)=u_i(\theta_j)k_i$
by 
(\ref{eq:uvk}). Applying
Lemma \ref{def:vi1a} to $\Phi^*$ we find
$v^*_j(\theta^*_0)=k^*_j$. 
Evaluating (\ref{eq:almostdone}) using
these comments we obtain
(\ref{eq:predu}).
\hfill $\Box $ \\

\begin{theorem}
\label{lem:uidual}
Let $\Phi$ denote the Leonard system
from 
Definition 
\ref{eq:ourstartingpt}. Let the polynomials
$u_i$ be as in Definition
\ref{def:ui1} and recall the $u^*_i$ are the 
corresponding polynomials for $\Phi^*$.
Let the scalars $\theta_i, \theta^*_i$ be as in
Definition \ref{def:evseq}.
Then
\begin{eqnarray}
u_i(\theta_j) = u^*_j(\theta^*_i) 
\qquad (0 \leq i,j\leq d).
\label{eq:uijus}
\end{eqnarray}
\end{theorem}
\noindent {\it Proof:}
Applying Lemma \ref{lem:predu} to $\Phi^*$
we find
\begin{equation}
\langle E_jv, E^*_iu\rangle = 
\nu^{-1}k^*_j k_i u^*_j(\theta^*_i)  \langle u,v \rangle 
\qquad  \qquad 
 (0 \leq i,j\leq d).
\label{eq:predudual}
\end{equation}
To finish the proof, compare
(\ref{eq:predu}), 
(\ref{eq:predudual}), and recall $\langle \,,\rangle$
is symmetric.
\hfill $\Box $ \\

\noindent In the following two theorems we show how 
(\ref{eq:uijus}) looks in terms of the polynomials
$v_i$ and $p_i$.

\begin{theorem}
\label{lem:vidual}
Let $\Phi$ denote the Leonard system
from 
Definition 
\ref{eq:ourstartingpt}.
With reference to 
Definition \ref{def:notconv},
 Definition
\ref{def:evseq},
and Definition
\ref{def:vi1},
\begin{eqnarray}
v_i(\theta_j)/k_i = v^*_j(\theta^*_i)/k^*_j 
\qquad (0 \leq i,j\leq d).
\label{eq:vijus}
\end{eqnarray}
\end{theorem}
\noindent {\it Proof:}
Evaluate (\ref{eq:uijus})
using
Lemma
\ref{def:ui4}.
\hfill $\Box $ \\

\begin{theorem}
\label{lem:pidual}
Let $\Phi$ denote the Leonard system
from 
Definition 
\ref{eq:ourstartingpt}.
With reference to 
Definition \ref{def:notconv},
 Definition
\ref{def:evseq},
and Definition
\ref{def:pi},
\begin{eqnarray}
\frac{p_i(\theta_j)}{p_i(\theta_0)}
= 
\frac{p^*_j(\theta^*_i)}{p^*_j(\theta^*_0)}
\qquad (0 \leq i,j\leq d).
\label{eq:pijus}
\end{eqnarray}
\end{theorem}
\noindent {\it Proof:}
Evaluate (\ref{eq:uijus}) 
using
Definition \ref{def:ui1}.
\hfill $\Box $ \\

\noindent The equations
(\ref{eq:uijus}),
(\ref{eq:vijus}),
(\ref{eq:pijus})
are often referred to as {\it Askey-Wilson duality}.

\section{The three-term recurrence and the difference equation}

In Lemma
\ref{def:ui2} we gave a three-term recurrence
for the polynomials $u_i$. This recurrence
is often expressed as follows.

\begin{theorem}
\label{cor:3term}
Let $\Phi$ denote the Leonard system
from 
Definition 
\ref{eq:ourstartingpt}
and let the polynomials
$u_i$ be as in
Definition 
\ref{def:ui1}. Let the scalars $\theta_i$
be as in Definition \ref{def:evseq}.
Then  for $0 \leq i,j\leq d$ we have 
\begin{eqnarray}
\theta_j u_i(\theta_j) = b_iu_{i+1}(\theta_j)
+ a_i u_i(\theta_j)
+ c_i u_{i-1}(\theta_j),
\label{eq:3termu}
\end{eqnarray}
where $u_{-1}=0$ and $u_{d+1}=0$.
\end{theorem}
\noindent {\it Proof:}
Apply Lemma
\ref{def:ui2}
(with $\lambda = \theta_j$)
and observe $p_{d+1}(\theta_j)=0$ by Theorem
\ref{thm:pseq}(ii).
\hfill $\Box $ \\

%The above equation is called the {\it three-term recurrence}
%for the polynomials $u_i$.
\noindent
Applying  
Theorem \ref{cor:3term}
to $\Phi^*$ and using
Theorem
\ref{lem:uidual} we routinely obtain the following.

\begin{theorem}
\label{lem:diffeq}
Let $\Phi$ denote the Leonard system
from 
Definition 
\ref{eq:ourstartingpt} and
let the polynomials $u_i$ be as in
 Definition 
\ref{def:ui1}.
Then for $0 \leq i,j\leq d$ we have
\begin{eqnarray}
\label{eq:diffeq}
\theta^*_i u_i(\theta_j) = 
b^*_ju_i(\theta_{j+1})
+
a^*_ju_i(\theta_j)
+
c^*_ju_i(\theta_{j-1}),
\end{eqnarray}
where $\theta_{-1}$, $\theta_{d+1}$ denote indeterminates.
\end{theorem}
%\noindent {\it Proof:}
%Apply 
%Corollary \ref{cor:3term}
%to $\Phi^*$ and evaluate the result using
%Theorem
%\ref{lem:uidual}.
%\hfill $\Box $ \\

\noindent We refer to 
(\ref{eq:diffeq}) as the {\it difference equation}
satisfied by the $u_i$.

\section{The orthogonality relations}

\noindent 
In this section we show that
each of the polynomials $p_i, u_i, v_i$
satisfy an orthogonality relation.
We begin with a lemma.

\begin{lemma}
\label{thm:Ptrans}
With reference to Definition
\ref{def:ip}, let $u$ denote a nonzero vector in 
$E_0V$ and let  $v$ denote a nonzero vector in $E^*_0V$.
Then for $0 \leq i \leq d$,  both 
\begin{eqnarray}
E^*_i u &=& 
\frac
{
\langle u,v 
\rangle
}
{
\Vert v 
\Vert^2
}
\sum_{j=0}^d v_i(\theta_j) E_jv,
%\qquad (0 \leq i \leq d),
\label{eq:trans1}
\\
E_i v &=& 
\frac
{
\langle u,v 
\rangle
}
{
\Vert u
\Vert^2
}
\sum_{j=0}^d v^*_i(\theta^*_j) E^*_ju.
%\qquad (0 \leq i \leq d).
\label{eq:trans2}
\end{eqnarray}
\end{lemma}
\noindent {\it Proof:}
We first show 
(\ref{eq:trans1}).
To do this we show
each side of 
(\ref{eq:trans1}) is equal to
$v_i(A)E^*_0u$.
By
Theorem \ref{def:vi3} we find
$v_i(A)E^*_0u$ 
is equal to the left-hand side of 
(\ref{eq:trans1}).
To see that 
$v_i(A)E^*_0u$ is equal to the right-hand side
of 
(\ref{eq:trans1}), multiply
$v_i(A)E^*_0u$ on the left by the identity $I$,
expand using
$I=\sum_{j=0}^d E_j$, and simplify the result using
$E_jA=\theta_jE_j$ $(0 \leq j\leq d)$ and
Lemma \ref{lem:anglebasic}(ii).
We have now proved 
(\ref{eq:trans1}).
Applying 
(\ref{eq:trans1}) to $\Phi^*$ we obtain
(\ref{eq:trans2}).
\hfill $\Box $ \\

\noindent We now display the othogonality relations
for the polynomials $v_i$.

\begin{theorem}
\label{thm:viorth}
Let $\Phi$ denote the Leonard system from Definition  
\ref{eq:ourstartingpt} and let
 the polynomials $v_i$ be as in
Definition
\ref{def:vi1}.
Then both 
\begin{eqnarray}
\sum_{r=0}^d v_i(\theta_r)v_j(\theta_r)k^*_r &=&
\delta_{ij} \nu k_i \qquad \qquad (0 \leq i,j\leq d),
\label{eq:vio1}
\\
\sum_{i=0}^d v_i(\theta_r)v_i(\theta_s)k^{-1}_i &=&
\delta_{rs} \nu k_r^{*-1}   \qquad \qquad (0 \leq r,s\leq d).
\label{eq:vio2}
\end{eqnarray}
\end{theorem}
\noindent {\it Proof:} 
Let $V$ denote an irreducible $\cal A$-module and 
let $v$ denote a nonzero vector in $E^*_0V$. 
Applying Theorem
\ref{thm:orthog} to $\Phi^*$ we find
$\langle E_rv, E_sv\rangle = \delta_{rs} k^*_r \nu^{-1} \Vert v\Vert^2
$ for $0 \leq r,s\leq d$.
To obtain
(\ref{eq:vio1}),
 in equation
(\ref{eq:vert}), eliminate 
each of $E^*_iu$, $E^*_ju$ using
(\ref{eq:trans1}), and simplify the result using
our preliminary comment and 
Lemma \ref{lem:anglebasic}.
To obtain
(\ref{eq:vio2}), apply
(\ref{eq:vio1}) to $\Phi^*$ and
use Askey-Wilson duality.
\hfill $\Box $ \\

\noindent 
We now 
turn to the polynomials $u_i$.

\begin{theorem}
\label{thm:uiorth}
Let $\Phi$ denote the Leonard system from 
Definition \ref{eq:ourstartingpt} and
let the polynomials $u_i$ be as in
Definition
\ref{def:ui1}.
Then both 
\begin{eqnarray*}
\sum_{r=0}^d u_i(\theta_r)u_j(\theta_r)k^*_r &=&
\delta_{ij}\nu k^{-1}_i \qquad \qquad (0 \leq i,j\leq d),
\\
\sum_{i=0}^d u_i(\theta_r)u_i(\theta_s)k_i &=&
\delta_{rs} \nu k_r^{*-1}   \qquad \qquad (0 \leq r,s\leq d).
\end{eqnarray*}
\end{theorem}
\noindent {\it Proof:}
Evaluate each of 
(\ref{eq:vio1}), 
(\ref{eq:vio2}) using
Lemma
\ref{def:ui4}.
\hfill $\Box $ \\

\noindent 
We now turn to the polynomials
$p_i$.

\begin{theorem}
\label{thm:piorth}
Let $\Phi$ denote the Leonard system from Definition  
\ref{eq:ourstartingpt}
and let the polynomials $p_i$ be as in
Definition
\ref{def:pi}.
Then both
\begin{eqnarray*}
 \sum_{r=0}^d p_i(\theta_r)p_j(\theta_r)m_r &=& \delta_{ij}
x_1x_2\cdots x_i \qquad \qquad (0 \leq i,j\leq d),
\\
\sum_{i=0}^d {{ p_i(\theta_r) p_i(\theta_s)}\over {x_1x_2\cdots x_i}}
&=& \delta_{rs}m_r^{-1} \qquad \qquad (0 \leq r,s\leq d).
\end{eqnarray*}
\end{theorem}
\noindent {\it Proof:}
Applying Definition
\ref{def:ki1} to $\Phi^*$ we find
$k^*_r=m_r\nu$ for $0 \leq r \leq d$.
Evaluate 
each of
(\ref{eq:vio1}), 
(\ref{eq:vio2}) using
this and 
Definition \ref{def:vi1},
Lemma \ref{def:bici}(i),
(\ref{eq:kibc}).
\hfill $\Box $ \\

\section{The matrix $P$}

In this section we express Lemma
\ref{thm:Ptrans} in matrix form and consider the consequences.

\begin{definition}
\label{def:pmat}
\rm
Let $\Phi$ denote the Leonard system
from
Definition \ref{eq:ourstartingpt}.
We define a matrix $P \in
\mbox{Mat}_{d+1}(\K)$  as follows.
For $0 \leq i,j\leq d$ the  entry
$P_{ij}=v_j(\theta_i)$, where
$\theta_i$ is from
Definition \ref{def:evseq} and
 $v_j$ is from
 Definition
 \ref{def:vi1}.
 \end{definition}
 
 \begin{theorem}
 \label{def:pmatpmats}
 Let $\Phi$ denote the Leonard system
 from
 Definition \ref{eq:ourstartingpt}.
 Let the matrix $P$ be as in
 Definition
 \ref{def:pmat} and recall $P^*$ is the corresponding
 matrix for $\Phi^*$.
 Then
 $P^*P=\nu I$,
 where $\nu$ is from
 Definition \ref{def:n}.
 \end{theorem}
 \noindent {\it Proof:}
 Compare
 (\ref{eq:trans1}), (\ref{eq:trans2}) and use
 Lemma
 \ref{lem:anglebasic}(iv).
 \hfill $\Box $ \\
 
 \begin{theorem}
 \label{lem:flatvssharp}
 Let $\Phi$ denote the Leonard system
 from
 Definition \ref{eq:ourstartingpt} and
 let the matrix $P$ be as in
 Definition
 \ref{def:pmat}.
Let  
 the  map $\flat : {\mathcal A}\rightarrow
 \mbox{Mat}_{d+1}(\K)$ be as in
 Definition
 \ref{def:flatcon}
 and let
  $\sharp : {\mathcal A}\rightarrow
  \mbox{Mat}_{d+1}(\K)$ denote the corresponding
  map for $\Phi^*$.
  Then for all $X \in {\mathcal A}$ we have
  \begin{eqnarray}
  X^\sharp P = P X^\flat.
  \end{eqnarray}
  \end{theorem}
  \noindent {\it Proof:}
  Let $V$ denote an irreducible $\mathcal A$-module.
  Let $u$ denote a nonzero vector in $E_0V$
  and recall $E^*_0u, E^*_1u, \ldots, E^*_du$ is
  a $\Phi$-standard basis for $V$.
  By Definition
  \ref{def:flatcon}, $X^\flat$ is the matrix in
  $\mbox{Mat}_{d+1}(\K)$
  that represents
  $X$ with respect to
  $E^*_0u, E^*_1u, \ldots, E^*_du$.
  Similarly for a nonzero $v \in E^*_0V$,
  $X^\sharp$ is the matrix in
  $\mbox{Mat}_{d+1}(\K)$
  that represents $X$ with respect to
  $E_0v, E_1v, \ldots, E_dv$.
  In view of (\ref{eq:trans1}),
  the
  transition matrix from
  $E_0v, E_1v, \ldots, E_dv$
  to
  $E^*_0u, E^*_1u, \ldots, E^*_du$ is
  a scalar multiple of $P$.
  The result follows
  from these comments and elementary linear
  algebra.
  \hfill $\Box $ \\

\section{The split decomposition}

\begin{notation}
\label{def:splitconsetupS99}
\rm
Throughout this section we 
let 
$\Phi=(A;A^*; \lbrace E_i\rbrace_{i=0}^d;  
\lbrace E^*_i\rbrace_{i=0}^d)$
denote a Leonard system in $\mathcal A$,
with eigenvalue sequence $\theta_0, \theta_1,\ldots, \theta_d$
and dual eigenvalue sequence 
 $\theta^*_0, \theta^*_1,\ldots, \theta^*_d$.
We let $V$ denote
an irreducible $\cal A$-module.
\end{notation}
\noindent With reference to Notation
\ref{def:splitconsetupS99}, 
by a {\it decomposition} of $V$ we mean
a sequence
$U_0, U_1, \ldots, U_d$ consisting of 1-dimensional
subspaces of $V$ such that 
\begin{eqnarray*}
V=U_0+U_1+\cdots + U_d \qquad \qquad \mbox{(direct sum)}.
\end{eqnarray*}
In this section we are concerned with the following type of 
decomposition.
\begin{definition}
\label{def:splitdec}
\rm
With reference to Notation
\ref{def:splitconsetupS99}, 
a decomposition
$U_0, U_1, \ldots, U_d$ of $V$
is said to be {\it $\Phi$-split} whenever
both
\begin{eqnarray}
(A-\theta_iI)U_i\subseteq U_{i+1} \quad (0 \leq i \leq d-1),
\qquad\quad 
(A-\theta_dI)U_d=0,
\label{eq:splitpart1}
\\
(A^*-\theta^*_iI)U_i\subseteq U_{i-1} \quad (1 \leq i \leq d),
\qquad\quad 
(A^*-\theta^*_0I)U_0=0.
\label{eq:splitpart2}
\end{eqnarray}
\end{definition}
Our goal in this section is to show there exists a unique
$\Phi$-split decomposition of $V$.
The following definition will be useful.

\begin{definition}
\label{def:VijS99}
\rm
With reference to 
Notation
\ref{def:splitconsetupS99}, 
we set 
\begin{equation}
V_{ij} =  \Biggl(\sum_{h=0}^i E^*_hV\Biggr)\cap
 \Biggl(\sum_{k=j}^d E_kV\Biggr)
\label{eq:defofvijS99}
\end{equation}
for all integers $i,j$. We interpret the sum on the left in
(\ref{eq:defofvijS99}) to be 0 (resp. $V$) if $i<0$  (resp. $i>d$).
Similarily, we interpret the sum on the right in
(\ref{eq:defofvijS99}) to be V (resp. $0$) if $j<0$  (resp. $j>d$).
\end{definition}

\begin{lemma}
\label{lem:thevijbasicfactsS99}
With reference to 
Notation
\ref{def:splitconsetupS99} and 
Definition \ref{def:VijS99}, we have
\begin{enumerate}
\item
$V_{i0} = 
 E^*_0V+E^*_1V+\cdots+E^*_iV \qquad \qquad  (0 \leq i \leq d),$
\item
$V_{dj} = 
 E_jV+E_{j+1}V+\cdots+E_dV \qquad \qquad (0 \leq j\leq d).$
\end{enumerate}

\end{lemma}

\noindent {\it Proof:}
To get (i), set $j=0$ in 
(\ref{eq:defofvijS99}), and apply 
(\ref{eq:VdecompS99}).
Line  (ii) is similarily obtained.
\hfill $\Box $ \\

\begin{lemma}
\label{lem:howaactsonvijS99}
With reference to 
Notation
\ref{def:splitconsetupS99} 
and Definition \ref{def:VijS99}, the following (i)--(iv) hold
for $0 \leq i, j\leq d$.
\begin{enumerate}
\item
$(A-\theta_jI)V_{ij} \subseteq V_{i+1,j+1}$,
\item
$AV_{ij} \subseteq V_{ij} + V_{i+1,j+1}$, 
\item
$(A^*-\theta^*_iI)V_{ij} \subseteq V_{i-1,j-1}$,
\item
$A^*V_{ij} \subseteq V_{ij} + V_{i-1,j-1}$.
\end{enumerate}
\end{lemma}

\noindent {\it Proof:}
(i) Using 
Definition \ref{def:defls}(v) we find
\begin{equation}
(A-\theta_jI)\sum_{h=0}^i E^*_hV \subseteq 
\sum_{h=0}^{i+1} E^*_hV.
\label{eq:aminusthetajAS99}
\end{equation}
Also observe
\begin{equation}
(A-\theta_jI)\sum_{k=j}^d E_kV = 
\sum_{k=j+1}^d E_kV.
\label{eq:aminusthetajBS99}
\end{equation}
Evaluating $(A-\theta_jI)V_{ij}$ using 
(\ref{eq:defofvijS99}),
(\ref{eq:aminusthetajAS99}), 
(\ref{eq:aminusthetajBS99}) 
 we
routinely find it is contained in $V_{i+1,j+1}$.
\\
\noindent (ii) Immediate from  (i) above.
\\
\noindent (iii) Similar to the proof of (i) above.
\\
\noindent (iv) Immediate from (iii) above.
\hfill $\Box $ \\

\begin{lemma}
\label{lem:manyvijzeroS99}
With reference to 
Definition
\ref{def:VijS99},  we have
\begin{equation}
V_{ij}= 0 \quad \hbox{if}\quad i<j,\qquad  \qquad (0 \leq i,j\leq d).
\label{eq:manyvijzeroS99}
\end{equation}
\end{lemma}

\noindent {\it Proof:}
We show the sum 
\begin{equation}
V_{0r}+V_{1,r+1}+\cdots +V_{d-r,d}
\label{eq:manyvijzeroAS99}
\end{equation}
is zero  for $0 <r \leq d$. Let $r$ be given, and let $W$ denote the sum
in 
(\ref{eq:manyvijzeroAS99}). Applying Lemma 
\ref{lem:howaactsonvijS99}(ii),(iv), we find
$AW\subseteq W$ and $A^*W\subseteq W$.
Now $W=0$ or $W=V$ in view of 
Corollary
\ref{cor:gensetlp}.
 By Definition
\ref{def:VijS99}, each term in 
(\ref{eq:manyvijzeroAS99}) is contained in 
\begin{equation}
E_rV+E_{r+1}V+\cdots +E_dV,
\label{eq:manyvijzeroBS99}
\end{equation}
so $W$ is contained in 
(\ref{eq:manyvijzeroBS99}). The sum 
(\ref{eq:manyvijzeroBS99}) is properly contained in $V$ by 
(\ref{eq:VdecompS99}), and since $r>0$.
Apparently $W\not=V$, so
$W=0$. We have now shown  
(\ref{eq:manyvijzeroAS99}) is zero for $0<r\leq d$, and 
(\ref{eq:manyvijzeroS99}) follows.
\hfill $\Box $ \\

\begin{theorem}
\label{thm:maincharls}
With reference to 
Notation
\ref{def:splitconsetupS99}, let 
$U_0, U_1,\ldots, U_d$ denote subspaces of $V$. 
Then the following (i)--(iii) are equivalent.
\begin{enumerate}
\item $U_i = (E^*_0V+E^*_1V+\cdots + E^*_iV)\cap (E_iV+E_{i+1}V+\cdots + E_dV) \qquad (0 \leq i \leq d)$.
\item  The sequence $U_0, U_1, \ldots, U_d $ is a $\Phi$-split
decomposition of $V$.
\item For $0 \leq i \leq d,$ both 
\begin{eqnarray}
U_i+U_{i+1}+\cdots +U_d &=& E_iV+E_{i+1}V +\cdots + E_dV,
\label{eq:vsumitod}
\\
U_0+U_1+\cdots +U_i &=& E^*_0V+E^*_1V +\cdots + E^*_iV. 
\label{eq:vsumzeroi}
\end{eqnarray}
\end{enumerate}
\end{theorem}

\noindent {\it Proof:}
$(i)\rightarrow (ii)$  
To get
(\ref{eq:splitpart1}) and
(\ref{eq:splitpart2}),
set $j=i$ in Lemma 
\ref{lem:howaactsonvijS99}(i),(iii),  and observe
$U_i = V_{ii}$.
We now show the sequence $U_0, U_1, \ldots, U_d$
is a decomposition of $V$.
Define $W=\sum_{i=0}^d U_i$.
Then $AW\subseteq W$ by 
(\ref{eq:splitpart1}) and
 $A^*W\subseteq W$ by 
(\ref{eq:splitpart2}).
Now 
$W=0$ or $W=V$  in view
of Corollary
\ref{cor:gensetlp}.
However $W$ contains 
 $U_{0}$, and $U_{0}=E^*_0V$ is nonzero,
 so
$W\not=0$.
It follows $W=V$, and in other words 
\begin{equation}
V= U_{0}+U_{1}+\cdots+U_{d}.
\label{eq:vijgivesdirsumCS99}
\end{equation}
We show the sum 
(\ref{eq:vijgivesdirsumCS99}) is direct.  To do this,
we show
\begin{eqnarray*}
 (U_{0}+U_{1}+\cdots +U_{i-1})\cap U_{i} =0
\end{eqnarray*}
for $1 \leq i \leq d$.
Let the integer $i$ be given. From the construction 
\begin{eqnarray*}
U_{j} \subseteq E^*_0V+E^*_1V+\cdots +E^*_{i-1}V
\end{eqnarray*}
for $0 \leq j \leq i-1$, and
\begin{eqnarray*}
U_{i}\subseteq E_iV+E_{i+1}V+\cdots +E_dV.
\end{eqnarray*}
It follows
\begin{eqnarray*}
 &&(U_{0}+U_{1}+\cdots +U_{i-1})\cap U_{i}
\\
 && \qquad  \subseteq
 (E^*_0V+E^*_1V+\cdots +E^*_{i-1}V)\cap
(E_iV+E_{i+1}V+\cdots +E_dV) \qquad \qquad 
\\
&& \qquad  = V_{i-1,i}
\\
&& \qquad  = 0
\end{eqnarray*}
in view of Lemma 
\ref{lem:manyvijzeroS99}.
We have now shown
the sum (\ref{eq:vijgivesdirsumCS99}) is direct.
We now show $U_i$ has dimension 1 for $0 \leq i \leq d$.
Since the sum
(\ref{eq:vijgivesdirsumCS99}) is direct, this will follow
if we can show $U_i \not=0$ for $0\leq i \leq d$.
Suppose there exists an integer $i$ $(0 \leq i \leq d)$ such
that $U_i=0$. We observe $i\not=0$, since $U_0=E^*_0V$
is nonzero, and $i\not=d$, since $U_d=E_dV$ is nonzero.
Set 
\begin{eqnarray*}
U=U_0+U_1+\cdots + U_{i-1},
\end{eqnarray*}
and observe $U\not=0$ and $U\not=V$ by our remarks above.
By Lemma
\ref{lem:howaactsonvijS99}(ii) and since $U_i=0$ we find
$AU\subseteq U$. By Lemma
\ref{lem:howaactsonvijS99}(iv) we find
$A^*U\subseteq U$. Now $U=0$ or $U=V$ in view of
Corollary
\ref{cor:gensetlp}, for a contradiction.
We conclude $U_i\not=0$ for $0 \leq i \leq d$ and it
follows
$U_0, U_1, \ldots, U_d$ is a decomposition of $V$.

\noindent 
 $(ii)\rightarrow (iii)$  
First consider 
(\ref{eq:vsumitod}).  
Let $i$ be given, and abbreviate
\begin{eqnarray*}
Z&=& E_iV+E_{i+1}V+\cdots + E_dV, \qquad \quad 
W=U_i+U_{i+1}+\cdots + U_d.
\end{eqnarray*}
We show $Z=W$. To obtain $Z\subseteq W$, set 
$X=\prod_{h=0}^{i-1} (A-\theta_hI)$, and observe
$Z=XV$.
 Using 
(\ref{eq:splitpart1})
we find $XU_j\subseteq W$ for $0 \leq j \leq d$. By this and since
$U_0, U_1, \ldots, U_d$ is a decomposition of $V$
we find
$XV\subseteq W$.
We now have 
$Z\subseteq W$. 
Each of $Z,W$ has dimension $d-i+1$ so $Z=W$. 
%To obtain $W\subseteq Z$, set 
%$Y=\prod_{h=i}^{d} (A-\theta_hI)$, and observe
%\begin{eqnarray}
%Z&=&\lbrace v \in V \;|\;Yv = 0\rbrace .
%\label{eq:Keraction}
%\end{eqnarray}
%Observe
%$YU_j=0$ for $i \leq j \leq d$, so $YW=0$. Combining 
%this with
%(\ref{eq:Keraction}), we find $W\subseteq Z$.
We now have
(\ref{eq:vsumitod}).
Line (\ref{eq:vsumzeroi})  is similarily obtained.

\noindent $(iii)\rightarrow (i)$ We first
show the sum $U_0+\cdots+U_d$ is direct.
 To do this, we show
\begin{equation} 
 (U_{0}+U_{1}+\cdots +U_{i-1})\cap U_{i}
\label{eq:linindepag}
\end{equation}
is zero for $1 \leq i \leq d$. Let $i$ be given. From
(\ref{eq:vsumitod}), 
 (\ref{eq:vsumzeroi}),
we find (\ref{eq:linindepag}) is contained in 
\begin{equation}
 (E^*_0V+E^*_1V+\cdots +E^*_{i-1}V)\cap
(E_iV+E_{i+1}V+\cdots +E_dV). \qquad \qquad 
\label{eq:thatintagain}
\end{equation}
The expression (\ref{eq:thatintagain}) equals $V_{i-1,i}$, and is hence zero
by Lemma 
\ref{lem:manyvijzeroS99}. It follows 
(\ref{eq:linindepag}) is zero, and we have now shown
the sum $U_0+\cdots+U_d$ is direct.
Combining  this with  
(\ref{eq:vsumitod}), 
 (\ref{eq:vsumzeroi}), we find
\begin{eqnarray*}
U_i &=& (U_0+U_1+\cdots + U_i)\cap (U_i+U_{i+1}+\cdots + U_d)
\\
&=& 
 (E^*_0V+E^*_1V+\cdots +E^*_iV)\cap
(E_iV+E_{i+1}V+\cdots +E_dV),
\end{eqnarray*}
as desired.
\hfill $\Box $ \\

\begin{corollary}
\label{cor:splitdecexists}
With reference to Notation 
\ref{def:splitconsetupS99},
there exists a unique $\Phi$-split decomposition of $V$.
\end{corollary}
\noindent {\it Proof:}
Immediate from Theorem
\ref{thm:maincharls}(i),(ii).
\hfill $\Box $ \\

\noindent We finish this section with a comment.

\begin{lemma}
\label{lem:splitstrongver}
With reference to Notation 
\ref{def:splitconsetupS99}, let $U_0, U_1, \ldots, U_d$
denote the $\Phi$-split decomposition of $V$.
Then the following (i), (ii) hold.
\begin{enumerate}
\item $(A-\theta_iI)U_i=U_{i+1} \qquad  (0 \leq i \leq d-1)$.
\item $(A^*-\theta^*_iI)U_i=U_{i-1} \qquad (1 \leq i \leq d)$.
\end{enumerate}
\end{lemma}
\noindent {\it Proof:}
(i) Let $i$ be given.
Recall $(A-\theta_iI)U_i$ is contained in $U_{i+1}$
by 
(\ref{eq:splitpart1}) and $U_{i+1}$ has dimension 1,
so it suffices to show
\begin{eqnarray}
(A-\theta_iI)U_i\not=0.
\label{eq:itszero}
\end{eqnarray}
Assume 
$(A-\theta_iI)U_i=0$, and set
$W=\sum_{h=0}^i U_h$.
Since $U_0, U_1, \ldots, U_d$ is a decomposition of $V$,
and since $0 \leq i \leq d-1$ we find $W\not=0$ and $W\not=V$.
Observe $AU_i\subseteq U_i$ by our above assumption;
combining this with
(\ref{eq:splitpart1})
we find
$AW\subseteq W$. By
(\ref{eq:splitpart2}) we find
$A^*W\subseteq W$. Now $W=0$ or $W=V$ in view of
Corollary
\ref{cor:gensetlp}, for a contradiction.
We conclude
(\ref{eq:itszero}) holds and the result follows.
\\
\noindent (ii) Similar to the proof of (i) above.
\hfill $\Box $ \\

\section{The split basis}

Let 
$\Phi=(A;A^*; \lbrace E_i\rbrace_{i=0}^d;  
\lbrace E^*_i\rbrace_{i=0}^d)$
denote a Leonard system in $\mathcal A$,
with eigenvalue sequence $\theta_0, \theta_1, \ldots, \theta_d$
and dual eigenvalue sequence
 $\theta^*_0, \theta^*_1, \ldots, \theta^*_d$.
Let $V$ denote
an irreducible $\cal A$-module and
let $U_0, U_1, \ldots, U_d$ denote the $\Phi$-split
decomposition of $V$ from
Definition
\ref{def:splitdec}.
Pick any integer $i$ $(1 \leq i\leq d)$. 
By Lemma
\ref{lem:splitstrongver} we have
$(A^*-\theta^*_iI)U_i=U_{i-1}$ 
and 
$(A-\theta_{i-1}I)U_{i-1}=U_{i}$.
Apparently $U_i$ is an eigenspace for 
$(A-\theta_{i-1}I)
(A^*-\theta^*_iI)$, and the corresponding eigenvalue
is a nonzero scalar in $\K$.
We denote this eigenvalue by $\varphi_i$.
We display a basis for $V$ that illuminates the 
significance of $\varphi_i$.
Setting $i=0$ in Theorem
\ref{thm:maincharls}(i) we find $U_0=E^*_0V$.
Combining this with 
Lemma
\ref{lem:splitstrongver}(i) we find
\begin{eqnarray}
U_i = (A-\theta_{i-1}I)\cdots (A-\theta_1I)(A-\theta_0I)E^*_0V
\qquad \qquad (0 \leq i \leq d).
\label{eq:uiinterp}
\end{eqnarray}
Let $v$ denote a nonzero vector in $E^*_0V$. From
(\ref{eq:uiinterp}) we find that for $0 \leq i \leq d$
the vector
$(A-\theta_{i-1}I)\cdots (A-\theta_0I)v$
is a basis for $U_i$. By this and since $U_0, U_1, \ldots, U_d$
is a decomposition of $V$ we find the sequence
\begin{eqnarray}
(A-\theta_{i-1}I)\cdots (A-\theta_1I)(A-\theta_0I)v
\qquad \qquad (0 \leq i \leq d)
\label{eq:splitbasisdef}
\end{eqnarray}
is a basis for $V$. With respect to this basis the
matrices representing $A$ and $A^*$ are
\begin{equation}
\left(
\begin{array}{c c c c c c}
\theta_0 & & & & & {\bf 0} \\
1 & \theta_1 &  & & & \\
& 1 & \theta_2 &  & & \\
& & \cdot & \cdot &  &  \\
& & & \cdot & \cdot &  \\
{\bf 0}& & & & 1 & \theta_d
\end{array}
\right),
\qquad  \quad 
%A^{*\heartsuit} = 
\left(
\begin{array}{c c c c c c}
\theta^*_0 &\varphi_1 & & & & {\bf 0} \\
 & \theta^*_1 & \varphi_2 & & & \\
&  & \theta^*_2 & \cdot & & \\
& &  & \cdot & \cdot &  \\
& & &  & \cdot & \varphi_d \\
{\bf 0}& & & &  & \theta^*_d
\end{array}
\right) 
\label{eq:matrepaastar}
\end{equation}
respectively.
By a {\it $\Phi$-split basis} for $V$ we mean
a sequence
 of the form 
(\ref{eq:splitbasisdef}), where $v$ is a nonzero vector
in $E^*_0V$.
We call
$ \varphi_1, \varphi_2,\ldots, \varphi_d$
the {\it first split sequence} of $\Phi$.
We let 
$ \phi_1, \phi_2,\ldots, \phi_d$
denote the first split sequence
of $\Phi^\Downarrow$ and call
this the 
 {\it second split sequence} of $\Phi$.
For notational convenience we define
$\varphi_0=0$, 
$\varphi_{d+1}=0$, 
$\phi_0=0$, 
$\phi_{d+1}=0$.

\section{The parameter array and the classifying space}

Our next goal is to describe the relationship between
the 
eigenvalue sequence,
the dual eigenvalue sequence,
the first split sequence,
and the second split sequence.
We will use the following concept.

\begin{definition} 
%\cite[Definition 10.1]{conform}
\label{def:pa}
\rm
Let $d$ denote a nonnegative integer.
By a {\it parameter array over $\K$ of diameter $d$}
we mean a sequence of scalars 
$
(\theta_i, \theta^*_i,i=0..d; \varphi_j, \phi_j,j=1..d)
$
taken from $\K$ that satisfy the following conditions
(PA1)--(PA5).
\begin{description}
\item[(PA1)]  
$
\theta_i \not=\theta_j, \qquad 
\theta^*_i \not=\theta^*_j \qquad \mbox{if}\quad i\not=j, \qquad \qquad
(0 \leq i,j\leq d)
$.
\\
\item[(PA2)]
$
\varphi_i \not=0, \qquad \phi_i \not=0 \qquad \qquad (1 \leq i \leq d).
$
\\
\item[(PA3)]
$\varphi_i = \phi_1 \sum_{h=0}^{i-1} 
\frac{\theta_h-\theta_{d-h}}{\theta_0-\theta_d}
+ (\theta^*_i-\theta^*_0)(\theta_{i-1}-\theta_d)
\qquad \qquad (1 \leq i \leq d)$.
\\
\item[(PA4)] 
$\phi_i = \varphi_1 \sum_{h=0}^{i-1} 
\frac{\theta_h-\theta_{d-h}}{\theta_0-\theta_d}
+ (\theta^*_i-\theta^*_0)(\theta_{d-i+1}-\theta_0)
\qquad \qquad (1 \leq i \leq d)$.
\\
\item[(PA5)]
The expressions
\begin{eqnarray} 
\frac{\theta_{i-2}-\theta_{i+1}}{\theta_{i-1}-\theta_i},
\qquad 
\frac{\theta^*_{i-2}-\theta^*_{i+1}}{\theta^*_{i-1}-\theta^*_i}
\label{eq:betaplusone}
\end{eqnarray}
are equal and independent of $i$ for $2 \leq i \leq d-1$.
\end{description}
\end{definition}

\begin{theorem} 
\cite[Theorem 1.9]{LS99}
\label{thm:ls}
Let $d$ denote a nonnegative  integer and let 
$(\theta_i, \theta^*_i,i=0..d;\varphi_j,\phi_j,j=1..d)$
denote a sequence of scalars taken from $\K$.
Then the following (i), (ii) are equivalent.
\begin{enumerate}
\item The sequence 
$(\theta_i, \theta^*_i,i=0..d;\varphi_j,\phi_j,j=1..d)$
is a parameter array over $\K$.
\item
There exists 
a Leonard system $\Phi$ over $\K$ that has
eigenvalue sequence 
$\theta_0, \theta_1, \ldots, \theta_d$, dual eigenvalue sequence
$\theta^*_0, \theta^*_1, \ldots, \theta^*_d$,
first split sequence 
$\varphi_1, \varphi_2, \ldots, \varphi_d$ and second split sequence
$\phi_1, \phi_2, $ $ \ldots, \phi_d$.
\end{enumerate}
Suppose (i), (ii) hold. Then $\Phi$ is 
unique up to isomorphism of Leonard systems.
\end{theorem}

\noindent Our proof of Theorem
\ref{thm:ls} is too long to be included in these notes.
A complete proof can be found in
\cite{LS99}.

\begin{definition}
\label{def:pararray}
\rm
Let $\Phi$ denote the Leonard system
from 
Definition
\ref{eq:ourstartingpt}.
By the {\it parameter array} of $\Phi$ we mean
the sequence
$(\theta_i,\theta^*_i,i=0..d;\varphi_j,\phi_j,j=1..d)$,
where  
 $\theta_0, \theta_1, \ldots, \theta_d$ 
  (resp. $\theta^*_0, \theta^*_1, \ldots, \theta^*_d$)
  is the eigenvalue sequence (resp. 
  dual eigenvalue sequence) of $\Phi$
and   
$ \varphi_1, \varphi_2,\ldots, \varphi_d$
(resp.
$ \phi_1, \phi_2,\ldots, \phi_d$)
is 
the first split sequence
(resp. second split sequence) of
 $\Phi$.
\end{definition}

\noindent
By Theorem \ref{thm:ls} the map which sends a given Leonard system
to its parameter array induces a bijection from the set of isomorphism
classes of Leonard systems over $\K$ to the set of parameter arrays over
$\K$.  Consequently we view the set of
parameter arrays over $\K$ as a ``classifying space'' for
the Leonard systems over $\K$.

\medskip
\noindent In the appendix to these notes we display all the
parameter arrays over $\K$.

\medskip
\noindent 
We now cite a result that shows
how the parameter arrays
 behave with respect to the $D_4$ action given in
 Section 4.

\begin{theorem} \cite[Theorem 1.11]{LS99}
\label{thm:phimod}
Let $\Phi$ denote a Leonard system
with parameter array
$(\theta_i,\theta^*_i,i=0..d; $ $\varphi_j,\phi_j,j=1..d)$.
 Then  (i)--(iii) hold below.
\begin{enumerate}
\item  The parameter array of $\Phi^*$
is
$(\theta^*_i,\theta_i,i=0..d;\varphi_j,\phi_{d-j+1},j=1..d)$.
\item 
 The parameter array of
 $\Phi^{\downarrow}$
is
$(\theta_i,\theta^*_{d-i},i=0..d;\phi_{d-j+1},\varphi_{d-j+1},j=1..d)$.
\item  The parameter array of $\Phi^\Downarrow$
is
$(\theta_{d-i},\theta^*_i,i=0..d;\phi_j,\varphi_j,j=1..d)$.
\end{enumerate}
\end{theorem}

\section{Everything in terms of the parameter array}

\noindent 
In this section we express  all the polynomials and
scalars that came up so far in the paper, in terms
the parameter array.
We will use the following notation.

\begin{definition}
\label{def:tau}
\rm
Suppose we are given an integer
$d\geq 0$ and two sequences of scalars
\begin{eqnarray*}
\theta_0, \theta_1, \ldots, \theta_d; \quad
\theta^*_0, \theta^*_1, \ldots, \theta^*_d
\end{eqnarray*}
taken from $\K$.
Then for $0 \leq i \leq d+1$ we let 
$\tau_i$, 
$\tau^*_i$, 
$\eta_i$,
$\eta^*_i$
denote the following polynomials in
$\K \lbrack \lambda \rbrack $.
\begin{eqnarray}
&&\tau_i = \prod_{h=0}^{i-1} (\lambda - \theta_h),
\qquad \qquad \;\;
\tau^*_i = \prod_{h=0}^{i-1} (\lambda - \theta^*_h),
\label{eq:tau}
\label{eq:taus}
\\
&&\eta_i = \prod_{h=0}^{i-1}(\lambda - \theta_{d-h}),
\qquad \qquad
\eta^*_i = \prod_{h=0}^{i-1}(\lambda - \theta^*_{d-h}).
\label{eq:eta}
\label{eq:etas}
\end{eqnarray}
%%%%%
%\begin{eqnarray}
%\tau_i &=& (\lambda - \theta_0)
%(\lambda - \theta_1) \cdots
%(\lambda - \theta_{i-1}),
%\label{eq:tau}
%\\
%\tau^*_i &=& (\lambda - \theta^*_0)
%(\lambda - \theta^*_1) \cdots
%(\lambda - \theta^*_{i-1}),
%\label{eq:taus}
%\\
%\eta_i &=& (\lambda - \theta_d)
%(\lambda - \theta_{d-1}) \cdots
%(\lambda - \theta_{d-i+1}),
%\label{eq:eta}
%\\
%\eta^*_i &=& (\lambda - \theta^*_d)
%(\lambda - \theta^*_{d-1}) \cdots
%(\lambda - \theta^*_{d-i+1}).
%\label{eq:etas}
%\end{eqnarray}
\noindent 
We observe that each of
$\tau_i$, 
$\tau^*_i$, 
$\eta_i$,
$\eta^*_i$
is monic with degree $i$.
\end{definition}

\begin{theorem}
\label{thm:uform}
Let $\Phi$ 
denote the Leonard system from Definition
\ref{eq:ourstartingpt}
and 
let 
$(\theta_i, \theta^*_i, i=0..d;
\varphi_j,\phi_j,j=1..d)$ denote the corresponding
parameter array.
Let the polynomials $u_i$ be as in
Definition \ref{def:ui1}.
Then
\begin{eqnarray}
\label{eq:uifinal}
u_i=\sum_{h=0}^i {{\tau^*_h(\theta^*_i)}\over 
{\varphi_1\varphi_2\cdots \varphi_h}} \tau_h 
\qquad \qquad 
(0 \leq i \leq d).
\end{eqnarray}
We are using the notation
 (\ref{eq:taus}).
\end{theorem}
\noindent {\it Proof:}
Let the integer $i$ be given. 
The polynomial $u_i$ has degree $i$ so
 there  exists 
 scalars $\alpha_0, \alpha_1, \ldots, \alpha_i$ 
in $\K$ such that
\begin{eqnarray}
\label{eq:ualphaform}
u_i = \sum_{h=0}^i \alpha_h \tau_h.
\end{eqnarray}
We show
\begin{eqnarray}
\alpha_h = \frac{\tau^*_h(\theta^*_i)}{\varphi_1\varphi_2 \cdots \varphi_h}
\qquad \qquad (0 \leq h \leq i).
\label{eq:alphagoal}
\end{eqnarray}
In order to do this we show
$\alpha_0=1$ and 
$\alpha_{h+1}\varphi_{h+1}=\alpha_h(\theta^*_i-\theta^*_h)$
for $0 \leq h \leq i-1$.
We now show $\alpha_0=1$.
We evaluate (\ref{eq:ualphaform})
at $\lambda=\theta_0$
and find $u_i(\theta_0)=\sum_{h=0}^i \alpha_h \tau_h(\theta_0)$.
Recall $u_i(\theta_0)=1$ by
(\ref{eq:unorm}).
Using 
(\ref{eq:tau})
we find $\tau_h(\theta_0)=1$  for $h=0$ and
$\tau_h(\theta_0)=0$  for
$1 \leq h \leq i$.
From these comments  
we find $\alpha_0=1$.
We  now show
$\alpha_{h+1} \varphi_{h+1} = 
\alpha_h(\theta^*_i-\theta^*_h)$
for $0 \leq h \leq i-1$.
Let $V$ denote an irreducible $\mathcal A$-module.
Let $v$ denote a nonzero vector in $E^*_0V$ and
define $e_i = \tau_i(A)v$ for $0 \leq i\leq d$.
Observe that the sequence $e_0, e_1, \ldots, e_d$
is the basis for $V$ from
(\ref{eq:splitbasisdef}). Using 
(\ref{eq:matrepaastar}) we find
$(A^*-\theta^*_jI)e_j = \varphi_je_{j-1}$ for
$1 \leq j \leq d$ and
$(A^*-\theta^*_0I)e_0 =0$.
By Theorem \ref{thm:psend} and
(\ref{uicl})
%, and since $v$ is a basis for $E^*_0V$
we find
$u_i(A)E^*_0V=E^*_iV$.
By this and since
 $v \in E^*_0V$ 
we find
$u_i(A)v \in E^*_iV$. 
Apparently 
$u_i(A)v $ is an eigenvector for $A^*$ with
eigenvalue $\theta^*_i$.
We may now argue
\begin{eqnarray*}
0  &=& (A^*-\theta^*_iI)u_i(A)v
\\
  &=& (A^*-\theta^*_iI)\sum_{h=0}^i \alpha_h \tau_h(A)v
\\
  &=& (A^*-\theta^*_iI)\sum_{h=0}^i \alpha_h e_h 
%\\
%  &=& \sum_{h=0}^i \alpha_h ((\theta^*_h-\theta^*_i)e_h +\varphi_h e_{h-1})
\\
  &=& \sum_{h=0}^{i-1} e_h 
  (\alpha_{h+1} \varphi_{h+1}-\alpha_{h}(\theta^*_i-\theta^*_h)).
\end{eqnarray*}
By this and since $e_0, e_1, \ldots, e_d$ are linearly independent we find
  $\alpha_{h+1} \varphi_{h+1}=\alpha_h(\theta^*_i-\theta^*_h)$ for
  $0 \leq h\leq i-1$.
Line
(\ref{eq:alphagoal}) follows and the theorem is proved.
\hfill $\Box $ \\

\begin{lemma}
\label{lem:p0}
Let $\Phi$ denote the Leonard system
from 
Definition 
\ref{eq:ourstartingpt}
and let
$(\theta_i, \theta^*_i, i=0..d;
\varphi_j,\phi_j,j=1..d)$ denote the corresponding
parameter array.
Let the polynomials 
$p_i$ be as in
Definition
\ref{def:pi}.
With reference to
Definition \ref{def:tau} we have
\begin{eqnarray}
\label{eq:pith0first}
p_i(\theta_0)=  \frac{\varphi_1\varphi_2\cdots \varphi_i}{\tau^*_i(\theta^*_i)}
\qquad \qquad (0 \leq i \leq d).
\end{eqnarray}
\end{lemma}
\noindent {\it Proof:}
In equation
(\ref{eq:uifinal}), each side is a polynomial
of degree $i$ in $\lambda $.
For the polynomial on the left 
in 
(\ref{eq:uifinal})
the coefficient of $\lambda^i$
is $p_i(\theta_0)^{-1}$
by 
(\ref{uicl}) and since  $p_i$ is monic.
For the polynomial
on the right
in 
(\ref{eq:uifinal})
the coefficient of
 $\lambda^i$ is
$\tau^*_i(\theta^*_i)(\varphi_1\varphi_2\cdots \varphi_i)^{-1}$.
Comparing these coefficients we obtain
the result.
\hfill $\Box $ \\

\begin{theorem}
\label{thm:pform}
Let $\Phi$ 
denote the Leonard system from Definition
\ref{eq:ourstartingpt}
and 
let 
$(\theta_i, \theta^*_i, i=0..d;
\varphi_j,\phi_j,j=1..d)$ denote the corresponding
parameter array.
Let the polynomials 
$p_i$ be as in
Definition \ref{def:pi}.
Then with reference to
Definition \ref{def:tau},
\begin{eqnarray*}
p_i=\sum_{h=0}^i \frac{\varphi_1\varphi_2\cdots \varphi_i}{
\varphi_1\varphi_2\cdots \varphi_h}
\frac{\tau^*_h(\theta^*_i)} 
{\tau^*_i(\theta^*_i)} \tau_h
\qquad \qquad 
(0 \leq i \leq d).
\end{eqnarray*}
%We are using the notation
%(\ref{eq:tau}), (\ref{eq:taus}).
\end{theorem}
\noindent {\it Proof:}
Observe $p_i = p_i(\theta_0)u_i$ by
(\ref{uicl}). In this equation we evaluate
$p_i(\theta_0)$ using
(\ref{eq:pith0first}) and
we evaluate $u_i$ using
(\ref{eq:uifinal}). The result follows.
\hfill $\Box $ \\

\begin{theorem}
\label{thm:bcform}
Let $\Phi$ 
denote the Leonard system from Definition
\ref{eq:ourstartingpt}
and 
let 
$(\theta_i, \theta^*_i, i=0..d;
\varphi_j,\phi_j,j=1..d)$ denote the corresponding
parameter array.
Let the scalars  $b_i, c_i$ be as in
Definition
\ref{def:sharpmp}.
Then with reference to
Definition \ref{def:tau}
the following (i), (ii) hold.
\begin{enumerate}
\item
$\displaystyle
{b_i = \varphi_{i+1} \frac{\tau^*_i(\theta^*_i)}{
\tau^*_{i+1}(\theta^*_{i+1})} \qquad \qquad 
(0 \leq i \leq d-1)
}$.
\\
\item
$
\displaystyle{
c_i = \phi_i \frac{\eta^*_{d-i}(\theta^*_i)}{\eta^*_{d-i+1}(\theta^*_{i-1})}
\qquad \qquad (1 \leq i \leq d)}$.
\end{enumerate}
\end{theorem}
\noindent {\it Proof:}
(i) Evaluate
(\ref{eqbi}) using
Lemma
\ref{lem:p0}.
\\
(ii) Using 
Definition
\ref{def:sharpmp}
we find, with reference to
Definition
\ref{def:notconv},
that $c_i=b^\downarrow_{d-i}$.
Applying
part (i) above to $\Phi^\downarrow$
and using
Theorem \ref{thm:phimod}(ii) we routinely
obtain the result.
\hfill $\Box $ \\

\noindent Let $\Phi$ denote the Leonard system
from 
 Definition
\ref{eq:ourstartingpt}
and let the scalars $a_i$ be as in
Definition
\ref{def:aixi}.
We  mention two formulae that give  
$a_i$ in terms of the parameter array of $\Phi$.
The first formula
is obtained using
Lemma
\ref{def:bici}(ii) and
Theorem \ref{thm:bcform}.
The second formula is given in the following theorem.

\begin{theorem}
%\cite[Lemma 5.1]{LS99}
\label{thm:ai}
Let $\Phi$ 
denote the Leonard system from Definition
\ref{eq:ourstartingpt}
and 
let 
$(\theta_i, \theta^*_i, i=0..d;
\varphi_j,\phi_j,j=1..d)$ denote the corresponding
parameter array.
Let the scalars $a_i$ be as in 
Definition
\ref{def:aixi}. Then
\begin{eqnarray}
a_i = \theta_i + \frac{\varphi_i}{\theta^*_i-\theta^*_{i-1}}
+ \frac{\varphi_{i+1}}{\theta^*_i-\theta^*_{i+1}} \qquad (0 \leq i \leq d),
\label{eq:aiform}
\end{eqnarray}
where we recall $\varphi_0=0$, $\varphi_{d+1}=0$, and where
$\theta^*_{-1}, \theta^*_{d+1}$ denote indeterminates.
\end{theorem}
\noindent {\it Proof:}
Let the polynomials $p_0, p_1, \ldots, p_{d+1}$ be as in
Definition
\ref{def:pi} and recall these polynomials are monic.
Let $i$ be given and consider the polynomial
\begin{eqnarray}
\lambda p_i-p_{i+1}.
\label{eq:picomb}
\end{eqnarray}
From (\ref{eq:prec}) 
we find
the polynomial (\ref{eq:picomb}) is equal to 
$a_ip_i+x_ip_{i-1}$.
Therefore 
the polynomial
(\ref{eq:picomb})
has degree $i$ and leading coefficient $a_i$.
In order to compute this leading coefficient,
in 
(\ref{eq:picomb})
we 
evaluate each of 
$p_i$, $p_{i+1}$
using
Theorem \ref{thm:pseq}(ii)
and 
Theorem
\ref{thm:pform}.
By this method we routinely obtain
(\ref{eq:aiform}).
\hfill $\Box $ \\

\begin{theorem}
\label{thm:xform}
Let $\Phi$ 
denote the Leonard system from Definition
\ref{eq:ourstartingpt}
and 
let 
$(\theta_i, \theta^*_i, i=0..d;
\varphi_j,\phi_j,j=1..d)$ denote the corresponding
parameter array.
Let the scalars $x_i$ be
as in
Definition
\ref{def:aixi}.
Then with reference to
Definition \ref{def:tau},
\begin{eqnarray}
\label{eq:xifinform}
x_i = 
 \varphi_{i}
\phi_i  
 \frac{\tau^*_{i-1}(\theta^*_{i-1})
\eta^*_{d-i}(\theta^*_i)}
{
\tau^*_{i}(\theta^*_{i})
\eta^*_{d-i+1}(\theta^*_{i-1})}
\qquad \qquad (1 \leq i \leq d).
\end{eqnarray}
\end{theorem}
\noindent {\it Proof:}
Use $x_i=b_{i-1}c_i$ and
Theorem
\ref{thm:bcform}.
\hfill $\Box $ \\

\begin{theorem}
\label{thm:nform}
Let $\Phi$ 
denote the Leonard system from Definition
\ref{eq:ourstartingpt}
and 
let 
$(\theta_i, \theta^*_i, i=0..d;
\varphi_j,\phi_j,j=1..d)$ denote the corresponding
parameter array.
Let the scalar $\nu$ be
as in
Definition
\ref{def:n}. 
Then with reference to
Definition \ref{def:tau},
\begin{eqnarray}
\nu = \frac{\eta_d(\theta_0) \eta^*_d(\theta^*_0)
}
{\phi_1 \phi_2 \cdots \phi_d 
}.
\label{eq:nfinalform}
\end{eqnarray}
\end{theorem}
\noindent {\it Proof:}
Evaluate 
(\ref{frame}) using
Theorem \ref{thm:bcform}(ii).
\hfill $\Box $ \\

\begin{theorem}
\label{thm:kiform}
Let $\Phi$ 
denote the Leonard system from Definition
\ref{eq:ourstartingpt}
and 
let 
$(\theta_i, \theta^*_i, i=0..d;
\varphi_j,\phi_j,j=1..d)$ denote the corresponding
parameter array.
Let the scalars $k_i$ be
as in
Definition
\ref{def:ki1}.
Then with reference to
Definition \ref{def:tau},
\begin{eqnarray}
\label{eq:kform2}
k_i = 
\frac{\varphi_1\varphi_2\cdots \varphi_i} 
{\phi_1\phi_2\cdots \phi_i}
\frac{\eta^*_d(\theta^*_0)}{\tau^*_i(\theta^*_i)\eta^*_{d-i}(\theta^*_i)}
\qquad \qquad (0 \leq i \leq d).
\end{eqnarray}
\end{theorem}
\noindent {\it Proof:}
Evaluate 
(\ref{eq:kibc})
using 
Theorem \ref{thm:bcform}.
\hfill $\Box $ \\

\begin{theorem}
\label{thm:miform}
Let $\Phi$ 
denote the Leonard system from Definition
\ref{eq:ourstartingpt}
and 
let 
$(\theta_i, \theta^*_i, i=0..d;
\varphi_j,\phi_j,j=1..d)$ denote the corresponding
parameter array.
Let the scalars $m_i$ be
as in
Definition
\ref{def:mid}.
Then with reference to
Definition \ref{def:tau},
\begin{eqnarray}
\label{eq:mform2}
m_i = 
\frac{\varphi_1\varphi_2\cdots \varphi_i 
\phi_1\phi_2\cdots \phi_{d-i}}
{\eta^*_d(\theta^*_0)\tau_i(\theta_i)\eta_{d-i}(\theta_i)}
\qquad \qquad (0 \leq i \leq d).
\end{eqnarray}
\end{theorem}
\noindent {\it Proof:}
Applying
Definition
\ref{def:ki1} to
$\Phi^*$ we find
$m_i=k^*_i\nu^{-1}$.
We compute $k^*_i$ using
Theorem \ref{thm:kiform}
and Theorem
\ref{thm:phimod}(i).
We compute 
$\nu$ using
Theorem \ref{thm:nform}.
The result follows.
\hfill $\Box $ \\

\section{The terminating branch of the Askey scheme}

Let $\Phi$ denote the Leonard system from
Definition
\ref{eq:ourstartingpt} and let the polynomials
$u_i$ be as in
Definition \ref{def:ui1}. In this section we discuss 
how the $u_i$ fit into the Askey scheme
\cite{KoeSwa}, 
 \cite[p260]{BanIto}. 
Our argument is summarized as follows.
In the appendix to these notes
we display all the parameter arrays over $\K$.
These parameter arrays fall into 13 families.
In (\ref{eq:uifinal})  the $u_i$ are expressed  as a sum involving
the parameter array of $\Phi$.
For each of the 13 families of parameter arrays
we
evaluate this sum.
We find the corresponding $u_i$  
form a class  
consisting of the 
$q$-Racah, $q$-Hahn, dual $q$-Hahn, 
$q$-Krawtchouk,
dual $q$-Krawtchouk,
quantum 
$q$-Krawtchouk,
affine 
$q$-Krawtchouk,
Racah, Hahn, dual Hahn, Krawtchouk,  Bannai/Ito, 
and 
orphan polynomials. 
This class coincides with the terminating branch of the 
Askey scheme. See the appendix for the details.
We remark the Bannai/Ito polynomials can be obtained
from the $q$-Racah polynomials by letting $q$ tend to $-1$ 
\cite[p260]{BanIto}. The orphan polynomials
exist for diameter $d=3$ and $\mbox{Char}(\K)=2$ only.

\medskip
\noindent
In this section we illustrate what is going on with
 some examples.
We will consider two families of parameter  arrays.
For the first family the corresponding 
$u_i$ will turn out to be some Krawtchouk
polynomials. For the  second family the corresponding
$u_i$ will turn out to
be the $q$-Racah polynomials.
%\noindent In \cite{TLT:array} 
%we displayed all the parameter arrays in parametric form.
%We showed that for these arrays the
%corresponding polynomials
%from the right-hand side of
%(\ref{eq:uifinal}) form a class consisting of
%the $q$-Racah, $q$-Hahn, dual $q$-Hahn, 
%$q$-Krawtchouk,
%dual $q$-Krawtchouk,
%quantum 
%$q$-Krawtchouk,
%affine 
%$q$-Krawtchouk,
%Racah, Hahn, dual Hahn, Krawtchouk,  Bannai/Ito, 
%and 
%orphan polynomials.
%This class coincides with  the terminating branch of the Askey scheme
%of orthogonal polynomials
%\cite{KoeSwa}, 
% \cite[p260]{BanIto}. 
%In this section we illustrate what is going
%on with some examples.
%We will consider two families of parameter  arrays.
%The first family will yield some Krawtchouk
%polynomials and the second family will yield the
%$q$-Racah polynomials.

\medskip
\noindent Our first example
is associated with the Leonard pair
(\ref{eq:fam1}). 
%For this first example the 
%polynomials $u_i$ 
%will turn out to be some Krawtchouk polynomials.
Let $d$ denote a nonnegative integer and 
consider the following elements of $\K$.
\begin{eqnarray}
&&\theta_i =  d-2i, \qquad \qquad \theta^*_i = d-2i 
\qquad \qquad (0 \leq i \leq d),
\label{eq:thsol}
\\
&&
\varphi_i = -2i(d-i+1), \qquad \qquad \phi_i = 2i(d-i+1) 
\qquad  \qquad (1 \leq i \leq d).
\label{eq:vpsol}
\end{eqnarray}
In order to avoid degenerate situations
we assume the characteristic
of $\K$ is zero or an odd prime greater than $d$.
It is routine to show 
(\ref{eq:thsol}), 
(\ref{eq:vpsol}) satisfy the conditions PA1--PA5 of 
Definition
\ref{def:pa}, so
$(\theta_i, \theta^*_i,i=0..d; \varphi_j, \phi_j,j=1..d)$
is a parameter array over $\K$.
By 
Theorem \ref{thm:ls}
there exists a Leonard system
$\Phi$
over $\K$ with this parameter array.
Let
the scalars $a_i$ for
  $\Phi$ be as in
  (\ref{eq:aitr}).
Applying Theorem
\ref{thm:ai} to $\Phi$ we find
\begin{eqnarray}
\label{eq:akraw}
a_i=0 \qquad (0 \leq i \leq d).
\end{eqnarray}
Let the scalars $b_i, c_i$ for $\Phi$
be as in
Definition
\ref{def:sharpmp}.
Applying
Theorem 
\ref{thm:bcform} to $\Phi$ we  find
\begin{eqnarray}
\label{eq:bckraw}
b_i = d-i, \qquad c_i = i \qquad \qquad (0 \leq i \leq d).
\end{eqnarray}
%Let the map $\flat$ for $\Phi$ be as in
%Definition \ref{def:flatcon}.
%Evaluating
%(\ref{eq:matrepls}) using
%(\ref{eq:akraw})
%and 
%(\ref{eq:bckraw})
%we find
%\begin{equation}
%\label{eq:matreplskraw}
%A^\flat =  \left(
%\begin{array}{ c c c c c c}
%0 & d  &      &      &   &{\bf 0} \\
%1 & 0  &  d-1   &      &   &  \\
%  & 2  &  \cdot    & \cdot  &   & \\
%  &   & \cdot     & \cdot  & \cdot   & \\
%  &   &           &  \cdot & \cdot & 1 \\
%{\bf 0} &   &   &   & d & 0  
%\end{array}
%\right).
%\end{equation} 
%Moreover
%\begin{eqnarray}
%\label{eq:asflatkraw}
%A^{*\flat} = \mbox{diag}(d, d-2, d-4, \ldots, -d)
%\end{eqnarray}
%by
%Lemma \ref{lem:firstcom}(ii).
%Apparently the Leonard pair $A,A^*$ is isomorphic
%to the Leonard pair in
%(\ref{eq:fam1}).
Pick any integers $i,j$ $(0 \leq i,j\leq d)$.
Applying 
Theorem
\ref{thm:uform}
%(\ref{eq:uifinal})
to $\Phi$
we find
\begin{eqnarray}
u_i(\theta_j)=\sum_{n=0}^d \frac{(-i)_n (-j)_n 2^n}{(-d)_n n! }, 
\label{eq:2F1expand}
\end{eqnarray}
 where
\begin{eqnarray*}
(a)_n:=a(a+1)(a+2)\cdots (a+n-1) \qquad \qquad n=0,1,2,\ldots
\end{eqnarray*}
Hypergeometric series are defined in \cite[p. 3]{GR}.
From this definition we find
the sum on the right in
(\ref{eq:2F1expand}) is the hypergeometric series
\begin{eqnarray}
{{}_2}F_1\Biggl({{-i, -j}\atop {-d}}\;\Bigg\vert \;2\Biggr).
\label{eq:2F1not}
\end{eqnarray}
A definition of the Krawtchouk polynomials can be found in 
 \cite{AAR} or 
\cite{KoeSwa}. Comparing this definition with 
(\ref{eq:2F1expand}),
(\ref{eq:2F1not})
we find the $u_i$  
are 
Krawtchouk polynomials but not the most general ones. 
%Let the scalars $x_i$  for $\Phi$ be as in
%(\ref{eq:xitr}). 
%Applying (\ref{eq:xifinform})
%to $\Phi$
% we find
%\begin{eqnarray*}
%x_i = i(d-i+1) \qquad \qquad  (1 \leq i \leq d).
%\end{eqnarray*}
Let the scalar $\nu$ for $\Phi$ be as in
Definition \ref{def:n}. 
Applying
Theorem
\ref{thm:nform}
%(\ref{eq:nfinalform})
to $\Phi$
 we find 
%\begin{eqnarray*}
$\nu=2^d. $
%\end{eqnarray*}
Let the scalars $k_i$ for $\Phi$ be as in
Definition \ref{def:ki1}.
Applying 
Theorem
\ref{thm:kiform}
%(\ref{eq:kform2}) 
to $\Phi$
we obtain  
a binomial coefficent
\begin{eqnarray*}
k_i=
\Biggl({{ d }\atop {i}}\Biggr)  \qquad \qquad (0 \leq i\leq d).
\end{eqnarray*}
Let the scalars $m_i$ for $\Phi$ be as in
Definition
\ref{def:mid}.
Applying
Theorem
\ref{thm:miform}
%(\ref{eq:mform2})
to $\Phi$
we find 
\begin{eqnarray*}
m_i = 
\Biggl({{ d }\atop {i}}\Biggr)2^{-d}  \qquad \qquad (0 \leq i \leq d).
\end{eqnarray*}

\medskip
\noindent
We now give our second example. For this
example the 
 polynomials $u_i$ will turn out to be the
$q$-Racah polynomials. To begin, 
let  $d$ denote a nonnegative integer  and consider the 
following elements in $\K$.
\begin{eqnarray}
\theta_i &=& \theta_0 + h(1-q^i)(1-sq^{i+1})/q^i,
\label{eq:thdefend}
\\
\theta^*_i &=& \theta^*_0 + h^*(1-q^i)(1-s^*q^{i+1})/q^i
\label{eq:thsdefend}
\end{eqnarray}
for $0 \leq i \leq d$, and
\begin{eqnarray}
\varphi_i &=& hh^*q^{1-2i}(1-q^i)(1-q^{i-d-1})(1-r_1q^i)(1-r_2q^i),
\label{eq:varphidefend}
\\
\phi_i &=& hh^*q^{1-2i}(1-q^i)(1-q^{i-d-1})(r_1-s^*q^i)(r_2-s^*q^i)/s^*
\label{eq:phidefend}
\end{eqnarray}
for $1 \leq i \leq d$. We assume 
$q, h, h^*, s, s^*, r_1, r_2$ are nonzero scalars
in the algebraic closure of $\K$, and that $r_1r_2 = s s^*q^{d+1}$.
To avoid degenerate situations
we assume
none of $q^i, r_1q^i, r_2q^i, s^*q^i/r_1, s^*q^i/r_2$ is
equal to 1 for $1 \leq i \leq d$
and neither of
$sq^i, s^*q^i$ is equal to 1 for $2 \leq i \leq 2d$.
It is routine to show 
(\ref{eq:thdefend})--(\ref{eq:phidefend})
satisfy the conditions PA1--PA5 of 
Definition
\ref{def:pa}, so
$(\theta_i, \theta^*_i,i=0..d; \varphi_j, \phi_j,j=1..d)$
is a parameter array over $\K$.
By 
Theorem
\ref{thm:ls}
there exists a Leonard system $\Phi$ over $\K$ with
this parameter array.
Let the scalars $b_i, c_i$ for $\Phi$ be as in
Definition
\ref{def:sharpmp}.
Applying
Theorem
\ref{thm:bcform} to $\Phi$ we find
\begin{eqnarray*}
b_0 &=& \frac{h(1-q^{-d})(1-r_1q)(1-r_2q)}
{1-s^*q^2},
\\
b_i &=& \frac{h(1-q^{i-d})(1-s^*q^{i+1})(1-r_1q^{i+1})(1-r_2q^{i+1})}
{(1-s^*q^{2i+1})(1-s^*q^{2i+2})} \qquad \quad (1 \leq i \leq d-1),
\\
c_i &=& \frac{h(1-q^i)(1-s^*q^{i+d+1})(r_1-s^*q^i)(r_2-s^*q^i)}
{s^*q^d(1-s^*q^{2i})(1-s^*q^{2i+1})} \qquad (1 \leq i \leq d-1),
\\
c_d &=& \frac{h(1-q^d)(r_1-s^*q^d)(r_2-s^*q^d)}
{s^*q^d(1-s^*q^{2d})}.
\end{eqnarray*}
Pick  integers $i,j$ $(0 \leq i,j\leq d)$.
Applying 
Theorem
\ref{thm:uform}
%(\ref{eq:uifinal})
to $\Phi$
we find
\begin{eqnarray}
u_i(\theta_j)=\sum_{n=0}^d \frac{(q^{-i};q)_n (s^*q^{i+1};q)_n 
(q^{-j};q)_n (sq^{j+1};q)_n q^n}
{(r_1q;q)_n(r_2q;q)_n (q^{-d};q)_n(q;q)_n},
\label{eq:uihyper}
\end{eqnarray}
where 
%%%%%%%%%%%%%%%%%%%%%%%%%%%%
\begin{eqnarray*}
(a;q)_n := (1-a)(1-aq)(1-aq^2)\cdots (1-aq^{n-1})\qquad \qquad n=0,1,2\ldots 
\end{eqnarray*}
Basic hypergeometric series are defined in \cite[p. 4]{GR}.
From that definition we find
the sum on the right in 
(\ref{eq:uihyper}) is the basic hypergeometric series
\begin{eqnarray}
 {}_4\phi_3 \Biggl({{q^{-i}, \;s^*q^{i+1},\;q^{-j},\;sq^{j+1}}\atop
{r_1q,\;\;r_2q,\;\;q^{-d}}}\;\Bigg\vert \; q,\;q\Biggr).
\label{eq:qrac}
\end{eqnarray}
A definition of the  $q$-Racah polynomials can be found in 
\cite{AWil} or 
\cite{KoeSwa}. Comparing this definition with 
(\ref{eq:uihyper}),
(\ref{eq:qrac})
 and
recalling $r_1r_2=s s^*q^{d+1}$,
we find  the $u_i$
are the $q$-Racah polynomials.
Let the scalar $\nu$ for $\Phi$ be as in
Definition \ref{def:n}. 
Applying 
Theorem
\ref{thm:nform}
%(\ref{eq:nfinalform})
to $\Phi$
 we find 
\begin{eqnarray*}
\nu = \frac{(sq^2;q)_d (s^*q^2;q)_d}{r^d_1q^d(sq/r_1;q)_d(s^*q/r_1;q)_d}. 
\end{eqnarray*}
Let
the scalars  $k_i$ for $\Phi$ be as in
Definition
\ref{def:ki1}.
Applying
Theorem
\ref{thm:kiform}
%(\ref{eq:kform2}) 
to $\Phi$
we obtain  
\begin{eqnarray*}
k_i = \frac{(r_1q;q)_i(r_2q;q)_i(q^{-d};q)_i(s^*q;q)_i(1-s^*q^{2i+1})}
{s^iq^i(q;q)_i(s^*q/r_1;q)_i(s^*q/r_2;q)_i(s^*q^{d+2};q)_i(1-s^*q)} 
\label{eq:qrack}
\qquad (0 \leq i \leq d).
\end{eqnarray*}
%provided $s^*q\not=1$.
Let the scalars $m_i$ for $\Phi$ be as in
Definition \ref{def:mid}.
Applying
Theorem
\ref{thm:miform}
%(\ref{eq:mform2})
to $\Phi$
we find
\begin{eqnarray*}
m_i = \frac{(r_1q;q)_i(r_2q;q)_i(q^{-d};q)_i(sq;q)_i(1-sq^{2i+1})}
{s^{*i}q^i(q;q)_i(sq/r_1;q)_i(sq/r_2;q)_i(sq^{d+2};q)_i(1-sq)\nu} 
\qquad (0 \leq i \leq d).
\end{eqnarray*}
%provided $sq\not=1$.

\section{A characterization of Leonard systems}

We are done describing the correspondence
between Leonard pairs and the terminating branch of the Askey scheme.
For the remainder of these notes we discuss applications and
related topics. We begin with a characterization of Leonard
systems.

\medskip
\noindent We recall some results from
 earlier in
the paper.
Let $\Phi$ denote the Leonard system from Definition 
\ref{eq:ourstartingpt}.
Let the polynomials
 $p_0, p_1, \ldots, p_{d+1}$
be as in Definition
\ref{def:pi} and recall 
$p^*_0, p^*_1, $ $ \ldots, p^*_{d+1} $
are the corresponding polynomials for $\Phi^*$.
For the purpose of this section, we call
 $p_0, p_1, \ldots, $ $ p_{d+1}$
the {\it monic polynomial sequence} (or {\it MPS}) of $\Phi$.
We call 
$p^*_0, p^*_1,\ldots, p^*_{d+1}$  the  {\it dual MPS} of  $\Phi$.
By Definition 
\ref{def:pi} 
we have
\begin{eqnarray}
&& \qquad p_0=1, \qquad \qquad p^*_0=1, 
\label{eq:rec1introS99}
\\
\lambda p_i &=& p_{i+1} + a_ip_i +x_ip_{i-1} \qquad \qquad (0 \leq i \leq d),
\label{eq:rec2introS99}
\\
\lambda p^*_i &=& p^*_{i+1} + a^*_ip^*_i +x^*_ip^*_{i-1} 
\qquad \qquad (0 \leq i \leq d),
\label{eq:rec3introS99}
\end{eqnarray}
where $\,x_0,\; x^*_0, \; 
p_{-1},\; 
p^*_{-1}\,$ are all zero,  and where 
\begin{eqnarray*} 
 a_i =\mbox{tr}(E^*_iA),\quad && \qquad    \quad   
a^*_i =\mbox{tr}(E_iA^*) \qquad \qquad \quad  (0 \leq i \leq d), \qquad 
\\
x_i =\mbox{tr}(E^*_iAE^*_{i-1}A),&& \qquad     
x^*_i =\mbox{tr}(E_iA^*E_{i-1}A^*) \qquad  \quad   (1 \leq i \leq d). 
\end{eqnarray*}
By 
Lemma \ref{lem:bmat}(iii) we have
\begin{equation}
x_i\not=0, \qquad \qquad 
x^*_i\not=0 \qquad \qquad (1 \leq i \leq d).
\label{eq:rec4introS99}
\end{equation}
%We call $p_0, p_1,\ldots, p_{d+1}$  the 
%{\it monic polynomial sequence} (or {\it MPS}) of $\Phi$.
%We call 
%$p^*_0, p^*_1,\ldots, p^*_{d+1}$  the  {\it dual MPS} of  $\Phi$.
Let $\theta_0,\theta_1,\ldots, \theta_d$ 
(resp. $\theta^*_0,\theta^*_1,\ldots,\theta^*_d$)
denote the eigenvalue sequence (resp. dual eigenvalue sequence)
of $\Phi$, and  recall
\begin{eqnarray}
&&\theta_i\not=\theta_j, \qquad \quad 
\theta^*_i\not=\theta^*_j \qquad \quad \hbox{if} \quad i\not=j,  
\qquad \;\;  (0 \leq i,j\leq d).
\label{eq:rec5introS99}
\end{eqnarray}
By 
Theorem
\ref{thm:pseq}(ii) we have
\begin{eqnarray}
&&\qquad p_{d+1}(\theta_i) = 0, \qquad \qquad 
p^*_{d+1}(\theta^*_i) = 0 \qquad  (0 \leq i \leq d). \qquad 
\label{eq:rec6introS99}
\end{eqnarray}
By Theorem 
\ref{def:bici3} we have
\begin{equation} 
p_i(\theta_0)\not=0, \qquad \qquad  
p^*_i(\theta^*_0)\not=0 \qquad \qquad  (0 \leq i \leq d).
\label{eq:rec7introS99}
\end{equation}
By Theorem
\ref{lem:pidual}
we have
\begin{equation}
{{p_i(\theta_j)}\over {
p_i(\theta_0)}}
= 
{{p^*_j(\theta^*_i)}\over {
p^*_j(\theta^*_0)}}\qquad \qquad (0 \leq i,j\leq d). 
\label{eq:rec8introS99}
\end{equation}

\noindent In the following theorem we show the
equations
(\ref{eq:rec1introS99})--(\ref{eq:rec8introS99})
characterize the Leonard systems.

\begin{theorem}
Let $d$ denote a nonnegative integer.
Given 
 polynomials
\begin{eqnarray}
 && p_0, p_1, \ldots, p_{d+1}
\label{eq:172L},
\\
&&p^*_0, p^*_1, \ldots, p^*_{d+1}
\label{eq:172R}
\end{eqnarray}
in 
 $\K\lbrack \lambda \rbrack $ satisfying  
(\ref{eq:rec1introS99})--(\ref{eq:rec4introS99}) 
and given   
 scalars 
\begin{eqnarray}
&&\theta_0, \theta_1, \ldots, \theta_d,
\label{eq:173L}
\\
&&\theta^*_0, \theta^*_1, \ldots, \theta^*_d 
\label{eq:173R}
\end{eqnarray}
in 
 $\K $ 
satisfying 
(\ref{eq:rec5introS99})--(\ref{eq:rec8introS99}),
 there exists  
 a Leonard system  $\Phi$ over $\K$
that has  MPS 
(\ref{eq:172L}),
 dual MPS
(\ref{eq:172R}),
  eigenvalue sequence 
(\ref{eq:173L})
and dual eigenvalue sequence (\ref{eq:173R}).
The system $\Phi $ is unique up to isomorphism of Leonard systems.   
\end{theorem}
\noindent {\it Proof:}
We abbreviate $V={\K}^{d+1}$.
Let $A$ and $A^*$ denote the following matrices in
$\mbox{Mat}_{d+1}(\K)$:
\begin{eqnarray*}
\label{eq:amonmat}
A:=\left(
\begin{array}{ c c c c c c}
a_0 & x_1  &      &      &   &{\bf 0} \\
1 & a_1  &  x_2   &      &   &  \\
  & 1  &  \cdot    & \cdot  &   & \\
  &   & \cdot     & \cdot  & \cdot   & \\
  &   &           &  \cdot & \cdot & x_d \\
{\bf 0} &   &   &   & 1 & a_d  
\end{array}
\right),
\qquad 
A^*:=\mbox{diag}(\theta^*_0, \theta^*_1, \ldots, \theta^*_d).
%\label{eq:asmonmat}
\end{eqnarray*} 
We show the pair
$A,A^*$ is a Leonard pair on $V$.
To do this we apply
Definition
\ref{def:lprecall}.
Observe that $A$ is irreducible tridiagonal and $A^*$
is diagonal. Therefore 
condition (i) of 
Definition
\ref{def:lprecall} is satisfied by the basis for
$V$ consisting of the columns of $I$, where $I$ denotes
the identity matrix in
$\mbox{Mat}_{d+1}(\K)$.
To verify condition (ii) of
Definition
\ref{def:lprecall},
we display an invertible matrix
$X$
such that $X^{-1}AX$ is diagonal and 
$X^{-1}A^*X$ is  irreducible tridiagonal.
Let 
$X$ denote the matrix in
$\mbox{Mat}_{d+1}(\K)$ 
that has entries
\begin{eqnarray}
X_{ij}&=&\frac{p_i(\theta_j)p^*_j(\theta^*_0)}
{x_1x_2\cdots x_i}
\label{eq:xver1}
\\
&=&\frac{p^*_j(\theta^*_i)p_i(\theta_0)}
{x_1x_2\cdots x_i}
\label{eq:xver2}
\end{eqnarray}
$0 \leq i ,j\leq d$. 
The matrix $X$ is invertible since it is essentially Vandermonde.
Using 
(\ref{eq:rec2introS99})
and 
(\ref{eq:xver1})
we find
$AX = XH$ where
$H=\mbox{diag}(\theta_0, \theta_1, \ldots, \theta_d)$.
Apparently $X^{-1}AX$ is equal to $H$ and is therefore diagonal.
Using (\ref{eq:rec3introS99})
 and
(\ref{eq:xver2})
we find
$A^*X=XH^*$ where
\begin{eqnarray*}
\label{eq:asmonmat}
H^*:=\left(
\begin{array}{ c c c c c c}
a^*_0 & x^*_1  &      &      &   &{\bf 0} \\
1 & a^*_1  &  x^*_2   &      &   &  \\
  & 1  &  \cdot    & \cdot  &   & \\
  &   & \cdot     & \cdot  & \cdot   & \\
  &   &           &  \cdot & \cdot & x^*_d \\
{\bf 0} &   &   &   & 1 & a^*_d  
\end{array}
\right).
\end{eqnarray*} 
Apparently $X^{-1}A^*X$ is equal to $H^*$ and is therefore irreducible
tridiagonal.
Now condition (ii) of Definition
\ref{def:lprecall} is satisfied by the basis for $V$ consisting
of the columns of $X$.
 We have now shown the pair $A,A^*$
 is a Leonard pair  on $V$.
Pick an integer $j$ $(0 \leq j \leq d)$.
Using $X^{-1}AX = H$ we find
 $\theta_j$ is the eigenvalue of 
$A$ associated with column $j$ of $X$.
From the definition of $A^*$ we find
 $\theta^*_j$ is the eigenvalue of 
$A^*$ associated with column $j$ of $I$.
Let $E_j$ (resp. $E^*_j$)
denote the primitive idempotent of $A$ (resp. $A^*$)
for $\theta_j$ (resp. $\theta^*_j$).
From our above comments the sequence
$\Phi:=(A; A^*; \lbrace E_i\rbrace_{i=0}^d; \lbrace E^*_i\rbrace_{i=0}^d)$
is a Leonard system. From the construction $\Phi$ is over $\K$.
We show
(\ref{eq:172L})
is the $MPS$
of $\Phi$.
To do this is suffices to show
 $a_i=\mbox{tr}(E^*_iA)$ for 
$0 \leq i \leq d$ and 
 $x_i=\mbox{tr}(E^*_iAE^*_{i-1}A)$ for $1 \leq i \leq d$.
Applying 
Lemma \ref{lem:bmat}(i),(ii) to $\Phi$
(with $v_i=\mbox{column  $i$ of   $I$}$, $B=A$)
we find
 $a_i=\mbox{tr}(E^*_iA)$ for 
$0 \leq i \leq d$ and 
 $x_i=\mbox{tr}(E^*_iAE^*_{i-1}A)$ for $1 \leq i \leq d$.
Therefore
(\ref{eq:172L})
is the $MPS$
of $\Phi$.
We show
(\ref{eq:172R})
is the dual $MPS$
of $\Phi$.
Applying
Lemma \ref{lem:bmat}(i),(ii) to $\Phi^*$
(with $v_i=\mbox{column  $i$ of $X$}$, $B=H^*$)
we find
 $a^*_i=\mbox{tr}(E_iA^*)$ for 
$0 \leq i \leq d$ and 
 $x^*_i=\mbox{tr}(E_iA^*E_{i-1}A^*)$ for $1 \leq i \leq d$.
Therefore 
(\ref{eq:172R})
is the dual $MPS$
of $\Phi$.
From the construction
we find (\ref{eq:173L}) (resp.  
(\ref{eq:173R}))
is the eigenvalue sequence
(resp.
dual eigenvalue sequence) of $\Phi$.
We show $\Phi$ is uniquely determined
by
(\ref{eq:172L})--(\ref{eq:173R})
up to isomorphism of Leonard systems.
Recall that $\Phi$ is determined
up to isomorphism of Leonard systems
by its own parameter array. We show the parameter array of 
$\Phi$ is determined by
(\ref{eq:172L})--(\ref{eq:173R}).
Recall the parameter array consists of
the eigenvalue sequence, the dual eigenvalue sequence,
the first split sequence and the second split sequence.
We mentioned earlier that the eigenvalue sequence of
$\Phi$ is
(\ref{eq:173L}) and the dual eigenvalue sequence
of $\Phi$ is 
(\ref{eq:173R}).
By Lemma \ref{lem:p0} the first split sequence of
$\Phi$ is determined by
(\ref{eq:172L})--(\ref{eq:173R}).
By this and
Theorem \ref{thm:xform} we find the second split sequence of 
$\Phi$ is determined by
(\ref{eq:172L})--(\ref{eq:173R}).
We have now shown the parameter array of $\Phi$ is
determined by 
(\ref{eq:172L})--(\ref{eq:173R}). We now see
that $\Phi$ 
is uniquely determined by
(\ref{eq:172L})--(\ref{eq:173R})
up to isomorphism of Leonard systems.
\hfill $\Box $ \\

\section{Leonard pairs $A,A^*$ with $A$ lower bidiagonal and
$A^*$ upper bidiagonal}

\noindent 
Let $A, A^*$ denote
matrices in $\hbox{Mat}_{d+1}(\K)$. Let us assume 
$A$ is lower bidiagonal 
and $A^*$ is upper  bidiagonal. 
We cite a necessary and sufficient condition
for $A,A^*$ to be  a Leonard pair.

\begin{theorem}
   \cite[Theorem 17.1]{conform}
\label{thm:lug}
Let $d$ denote a nonnegative integer and let $A, A^*$ denote
matrices in $\hbox{Mat}_{d+1}(\K)$. Assume $A$  
lower bidiagonal and $A^*$ is upper bidiagonal.
Then the following (i), (ii) are equivalent.
\begin{enumerate}
\item
The pair $A,A^*$ is a Leonard pair.
\item
There exists a parameter array
$(\theta_i, \theta^*_i, i=0..d;  \varphi_j, \phi_j, j=1..d)$
over $\K$ 
such that
\begin{eqnarray*}
&&A_{ii} =\theta_i,
\qquad \qquad 
A^*_{ii} =\theta^*_i \qquad \qquad (0 \leq i \leq d),
\label{eq:comb1}
\\
&&\qquad A_{i,i-1}A^*_{i-1,i} = \varphi_i \qquad \qquad (1 \leq i \leq d).
\label{eq:comb2}
\end{eqnarray*}
\end{enumerate}
\end{theorem}

\section{Leonard pairs $A,A^*$ with $A$ tridiagonal and $A^*$ diagonal}

\noindent 
Let $A, A^*$ denote
matrices in $\hbox{Mat}_{d+1}(\K)$. Let us assume 
$A$ is tridiagonal 
and $A^*$ is diagonal. 
We cite  a necessary and sufficient condition for
 $A,A^*$ to be a Leonard pair. 

\begin{theorem}
   \cite[Theorem 25.1]{conform}
\label{thm:tdcrit}
Let $d$ denote a nonnegative integer and let $A, A^*$ denote
matrices in $\hbox{Mat}_{d+1}(\K)$. Assume $A$ is
tridiagonal and $A^*$ is diagonal.
Then the following (i), (ii) are equivalent.
\begin{enumerate}
\item
The pair $A,A^*$ is a Leonard pair.
\item
There exists
a parameter array 
$(\theta_i, \theta^*_i, i=0..d;  \varphi_j, \phi_j, j=1..d)$
over $\K$
such that
\begin{eqnarray*}
A_{ii} &=&
\theta_i + \frac{\varphi_i}{\theta^*_i-\theta^*_{i-1}}
+
 \frac{\varphi_{i+1}}{\theta^*_i-\theta^*_{i+1}}
\qquad \qquad  
(0 \leq i \leq d), 
\label{eq:aiiform}
\\
A_{i,i-1}A_{i-1,i}&=&
\varphi_i\phi_i \frac{\prod_{h=0}^{i-2}(\theta^*_{i-1}-\theta^*_h)
}
{\prod_{h=0}^{i-1}(\theta^*_{i}-\theta^*_h)
}
\, \frac{\prod_{h=i+1}^d (\theta^*_i-\theta^*_h)
}
{\prod_{h=i}^d(\theta^*_{i-1}-\theta^*_h)
}
\qquad  (1 \leq i \leq d),
\label{eq:crossprod}
\\
A^*_{ii} &=&\theta^*_i \qquad \qquad (0 \leq i \leq d).
\label{eq:ths}
\end{eqnarray*}
\end{enumerate}
\end{theorem}

\section{A characterization of the parameter arrays I}

In this section we cite a characterization of the parameter
arrays in terms of bidiagonal matrices.
We will refer to the following set-up.

\begin{definition}
\label{def:setup}
Let $d$ denote a nonnegative integer and let
$(\theta_i, \theta^*_i,i=0..d; \varphi_j, \phi_j,j=1..d)
$
denote a sequence of  scalars taken from $\K$. We assume
this sequence satisfies 
PA1 and PA2.
\end{definition}

\begin{theorem}
 \cite[Theorem 3.2]{TLT:array}
\label{thm:gmat}
With reference to Definition \ref{def:setup},
 the following (i), (ii) are equivalent.
\begin{enumerate}
\item
The sequence
$(\theta_i, \theta^*_i,i=0..d;\varphi_j,\phi_j,j=1..d)$
satisfies PA3--PA5.
\item
There exists an invertible matrix 
$G \in \mbox{Mat}_{d+1}(\K)$ such that both
\begin{eqnarray*}
G^{-1}
\left(
\begin{array}{c c c c c c}
\theta_0 & & & & & {\mathbf 0} \\
1 & \theta_{1} &  & & & \\
& 1 & \theta_{2} &  & & \\
& & \cdot & \cdot &  &  \\
& & & \cdot & \cdot &  \\
{\mathbf 0}& & & & 1 & \theta_d
\end{array}
\right) G
&=&\left(
\begin{array}{c c c c c c}
\theta_d & & & & & {\mathbf 0} \\
1 & \theta_{d-1} &  & & & \\
& 1 & \theta_{d-2} &  & & \\
& & \cdot & \cdot &  &  \\
& & & \cdot & \cdot &  \\
{\mathbf 0}& & & & 1 & \theta_0
\end{array}
\right),
\label{eq:g1}
\\
G^{-1}
\left(
\begin{array}{c c c c c c}
\theta^*_0 &\varphi_1 & & & & {\mathbf 0} \\
 & \theta^*_1 & \varphi_2 & & & \\
&  & \theta^*_2 & \cdot & & \\
& &  & \cdot & \cdot &  \\
& & &  & \cdot & \varphi_d \\
{\mathbf 0}& & & &  & \theta^*_d
\end{array}
\right) G
&=&
\left(
\begin{array}{c c c c c c}
\theta^*_0 &\phi_1 & & & & {\mathbf 0} \\
 & \theta^*_1 & \phi_2 & & & \\
&  & \theta^*_2 & \cdot & & \\
& &  & \cdot & \cdot &  \\
& & &  & \cdot & \phi_d \\
{\mathbf 0}& & & &  & \theta^*_d
\end{array}
\right).
\label{eq:g2}
\end{eqnarray*}
\end{enumerate}
\end{theorem}

\section{A characterization of the parameter arrays II}

\noindent In this section we cite a characterization of the parameter
arrays in terms of polynomials.
%We will use the following notation.
%Let $\lambda $ denote an indeterminate, and let
%$\K \lbrack \lambda \rbrack $ denote the $\K$-algebra
%consisting of all polynomials in $\lambda $  which have
%coefficients in $\K$. From now on
%all polynomials which we discuss
%are assumed to lie in $\K\lbrack \lambda \rbrack $. 

\begin{theorem}
 \cite[Theorem 4.1]{TLT:array}
\label{eq:mth}
With reference to Definition \ref{def:setup},
the following (i), (ii) are equivalent.
\begin{enumerate}
\item
The sequence
 $(\theta_i, \theta^*_i,i=0..d; \varphi_j, \phi_j,j=1..d)$
satisfies PA3--PA5.
\item For $0 \leq i \leq d$ the 
polynomial 
\begin{equation}
\sum_{n=0}^i \frac{
(\lambda-\theta_0)
(\lambda-\theta_1) \cdots
(\lambda-\theta_{n-1})
(\theta^*_i-\theta^*_0)
(\theta^*_i-\theta^*_1) \cdots
(\theta^*_i-\theta^*_{n-1})
}
{\varphi_1\varphi_2\cdots \varphi_n}
\label{eq:poly1}
\end{equation}
is a scalar multiple of the polynomial 
\begin{eqnarray*}
\sum_{n=0}^i \frac{
(\lambda-\theta_d)
(\lambda-\theta_{d-1}) \cdots
(\lambda-\theta_{d-n+1})
(\theta^*_i-\theta^*_0)
(\theta^*_i-\theta^*_1) \cdots
(\theta^*_i-\theta^*_{n-1})
}
{\phi_1\phi_2\cdots \phi_n}.
\label{eq:poly2}
\end{eqnarray*}
\end{enumerate}
\end{theorem}

%\section{Some orthogonal polynomials of the Askey scheme}
%
%There is a natural correspondence between 
%Leonard systems and a class of orthogonal polynomials 
%consisting of the $q$-Racah polynomials and some related
%polynomials of the Askey scheme. This correspondence
%is described as follows
% \cite{TLT:array}.
%
%\medskip
%\noindent
%Let 
%$\Phi$
%denote a Leonard system.  
%By the {\it polynomials which correspond to $\Phi$} we mean 
%the polynomials 
%from  
%(\ref{eq:poly1}), where the parameter array from
%that line is the one associated with 
%$\Phi$.
%The polynomials which correspond to a Leonard system
%are listed below:
%
%\medskip
%\noindent
%$q$-Racah, $q$-Hahn, dual $q$-Hahn,
%$q$-Krawtchouk,
%dual $q$-Krawtchouk,
%quantum $q$-Krawtchouk, 
%affine $q$-Krawtchouk, 
%Racah, Hahn, dual-Hahn, Krawtchouk, Bannai/Ito, and  orphan polynomials.
%
%\medskip
%\noindent
%The Bannai/Ito polynomials can be obtained from the $q$-Racah polynomials by letting $q$ 
%tend to $-1$. The orphan polynomials have maximal degree 3 and 
%exist for $\mbox{char}(\K)=2$  only. 
%
%\medskip
%\noindent See 
% \cite{KoeSwa} for information on the Askey scheme.

\section{The Askey-Wilson relations}

We turn our attention to the representation theoretic aspects
of Leonard pairs.

\begin{theorem}
\cite[Theorem 1.5]{aw}
\label{lptheorem}
Let $V$ denote a vector space over $\K$ with finite positive
dimension. Let $A,A^*$ denote
a Leonard pair on $V$. Then there
exists a sequence of scalars
$\beta,\gamma,\gamma^*,\varrho,\varrho^*$, $\omega,\eta,\eta^*$
taken from $\K$ such that both
\begin{eqnarray}  \label{askwil1}
A^2 A^*-\beta A A^*\!A+A^*\!A^2-\gamma\left( A A^*\!+\!A^*\!A
\right)-\varrho\,A^* &=& \gamma^*\!A^2+\omega A+\eta\,I,\\
\label{askwil2} A^*{}^2\!A-\beta A^*\!AA^*\!+AA^*{}^2-
\gamma^*\!\left(A^*\!A\!+\!A A^*\right)-\varrho^*\!A &=&
\gamma A^*{}^2+\omega A^*\!+\eta^*I.
\end{eqnarray}
The sequence is uniquely determined by the pair $A,A^*$ provided
the diameter $d\geq 3$.
\end{theorem}

\noindent We refer to 
(\ref{askwil1}), 
(\ref{askwil2})  as the {\it Askey-Wilson relations}.
As far as we know these relations first appeared in
\cite{Zhidd}.

\medskip
\noindent 
Our next result is a kind of converse to Theorem 
 \ref{lptheorem}.
 
% In order to state the result we make a definition.
%\begin{definition}
%\label{def:awalg}
%\rm
%Given a sequence of scalars
%$\beta,\gamma,\gamma^*,\varrho,\varrho^*, \omega,\eta,\eta^*$
%taken from $\K$, we let $W$ denote the unital associative $\K$-algebra
%generated by two symbols $A,A^*$ subject to the relations
%(\ref{askwil1}), 
%(\ref{askwil2}). We call $W$ the {\it Askey-Wilson algebra}.
%We refer to $A$ and $A^*$ as the {\it standard generators} of $W$.
%\end{definition}

\begin{theorem}
\label{th:awgiveslp}
\cite[Theorem 6.2]{aw}
Given a sequence of scalars
$\beta,\gamma,\gamma^*,\varrho,\varrho^*, \omega,\eta,\eta^*$
taken from $\K$, we let $A_w$ denote the unital associative $\K$-algebra
generated by two symbols $A,A^*$ subject to the relations
(\ref{askwil1}), 
(\ref{askwil2}). 
Let $V$ denote 
 a finite dimensional irreducible $A_w$-module and assume
 each of $A,A^*$ is multiplicity-free on $V$.
Then $A,A^*$ act on $V$ as a Leonard pair provided $q$ is not
a root of unity, where $q+q^{-1}=\beta$.
\end{theorem}

\noindent
The algebra $A_w$ in
Theorem
\ref{th:awgiveslp}
is called the {\it Askey-Wilson algebra}
\cite{Zhidd}.
%\cite{GYZnature},
% \cite{GYLZmut},
%\cite{GYZTwisted},
%\cite{GYZlinear},
%\cite{GYZspherical},
%\cite{Zhidd}, 
%\cite{ZheCart},
%\cite{Zhidden}.

\medskip
\noindent We finish this section with an open problem.

\begin{problem}
\label{pr:awsolve}
\rm
Let $A_w$ denote the Askey-Wilson algebra from Theorem
\ref{th:awgiveslp}.
Let $V$ denote an irreducible $A_w$-module with either finite
or countably infinite dimension.
 We say $V$ has {\it polynomial type} whenever
there exists a basis $v_0, v_1, \ldots$ for $V$
with respect to which
the matrix representing $A$ is irreducible tridiagonal
and the matrix representing $A^*$ is diagonal.
Determine up to isomorphism
the irreducible $A_w$-modules of polynomial type.
We expect that the solutions correspond to the entire Askey scheme
of orthogonal polynomials.  
\end{problem}

\begin{remark}
\rm
The papers 
\cite{GYZnature},
 \cite{GYLZmut},
\cite{GYZTwisted},
\cite{GYZlinear},
\cite{GYZspherical},
\cite{qSerre},
\cite{Zhidd}, 
\cite{ZheCart},
\cite{Zhidden}
contain some results related to Problem
\ref{pr:awsolve},  
but a
complete and rigorous treatment
has yet to be carried out.
See also the work of 
Gr\"unbaum and Haine on the 
 ``bispectral problem''
\cite{GH4},
\cite{GH5},
\cite{GH7},
\cite{GH6},
\cite{GH1}, 
\cite{GH3},
\cite{GH2} 
as well as
\cite{atak},
\cite{atak2},
\cite{atak3},
\cite[p. 263]{BanIto},
\cite{Marco},
\cite{noumi1}, \cite{noumi2},
\cite{noumi3}, \cite{noumi4}.
\end{remark}

\begin{remark}
\rm 
Referring to Theorem
\ref{th:awgiveslp}, for the special case 
$\beta=q+q^{-1}$, $\gamma=\gamma^*=0$, $\omega=0$, $\eta=\eta^*=0$
the Askey-Wilson algebra
is related to the quantum groups $U_q(\mbox{su}_2)$, $U_q(\mbox{so}_3)$
\cite[Theorem 8.10]{ciccoli}, \cite{fairlie},
\cite{havlicek}, \cite{odesskii}
as well as the bipartite 2-homogeneous distance-regular graphs
\cite[Lemma 3.3]{curt2hom},
\cite{CurNom}, \cite[p.~427]{go}.
\end{remark}

%\begin{theorem}
%\cite{aw}
%Let $V$ denote a vector space over $\K$ with finite positive
%dimension. Let $A:V\to V$ and $A^*:V\to V$ denote linear
%transformations. Suppose that:
%\begin{itemize}
%\item There exists a sequence of scalars
%$\beta,\gamma,\gamma^*,\varrho,\varrho^*,\omega,\eta,\eta^*$ taken
%from $\K$ which satisfies {\rm (\ref{askwil1})}, {\rm
%(\ref{askwil2})}.
%\item $q$ is not a root of unity, where $q+q^{-1}=\beta$.
%\item Each of $A$ and $A^*$ is multiplicity-free.
%\item There does not exist a subspace $W\subseteq V$ such that
%$W\not=0$,
%$W\not=V$,
%$AW\subseteq W$,
%$A^*W\subseteq W$.
%\end{itemize}
%Then $A,A^*$ is a Leonard pair on $V$.
%\end{theorem}

\section{Leonard pairs and the Lie algebra $sl_2$}

\medskip
\noindent
In this section we assume
the field $\K$ is algebraically closed 
with characteristic zero.

\medskip
\noindent We recall the Lie algebra $sl_2=sl_2(\K)$. 
This algebra has a basis  $e,f,h$ satisfying
\begin{eqnarray*}
&&\lbrack h,e \rbrack  = 2e, \qquad  
\lbrack h,f \rbrack  = -2f, 
\qquad
\lbrack e,f \rbrack  = h,
\end{eqnarray*}
where $\lbrack\,,\, \rbrack $ denotes the Lie bracket.

\medskip
\noindent 
We recall the irreducible finite dimensional modules for $sl_2$.
\begin{lemma} \cite[p. 102]{Kassel} 
\label{lemvd}
There exists a family 
\begin{equation}
\label{sl2modlist}
V_d  \qquad \qquad d = 0,1,2\ldots
\end{equation}
of irreducible finite dimensional $sl_2$-modules with the following
properties. The module $V_d$ has a basis $v_0, v_1, \ldots, v_d$
satisfying $h v_i =(d-2i)v_i$ for $0 \leq i \leq d$,
$fv_i = (i+1)v_{i+1}$ for $0 \leq i \leq d-1$, $fv_d=0$, 
$ev_i = (d-i+1)v_{i-1}$ for $1 \leq i \leq d$,  $ev_0=0$.
Every irreducible finite dimensional $sl_2$-module
is isomorphic to exactly one of the modules 
in line (\ref{sl2modlist}).
\end{lemma}

\begin{example}
\label{ex:sl2}
Let $A$ and $A^*$ denote the following
elements of $sl_2$.
$$
A = e+f,  \qquad \qquad A^*= h.
$$
Let $d$ denote a nonnegative integer and consider
the action of $A$, $A^*$ on the module $V_d$.
With respect to the  basis $v_0, v_1, \ldots, v_d$ from
Lemma
\ref{lemvd},
the matrices representing $A$ and $A^*$ are
\begin{eqnarray*}
A: \; 
\left(
\begin{array}{ c c c c c c}
0 & d  &      &      &   &{\bf 0} \\
1 & 0  &  d-1   &      &   &  \\
  & 2  &  \cdot    & \cdot  &   & \\
    &   & \cdot     & \cdot  & \cdot   & \\
      &   &           &  \cdot & \cdot & 1 \\
      {\bf 0} &   &   &   & d & 0
      \end{array}
      \right),
      \qquad A^*:\; \hbox{diag}(d, d-2, d-4, \ldots, -d).
      \end{eqnarray*}
      The pair $A, A^*$ acts on
      $V_d$ as a Leonard pair. 
      The resulting
      Leonard pair is isomorphic to the one in
(\ref{eq:fam1}).
      \end{example}
      
      \noindent
      The Leonard pairs in
      Example \ref{ex:sl2}
      are not the
      only ones associated with $sl_2$.
      To get more Leonard pairs
       we replace $A$ and $A^*$ by more general
       elements in $sl_2$.
       Our result is the following.

\begin{theorem} 
\label{ex:slgen}
\cite[Ex. 1.5]{TD00}
Let $A$ and $A^*$ denote semi-simple 
elements in $sl_2$ and assume $sl_2$ is generated by
these elements. 
Let $V$ denote an irreducible finite dimensional
 module for $sl_2$. Then
the pair $A, A^*$ acts on $V$ as a Leonard  pair.
\end{theorem}

\noindent We remark the Leonard pairs in Theorem 
\ref{ex:slgen} correspond to the  Krawtchouk
polynomials
 \cite{KoeSwa}.

\section{Leonard pairs and the quantum algebra $U_q(sl_2)$}

\medskip
\noindent In this section we assume $\K$ is algebraically closed.
We fix a nonzero scalar $q\in \K$ that is not a root of unity. 
We recall the quantum algebra
 $U_q(sl_2)$.

\begin{definition} 
\cite[p.122]{Kassel}
\label{def:uqsl2}
Let $U_{q}(sl_2)$ denote the unital associative $\K$-algebra 
with generators  $ e,  f, k, k^{-1}$ and 
 relations
\begin{eqnarray*}
kk^{-1} = k^{-1}k= 1,
\end{eqnarray*}
\begin{eqnarray*}
 ke = q^2 ek,\qquad\qquad  kf = q^{-2}fk,
\end{eqnarray*}
\begin{eqnarray*}
ef - fe  = {{k-k^{-1}}\over {q - q^{-1}}}.
\end{eqnarray*}

\end{definition}
\noindent We recall the irreducible finite dimensional modules
for 
 $U_{q}(sl_2)$. We use the following notation.
\begin{eqnarray*}
\lbrack n \rbrack_q = 
{{q^n - q^{-n}}\over {q-q^{-1}}} \qquad 
\qquad n \in \Z .
\end{eqnarray*}

\begin{lemma} \cite[p. 128]{Kassel}
\label{lem:uqmods}
With reference to Definition \ref{def:uqsl2},
there exists a family 
\begin{eqnarray}
V_{\varepsilon,d} \qquad \quad 
\varepsilon \in \lbrace 1,-1\rbrace, \qquad \quad  d = 0,1,2\ldots
\label{eq:uqmods}
\end{eqnarray}
of 
 irreducible finite dimensional 
 $U_{q}(sl_2)$-modules
 with the following properties.
The module 
$V_{\varepsilon,d}$ has a basis $u_0, u_1, \ldots, u_d$
satisfying $ku_i=\varepsilon q^{d-2i}u_i$ for
$0 \leq i \leq d$, $fu_i=\lbrack i+1\rbrack_q u_{i+1}$
for $0 \leq i \leq d-1$, $fu_d=0$, 
 $eu_i=\varepsilon\lbrack d-i+1\rbrack_q u_{i-1}$
for $1 \leq i \leq d$, $eu_0=0$. 
Every irreducible finite dimensional
 $U_{q}(sl_2)$-module is isomorphic to exactly one of  
the modules $V_{\varepsilon,d}$.
(Referring to line
(\ref{eq:uqmods}), if 
$\K$ has characteristic 2
we interpret the set $\lbrace 1,-1 \rbrace $ as having 
a single element.) 
\end{lemma}

%\medskip
%\noindent
%We use the above $U_q(sl_2)$-modules to get Leonard pairs.
%The following result was proved by the present author in
%\cite{Terint}
%and is implicit in the results of Koelink and Van der Jeugt
%\cite{cite38}, \cite{cite39}. For related results see
%\cite{cite37},
%\cite{cite40},
%\cite{Rosengren},
%\cite{qSerre}. 

\begin{theorem}
\cite{cite38}, \cite{cite39},
\cite{Terint}
\label{ex:uqex}
Referring to Definition
\ref{def:uqsl2} and Lemma 
\ref{lem:uqmods},
let  $\alpha, \beta $ denote nonzero scalars in  $\K $ and 
 define $A$, $A^*$ as follows.
\begin{eqnarray*}
A= \alpha f + {{ k}\over {q-q^{-1}}},
\qquad \qquad 
A^*=\beta e + {{k^{-1}}\over {q-q^{-1}}}.
\label{eq:abdef}
\end{eqnarray*}
Let $d$ denote a nonnegative integer and choose 
$\varepsilon \in \lbrace 1, -1 \rbrace $.
 Then the pair $A, A^*$ acts on $V_{\varepsilon,d}$ as a Leonard pair
provided $\varepsilon \alpha \beta $ is not among
$q^{d-1}, q^{d-3},\ldots, q^{1-d}$. 
\end{theorem}

\noindent 
We remark  
the Leonard pairs in 
Theorem
\ref{ex:uqex}
correspond to the quantum $q$-Krawtchouk
polynomials \cite{KoeSwa}, 
\cite{cite37}. 

\section{Leonard pairs in  combinatorics}

Leonard pairs arise in many branches of combinatorics.
For instance they arise in
the theory of partially ordered sets (posets).
We illustrate this with a poset called the  
subspace lattice $L_n(q)$.

\medskip
\noindent 
In this section we assume our field $\K$ is the field $\C$ of
complex numbers.

\medskip
\noindent To define the subspace lattice
 we introduce a second field. 
Let $GF(q)$ denote a finite field of order $q$.
Let $n$ denote a positive integer
and let $W$ denote an $n$-dimensional vector space over $GF(q)$.
Let $P$ denote the set consisting of all subspaces of $W$.
The set $P$, together with the containment relation, is a poset
called $L_n(q)$. 

\medskip
\noindent
Using $L_n(q)$ we obtain a family of Leonard pairs as follows.
Let  $\C P$  denote
the vector space over $\C$ consisting of all formal $\C$-linear combinations
of elements of $P$. 
We observe $P$ is a basis for  
 $\C P$ so
the dimension of $\C P$ is equal to the cardinality of P.

\medskip
\noindent
We define  three linear transformations on $\C P$. We call these
$K$, $R$ (for ``raising''), $L$ (for ``lowering'').

\medskip
\noindent
We begin with $K$. For all $x \in P$,
\begin{eqnarray*} 
	  K x =   q^{{n/2} - dim\, x} x.
	 % \qquad \qquad    (x \in P).
\end{eqnarray*}
Apparently each element of $P$ is an eigenvector for $K$.

\medskip
\noindent
To define $R$ and $L$ we use the following notation. 
For $x, y \in P$ we say $y$ {\it covers}  $x$
whenever (i) $x \subseteq y$ and  (ii) $\hbox{dim}\, y= 1+\hbox{dim}\,x$.

\medskip
\noindent The maps 
$R$ and $L$ are defined as follows. For all $x \in  P$, 
\begin{eqnarray*}
       %    R x = \sum_{{y \in P}\atop{y \;covers\; x}} y.
           R x = \sum_{y \;covers \; x} y.
\end{eqnarray*}  
%where the sum is over all elements y of P that cover x.
Similarly
%for all $x \in P$, 
 \begin{eqnarray*} 
         %  L x =q^{(1-n)/2} \sum_{{y \in P}\atop{x \;\hbox{covers}\; y}} y.
           L x =q^{(1-n)/2} \sum_{x \; covers \; y} y.
\end{eqnarray*}
%where the sum is over all elements y of P that are covered by x.
(The scalar $q^{(1-n)/2}$ is included for aesthetic reasons.)

\medskip
\noindent We consider the properties of $K, R, L$.
From the construction we find $K^{-1}$ exists.
By combinatorial counting we  verify 
\begin{eqnarray*}
     &&     K L = q L K,   \qquad \qquad \qquad   K R = q^{-1} R K,
     \\
       && \qquad          L R - R L = \frac{K - K^{-1}}{q^{1/2} - q^{-1/2}}.
\end{eqnarray*}
We recognize these equations. They are the defining relations for 
$U_{q^{1/2}}(sl_2)$.  Apparently $K$, $R$, $L$ turn  $\C P$ into
a module for $U_{q^{1/2}}(sl_2)$. 

\medskip
\noindent
We now see how to get Leonard pairs from $L_n(q)$.
Let  $\alpha, \beta$ denote nonzero complex scalars and 
define $A$, $A^*$ as follows.
\begin{eqnarray*}
                  A =  \alpha R    + \frac{K}{q^{1/2}-q^{-1/2}},
        \qquad \qquad  A^* =  \beta  L    + \frac{K^{-1}}{q^{1/2}-q^{-1/2}}.
\end{eqnarray*}
To avoid 
degenerate situations we assume  
$\alpha \beta $  is not among $q^{(n-1)/2}, q^{(n-3)/2},\ldots, q^{(1-n)/2}$.

\medskip
\noindent The  $U_{q^{1/2}}(sl_2)$-module $\C P$ is 
completely reducible \cite[p. 144]{Kassel}.
In other words $\C P$ 
is a direct sum of irreducible
$U_{q^{1/2}}(sl_2)$-modules. On
each irreducible module in this sum
the pair $A, A^*$ acts as a Leonard pair. This follows from
Theorem 
\ref{ex:uqex}.

\medskip
\noindent  
We just saw how the subspace lattice gives Leonard pairs.
We expect that some other classical posets,
such as
the polar spaces and attenuated spaces
\cite{uniform},
give
Leonard pairs in a similar fashion.
However
the details remain to be worked out. See
\cite{uniform}
for more information on this topic.

\medskip
\noindent 
Another combinatorial object that gives Leonard pairs is
 a  $P$- and $Q$-polynomial
association scheme \cite{BanIto}, \cite{bcn}, 
 \cite{TersubI}.
Leonard pairs have been used to describe
certain irreducible modules for the subconstituent algebra of
 these schemes
\cite{Cau}, \cite{CurNom}, \cite{Curspin}, 
\cite{TD00},
 \cite{TersubI}.

\section{Tridiagonal pairs}
\noindent
There is a mild generalization of a Leonard
pair called a 
{\it tridiagonal pair}
\cite{TD00}, 
\cite{shape},
\cite{tdanduq},
\cite{qSerre}.
In order to define this, we use the following terms. 
 Let $V$  denote
 a vector space over $\K$ with finite positive dimension.
  Let $A:V\rightarrow V$ denote a linear transformation and
   let $W$ denote a subspace of $V$. 
   We call $W$ an {\it eigenspace} of $A$ whenever
   $W\not=0$ and there exists $\theta \in \K$ such that
   \begin{eqnarray*}
   W=\lbrace v \in V \;\vert \;Av = \theta v\rbrace.
   \end{eqnarray*}
   We say $A$ is {\it diagonalizable} whenever
   $V$ is spanned by the eigenspaces of $A$.

   \begin{definition}
   \cite[Definition 1.1]{TD00}
   \label{def:tdp}
   \rm
   Let $V$ denote a vector space over $\K$ with finite positive dimension.
   By a {\it tridiagonal pair}  on $V$,
   we mean an ordered pair of linear transformations
   $A:V\rightarrow V$ and
   $A^*:V\rightarrow V$
     that satisfy
    the following four conditions.
    \begin{enumerate}
    \item Each of $A,A^*$ is diagonalizable.
    \item There exists an ordering $V_0, V_1,\ldots, V_d$ of the
    eigenspaces of $A$ such that
    \begin{equation}
    A^* V_i \subseteq V_{i-1} + V_i+ V_{i+1} \qquad \qquad (0 \leq i \leq d),
    \label{eq:t1}
    \end{equation}
    where $V_{-1} = 0$, $V_{d+1}= 0$.
    \item There exists an ordering $V^*_0, V^*_1,\ldots, V^*_\delta$ of
    the
    eigenspaces of $A^*$ such that
    \begin{equation}
    A V^*_i \subseteq V^*_{i-1} + V^*_i+ V^*_{i+1} \qquad \qquad
    (0 \leq i \leq 
    \delta),
    \label{eq:t2}
    \end{equation}
    where $V^*_{-1} = 0$, $V^*_{\delta+1}= 0$.
    \item There does not exist a 
    subspace $W$ of $V$ such  that $AW\subseteq W$,
    $A^*W\subseteq W$, $W\not=0$, $W\not=V$.
    \end{enumerate}
    \end{definition}

\noindent The following problem is open.

\begin{problem}
\label{pr:classtd}
\rm
Classify the tridiagonal pairs.
\end{problem}

\noindent For the rest of this section we 
discuss what is known about tridiagonal pairs,
and give some conjectures.

\medskip
\noindent 
Let $A,A^*$ denote a tridiagonal pair on $V$
and let the integers $d$, $\delta$ be as in
Definition
   \ref{def:tdp}(ii), (iii) respectively.
By \cite[Lemma 4.5]{TD00}
we have $d=\delta$; we call this common value
the {\it diameter} of the pair.
An ordering of the eigenspaces of $A$ (resp. $A^*$)
will be called
{\it standard} whenever it satisfies
    (\ref{eq:t1})
   (resp.  (\ref{eq:t2})).
We comment on the uniqueness of the standard ordering.
Let
 $V_0, V_1,\ldots, V_d$ denote a standard ordering
of the eigenspaces of $A$.
Then the ordering
 $V_d, V_{d-1},\ldots, V_0$
 is standard and no other ordering is standard.
 A similar result holds for the
 eigenspaces of $A^*$.
 Let
 $V_0, V_1,\ldots, V_d$
 (resp.  $V^*_0, V^*_1,\ldots, V^*_d$)
 denote a standard ordering
 of the eigenspaces of $A$ (resp. $A^*$).
 By \cite[Corollary 5.7]{TD00},  for $0 \leq i \leq d$
the spaces $V_i$,  $V^*_i$ have the same dimension;
we denote this common dimension by $\rho_i$.
 By the construction $\rho_i\not=0$.
 By  \cite[Corollary 5.7]{TD00}
 and \cite[Corollary 6.6]{TD00},
  the sequence $\rho_0, \rho_1, \ldots, \rho_d$
 is symmetric and unimodal; that is
$ \rho_i =\rho_{d-i}$ for $0 \leq i \leq d$
 and
 $\rho_{i-1} \leq \rho_{i}$ for  $1 \leq i \leq d/2$.
 We refer to the sequence
 $(\rho_0, \rho_1, \ldots, \rho_d)$ as the {\it shape vector} of $A,A^*$.
      A Leonard pair is the same thing as a tridiagonal pair
      that has shape vector $(1,1,\ldots, 1)$.

\begin{conjecture}
\label{conj:shape}
\cite[Conjecture 13.5]{TD00}
Referring to Definition 
\ref{def:tdp}, 
assume $\K$ is algebraically closed 
and let 
$(\rho_0, \rho_1, \ldots, \rho_d)$
denote the 
shape vector for $A,A^*$.
Then the entries in this shape vector are bounded
above by binomial coefficients as follows:
\begin{eqnarray*}
\rho_i \;\leq \;\biggl( {{d}\atop {i}} \biggr) \qquad \qquad (0 \leq i \leq d).
\end{eqnarray*}
\end{conjecture}

\noindent 
See \cite{shape} for some partial results on
Conjecture
\ref{conj:shape}.
We now give some examples of tridiagonal
pairs.

\begin{example}
\label{ex:onsag}
\cite[Example 1.6]{TD00}
\rm
Assume $\K$ is algebraically closed with characteristic
0. Let $b, b^*$ denote
nonzero scalars in $\K$. Let $O$ denote the Lie
algebra over $\K$ generated by symbols
$A,A^*$ subject to the relations
\begin{eqnarray}
\lbrack A,\lbrack A,\lbrack A,A^*\rbrack \rbrack\rbrack 
&=& b^2\lbrack A,A^* \rbrack,
\label{eq:td1dg}
 \label{eq:DG1}
\\
\lbrack A^*,\lbrack A^*,\lbrack A^*,A\rbrack \rbrack\rbrack 
&=& b^{*2}\lbrack A^*,A \rbrack.
\label{eq:td2dg}
 \label{eq:DG2}
\end{eqnarray}
Let $V$ denote a finite dimensional irreducible
$O$-module. Then $A,A^*$ act on $V$
as a tridiagonal pair.  
\end{example}

\begin{remark}
\rm
The algebra $O$ from 
Example
\ref{ex:onsag}
is called the 
 {\it Onsager algebra}.
%\cite{Davfirst}.
It
first appeared in the seminal paper
by  Onsager \cite{Onsager} in which the
free energy of the two dimensional Ising
model was computed exactly.
Onsager presented his algebra by displaying a basis;
the above presentation 
using generators and relations
(\ref{eq:td1dg}), 
(\ref{eq:td2dg}) 
was  established by Perk
\cite{perk}. The relations themselves
first appeared 
in work of Dolan and Grady \cite{Dolgra}.
A few years later they were used 
by von Gehlen and Rittenberg
\cite{vongehlen}
to describe
the superintegrable chiral Potts model.
In
%Some 47 years after its discovery and upon
%the invention of Kac-Moody algebras, 
\cite{roan}
Roan
observed that  
$O$ is isomorphic
to the invariant subalgebra of the loop algebra
$\K\lbrack t,t^{-1}\rbrack \otimes sl_2$ 
by an involution.
Of course this last
result was not available to Onsager since
his discovery predates the invention
of 
Kac-Moody algebras by some 25 years.
See 
\cite{CKOnsn},
\cite{Albert},
\cite{agmpy},
\cite{auyangandperk},
\cite{auyang2},
\cite{Baz},
\cite{Davfirst},
\cite{McCoyd},
%\cite{vongehlen},
\cite{McCoy1},
\cite{uglov1}
for recent work involving the
Onsager algebra and integrable lattice models.
The equations
(\ref{eq:td1dg}), 
(\ref{eq:td2dg}) are called the {\it Dolan-Grady} relations
\cite{DateRoan2},
\cite{Dav},
%\cite{Dolgra},
\cite{Klish1},
\cite{Klish2},
\cite{Klish3},
\cite{Lee}.
\end{remark}

%\begin{remark}
%\rm
%The algebra $O$ from 
%Example
%\ref{ex:onsag}
%is called the 
% {\it Onsager algebra}
%\cite{Davfirst}. It
%first appeared in the seminal paper
%by  Onsager \cite{Onsager} in which the
%free energy of the two dimensional Ising
%model was computed exactly.
%Some 25 years after its discovery and upon
%the invention of Kac-Moody algebras, Roan \cite{roan} 
%observed that  
%$O$ is isomorphic
%to the invariant subalgebra of the loop algebra
%$\K\lbrack t,t^{-1}\rbrack \otimes sl_2$ 
%by an involution.
%See 
%\cite{CKOnsn},
%\cite{agmpy},
%\cite{auyangandperk},
%\cite{McCoyd},
%\cite{vongehlen},
%\cite{McCoy1},
%\cite{perk}
%for recent work involving the
%Onsager algebra and integrable lattice models.
%The equations
%(\ref{eq:td1dg}), 
%(\ref{eq:td2dg}) are called the {\it Dolan-Grady} relations
%\cite{Albert},
%\cite{Baz},
%\cite{DateRoan2},
%\cite{Dav},
%\cite{Dolgra},
%\cite{Klish1},
%\cite{Klish2},
%\cite{Klish3},
%\cite{Lee}.
%\end{remark}

\begin{example} \cite[Example 1.7]{TD00}
\label{ex:qserre}
\rm
Assume $\K$ is algebraically closed, and let
$q$ denote a nonzero scalar in $\K$ that is not a root of
unity. Let 
$U_q({\widehat{sl}}_2)^{>0}$
 denote the unital associative $\K$-algebra
generated by symbols $A,A^*$ 
subject to the relations
\begin{eqnarray}
0 &=&A^3A^*-\lbrack 3 \rbrack_q A^2A^*A +\lbrack 3 \rbrack_q AA^*A^2
- A^*A^3,
\label{eq:td1s} 
     \label{eq:qserre1}
\\
\label{eq:td2s}
     \label{eq:qserre2}
0 &=&A^{*3}A-\lbrack 3 \rbrack_q A^{*2}AA^* +\lbrack 3 \rbrack_q A^*AA^{*2}
- AA^{*3}.
\end{eqnarray}
Let $V$ denote a finite dimensional irreducible 
$U_q({\widehat{sl}}_2)^{>0}$-module
and assume neither of $A,A^*$ is nilpotent on $V$.
Then $A,A^*$ act on $V$ as a tridiagonal pair.
\end{example}

\begin{remark}
\rm
The equations
(\ref{eq:td1s}),  
(\ref{eq:td2s}) are known as the {\it $q$-Serre relations},
and are among
the defining relations for the quantum affine algebra
$U_q({\widehat{sl}}_2)$
\cite{charp}, \cite{uglov2}.
The algebra $U_q({\widehat{sl}}_2)^{>0}$
is called the 
{\it positive part} of
$U_q({\widehat{sl}}_2)$.
The tridiagonal pairs from
Example
\ref{ex:qserre}
are said to have {\it $q$-geometric type}.
%See 
%\cite{shape},
%\cite{tdanduq} for more on the connection
%to 
%$U_q({\widehat{sl}}_2)$. 
\end{remark}

\noindent In order to get the most general tridiagonal
pairs, we consider a pair of relations that
generalize both the
Dolan-Grady relations and the $q$-Serre relations.
We call these the {\it tridiagonal relations}.
These relations are given as follows.

\begin{theorem} \cite[Theorem 10.1]{TD00}
\label{thm:tdrels}
Let $V$ denote a vector space over $\K$ with finite positive
dimension and let
 $A,A^*$ denote a tridiagonal pair on $V$.
Then
there exists a sequence of scalars
$\beta, \gamma, \gamma^*, \varrho, \varrho^*$ taken from $\K$
such that both
\begin{eqnarray}
0 &=&\lbrack A,A^2A^*-\beta AA^*A+A^*A^2-\gamma (AA^*+A^*A)-\varrho A^*\rbrack
\label{eq:td1}
\label{eq:TD1}
\\
0 &=&\lbrack A^*,A^{*2}A-\beta A^*AA^*+AA^{*2}-\gamma^*(AA^*+A^*A)-\varrho^*
A\rbrack,
\label{eq:td2}
\label{eq:TD2}
\end{eqnarray}
where $\lbrack r,s\rbrack$ means $rs-sr$. The sequence is unique if
the diameter $d\geq 3$.
\end{theorem}

\noindent 
We call 
(\ref{eq:td1}),
(\ref{eq:td2}) the {\it tridiagonal relations} \cite{qSerre}.
As far as we know these relations first appeared in
\cite[Lemma 5.4]{tersub3}.

\begin{remark}
\label{rem:dg}
 \rm
 The Dolan-Grady relations (\ref{eq:DG1}), (\ref{eq:DG2}) are
   the tridiagonal relations with parameters
   $\beta  = 2$, $\gamma = \gamma^*=0$,
   $\varrho = b^2$,
   $\varrho^* = b^{*2}$, if we interpret the bracket
    in (\ref{eq:DG1}), (\ref{eq:DG2}) as
     $\lbrack r,s\rbrack = rs-sr$.
     The $q$-Serre relations
     (\ref{eq:qserre1}),
     (\ref{eq:qserre2})
     are the tridiagonal relations with parameters
     $\beta = q^2+q^{-2}$, $\gamma=\gamma^*=0$,
     $\varrho=\varrho^*=0$.
     \end{remark}

\noindent Our next result is a  kind of converse to
Theorem
\ref{thm:tdrels}.

     \begin{theorem}
     \cite[Theorem 3.10]{qSerre}
     \label{def:tdalg}
\rm
Let
     $\beta, \gamma, \gamma^*, \varrho, \varrho^* $ denote a sequence
     of scalars taken from $\K$.
     Let
      $T$
       denote the unital associative $\K$-algebra
      generated by symbols $A$, $A^*$
      subject to the tridiagonal relations
  (\ref{eq:TD1}),
      (\ref{eq:TD2}).
%
%      We call $T$ the {\it Tridiagonal algebra}.
%      We refer to $A $ and $A^*$ as the {\it standard generators} of $T$.
	  Let $V$ denote an irreducible finite dimensional $T$-module
	  and assume each of $A, A^*$ is diagonalizable on $V$.
	  Then $A, A^*$ act on $V$ as a tridiagonal pair provided
           $q$ is not a root of unity, where $q+q^{-1}=\beta$. 
\end{theorem}

\begin{remark}
\rm The algebra $T$ in 
  Theorem   \ref{def:tdalg} is called the
 {\it tridiagonal algebra} 
\cite{TD00},
\cite{LS99},
\cite{qSerre}.
\end{remark}

\noindent So far in our research on tridiagonal pairs, our 
strongest result concerns the case of
 $q$-geometric type.
In order to describe this result
we define one more algebra.
In what follows $\Z_4 = \Z/4\Z$ denotes the cyclic group of order 4.

\begin{definition} 
\label{def:qtet}
\rm \cite{qtet}
Let $q$ denote a nonzero element of $\K$ such that $q^2\not=1$.
Let $\boxtimes_q$ denote the unital associative $\K$-algebra that has
generators
\begin{eqnarray*}
\lbrace x_{ij}\;|\; i,j \in \Z_4,\;j-i=1 \;\mbox{or} \;j-i=2\rbrace
\end{eqnarray*}
and the following relations:
\begin{enumerate}
\item For $i,j\in \Z_4$ such that $j-i=2$,
\begin{eqnarray*}
x_{ij}x_{ji} = 1.
\label{eq:qrel0}
\end{eqnarray*}
\item For $i,j,k\in \Z_4$ such that the pair $(j-i,k-j)$ is one of
$(1,1), (1,2), (2,1)$,
\begin{eqnarray*}
\frac{qx_{ij}x_{jk}-q^{-1}x_{jk}x_{ij}}{q-q^{-1}}=1.
\label{eq:qrel1}
\end{eqnarray*}
\item For $i,j,k,\ell\in \Z_4$ such that $j-i=k-j=\ell-k=1$,
\begin{eqnarray*}
\label{eq:qserre}
x_{ij}^3x_{k\ell} -
\lbrack 3 \rbrack_q
x_{ij}^2x_{k\ell}x_{ij} +
\lbrack 3 \rbrack_q
x_{ij}x_{k\ell}x_{ij}^2-
x_{k\ell}x_{ij}^3=0.
\end{eqnarray*}
\end{enumerate}
We call $\boxtimes_q$ the {\it $q$-tetrahedron algebra}.
\end{definition}

%%%what follows is the original published version. Changed dec10 06
%\begin{definition}
%\label{conj:box}
%Assume $\K$ is algebraically closed.
%Let $q$ denote a nonzero scalar in $\K$ that is
%not a root of unity.
%We let $\boxtimes_q$ denote the unital associative $\K$-algebra
%with generators $x_i,k_i$ $(i=0,1,2,3)$
%and relations
%\begin{eqnarray*}
%&&\frac{qx_i k_{i}-q^{-1}k_{i}x_i}{q-q^{-1}}=1,
%\qquad \qquad
%\frac{qx_i x_{i+1}-q^{-1}x_{i+1}x_i}{q-q^{-1}}=1,
%\\
%&&\frac{qk_{i+1}x_i-q^{-1}x_ik_{i+1}}{q-q^{-1}}=1,
%\qquad \qquad
%k_ik_{i+2}=1,
%\\
%&&
%x_i^3x_{i+2}-
%\lbrack 3 \rbrack_q x_i^2x_{i+2}x_i
%+
%\lbrack 3 \rbrack_q x_ix_{i+2}x^2_i
%-
%x_{i+2}x^3_i=0
%\end{eqnarray*}
%for $0 \leq i \leq 3$,
%where all subscripts are computed modulo 4.
%\end{definition}
%
\begin{remark} 
\label{rem:uq}
\rm
The 
algebra $\boxtimes_q$ is closely related to
the quantum affine algebra $U_q({\widehat{sl}}_2)$;
see \cite{qtet} for the details.
%Indeed there exists a homomorphism
%of $\K$-algebras from 
%$U_q({\widehat{sl}}_2)$ into $\Box_q$
%\cite[Theorem 2.1]{tdanduq}.
%This homomorphism
%induces on each  $\Box_q$-module the structure of
%a $U_q({\widehat{sl}}_2)$-module.
%The 
%$\Box_q$-module structure is irreducible 
%if and only if 
%the 
%$U_q({\widehat{sl}}_2)$-module structure is
%irreducible.
\end{remark}

\begin{theorem}
\cite[Section 10]{qtet}
\label{thm:box}
Assume $\K$ is algebraically closed.
Let $q$ denote a nonzero scalar in $\K$
that is not a root of unity.
 Let $V$ denote a vector space over $\K$ with finite
positive dimension. 
Let $A,A^*$ denote a tridiagonal pair on
$V$ that has $q$-geometric type.
Then there exists an irreducible $\boxtimes_q$-module
structure on $V$ such that $A$ acts as a scalar
multiple of $x_{01}$ and $A^*$ acts as a scalar multiple
of $x_{23}$.
Conversely, let $V$ denote a finite dimensional
irreducible $\boxtimes_q$-module. Then
the generators $x_{01}, x_{23}$ act on $V$ as a tridiagonal
pair of $q$-geometric type.
\end{theorem}

%\begin{note}
%\rm Determine up to isomorphism the tridiagonal
%pairs of $q$-geometric type. This can probably
%be done by
%using Remark \ref{rem:uq}
%and the
%classification of finite dimensional
%irreducible $U_q({\widehat{sl}}_2)$-modules \cite{charp}.
%\end{problem}

\noindent We end this section with a conjecture.

\begin{conjecture}
\label{conj:tdgen}
Assume $\K$ is algebraically closed.
Let $V$ denote a vector space over $\K$ with finite positive
dimension and
let $A,A^*$ denote a tridiagonal pair
on $V$.
To avoid degenerate situations we assume
$q$ is not a root of unity, where
$\beta=q^2+q^{-2}$, 
and where
$\beta$ is from
Theorem
\ref{thm:tdrels}.
Then referring to Definition
\ref{def:qtet},
there exists an irreducible $\boxtimes_q$-module
structure on $V$ such that
$A$ acts as a linear combination of $x_{01}, x_{12}, I$
and
$A^*$ acts as a linear combination of $x_{23}, x_{30}, I$.
\end{conjecture}

%\noindent We mention some results on the eigenvalues
%of a tridiagonal pair.
%For $0 \leq i \leq d$ 
%let $\theta_i$ (resp. $\theta^*_i$) denote the eigenvalue
%of $A$ (resp. $A^*$) associated with $V_i$ (resp. $V^*_i$).
%Then both expressions
%\begin{eqnarray*}
%\frac{\theta_{i-2}-\theta_{i+1}}{\theta_{i-1}-\theta_i},
%\qquad \qquad 
%\frac{\theta^*_{i-2}-\theta^*_{i+1}}{\theta^*_{i-1}-\theta^*_i}
%\end{eqnarray*}
%equal $\beta+1$ for  $2 \leq i \leq d-1$.

%The equations
%(\ref{eq:td1dg}), 
%(\ref{eq:td2dg}) are called the {\it Dolan-Grady} relations
%\cite{CKOnsn},
%\cite{DateRoan2},
%\cite{Dav},
%\cite{Dolgra}. We note that the Lie algebra  over
%$\C$ generated by two symbols $A,A^*$ subject to
%(\ref{eq:td1dg}), 
%(\ref{eq:td2dg}) (where we interpret $\lbrack \,,\,\rbrack $
%as the Lie bracket) is infinite dimensional and is known as
%the {\it Onsager algebra}
%\cite{CKOnsn}.
%\end{remark}

%The reference \cite{} contains some basic information 
%on tridiagonal pairs. In \cite{} 
%%we show that  certain tridiagonal pairs are related to
%irreducible modules for the quantum affine algebra
%$U_q({\widehat{sl}}_2)$ \cite{}.

\section{Appendix: List of parameter arrays}

In this section we display all the parameter arrays over $\K$.
We will use the following notatation.

%\begin{definition}
%\label{def:base}
%Let
%$p=(\theta_i, \theta^*_i, i=0...d;  \varphi_j, \phi_j, j=1...d)$
%denote a parameter array over $\K$.  By a {\it base} for
%$p$, we mean a nonzero scalar $q$ in the algebraic closure
%of $\K$  such that
%$q+q^{-1}+1$ is equal to the common value of 
%(\ref{eq:betaplusone}) for $2 \leq i \leq d-1$.
%We remark on the uniqueness of the base.
%Suppose $d\geq 3$.
%If $q$ is a base for $p$ then so is $q^{-1}$ and 
%$p$ has no other base. Suppose $d<3$.
%Then any nonzero scalar in the algebraic closure of 
%$\K$ is a base for $p$.
%\end{definition}

\begin{definition}
\label{def:assocpoly}
\rm
Let
$p=(\theta_i, \theta^*_i, i=0..d;  \varphi_j, \phi_j, j=1..d)$
denote a parameter array over $\K$. 
For $0 \leq i\leq d$ we let $u_i$ denote the following polynomial
in
$\K\lbrack \lambda \rbrack$.
\begin{eqnarray}
u_i = \sum_{n=0}^i \frac{
(\lambda-\theta_0)
(\lambda-\theta_1) \cdots
(\lambda-\theta_{n-1})
(\theta^*_i-\theta^*_0)
(\theta^*_i-\theta^*_1) \cdots
(\theta^*_i-\theta^*_{n-1})
}
{\varphi_1\varphi_2\cdots \varphi_n}.
\label{eq:fipoly}
\end{eqnarray}
We call $u_0, u_1, \ldots, u_d$ {\it the polynomials that 
correspond to $p$}.
\end{definition}
\noindent 
We now display all the  parameter arrays over $\K$.
For each displayed array 
$(\theta_i, \theta^*_i, i=0..d;  \varphi_j, \phi_j, j=1..d)$
we present  $u_i(\theta_j)$
for $0 \leq i,j\leq d$, where  $u_0, u_1, \ldots, u_d$ are the corresponding
polynomials.
Our  
presentation is organized as follows.
In each of Example 
\ref{ex:pa2}--\ref{ex:orphan} below we give a family of
parameter arrays over $\K$. 
In Theorem
\ref{thm:thatsit}
 we show every parameter array over $\K$ is contained in
 at least one of these families.

\medskip
\noindent 
In each of
 Example 
\ref{ex:pa2}--\ref{ex:orphan} below 
the following implicit assumptions apply:
$d$ denotes a nonnegative integer,
the scalars
$(\theta_i, \theta^*_i, i=0..d;  \varphi_j, \phi_j, j=1..d)$
are contained in  $\K$, and the scalars $q,h,h^*\ldots $
are contained in the algebraic closure of $\K$.

\begin{example}($q$-Racah)
\label{ex:pa2}
\label{ex:qracah}
Assume
%\begin{eqnarray}
\begin{eqnarray}
\theta_i &=& \theta_0+h(1-q^i)(1-sq^{i+1})q^{-i},
\label{eq:thrac}
\\
\theta^*_i &=& \theta^*_0+h^*(1-q^i)(1-s^*q^{i+1})q^{-i}
\label{eq:thsrac}
\end{eqnarray}
for $0 \leq i \leq d$ and
\begin{eqnarray}
\varphi_i &=& hh^*q^{1-2i}(1-q^i)(1-q^{i-d-1})(1-r_1q^i)(1-r_2q^i),
\label{lem:varphiphi1}
\\
\phi_i &=& hh^*q^{1-2i}(1-q^i)(1-q^{i-d-1})(r_1-s^*q^i)(r_2-s^*q^i)/s^*
\label{lem:varphiphi2}
\end{eqnarray}
for $1 \leq i \leq d$.
Assume 
$h, h^*, q, s, s^*, r_1, r_2$ are nonzero
 and  $r_1r_2=ss^*q^{d+1}$.
Assume
none of
$q^i,
r_1q^i, r_2q^i,
s^*q^i/r_1,$ $ s^*q^i/r_2$
is equal to $1$ for $1 \leq i \leq d$ and 
that neither of $sq^i, s^*q^i$ is equal to $1$ for $2 \leq i \leq 2d$.
Then 
$(\theta_i, \theta^*_i, i=0..d;  \varphi_j, \phi_j, j=1..d)$
is a parameter array over $\K$.
The corresponding polynomials
$u_i$  satisfy
\begin{eqnarray*}
u_i(\theta_j)=
{}_4\phi_3 \Biggl({{q^{-i}, \;s^*q^{i+1},\;q^{-j},\;sq^{j+1}}\atop
{r_1q,\;\;r_2q,\;\;q^{-d}}}\;\Bigg\vert \; q,\;q\Biggr)
\label{eq:fqrac}
\end{eqnarray*}
for $0 \leq i,j\leq d$. These $u_i$ are the $q$-Racah polynomials.
\end{example}

\begin{example}
($q$-Hahn)
\label{ex:qhahn}
Assume
\begin{eqnarray*}
\theta_i &=&\theta_0+h(1-q^i)q^{-i}, 
\\
\theta^*_i &=& \theta^*_0+h^*(1-q^i)(1-s^*q^{i+1})q^{-i}
\end{eqnarray*}
for $0 \leq i \leq d$ and
\begin{eqnarray*}
\varphi_i &=& hh^*q^{1-2i}(1-q^i)(1-q^{i-d-1})(1-rq^i),
\\
\phi_i &=& -hh^*q^{1-i}(1-q^i)(1-q^{i-d-1})
(r-s^*q^i)
\end{eqnarray*}
for $1 \leq i \leq d$.
Assume 
$h, h^*, q, s^*, r$ are nonzero.
Assume none of
$q^i,
rq^i,
s^*q^i/r$
is equal to $1$ for $1 \leq i \leq d$ and 
that $s^*q^i\not=1$ for $2 \leq i \leq 2d$.
Then the sequence
$(\theta_i, \theta^*_i, i=0..d;  \varphi_j, \phi_j, j=1..d)$
is a parameter array over $\K$.
The corresponding polynomials
$u_i$
satisfy
\begin{eqnarray*}
u_i(\theta_j)=
{}_3\phi_2 \Biggl({{q^{-i},\;s^*q^{i+1},\;q^{-j}}\atop
{rq,\;\;q^{-d}}}\;\Bigg\vert \; q,\;q\Biggr)
%\label{eq:qrac}
\end{eqnarray*}
for $0 \leq i,j\leq d$. These $u_i$ are the $q$-Hahn polynomials.
\end{example}

\begin{example}
\label{ex:dual qhahn}
(Dual $q$-Hahn)
Assume
\begin{eqnarray*}
\theta_i &=& \theta_0+h(1-q^i)(1-sq^{i+1})q^{-i},
\\
\theta^*_i &=&\theta^*_0+h^*(1-q^i)q^{-i} 
%\label{lem:eigdeig}
\end{eqnarray*}
for $0 \leq i \leq d$ and
\begin{eqnarray*}
\varphi_i &=& hh^*q^{1-2i}(1-q^i)(1-q^{i-d-1})(1-rq^i),
%\label{lem:varphiphi1}
\\
\phi_i &=& hh^*q^{d+2-2i}(1-q^i)(1-q^{i-d-1})(s-rq^{i-d-1})
%\label{lem:varphiphi2}
\end{eqnarray*}
for $1 \leq i \leq d$.
Assume $h,h^*, q, r,s$  are nonzero.
Assume
none of
$q^i,
rq^i,
sq^i/r$
is equal to $1$ for $1 \leq i \leq d$ and 
that $sq^i\not=1$ for $2 \leq i \leq 2d$.
Then the sequence
$(\theta_i, \theta^*_i, i=0..d;  \varphi_j, \phi_j, j=1..d)$
is a parameter array over $\K$.
The corresponding polynomials
$u_i$
satisfy
\begin{eqnarray*}
u_i(\theta_j)=
{}_3\phi_2 \Biggl({{q^{-i},\;q^{-j},\;sq^{j+1}}\atop
{rq,\;\;q^{-d}}}\;\Bigg\vert \; q,\;q\Biggr)
%\label{eq:qrac}
\end{eqnarray*}
for $0 \leq i,j\leq d$. These $u_i$ are the dual $q$-Hahn polynomials.
\end{example}

\begin{example}
(Quantum $q$-Krawtchouk)
\label{ex:qkrawquantum}
Assume
\begin{eqnarray*}
\theta_i &=& \theta_0-sq(1-q^i),
\\
\theta^*_i &=& \theta^*_0+h^*(1-q^i)q^{-i} 
\end{eqnarray*}
for 
 $0 \leq i \leq d$  and
\begin{eqnarray*}
\varphi_i &=& -rh^*q^{1-i}(1-q^i)(1-q^{i-d-1}),
\\
\phi_i &=& h^*q^{d+2-2i}(1-q^i)(1-q^{i-d-1})(s-rq^{i-d-1})
\end{eqnarray*}
for $1 \leq i \leq d$.
Assume $ h^*, q, r, s$ are nonzero.
Assume
neither of
$q^i,
sq^i/r$
is equal to $1$ for $1 \leq i \leq d$.
Then the sequence
$(\theta_i, \theta^*_i, i=0..d;  \varphi_j, \phi_j, j=1..d)$
is a parameter array over $\K$.
The corresponding polynomials
$u_i$
satisfy
\begin{eqnarray*}
u_i(\theta_j)=
{}_2\phi_1 \Biggl({{q^{-i}, \;q^{-j}}\atop
{q^{-d}}}\;\Bigg\vert \; q,\;sr^{-1}q^{j+1}\Biggr)
%\label{eq:qrac}
\end{eqnarray*}
for $0 \leq i,j\leq d$. These $u_i$ are the
quantum $q$-Krawtchouk polynomials.
\end{example}

\begin{example}
($q$-Krawtchouk)
\label{ex:sbutnotss2}
\label{ex:qkraw}
Assume
\begin{eqnarray*}
\theta_i &=& 
\theta_0+h(1-q^i)q^{-i},
\\
\theta^*_i &=&  \theta^*_0+h^*(1-q^i)(1-s^*q^{i+1})q^{-i}
\end{eqnarray*}
for $0 \leq i \leq d$ and
\begin{eqnarray*}
\varphi_i &=& hh^*q^{1-2i}(1-q^i)(1-q^{i-d-1}),
%\label{lem:varphiphi1}
\\
\phi_i &=& hh^*s^*q(1-q^i)(1-q^{i-d-1})
%\label{lem:varphiphi2}
\end{eqnarray*}
for $1 \leq i \leq d$.
Assume $h,h^*, q, s^*$ are nonzero.
Assume  
$q^i\not=1$
for $1 \leq i \leq d$ and 
that $s^*q^i\not=1$ for $2 \leq i \leq 2d$.
Then the sequence
$(\theta_i, \theta^*_i, i=0..d;  \varphi_j, \phi_j, j=1..d)$
is a parameter array over $\K$.
The corresponding polynomials
$u_i$
  satisfy
\begin{eqnarray*}
u_i(\theta_j)=
{}_3\phi_2 \Biggl({{q^{-i}, \;s^*q^{i+1},\;q^{-j}}\atop
{0,\;\;q^{-d}}}\;\Bigg\vert \; q,\;q\Biggr)
%\label{eq:qrac}
\end{eqnarray*}
for $0 \leq i,j\leq d$. These $u_i$ are the $q$-Krawtchouk polynomials.
\end{example}

\begin{example}(Affine $q$-Krawtchouk)
\label{ex:qkrawaffine}
Assume
\begin{eqnarray*}
\theta_i &=& \theta_0+h(1-q^i)q^{-i},
\\
\theta^*_i &=& \theta^*_0+h^*(1-q^i)q^{-i}
\end{eqnarray*}
for 
 $0 \leq i \leq d$  and
\begin{eqnarray*}
\varphi_i &=& hh^*q^{1-2i}(1-q^i)(1-q^{i-d-1})(1-rq^i),
\\
\phi_i &=& -hh^*rq^{1-i}(1-q^i)(1-q^{i-d-1})
\end{eqnarray*}
for $1 \leq i \leq d$.
Assume $h,h^*, q, r$ are nonzero.
Assume
neither of
$q^i,
rq^i$
is equal to $1$ for $1 \leq i \leq d$.
Then the sequence
$(\theta_i, \theta^*_i, i=0..d;  \varphi_j, \phi_j, j=1..d)$
is a parameter array over $\K$.
The corresponding  polynomials 
$u_i$
 satisfy
\begin{eqnarray*}
u_i(\theta_j)=
{}_3\phi_2 \Biggl({{q^{-i}, \;0,\;q^{-j}}\atop
{rq,\;\;q^{-d}}}\;\Bigg\vert \; q,\;q\Biggr)
\end{eqnarray*}
for $0 \leq i,j\leq d$. These $u_i$ are the affine $q$-Krawtchouk polynomials.
\end{example}

\begin{example}
(Dual $q$-Krawtchouk)
\label{ex:sbutnotss}
Assume 
\begin{eqnarray*}
\theta_i &=& \theta_0+h(1-q^i)(1-sq^{i+1})q^{-i},
\\
\theta^*_i &=& \theta^*_0+h^*(1-q^i)q^{-i}
%\label{lem:eigdeig}
\end{eqnarray*}
for $0 \leq i \leq d$ and
\begin{eqnarray*}
\varphi_i &=& hh^*q^{1-2i}(1-q^i)(1-q^{i-d-1}),
%\label{lem:varphiphi1}
\\
\phi_i &=& hh^*sq^{d+2-2i}(1-q^i)(1-q^{i-d-1})
%\label{lem:varphiphi2}
\end{eqnarray*}
for $1 \leq i \leq d$.
Assume $h,h^*,q, s$ are nonzero.
Assume 
$q^i\not=1$
for $1 \leq i \leq d$ and 
$sq^i\not=1$ for $2 \leq i \leq 2d$.
Then the sequence
$(\theta_i, \theta^*_i, i=0..d;  \varphi_j, \phi_j, j=1..d)$
is a parameter array over $\K$.
The corresponding polynomials
$u_i$
 satisfy
\begin{eqnarray*}
u_i(\theta_j)=
{}_3\phi_2 \Biggl({{q^{-i}, \;q^{-j},sq^{j+1}}\atop
{0,\;\;q^{-d}}}\;\Bigg\vert \; q,\;q\Biggr)
\end{eqnarray*}
for $0 \leq i,j\leq d$. These $u_i$ are the dual $q$-Krawtchouk polynomials.
\end{example}

\begin{example}(Racah)
\label{ex:racah}
Assume 
\begin{eqnarray}
\theta_i &=& \theta_0+hi(i+1+s),
\label{eq:thracII}
\\
\theta^*_i &=& \theta^*_0+h^*i(i+1+s^*)
\label{eq:thsracII}
\end{eqnarray}
for $0 \leq i \leq d$ and
\begin{eqnarray}
\varphi_i &=& hh^*i(i-d-1)(i+r_1)(i+r_2),
\label{lem:varphiphi1II}
\\
\phi_i &=& hh^*i(i-d-1)(i+s^*-r_1)(i+s^*-r_2)
\label{lem:varphiphi2II}
\end{eqnarray}
for $1 \leq i \leq d$.
Assume $h, h^*$ are nonzero 
and that $r_1+r_2=s+s^*+d+1$.
Assume the characteristic of $\K$ is $0$ or a prime 
greater than $d$.
Assume none of
$
r_1, r_2,
s^*-r_1$, $ s^*-r_2$
is equal to $-i$ for $1 \leq i \leq d$ and 
that neither of $s, s^*$ is equal to $-i$ for $2 \leq i \leq 2d$.
Then the sequence
$(\theta_i, \theta^*_i, i=0..d;  \varphi_j, \phi_j, j=1..d)$
is a parameter array over $\K$.
The corresponding polynomials
$u_i$
 satisfy
\begin{eqnarray*}
u_i(\theta_j)=
{}_4F_3 \Biggl({{-i, \;i+1+s^*,\;-j,\;j+1+s}\atop
{r_1+1,\;\;r_2+1,\;\;-d}}\;\Bigg\vert \; 1\Biggr)
%\label{eq:qrac}
\end{eqnarray*}
for $0 \leq i,j\leq d$. These $u_i$ are the Racah polynomials.
\end{example}

\begin{example}(Hahn)
\label{ex:Hahn}
Assume
\begin{eqnarray*}
\theta_i &=& \theta_0+si,
\\
\theta^*_i &=& \theta^*_0+h^*i(i+1+s^*)
\end{eqnarray*}
for $0 \leq i \leq d$ and
\begin{eqnarray*}
\varphi_i &=& h^*si(i-d-1)(i+r),
%\label{lem:varphiphi1}
\\
\phi_i &=& -h^*si(i-d-1)(i+s^*-r)
%\label{lem:varphiphi2}
\end{eqnarray*}
for $1 \leq i \leq d$.
Assume $h^*, s$ are nonzero.
Assume the characteristic of $\K$ is $0$ or a prime greater than $d$.
Assume neither of
$
r, 
s^*-r$
is equal to $-i$ for $1 \leq i \leq d$ and 
that  $s^*\not=-i$ for $2 \leq i \leq 2d$.
Then the sequence
$(\theta_i, \theta^*_i, i=0..d;  \varphi_j, \phi_j, j=1..d)$
is a parameter array over $\K$.
The corresponding polynomials
$u_i$
 satisfy
\begin{eqnarray*}
u_i(\theta_j)=
{}_3F_2 \Biggl({{-i, \;i+1+s^*,\;-j}\atop
{r+1,\;\;-d}}\;\Bigg\vert \; 1\Biggr)
%\label{eq:qrac}
\end{eqnarray*}
for $0 \leq i,j\leq d$. These $u_i$ are the Hahn polynomials.
\end{example}

\begin{example}(Dual Hahn)
\label{ex:dHahn}
Assume
\begin{eqnarray*}
\theta_i &=& \theta_0+hi(i+1+s),
\\
\theta^*_i &=& \theta^*_0+s^*i
\end{eqnarray*}
for $0 \leq i \leq d$ and
\begin{eqnarray*}
\varphi_i &=& hs^*i(i-d-1)(i+r),
%\label{lem:varphiphi1}
\\
\phi_i &=& hs^*i(i-d-1)(i+r-s-d-1)
%\label{lem:varphiphi2}
\end{eqnarray*}
for $1 \leq i \leq d$.
Assume $h, s^*$ are nonzero.
Assume the characteristic of $\K$ is $0$ or a prime greater than $d$.
Assume neither of
$
r, 
s-r$
is equal to $-i$ for $1 \leq i \leq d$ and 
that  $s\not=-i$ for $2 \leq i \leq 2d$.
Then the sequence
$(\theta_i, \theta^*_i, i=0..d;  \varphi_j, \phi_j, j=1..d)$
is a parameter array over $\K$.
The corresponding polynomials 
$u_i$
 satisfy
\begin{eqnarray*}
u_i(\theta_j)=
{}_3F_2 \Biggl({{-i, \;-j,\;j+1+s}\atop
{r+1,\;\;-d}}\;\Bigg\vert \; 1\Biggr)
%\label{eq:qrac}
\end{eqnarray*}
for $0 \leq i,j\leq d$. These $u_i$ are the dual Hahn polynomials.
\end{example}

\begin{example}(Krawtchouk)
\label{ex:kraw}
Assume
\begin{eqnarray*}
\theta_i &=& \theta_0+si,
\\
\theta^*_i &=& \theta^*_0+s^*i
\end{eqnarray*}
for $0 \leq i \leq d$ and
\begin{eqnarray*}
\varphi_i &=& ri(i-d-1)
%\label{lem:varphiphi1}
\\
\phi_i &=& (r-ss^*)i(i-d-1)
%\label{lem:varphiphi2}
\end{eqnarray*}
for $1 \leq i \leq d$.
Assume $r,s,s^*$ are nonzero.
Assume the characteristic of $\K$ is $0$ or a prime greater than $d$.
Assume
$
r\not=ss^*$. 
Then the sequence
$(\theta_i, \theta^*_i, i=0..d;  \varphi_j, \phi_j, j=1..d)$
is a parameter array over $\K$.
The corresponding polynomials 
$u_i$
 satisfy
\begin{eqnarray*}
u_i(\theta_j)=
{}_2F_1 \Biggl({{-i, \;-j}\atop
{-d}}\;\Bigg\vert \; r^{-1}ss^*\Biggr)
%\label{eq:qrac}
\end{eqnarray*}
for $0 \leq i,j\leq d$. These $u_i$ are the Krawtchouk polynomials.
\end{example}

\begin{example}(Bannai/Ito)
\label{ex:bi}
Assume
\begin{eqnarray}
\theta_i &=& \theta_0+h(s-1+(1-s+2i)(-1)^i),
\label{eq:thracIII}
\\
\theta^*_i &=& \theta^*_0+h^*(s^*-1+(1-s^*+2i)(-1)^i)
\label{eq:thsracIII}
\end{eqnarray}
for $0 \leq i \leq d$ and
\begin{eqnarray}
\varphi_i &=& 
\cases{
-4hh^*i(i+r_1),    &if $\;i$ even, $\;d$ even;\cr
-4hh^*(i-d-1)(i+r_2),    &if $\;i$ odd, $\;d$ even;\cr
-4hh^*i(i-d-1),    &if $\;i$ even, $\;d$ odd;\cr
-4hh^*(i+r_1)(i+r_2),    &if $\;i$ odd, $\;d$ odd,
}
\label{lem:varphiphi1III}
\\
\phi_i &=& 
\cases{
4hh^*i(i-s^*-r_1),    &if $\;i$ even, $\;d$ even;\cr
4hh^*(i-d-1)(i-s^*-r_2),    &if $\;i$ odd, $\;d$ even;\cr
-4hh^*i(i-d-1),    &if $\;i$ even, $\;d$ odd;\cr
-4hh^*(i-s^*-r_1)(i-s^*-r_2),    &if $\;i$ odd, $\;d$ odd
}
\label{lem:varphiphi2III}
\end{eqnarray}
for $1 \leq i \leq d$.
Assume $h, h^*$ are nonzero and that
$r_1+r_2=-s-s^*+d+1$.
Assume the characteristic of $\K$ is either $0$
or an odd prime greater
than $d/2$.
Assume neither of
$
r_1, -s^*-r_1$
is equal to $-i$ for $1 \leq i \leq d$, $d-i$ even.
Assume neither of 
$r_2$, $ -s^*-r_2$
is equal to $-i$ for $1 \leq i \leq d$, $i$ odd.
Assume neither of $s, s^*$ is equal to $2i$ for $1 \leq i \leq d$.
Then the sequence
$(\theta_i, \theta^*_i, i=0..d;  \varphi_j, \phi_j, j=1..d)$
is a parameter array over $\K$.
We call the corresponding polynomials  from
Definition \ref{def:assocpoly}
the Bannai/Ito polynomials 
\cite[p. 260]{BanIto}.
\end{example}

\begin{example} (Orphan) 
\label{ex:orphan}
For this example assume
$\K$ has characteristic 2.
%Put $d=3$.
For notational convenience we
define some scalars $\gamma_0, 
\gamma_1, 
\gamma_2, 
\gamma_3$  
in $ \K$.
We define
$\gamma_i=0$ for $i \in \lbrace 0,3\rbrace $
 and $\gamma_i=1$ for $i\in \lbrace 1,2\rbrace$.
Assume
\begin{eqnarray}
\theta_i &=& \theta_0 + h(si+ \gamma_i),
\label{eq:thracIV}
\\
\theta^*_i &=& \theta^*_0 + h^*(s^*i+ \gamma_i)
\label{eq:thsracIV}
\end{eqnarray}
for $0 \leq i \leq 3$.
Assume $\varphi_1 = hh^*r$, 
$\varphi_2 = hh^*$, 
$\varphi_3 = hh^*(r+s+s^*)$ 
and 
$\phi_1 = hh^*(r+s(1+s^*))$, 
$\phi_2 = hh^*$, 
$\phi_3 = hh^*(r+s^*(1+s))$.
Assume each of $h,h^*, s,s^*,r$ is nonzero.
Assume neither of $s, s^*$ is equal to $1$
and that $r$ is equal to none of 
$s+s^*$, 
$s(1+s^*)$, 
$s^*(1+s)$.
Then the sequence
$(\theta_i, \theta^*_i, i=0..3;  \varphi_j, \phi_j, j=1..3)$
is a parameter array over $\K$ which has
diameter $3$.
We call the corresponding polynomials  from 
Definition \ref{def:assocpoly}
the orphan polynomials.
\end{example}

%\begin{remark}
%For $d\geq 3$ the families of parameter arrays
%given in Examples \ref{ex:qracah}--\ref{ex:orphan}
%are mutually disjoint.
%\end{remark}

\begin{theorem}
\label{thm:thatsit}
Every parameter array over $\K$
is listed in at least one of the Examples
\ref{ex:qracah}--\ref{ex:orphan}.
\end{theorem}
\noindent {\it Proof:}
Let 
$p:=(\theta_i, \theta^*_i, i=0..d;  \varphi_j, \phi_j, j=1..d)$
denote a parameter array over $\K$.
We show this array is given in
at least one of the Examples \ref{ex:qracah}--\ref{ex:orphan}.
We assume $d\geq 1$; otherwise the result is trivial. 
For notational convenience let ${\tilde K}$ denote the 
algebraic closure of $\K$.
We define a scalar $q \in {\tilde \K}$ as follows.
For $d\geq 3$, we let $q$ denote a nonzero scalar in
$\tilde \K$ such that $q+q^{-1}+1$ is equal to the
common value of
(\ref{eq:betaplusone}). For $d<3$ we let $q$
denote a nonzero scalar in
$\tilde \K$ such that
$q\not=1$ and $q \not=-1$.
By PA5,
both
\begin{eqnarray}
%\theta_{i-2}- \lbrack 3 \rbrack_q\theta_{i-1}+
%\lbrack 3 \rbrack_q\theta_{i} -\theta_{i+1}&=&0,
\theta_{i-2}- \xi \theta_{i-1}+
\xi \theta_{i} -\theta_{i+1}&=&0,
\label{eq:lin1}
\\
%\theta^*_{i-2}- \lbrack 3 \rbrack_q\theta^*_{i-1}+
%\lbrack 3 \rbrack_q\theta^*_{i} -\theta^*_{i+1}&=&0
\theta^*_{i-2}- \xi\theta^*_{i-1}+
\xi\theta^*_{i} -\theta^*_{i+1}&=&0
\label{eq:lin2}
\end{eqnarray}
for $2 \leq i \leq d-1$, where $\xi = q+q^{-1}+1$.
We divide the argument into the following four cases.
(I) $q \not=1$, $q \not=-1$;
(II) $q = 1$ and $\mbox{char}(\K) \not=2$;
(III) $q = -1$ and $\mbox{char}(\K) \not=2$;
(IV) $q= 1 $ and 
 $\mbox{char}(\K)=2$.

\medskip
\noindent Case I: $q\not=1$, $q\not=-1$.  \\
\noindent By (\ref{eq:lin1})  there exist scalars $\eta, \mu, h$ in
${\tilde \K}$ such that
\begin{eqnarray}
\theta_i &=& \eta + \mu q^i+ h q^{-i}  \qquad \qquad 
(0 \leq i\leq d).
\label{eq:thI}
\end{eqnarray}
By 
(\ref{eq:lin2}) there exist 
 scalars $\eta^*, \mu^*,h^*$ in
${\tilde \K}$ such that
\begin{eqnarray}
\theta^*_i &=& \eta^* + \mu^* q^i+ h^*q^{-i}  \qquad \qquad 
(0 \leq i\leq d).
\label{eq:thsI}
\end{eqnarray}
Observe $\mu, h$ are not both 0; 
otherwise $\theta_1=\theta_0$ by
(\ref{eq:thI}).
Similarly 
$\mu^*, h^*$ are not both 0.
For $1 \leq i \leq d$ we have $q^i \not=1$; otherwise
$\theta_i=\theta_0$ by
(\ref{eq:thI}).
Setting $i=0$ 
in 
(\ref{eq:thI}), (\ref{eq:thsI}) we obtain
\begin{eqnarray}
\theta_0 &=& \eta + \mu +h,
\label{eq:thI0}
\\
\theta^*_0 &=& \eta^* + \mu^*+h^*.
\label{eq:thsI0}
\end{eqnarray}
%Using 
%(\ref{eq:thI}) we obtain
%\begin{eqnarray}
%\sum_{n=0}^{i-1} \frac{\theta_n-\theta_{d-n}}{\theta_0-\theta_d}
%= \frac
%{
%(q^i-1)(q^{d-i+1}-1)}
%{
%(q-1)(q^d-1)
%}
%\qquad \qquad (1 \leq i \leq d).
%\label{eq:sumi}
%\end{eqnarray}
We claim there exists $\tau \in {\tilde \K}$ such that both
\begin{eqnarray}
\varphi_i &=&
(q^i-1)(q^{d-i+1}-1)(\tau -\mu \mu^* q^{i-1} - h h^* q^{-i-d}),
\label{eq:varphiIfirst}
\\
\phi_i &=&
(q^i-1)(q^{d-i+1}-1)(\tau -h \mu^* q^{i-d-1} - \mu h^* q^{-i})
\label{eq:phiIfirst}
\end{eqnarray}
for $1 \leq i \leq d$.
Since $q\not=1$ and $q^d\not=1$ there exists $\tau \in {\tilde \K}$
such that 
(\ref{eq:varphiIfirst})
holds for $i=1$.
In the equation of PA4, we eliminate
$\varphi_1$  using  
(\ref{eq:varphiIfirst}) at $i=1$, and evaluate the result using
(\ref{eq:thI}), 
(\ref{eq:thsI})
in order to obtain 
(\ref{eq:phiIfirst}) for $1 \leq i \leq d$.
In the equation of PA3, we
eliminate $\phi_1$ using  
(\ref{eq:phiIfirst}) at $i=1$, and evaluate the result using
(\ref{eq:thI}), 
(\ref{eq:thsI})
in order to obtain
(\ref{eq:varphiIfirst})  for $1 \leq i \leq d$.
We have now proved the claim. 
We now break the argument into subcases.
For each subcase our argument is similar.
We will discuss the first subcase in detail in order to
give the idea; for the remaining subcases 
we give the essentials only.

\medskip
\noindent Subcase $q$-Racah: 
$\mu\not=0, \mu^*\not=0, h\not=0, h^*\not=0$.  
We show $p$ is listed in Example
\ref{ex:qracah}.
Define 
\begin{eqnarray}
s:=\mu h^{-1} q^{-1}, \qquad   
s^*:=\mu^* h^{*-1} q^{-1}.
\label{eq:hdef}
\end{eqnarray}
 Eliminating
$\eta $ in 
(\ref{eq:thI})  using
(\ref{eq:thI0}) and eliminating $\mu $ in the result
using the equation on the left in  
(\ref{eq:hdef}),
we obtain
(\ref{eq:thrac}) for $0 \leq i\leq d$.
%Eliminating
%$\eta^* $ in 
%(\ref{eq:thsI})  using
%(\ref{eq:thsI0}) and eliminating $\mu^*$ in the result
%using the equation on the right in 
%(\ref{eq:hdef})
Similarly we obtain
(\ref{eq:thsrac}) for $0 \leq i\leq d$.
Since $\tilde \K $ is algebraically closed it contains
scalars $r_1, r_2$
such that both
\begin{eqnarray}
r_1r_2 = s s^*q^{d+1}, \qquad \qquad 
r_1+r_2=\tau h^{-1} h^{*-1} q^d.
\label{eq:rdef}
\end{eqnarray}
Eliminating $\mu, \mu^*, \tau$ in
(\ref{eq:varphiIfirst}), 
(\ref{eq:phiIfirst}) using
(\ref{eq:hdef}) and the equation on the right in
(\ref{eq:rdef}), and evaluating the result using
the equation on the left in 
(\ref{eq:rdef}), 
we obtain 
(\ref{lem:varphiphi1}), 
(\ref{lem:varphiphi2}) 
for $1 \leq i \leq d$.
By the construction 
each of $h, h^*, q, s, s^*$ is nonzero.
Each of $r_1, r_2$ is nonzero by the equation on the left in
(\ref{eq:rdef}). The remaining inequalities mentioned below
(\ref{lem:varphiphi2}) follow
from PA1, PA2 and 
(\ref{eq:thrac})--(\ref{lem:varphiphi2}).
We have now shown $p$ is listed in Example
\ref{ex:qracah}.
%To obtain
%(\ref{eq:fqrac}) set $\lambda=\theta_j$ in
%(\ref{eq:fipoly}) and evaluate the result using
%(\ref{eq:thrac}), (\ref{eq:thsrac}),  (\ref{lem:varphiphi1}).

\medskip
\noindent
We now give the remaining subcases of Case I.
We list the essentials only. 

\medskip
\noindent Subcase $q$-Hahn:
$\mu=0, \mu^*\not=0, h\not=0, h^*\not=0, \tau \not=0$.  
Definitions:
\begin{eqnarray*}
s^*:=\mu^* h^{*-1} q^{-1},\qquad
r:=\tau h^{-1} h^{*-1} q^d.
\end{eqnarray*}
\noindent Subcase dual $q$-Hahn:
$\mu\not=0, \mu^*=0, h\not=0, h^*\not=0, \tau \not=0$.  
Definitions:
\begin{eqnarray*}
 s:=\mu h^{-1} q^{-1},\qquad
r:=\tau h^{-1} h^{*-1} q^d.
\end{eqnarray*}
\noindent Subcase quantum $q$-Krawtchouk: 
$\mu\not=0, \mu^*=0, h=0, h^*\not=0, \tau \not=0$. 
Definitions: 
\begin{eqnarray*}
s:=\mu  q^{-1}, \qquad   
r :=\tau h^{*-1}q^d.
\end{eqnarray*}
\noindent Subcase  $q$-Krawtchouk: 
$\mu=0, \mu^*\not=0, h\not=0, h^*\not=0, \tau =0$. 
Definition: 
\begin{eqnarray*}
s^*:=
 \mu^* h^{*-1} q^{-1}.
\end{eqnarray*}
\noindent Subcase affine $q$-Krawtchouk:
$\mu=0, \mu^*=0, h\not=0, h^*\not=0, \tau \not=0$. 
Definition:
\begin{eqnarray*}
r:=\tau h^{-1} h^{*-1} q^d.
\end{eqnarray*}
\noindent Subcase dual $q$-Krawtchouk:
$\mu \not=0, \mu^*=0, h\not=0, h^*\not=0, \tau =0$. 
Definition:
\begin{eqnarray*}
s:=
 \mu h^{-1} q^{-1}.
\end{eqnarray*}
\noindent 
We have a few more comments concerning Case I. 
Earlier we mentioned that $\mu,h$ are not both 0 and 
that $\mu^*,h^*$ are not both 0.
Suppose one of 
$\mu,h$ is 0 
and one of  $\mu^*,h^*$ is 0.
Then $\tau\not=0$; otherwise $\varphi_1=0$ by 
(\ref{eq:varphiIfirst}) or 
$\phi_1=0$ by 
(\ref{eq:phiIfirst}).
Suppose  $\mu^*\not=0$,  $h^*=0$.
Replacing $q$ by $q^{-1}$
we obtain  
$\mu^* =0$, $h^*\not=0$. 
Suppose 
 $\mu^* \not=0$, $h^*\not=0$,  
 $\mu\not=0$,  $h=0$. Replacing $q$ by $q^{-1}$
we obtain
 $\mu^* \not=0$, $h^*\not=0$,  
 $\mu=0$,  $h\not=0$.
By these comments we find that after 
replacing $q$ by $q^{-1}$ if necessary,
one of the above subcases holds.
This completes our argument for Case I.

\medskip
\noindent Case II: $q=1$ and $\mbox{char}(\K)\not=2$.
\\
\noindent By (\ref{eq:lin1}) and since
$\mbox{char}(\K)\not=2$, 
there exist scalars $\eta, \mu, h$ in
${\tilde \K}$ such that
\begin{eqnarray}
\theta_i &=& \eta + (\mu + h) i+ h i^2  \qquad \qquad 
(0 \leq i\leq d).
\label{eq:thII}
\end{eqnarray}
Similarly  
there exist 
 scalars $\eta^*, \mu^*,h^*$ in
${\tilde \K}$ such that
\begin{eqnarray}
\theta^*_i &=& \eta^* + (\mu^* +h^*)i+ h^*i^2  \qquad \qquad 
(0 \leq i\leq d).
\label{eq:thsII}
\end{eqnarray}
Observe $\mu, h$ are not both 0; 
otherwise $\theta_1=\theta_0$.
Similarly 
$\mu^*, h^*$ are not both 0.
For any prime $i$ such that $ i \leq d$ we have
$\mbox{char}(\K)\not=i$;
otherwise
$\theta_i=\theta_0$ by
(\ref{eq:thII}).
Therefore
$\mbox{char}(\K)$ is 0 or a prime greater than $d$.
Setting $i=0$ 
in 
(\ref{eq:thII}), (\ref{eq:thsII}) we obtain
\begin{eqnarray}
\theta_0 = \eta, \qquad \qquad 
\theta^*_0 = \eta^*.
\label{eq:thII0}
\end{eqnarray}
We claim there exists $\tau \in {\tilde \K}$ such that both
\begin{eqnarray}
%\varphi_i &=&
%i(d-i+1)(\tau -\mu \mu^*-( \mu h^*+h \mu^*) i - h h^* i(i+d+1)),
\varphi_i &=&
i(d-i+1)(\tau -( \mu h^*+h \mu^*) i - h h^* i(i+d+1)),
\label{eq:varphiIIfirst}
\\
%\phi_i &=&
%i(d-i+1)(\tau +h \mu^*(1+d)  +(\mu h^*-h \mu^*) i + h h^* i(d-i+1))
\phi_i &=&
i(d-i+1)(\tau +\mu \mu^*+ h \mu^*(1+d)  +(\mu h^*-h \mu^*) i + h h^* i(d-i+1))
\label{eq:phiIIfirst}
\end{eqnarray}
for $1 \leq i \leq d$.
There exists $\tau \in {\tilde \K}$
such that 
(\ref{eq:varphiIIfirst})
holds for $i=1$.
In the equation of PA4, we eliminate
$\varphi_1$  using  
(\ref{eq:varphiIIfirst}) at $i=1$, and evaluate the result using
(\ref{eq:thII}), 
(\ref{eq:thsII})
 in order to obtain 
(\ref{eq:phiIIfirst}) for $1 \leq i \leq d$.
In the equation of PA3, we
eliminate $\phi_1$ using  
(\ref{eq:phiIIfirst}) at $i=1$, and evaluate the result using
(\ref{eq:thII}), 
(\ref{eq:thsII})
in order to obtain
(\ref{eq:varphiIIfirst})  for $1 \leq i \leq d$.
We have now proved the claim. 
We now break the argument into subcases.
%For each subcase our argument is similar.
%We will discuss the first one in detail in order to
%give the idea, and for the others give the essentials only.

\medskip
\noindent Subcase Racah: 
$h\not=0$, $h^*\not=0$.  
We show $p$ is listed in Example
\ref{ex:racah}.
Define 
\begin{eqnarray}
s:=\mu h^{-1} , \qquad   
s^*:=\mu^* h^{*-1}.
\label{eq:hdefII}
\end{eqnarray}
 Eliminating
$\eta, \mu $ in 
(\ref{eq:thII})  using
(\ref{eq:thII0}),
(\ref{eq:hdefII})
we obtain
(\ref{eq:thracII}) for $0 \leq i\leq d$.
Eliminating
$\eta^*, \mu^*$  in 
(\ref{eq:thsII})  using
(\ref{eq:thII0}),
(\ref{eq:hdefII})
we obtain
(\ref{eq:thsracII}) for $0 \leq i\leq d$.
Since $\tilde \K $ is algebraically closed it contains
scalars $r_1, r_2$
such that both
\begin{eqnarray}
r_1r_2 =  - \tau h^{-1} h^{*-1}, \qquad \qquad 
r_1+r_2=s + s^*+d+1.
\label{eq:rdefII}
\end{eqnarray}
Eliminating $\mu, \mu^*, \tau$ in
(\ref{eq:varphiIIfirst}), 
(\ref{eq:phiIIfirst}) using
(\ref{eq:hdefII}) and the equation on the left in
(\ref{eq:rdefII})
we obtain 
(\ref{lem:varphiphi1II}), 
(\ref{lem:varphiphi2II}) 
for $1 \leq i \leq d$.
By the construction 
each of $h, h^*$ is nonzero.
The remaining inequalities mentioned below
(\ref{lem:varphiphi2II}) follow
from PA1, PA2 and 
(\ref{eq:thracII})--(\ref{lem:varphiphi2II}).
We have now shown $p$ is listed in Example
\ref{ex:racah}.

\medskip
\noindent We now give the remaining subcases of Case II. We list the
essentials only.

\medskip
\noindent Subcase Hahn:
$ h=0$, $h^*\not=0 $.  
Definitions:
\begin{eqnarray*}
s= \mu, 
\qquad s^*:=\mu^* h^{*-1}, \qquad 
r:=-\tau \mu^{-1} h^{*-1}.
\end{eqnarray*}
\noindent Subcase dual Hahn:
$ h\not=0, h^*=0$.  
Definitions:
\begin{eqnarray*}
 s:=\mu h^{-1},\qquad
s^*= \mu^*, \qquad 
r:=- \tau 
h^{-1}
\mu^{*-1} .
\end{eqnarray*}
\noindent Subcase Krawtchouk: 
$ h=0, h^*=0$. 
Definitions: 
\begin{eqnarray*}
s:=\mu , \qquad   
s^*:=\mu^* , \qquad   
r :=- \tau.
\end{eqnarray*}

\medskip
\noindent Case III: $q = -1$ and $\mbox{char}(\K)\not=2$.
\\
\noindent
We show $p$ is listed in
Example \ref{ex:bi}.
By (\ref{eq:lin1})  and since
$\mbox{char}(\K)\not=2$, 
there exist scalars $\eta, \mu, h$ in
${\tilde \K}$ such that
\begin{eqnarray}
\theta_i &=& \eta + \mu (-1)^i+ 2h i(-1)^i  \qquad \qquad 
(0 \leq i\leq d).
\label{eq:thIII}
\end{eqnarray}
Similarly  
there exist 
 scalars $\eta^*, \mu^*,h^*$ in
${\tilde \K}$ such that
\begin{eqnarray}
\theta^*_i &=& \eta^* + \mu^*(-1)^i+ 2 h^*i(-1)^i  \qquad \qquad 
(0 \leq i\leq d).
\label{eq:thsIII}
\end{eqnarray}
Observe $h\not=0$;
otherwise $\theta_2=\theta_0$ by
(\ref{eq:thIII}).
Similarly 
$h^*\not=0$.
For any prime $i$ such that $ i \leq d/2$ we have
$\mbox{char}(\K)\not=i$;
otherwise
$\theta_{2i}=\theta_0$ by
(\ref{eq:thIII}).
By this and since 
$\mbox{char}(\K)\not=2$
we find $\mbox{char}(\K)$ is either
0 or an odd prime greater than $d/2$.
Setting $i=0$ 
in 
(\ref{eq:thIII}), (\ref{eq:thsIII}) we obtain
\begin{eqnarray}
\theta_0 &=& \eta+\mu, 
\qquad \qquad 
\theta^*_0 = \eta^*+\mu^*.
\label{eq:thIII0}
\end{eqnarray}
We define 
\begin{eqnarray}
s:=1-\mu h^{-1}, \qquad \qquad  s^*=1-\mu^*h^{*-1}. 
\label{eq:sssIII}
\end{eqnarray}
Eliminating $\eta$ in 
(\ref{eq:thIII}) using
(\ref{eq:thIII0}) and eliminating $\mu$ in the result using
(\ref{eq:sssIII}) we find 
(\ref{eq:thracIII}) holds for $0 \leq i \leq d$.
Similarly we find
(\ref{eq:thsracIII}) holds for $0 \leq i \leq d$.
We now define $r_1, r_2$. First 
assume 
$d$ is odd. 
Since $\tilde \K$ is algebraically closed it contains 
$r_1, r_2$ such that
\begin{eqnarray}
r_1+r_2=-s-s^*+d+1
\label{r1pr2III}
\end{eqnarray}
and such that
\begin{eqnarray}
4hh^*(1+r_1)(1+r_2)= -\varphi_1.
\label{eq:rprodIII}
\end{eqnarray}
Next assume $d$ is even. Define
\begin{eqnarray}
r_2 := -1 + \frac{\varphi_1}{4hh^*d}
\label{eq:deveIII}
\end{eqnarray}
and define $r_1$ so that 
(\ref{r1pr2III}) holds.
We have now defined $r_1, r_2$ for either parity of $d$.
In the equation of PA4, we eliminate
$\varphi_1$  using 
(\ref{eq:rprodIII}) or 
(\ref{eq:deveIII}), and evaluate the result using
(\ref{eq:thracIII}), 
(\ref{eq:thsracIII}) 
 in order to obtain
(\ref{lem:varphiphi2III}) for $1 \leq i\leq d$.
In the equation of PA3, we
eliminate $\phi_1$ using  
(\ref{lem:varphiphi2III}) 
at $i=1$, and evaluate the result using
(\ref{eq:thracIII}), 
(\ref{eq:thsracIII})
in order to obtain
(\ref{lem:varphiphi1III}) for $1 \leq i\leq d$.
We mentioned each of $h, h^*$ is nonzero.
The remaining inequalities mentioned below
(\ref{lem:varphiphi2III}) 
follow from PA1, PA2 and
(\ref{eq:thracIII})--(\ref{lem:varphiphi2III}).
We have now shown $p$  is listed in
Example \ref{ex:bi}.

\medskip
\noindent Case IV: $q = 1$ and $\mbox{char}(\K)=2$.
\\
\noindent 
We show $p$ is listed in
  Example
\ref{ex:orphan}.
We first show $d=3$.
Recall $d\geq 3$
since $q=1$.
Suppose $d\geq 4$.
By 
(\ref{eq:lin1}) we have $\sum_{j=0}^3 \theta_j=0$
and 
$\sum_{j=1}^4 \theta_j=0$.
Adding 
these sums
we find $\theta_0=\theta_4$ which contradicts
PA1.
Therefore $d=3$. 
We claim there exist nonzero scalars 
$h,s$ in $\K$ such that 
(\ref{eq:thracIV}) holds for $0 \leq i \leq 3$.
Define $h=\theta_0+\theta_2$. Observe $h\not=0$; otherwise
$\theta_0=\theta_2$.
Define $s=(\theta_0+\theta_3)h^{-1}$.
Observe
$s \not=0$; otherwise $\theta_0=\theta_3$.
Using these values for $h,s$ we find 
(\ref{eq:thracIV}) holds for $i=0,2,3$.
By this and $\sum_{j=0}^3 \theta_j=0$ we find 
(\ref{eq:thracIV}) holds for $i=1$.
We have now proved our claim.
Similarly there exist nonzero scalars
$h^*,s^*$ in $\K$ such that 
(\ref{eq:thsracIV}) holds for $0 \leq i \leq 3$.
Define 
$r:=\varphi_1h^{-1}h^{*-1}$.
Observe $r\not=0$ and that
$\varphi_1=h h^*r$.
In the equation of PA4, we eliminate
$\varphi_1$  using  
$\varphi_1=h h^*r$
 and evaluate the result using
(\ref{eq:thracIV}),
(\ref{eq:thsracIV})
 in order to obtain 
$\phi_1 = hh^*(r+s(1+s^*))$, 
$\phi_2 = hh^*$, 
$\phi_3 = hh^*(r+s^*(1+s))$.
In the equation of PA3, we
eliminate $\phi_1$ using  
$\phi_1 = hh^*(r+s(1+s^*))$ 
 and evaluate the result using
(\ref{eq:thracIV}),
(\ref{eq:thsracIV})
 in order to obtain
$\varphi_2 = hh^*$, 
$\varphi_3 = hh^*(r+s+s^*)$.
We mentioned each of $h, h^*, s,s^*, r$ is nonzero.
Observe $s\not=1$; otherwise $\theta_1=\theta_0$.
Similarly 
 $s^*\not=1$.
Observe $r\not=s+s^*$; otherwise $\varphi_3=0$.
Observe $r\not=s(1+s^*)$; otherwise $\phi_1=0$.
Observe $r\not=s^*(1+s)$; otherwise $\phi_3=0$.
We have now shown  $p$ is
listed in Example
\ref{ex:orphan}. 
We are done with Case IV and the proof is complete.
\hfill $\Box $ \\

\section{Suggestions for further research}

\noindent In this section we give some suggestions for further
research.

% solved fall 06 in the paper affine...by nomura and terwilliger%%%
%\begin{problem}
%\rm
%Let 
%$\Phi=(A;A^*; \lbrace E_i\rbrace_{i=0}^d;  
%\lbrace E^*_i\rbrace_{i=0}^d)$ denote a Leonard system
%in $\mathcal A$. Let $\alpha, \alpha^*, \beta, \beta^*$
%denote scalars in $\K$ such that $\alpha\not=0$ and $\alpha^*\not=0$.
%We mentioned in Section 4 that the sequence
%$(\alpha A+\beta I;\alpha^*A^*+\beta^* I; \lbrace E_i\rbrace_{i=0}^d;  
%\lbrace E^*_i\rbrace_{i=0}^d)$ is a Leonard system
%in $\mathcal A$. In some cases this system is isomorphic to
%a relative of $\Phi$; determine all the cases where this occurs.
%\end{problem}

\begin{problem} \rm Let $V$ denote a vector space over $\K$ 
with finite positive dimension and 
let $A,A^*$ denote a tridiagonal pair
on $V$.
 Let $\alpha, \alpha^*, \beta, \beta^*$
denote scalars in $\K$ with $\alpha, \alpha^*$ nonzero,
and note that the pair $\alpha A+\beta I$, $\alpha^* A^*+\beta^*I$
is a tridiagonal pair on $V$. 
Find necessary and sufficient conditions 
for this tridiagonal pair to be isomorphic to the tridiagonal
pair $A,A^*$.
Also, find necessary and sufficient conditions 
for this tridiagonal pair to be isomorphic to the tridiagonal
pair $A^*,A$. This problem has been solved for Leonard pairs
 \cite{NT:aff}. 
%Vidar is considering the problem for $d=2$.
\end{problem}

%%%%%%%solved fall 06 nomura and ter 
%\begin{conjecture}
%Let 
%$(A;A^*; \lbrace E_i\rbrace_{i=0}^d;  
%\lbrace E^*_i\rbrace_{i=0}^d)$ 
%and 
%$(B;B^*; \lbrace E_i\rbrace_{i=0}^d;  
%\lbrace E^*_i\rbrace_{i=0}^d)$ 
%denote  Leonard systems
%in $\mathcal A$. Then there exist scalars
%$\alpha, \alpha^*, \beta, \beta^*$ in
%$\K$ such that $\alpha \not=0, \alpha^*\not=0$ and
%\begin{eqnarray*}
%B= 
%\alpha A+\beta I,
%\qquad \qquad 
%B^*=\alpha^*A^*+\beta^* I.
%\end{eqnarray*}
%\end{conjecture}
%%%%%%%

\begin{problem} 
\rm
Assume
 $\K=\R$.
With reference to
Definition \ref{def:ip},
find a necessary and sufficient condition on the parameter
array of $\Phi$, for
the bilinear form $\langle \,,\,\rangle $ to be positive definite.
By definition the form 
$\langle \,,\,\rangle $ is positive definite whenever
$\Vert u \Vert^2 > 0 $ for all nonzero $u \in V$.
\end{problem}

%extra below
%\begin{problem}
%\rm
%Assume
% $\K=\R$ and
%let $\Phi$ denote the Leonard system
%from Definition \ref{eq:ourstartingpt}.
%For $0 \leq i\leq d$ we define $A_i=v_i(A)$, where
%the polynomial $v_i$ is from Definition
%\ref{def:vi1}. Observe there exist real scalars $p^h_{ij}$
%$(0 \leq h,i,j\leq d)$ such that 
%\begin{eqnarray*}
%A_iA_j = \sum_{h=0}^d p^h_{ij}A_h
%\qquad (0 \leq i,j\leq d).
%\end{eqnarray*}
%Determine those $\Phi$ for which 
%$p^h_{ij}\geq 0$ for $0 \leq h,i,j\leq d$.
%\end{problem}
%%extra above

\noindent In order to motivate the next problem
we make a definition.

\begin{definition}
\label{def:int}
\rm
Let $\Phi$ denote the Leonard system
from Definition \ref{eq:ourstartingpt}.
For $0 \leq i\leq d$ we define $A_i=v_i(A)$, where
the polynomial $v_i$ is from Definition
\ref{def:vi1}. Observe that there exist 
scalars $p^h_{ij} \in \K$
$(0 \leq h,i,j\leq d)$ such that 
\begin{eqnarray*}
A_iA_j = \sum_{h=0}^d p^h_{ij}A_h
\qquad (0 \leq i,j\leq d).
\end{eqnarray*}
We call the $p^h_{ij}$ the {\it intersection numbers}
of $\Phi$.
%By the  {\it dual intersection numbers}
%of $\Phi$ we mean the intersection numbers for $\Phi^*$.
\end{definition}

\begin{problem}
\rm
%Assume
% $\K=\R$ and
Let $\Phi$ denote the Leonard system
from Definition \ref{eq:ourstartingpt}.
For each of the Examples
\ref{ex:pa2}--\ref{ex:orphan}, if possible
express each intersection number as a hypergeometric
series or a basic hypergeometric series.
Also for  $\K=\R$,
determine those $\Phi$ for which 
the intersection numbers are all nonnegative.
\end{problem}

\begin{problem}
\rm
Assume $\K=\R$ and let 
 $\Phi$ denote the Leonard system
from Definition \ref{eq:ourstartingpt}.
Determine those $\Phi$ for which the 
intersection numbers of each of $\Phi$, $\Phi^\downarrow$,
 $\Phi^\Downarrow$,
 $\Phi^{\downarrow \Downarrow}$ are all nonnegative.
Also, determine those $\Phi$ for which the intersection
numbers of each relative of $\Phi$ are all nonnegative.
\end{problem}

\begin{problem}
\rm
Assume
 $\K=\R$ and
let $\Phi$ denote the Leonard system
from Definition \ref{eq:ourstartingpt}.
Assume that for each of $\Phi$,$\Phi^*$ the intersection
numbers are nonnegative. Show that the scalars
\begin{eqnarray*}
k_0+k_1+\cdots+k_i-k^*_0-k^*_1-\cdots-k^*_i
\qquad \qquad (0 \leq i \leq d)
\end{eqnarray*}
are all nonnegative or all nonpositive.
\end{problem}

\begin{problem}
\label{prob:conf}
\rm
Assume
 $\K=\R$ 
and let $\Phi$ denote the Leonard system
from Definition \ref{eq:ourstartingpt}. Define
\begin{eqnarray*}
N^h_{ij}:=p^h_{ij} \sqrt{\frac{k_h}{k_ik_j}}
\qquad \qquad (0 \leq h,i,j\leq d),
\end{eqnarray*}
where the $p^h_{ij}$ are from 
Definition \ref{def:int}
 and $k_0, \ldots, k_d$ are from
Definition \ref{def:ki1}. Determine those $\Phi$ for
which (i) $\Phi$ is isomorphic to $\Phi^*$;
and (ii) $N^h_{ij}$ is a nonnegative integer for $0 \leq h,i,j\leq d$.
We remark that such $\Phi$ arise in conformal field theory
\cite{gan}.
\end{problem}

\begin{example}
\rm
Assume $\K=\R$.
Let $\Phi$ denote the Leonard system
from Definition \ref{eq:ourstartingpt} and
let 
 $(\theta_i, \theta^*_i,i=0..d; \varphi_j, \phi_j,j=1..d)$
denote the corresponding parameter array. Assume
\begin{eqnarray*}
\theta_i &=& q^i-q^{d-i},
\\
\theta^*_i &=& q^i-q^{d-i}
\end{eqnarray*}
for $0 \leq i \leq d$ and
\begin{eqnarray}
\varphi_i &=& -(1-q^i)(1-q^{d-i+1})(q^{i-1}+q^{d-i}),
\\
\phi_i &=& (1-q^i)(1-q^{d-i+1})(q^{i-1}+q^{d-i})
\end{eqnarray}
for $1 \leq i \leq d$, where $q$ is a primitive $(2d+4)$th
root of 1.
Then $\Phi$ satisfies the conditions
(i), (ii) of
Problem
\ref{prob:conf}. Moreover $N^i_{i-1,1}=1$ for $1 \leq i \leq d$.
This example is related to the modular data for 
the affine
Kac-Moody algebra $A^{(1)}_1$ at level $d$
\cite[p.~223]{gan}.
\end{example}

\begin{problem} 
\rm
Assume
 $\K=\R$.
Let $\Phi$ denote the Leonard system
from Definition \ref{eq:ourstartingpt} and
let $\theta_0, \theta_1, \ldots,\theta_d$
denote the corresponding eigenvalue sequence.
 Consider the permutation
$\sigma$ of $0,1,\ldots, d$ such that
%\begin{eqnarray*}
$\theta_{\sigma(0)}>
\theta_{\sigma(1)}>
\cdots >
\theta_{\sigma(d)}.$
%\end{eqnarray*}
What are the possibilities for $\sigma$?
\end{problem}

\begin{problem}
\label{prob:abs}
\rm
Assume
 $\K=\R$.
Let $\Phi$ denote the Leonard system
from Definition \ref{eq:ourstartingpt}
 and
let $\theta_0, \theta_1, \ldots,\theta_d$
denote the corresponding eigenvalue sequence.
Let the polynomials
$u_i$ be as in
Definition
\ref{def:ui1}. Find a necessary and sufficient condition
on the parameter array of $\Phi$, so that 
the absolute value
$|u_i(\theta_j)| \leq 1$ for $0 \leq i,j\leq d$.
See \cite[Conjecture 2]{kresh} for an application.
\end{problem}

\begin{problem} 
\rm
Find a short direct proof of
Theorem \ref{eq:mth}.
Such a proof is likely
to lead to an improved proof of 
Theorem
\ref{thm:ls}.
The current proof 
of Theorem
\ref{thm:ls}
is
in  
 \cite{LS99}.
\end{problem}

\begin{problem}
\rm
Let $\Phi$ denote the Leonard system
from Definition
\ref{eq:ourstartingpt} and let $V$ denote an irreducible
$\mathcal A$-module. Let $\langle\,,\,\rangle$  denote
the bilinear form on $V$ from Definition
\ref{def:ip}. We recall $\langle \,,\,\rangle $ is nondegenerate.
What is the Witt index of $\langle \,,\,\rangle $?
The definition of the Witt index is given in
\cite{grove}.
\end{problem}

\begin{problem}
\rm
Let $A_w$ denote the 
Askey-Wilson algebra from Theorem
\ref{th:awgiveslp}.
 An element of $A_w$ is called {\it central}
whenever it commutes with every element of $A_w$. By definition
the {\it center} of $A_w$ is the $\K$-subalgebra of $A_w$
consisting
of the central elements of $A_w$. Describe the center of $A_w$.
Find a generating set for this center. 
By \cite[p.~6]{GYLZmut} the following element of $A_w$ is central:
\begin{eqnarray*}
&&AA^*AA^*-
\beta AA^{*2}A+
A^*AA^*A
-\gamma^* AA^*A
-\gamma (1+\beta) A^*AA^*
-\beta \varrho A^{*2}
\\
&&\quad -(\omega +\gamma \gamma^*)(AA^*+A^*A)
-(\eta^*+\gamma \varrho^*)A
-(\eta(1+\beta)+\gamma^*\varrho)A^*.
\end{eqnarray*}
This can be verified using the Askey-Wilson
relations. Does this element generate the center of 
$A_w$?
\end{problem}

 \begin{problem}
\label{prob:bi-tri}
\rm 
 Let $d$ denote a nonnegative integer.
 Find all Leonard pairs
 $A,A^*$ in $\hbox{Mat}_{d+1}(\K)$ that satisfy
 the following two conditions:
 (i) $A$ is
 irreducible tridiagonal;
 (ii)  $A^*$ is lower bidiagonal with
 $A^*_{i,i-1}=1$ for $1 \leq i \leq d$.
 \end{problem}

 \begin{problem}
\rm 
 Let $d$ denote a nonnegative integer.
 Find all Leonard pairs $A,A^*$  in $\hbox{Mat}_{d+1}(\K)$
 such that each of $A, A^*$ is
 irreducible tridiagonal.
 \end{problem}

 \begin{problem}
\label{prob:diagzero}
\rm 
 Let $d$ denote a nonnegative integer.
 Find all Leonard pairs $A,A^*$  in $\hbox{Mat}_{d+1}(\K)$
 such that 
 each of $A, A^*$ is
 irreducible tridiagonal with all diagonal entries 0. 
Note that in this case the Leonard pair $A,A^*$ is isomorphic to
the Leonard pair $-A,-A^*$. 
%Brown
\end{problem}

\begin{problem}
\label{prob:bipdbip}
\rm
Let $A,A^*$ denote a Leonard pair of diameter $d$,
such that
$a_i=0$ and $a^*_i=0$ for $0 \leq i \leq d$.
Find all the bases for the underlying vector space,
with respect to which each of $A,A^*$ is tridiagonal
with all diagonal entries zero.
Show that such a basis is an eigenbasis for
$qAA^*-q^{-1}A^*A$ or
$q^{-1}AA^*-qA^*A$, where
$q^2+q^{-2}+1$ is the common value of
(\ref{eq:betaplusone}).
%Brown
\end{problem}

 \begin{problem}
\label{prob:lt}
 \rm
 Let $V$ denote a vector space over $\K$ with finite positive
 dimension. By a {\it Leonard triple} on $V$, we mean a three-tuple
 of linear transformations
 $A:V\rightarrow V$, $A^*:V\rightarrow V$,
  $A^\varepsilon:V\rightarrow V$
 that 
  satisfy conditions (i)--(iii) below.
   \begin{enumerate}
   \item There exists a basis for $V$ with respect to which
   the matrix representing $A$ is
   diagonal and 
   the matrices representing $A^*$ and $A^{\varepsilon}$ are each irreducible
   tridiagonal.
   \item There exists a basis for $V$ with respect to which
   the matrix representing $A^*$ is
   diagonal and 
   the matrices representing $A^\varepsilon $ and $A$ are each irreducible
   tridiagonal.
  \item There exists a basis for $V$ with respect to which
   the matrix representing $A^{\varepsilon}$ is
   diagonal and 
   the matrices representing $A$ and $A^*$ are each irreducible
   tridiagonal.
  \end{enumerate}
   Find all the Leonard triples. See \cite{mlt} for a connection
   between Leonard triples and spin models.
   \end{problem}

\begin{conjecture}
 Let $V$ denote a vector space over $\K$ with finite
positive dimension. 
Given three linear transformations
 $A:V\rightarrow V$, $A^*:V\rightarrow V$,
 and  $A^\varepsilon:V\rightarrow V$,
if each two-element subset of $A,A^*,A^\varepsilon$ is
a Leonard pair on $V$ then
 $A,A^*,A^\varepsilon$ is a Leonard triple on $V$.
\end{conjecture}

\begin{problem}
\rm Let $V$ denote a vector space over $\K$ with finite
positive dimension. Let $\mbox{End}(V)$
denote the $\K$-algebra consisting of all linear transformations
from $V$ to $V$.
Let $A,A^*,A^\varepsilon$ denote
a Leonard triple on $V$. Each of the pairs
 $A,A^*$; 
 $A,A^\varepsilon$; 
 $A^*,A^\varepsilon$  is a Leonard pair on $V$; let 
 $r,s,t$ denote the corresponding antiautomorphisms
of $\mbox{End}(V)$ from Definition
\ref{def:dag}. Determine the subgroup of $\mbox{GL(End}(V))$
generated by $r,s,t$. Since
$r^2=s^2=t^2=1$, it is conceivable that this subgroup is
a Coxeter group. For which Leonard triples is this the case?
\end{problem}

\begin{problem}
\rm
Let $V$ denote a vector space over $\K$ with
finite positive dimension
and let $A,A^*,A^\varepsilon$ denote a Leonard triple on $V$.
Show that for any permutation $x,y,z$ of 
$A,A^*,A^\varepsilon$ there exists an antiautomorphism
$\sigma $
of $\mbox{End}(V)$ such that $x^\sigma=x$ and
each of $\lbrack x,y\rbrack^\sigma, \lbrack x,z\rbrack^\sigma$ is
a scalar multiple of the other. Here
$\lbrack r,s\rbrack$ means $rs-sr$.
\end{problem}

\begin{problem}
\label{prob:krawtriple}
\rm 
Assume $\K$ is algebraically closed with characteristic 0.
Let $d$ denote a nonnegative integer and let
 $A,A^*$ denote the Leonard pair on $\K^{d+1}$
given in 
(\ref{eq:fam1}).
Find all the matrices $A^{\varepsilon}$ such that
$A,A^*, A^\varepsilon$ is a Leonard triple on
$\K^{d+1}$. Given a solution
$A^\varepsilon$, show that each of
\begin{eqnarray*}
\lbrack A,A^*\rbrack,
\qquad \lbrack A^*,A^\varepsilon\rbrack,
\qquad \lbrack A^\varepsilon,A\rbrack
\end{eqnarray*}
is contained in the $\K$-linear span of
$I, A,A^*, A^\varepsilon$.
Here
$\lbrack r,s\rbrack$ means $rs-sr$.
\end{problem}

\begin{problem}
\label{prob:qc}
\rm
Let $V$ denote a vector space over $\K$ with finite
positive dimension and let $A,A^*,A^\varepsilon$ denote
a Leonard triple on $V$. Show that there exists a nonzero
scalar $q \in \K$ such that each of
\begin{eqnarray*}
AA^*-qA^*A, 
\qquad \quad 
A^*A^\varepsilon-qA^\varepsilon A^*, 
\qquad \quad 
A^\varepsilon A-qAA^\varepsilon 
\end{eqnarray*}
is contained in the $\K$-linear span of
$I, A,A^*, A^\varepsilon $.
\end{problem}

\begin{problem}
\label{prob:balg}
\rm
Assume $\K$ is algebraically closed.
Let $q$ denote a nonzero scalar 
in $\K$ that is not a root of
unity.
Let $B$ denote the unital associative $\K$-algebra 
with generators $x,y,z$ 
and the following relations. The relations are that
each of 
\begin{eqnarray*}
\frac{q^{-1}xy-qyx}{q^2-q^{-2}}-z,
\qquad \qquad 
\frac{q^{-1}yz-qzy}{q^2-q^{-2}}-x,
\qquad \qquad 
\frac{q^{-1}zx-qxz}{q^2-q^{-2}}-y
\end{eqnarray*} is central in $B$.
Let $V$ denote a
finite dimensional irreducible $B$-module
on which each of $x,y,z$ is multiplicity-free.
Show that $x, y, z$ act on $V$ as a Leonard triple.
Determine all the $B$-modules of this type, up to isomorphism.
\end{problem}

\begin{problem}
\rm
Classify up to isomorphism the
finite dimensional irreducible $B$-modules,
where the algebra $B$ is from
Problem
\ref{prob:balg}.
This problem is
closely related to Problem
\ref{pr:awsolve}.
\end{problem}

%%%%%%%%%problem below solved spring 06 by Nomura and Ter
%   \begin{conjecture}
%   \label{conj:w1111}
%   Let $\Phi$ denote the Leonard system
%   from
%   Definition \ref{eq:ourstartingpt} and let $I$ denote
%   the identity element of $\mathcal A$. Then
%    for all
%    $X \in {\mathcal A}$ the following are equivalent:
%    (i) both
%    \begin{eqnarray*}
%    &&E_iX E_j = 0  \quad
%    \hbox{if} \quad |i-j|>1,
%    \qquad \qquad (0 \leq i,j\leq d),
%    \\
%    &&E^*_iX E^*_j = 0 \quad  \hbox{if} \quad |i-j|>1,
%    \qquad \qquad (0 \leq i,j\leq d);
%    \end{eqnarray*}
%    (ii) $X$ is a $\K$-linear combination of
%    $I, A, A^*, AA^*, A^*A.$
%    \end{conjecture}
%%%%%%%%%%%%%%%%%%%%%%%%%%%%%%%%%
\begin{problem}
\rm
\label{conj:span}
Referring to the tridiagonal pair $A,A^*$ in Definition
   \ref{def:tdp}, consider the space of all linear transformations $X:V\to V$
such that both
\begin{eqnarray*}
X V_i &\subseteq & V_{i-1}+V_i+V_{i+1} \qquad \qquad (0 \leq i \leq d),
\\
X V^*_i &\subseteq & V^*_{i-1}+V^*_i+V^*_{i+1} \qquad \qquad (0 \leq i \leq d).
\end{eqnarray*}
Find a basis for this space.
\end{problem}

\begin{note} \rm Problem
\ref{conj:span}
is solved if $A,A^*$ is a Leonard pair \cite{NT:span}.
\end{note}

    \begin{conjecture}
\label{conj:genset}
Let $\Phi$ denote the Leonard system from
    Definition \ref{eq:ourstartingpt}.
    Then for $0 \leq  r \leq d$
    the elements
    \begin{eqnarray*}
    E^*_0, E^*_1, \ldots, E^*_r, E_r, E_{r+1}, \ldots, E_d
    \end{eqnarray*}
    together generate $\mathcal A$.
\end{conjecture}

\begin{remark}
\rm
Conjecture \ref{conj:genset} holds for $r=0$
by Corollary
\ref{cor:genset}, and since 
$A$ is a linear combination
of $E_0, E_1, \ldots, E_d$.
Similarly 
Conjecture \ref{conj:genset} holds for $r=d$.
\end{remark}

\begin{problem} 
\label{prob:psipre}
\rm
Referring to the tridiagonal pair $A,A^*$ in Definition
   \ref{def:tdp},  
find all the linear transformations $X:V\to V$
such that both
\begin{eqnarray*}
X V^*_i &\subseteq & V^*_0+V^*_1+\cdots+V^*_{i-1} \qquad \qquad (0 \leq i \leq d),
\\
X V_i &\subseteq & V_{i-1}+V_i+V_{i+1} \qquad \qquad (0 \leq i \leq d).
\end{eqnarray*}
%Vidar is considering for $d=2$.
\end{problem}

\begin{problem}
\label{prob:psi}
\rm
Referring to the tridiagonal pair $A,A^*$ in Definition
   \ref{def:tdp},  
find all the linear transformations $X:V\to V$
such that both
\begin{eqnarray*}
X V^*_i &\subseteq & V^*_0+V^*_1+\cdots+V^*_i \qquad \qquad (0 \leq i \leq d),
\\
X V_i &\subseteq & V_{i-1}+V_i+V_{i+1} \qquad \qquad (0 \leq i \leq d).
\end{eqnarray*}
%Vidar is considering for $d=2$.
\end{problem}

\begin{problem}
\rm
Recall the algebra
$U_q({\widehat{sl}}_2)^{>0}$
from Example
\ref{ex:qserre}. Describe the center of
$U_q({\widehat{sl}}_2)^{>0}$.
Find a generating set for this center.
We remark that
$U_q({\widehat{sl}}_2)^{>0}$
has infinite dimension as
a vector space over $\K$. A basis for this vector space
is given in
\cite[Theorem 2.29]{shape}.
\end{problem}

\begin{problem}
\rm
Assume $\K$ is algebraically closed.
Let $V$ denote a vector space over $\K$ with 
finite positive dimension
and let $A,A^*$ denote a tridiagonal pair on $V$.
Compute the Jordan Canonical form for 
$q^{-1}AA^*-qA^*A$, where $q^2+q^{-2}=\beta$ and
$\beta $ is from Theorem 
\ref{thm:tdrels}.
\end{problem}

%%%%%conjecture below solved by Nomura Summer 05
%%the determinant ..
%\begin{conjecture}
%\label{conj:det}
%Let 
%$\Phi=(A;A^*; \lbrace E_i\rbrace_{i=0}^d;  
%\lbrace E^*_i\rbrace_{i=0}^d)$ denote a Leonard system
%and let $(\theta_i, \theta^*_i,i=0..d; \varphi_j, \phi_j,j=1..d)$
%denote the corresponding parameter array.
%Assume $q\not=1$, $q\not=-1$,
%where $q+q^{-1}+1$ is the common value of
%(\ref{eq:betaplusone}).
%For $d$ odd, 
%\begin{eqnarray*}
%det(AA^*-A^*A)= \prod_{{1 \leq i \leq d} \atop{i \; odd}}
%\varphi_i\phi_i\,\frac{q-1}{1-q^i}\, \frac{1-q^{-1}}{1-q^{-i}}.
%\end{eqnarray*}
%\end{conjecture}
%%%%%%%%%%%%%%%
%\begin{remark}
%\label{rem:deven}
%\rm
%Referring to 
%Conjecture \ref{conj:det},
%for $d$ even we have
%$\mbox{det}(AA^*-A^*A)=0$.
%This is because
%$\mbox{det}(X^\dagger)=
% \mbox{det}(X)$ and 
%$X^\dagger=-X$, where
%$X=AA^*-A^*A$ and where
%$\dagger$ is the antiautomorphism from
%Definition
%\ref{def:dag}.
%\end{remark}
%%%%%%%%%%%%%%%
%problem below solved by Nomura Summer 05
%%the determinant..
%\begin{problem}
%\rm
%Referring to
%Conjecture \ref{conj:det} and
%Remark
%\ref{rem:deven}, show that for $d$ even
%the null space of $AA^*-A^*A$ has dimension 1. 
%Find a basis vector for this null space.
%Express this basis vector in terms of
%a $\Phi$-standard basis and a $\Phi^*$-standard 
%basis.
%\end{problem}
%%%%%%%%%%%%%%%

\begin{problem} 
\label{prob:invpair}
\rm
Let $V$ denote a vector space over $\K$ with finite
positive dimension. By an {\it inverting pair} on $V$
we mean an ordered pair of invertible linear transformations
$K:V\to V$ and $K^*:V\to V$ that satisfy both (i), (ii) below.
\begin{enumerate}
\item There exists a basis for $V$ with respect to which
the matrix
representing $K$ has all entries 0 above the superdiagonal,
the matrix representing $K^{-1}$ has all entries 0 below
the subdiagonal, and
the  matrix representing $K^*$ is diagonal.
\item There exists a basis for $V$ with respect to which
the matrix
representing $K^*$ has all entries 0 above the superdiagonal,
the matrix representing $K^{*-1}$ has all entries 0 below
the subdiagonal, and
the  matrix representing $K$ is diagonal.
\end{enumerate}
Find all the inverting pairs. See  
\cite[Section 11]{tdanduq} for a connection between inverting pairs and
 Leonard pairs of $q$-geometric type.
\end{problem}

%%Nomura solved this Spring 05, but it turned out
%to be already solved 
%\begin{problem}
%\rm
%For an integer $d\geq 0$, 
%find all the invertible matrices
%$K \in 
%\hbox{Mat}_{d+1}(\K)$ such that (i) $K$ has 
%all entries 0 above the superdiagonal,
%(ii) $K^{-1}$ has all entries 0 below the subdiagonal.
%\end{problem}
%%%%%%%%%%%%%%%%%%%

\begin{problem}
\label{prob:kki}
\rm Let $V$ denote a vector space over $\K$ with finite
positive dimension and let $A,A^*$ denote a Leonard pair
on $V$. To avoid degenerate situations
assume $q$ is not a root of unity, where $q+q^{-1}=\beta$
and $\beta $ is from 
Theorem
\ref{thm:tdrels}.
Find all the invertible linear transformations
$K:V\rightarrow V$ such that (i) $A,K$ is
a Leonard pair, and (ii) $A^*, K^{-1}$ is a Leonard pair.
\end{problem}

\begin{problem}
\label{prob:raise}
\rm 
Let $d$ denote a nonnegative integer and 
let $\theta_0, \theta_1, \ldots, \theta_d$ denote a sequence of
mutually distinct scalars in $\K$.
Let $\lambda $ denote an indeterminate and
let $V$ denote the vector space over $\K$ consisting
of all polynomials in $\lambda$ that have degree at most $d$.
Define
a polynomial $\tau_i=\prod_{h=0}^{i-1}(\lambda-\theta_h)$
for $0 \leq i \leq d$ and observe
$\tau_0, \tau_1, \ldots, \tau_d$ is a basis for $V$.
Define
$\rho_i=\prod_{h=0}^{i-1}(\lambda-\theta_{d-h})$
for $0 \leq i \leq d$ and observe
$\rho_0, \rho_1, \ldots, \rho_d$ is a basis for $V$.
By a {\it lowering} map on $V$ we mean
a linear transformation $\Psi: V\to V$ that satisfies both
\begin{eqnarray*}
&&\Psi \tau_i \in \mbox{span}(\tau_{i-1}) \qquad (1 \leq i \leq d),
\qquad \qquad \Psi\tau_0=0,
\\
&&\Psi \rho_i \in \mbox{span}(\rho_{i-1}) \qquad (1 \leq i \leq d),
\qquad \qquad \Psi\rho_0=0.
\end{eqnarray*}
Show that there exists a nonzero lowering map on $V$ provided 
$(\theta_{i-2}-\theta_{i+1})(\theta_{i-1}-\theta_i)^{-1}$
is independent of $i$ for $2 \leq i \leq d-1$. 
\end{problem}

\begin{remark} 
\rm (Vidunas 2007)
Referring to Problem \ref{prob:raise},
the existence of a nonzero lowering map does not imply
that 
$(\theta_{i-2}-\theta_{i+1})(\theta_{i-1}-\theta_i)^{-1}$
is independent of $i$ for $2 \leq i \leq d-1$. 
Here is an example with $d=5$.
For $b,c \in \K$ define
\begin{eqnarray*}
\theta_0=0,
\quad
\theta_1=b+c,
\quad
\theta_2=b,
\quad
\theta_3=c,
\quad
\theta_4=b+c-1,
\quad
\theta_5=1,
\end{eqnarray*}
with $b,c$ chosen so that
$\lbrace \theta_i \rbrace_{i=0}^5$ are mutually distinct.
Note that
\begin{eqnarray*}
\frac{\theta_0-\theta_3}{\theta_1-\theta_2}=-1,
\qquad 
\frac{\theta_1-\theta_4}{\theta_2-\theta_3}=\frac{-1}{c-b},
\qquad 
\frac{\theta_2-\theta_5}{\theta_3-\theta_4}=-1,
\end{eqnarray*}
and that $c-b\not=1$ in general. Note also that
there exists a lowering map $\Psi$ on $V$ that sends
\begin{eqnarray*}
&& \tau_0 \mapsto 0,
\quad 
\tau_1 \mapsto \tau_0,
\quad 
\tau_2 \mapsto 0,
\quad 
\tau_3 \mapsto \tau_2,
\quad 
\tau_4 \mapsto 0,
\quad 
\tau_5 \mapsto \tau_4,
\\
&&\rho_0 \mapsto 0,
\quad
\rho_1 \mapsto \rho_0,
\quad
\rho_2 \mapsto 0,
\quad
\rho_3 \mapsto \rho_2,
\quad
\rho_4 \mapsto 0,
\quad
\rho_5 \mapsto \rho_4.
\end{eqnarray*}
The transition matrix from the basis 
$\lbrace \tau_i\rbrace_{i=0}^5$ to
the basis $\lbrace \rho_i\rbrace_{i=0}^5$ is
\begin{eqnarray*}
\left(
\begin{array}{ c c c c c c }
1 & -1  &  a    & -ac & abc & -a(a+1)bc  \\
0 & 1  &  0   &  a & 0 &abc    \\
0 &0 & 1&  b-c  & a   & -ac \\
0 & 0 & 0  & 1  & 0  & a \\
0 & 0 & 0  & 0  & 1  & -1 \\
0 & 0 & 0  & 0  & 0  & 1 
\end{array}
\right),
\end{eqnarray*}
where $a=b+c-1$.
%for example $\rho_1 = \tau_1-\tau_0$.
\end{remark}

\begin{problem}
\rm
By a {\it generalized Leonard system} in $\mathcal A$
we mean a sequence \\
$(A;A^*;\lbrace \theta_i\rbrace_{i=0}^d;\lbrace \theta^*_i \rbrace_{i=0}^d)$
that satisifies  (i)--(v) below.
\begin{enumerate}
\item $A,A^* \in {\mathcal A}$ and 
$\theta_i, \theta^*_i \in \K$ for $0 \leq i \leq d$.
\item $\theta_0, \theta_1,\ldots,\theta_d$ is an ordering of the 
roots of the characteristic polynomial of $A$.
\item $\theta^*_0, \theta^*_1,\ldots,\theta^*_d$ is an ordering of the 
roots of the characteristic polynomial of $A^*$.
\item For $0 \leq i ,j\leq d$,
\begin{eqnarray*}
\tau_i(A)A^*\eta_{d-j}(A) &=& \cases{0, &if $\;i-j > 1$;\cr
\not=0, &if $\; i-j  = 1$,\cr}
\\
\eta_{d-i}(A)A^*\tau_j(A) &=& \cases{0, &if $\;j-i > 1$;\cr
\not=0, &if $\;j-i  = 1$.\cr}
\end{eqnarray*}
\item For $0 \leq i ,j\leq d$,
\begin{eqnarray*}
\tau^*_i(A^*)A\eta^*_{d-j}(A^*) &=& \cases{0, &if $\;i-j > 1$;\cr
\not=0, &if $\; i-j  = 1$,\cr}
\\
\eta^*_{d-i}(A^*)A\tau^*_j(A^*) &=& \cases{0, &if $\;j-i > 1$;\cr
\not=0, &if $\; j-i  = 1$.\cr}
\end{eqnarray*}
\end{enumerate}
(We are using the notation
 (\ref{eq:taus}),
 (\ref{eq:etas})).
We are not assuming $\theta_0, \theta_1, \ldots,\theta_d$
are mutually distinct or that
 $\theta^*_0, \theta^*_1, \ldots,\theta^*_d$
are mutually distinct.
Classify the generalized Leonard systems.
Extend the theory of Leonard systems to the level
of generalized Leonard systems.
\end{problem}

\begin{problem}
\label{prob:4vand}
\rm
For an integer $d\geq 0$
and for $X \in 
\hbox{Mat}_{d+1}(\K)$, we define $X$ to be {\it north Vandermonde}
whenever the entries
$X_{ij}=X_{0j}f_i(\theta_j)$ for $0 \leq i,j\leq d$,
where $\theta_0, \theta_1, \ldots, \theta_d$
are mutually distinct scalars in $\K$ and $f_i \in \K\lbrack \lambda 
\rbrack $ has degree  $i$ for $0 \leq i \leq d$.
Let $X' \in
\hbox{Mat}_{d+1}(\K)$
 denote the matrix obtained 
by rotating $X$ counterclockwise 90 degrees.
We define $X$ to be {\it east Vandermonde}
(resp. {\it south Vandermonde})
(resp. {\it west Vandermonde})
whenever $X'$ (resp. $X''$)  (resp. $X'''$) is north Vandermonde.
Find all the matrices in
$\hbox{Mat}_{d+1}(\K)$ that are simultaneously
north, south, east, and west Vandermonde.
\end{problem}

\begin{definition}
\label{def:cycliclp}
\rm
Let $V$ denote a vector space over $\K$ with finite
positive dimension $n$. By a {\it cyclic Leonard pair}
on $V$, we mean an ordered pair of linear transformations
$A: V\to V$ and $A^*:V\to V$ that satisfy (i)--(iv) below.
    \begin{enumerate}
    \item Each of $A,A^*$ is multiplicity free. 
    \item There exists a bijection
    $i \to V_i$ from the cyclic group $\Z / n\Z$
    to the set of eigenspaces of $A$ such that
    \begin{eqnarray*}
    A^* V_i \subseteq V_{i-1} + V_i+ V_{i+1} \qquad 
    \qquad (\forall i \in \Z /n \Z).
    \end{eqnarray*}
    \item There exists a bijection
    $i \to V^*_i$ from  $\Z /n\Z$
    to the set of eigenspaces of $A^*$ such that
    \begin{eqnarray*}
    A V^*_i \subseteq V^*_{i-1} + V^*_i+ V^*_{i+1} \qquad 
    \qquad (\forall i \in \Z /n \Z).
    \end{eqnarray*}
    \item There does not exist a 
    subspace $W$ of $V$ such  that $AW\subseteq W$,
    $A^*W\subseteq W$, $W\not=0$, $W\not=V$.
    \end{enumerate}
\end{definition}

\begin{example}
\label{ex:cycliclp}
\rm
For an integer $d\geq 3$ 
let $\lbrace \theta_i\rbrace_{i=0}^d$
and 
$\lbrace \theta^*_i\rbrace_{i=0}^d$ denote scalars in $\K$ that
satisfy the conditions PA1 and PA5 of
Definition
\ref{def:pa}.
 Assume $q$ is a primitive $n^{th}$ root of unity,
 where $q+q^{-1}+1$ is the common value of
(\ref{eq:betaplusone}) and $n=d+1$. Pick distinct $b,c\in \K$ 
and define
\begin{eqnarray*}
\varphi_i = (q^i-1)(b-cq^{-i}) + 
(\theta^*_i-\theta^*_0)
(\theta_{i-1}-\theta_d)
\qquad (1 \leq i \leq d),
\end{eqnarray*}
with $b,c$ chosen such that $\varphi_i \not=0$ for
$1 \leq i \leq d$. 
Let $A$ (resp. $A^*$)
denote the matrix in 
$\hbox{Mat}_{d+1}(\K)$
with entries $A_{ii} = \theta_i$
(resp. $A^*_{ii} = \theta^*_i$)
for $0 \leq i \leq d$,
$A_{i,i-1}=1$ 
(resp. $A^*_{i-1,i}=\varphi_i$)
for $1 \leq i \leq d$, and all
other entries 0.
Show that the pair $A,A^*$
acts on 
$V=\K^n$
as a cyclic Leonard pair.
Show that $A,A^*$ satisfy the Askey-Wilson relations
(\ref{askwil1}),
(\ref{askwil2}).
\end{example}

\begin{problem}
\label{prob:cyclicaw}
\rm Recall the tridiagonal algebra $T$ from
Theorem
     \ref{def:tdalg}.
Let $V$ denote a finite-dimensional irreducible $T$-module
on which the generators $A$, $A^*$ are multiplicity-free.
Show that the pair $A,A^*$ acts on $V$ as a cyclic Leonard pair.
\end{problem}

\begin{problem}
\label{prob:cycliclp}
\rm
 \noindent Classify the cyclic Leonard pairs. Extend the theory
  of Leonard pairs to the level of cyclic Leonard pairs.
Does a cyclic Leonard pair satisfy the Askey-Wilson relations
(\ref{askwil1}),
(\ref{askwil2})?
If not, does it satisfy the tridiagonal relations
(\ref{eq:td1}), 
(\ref{eq:td2})? 
\end{problem}

\begin{problem}
\label{prob:tdalgdiag}
\rm 
Referring to the tridiagonal algebra
$T$ from
Theorem
     \ref{def:tdalg}, give examples of finite-dimensional
irreducible $T$-modules on which the generators $A,A^*$
are not diagonalizable. If possible, classify all such modules.
\end{problem}

\begin{problem}
\label{prob:tdalgd2}
\rm
Referring to the tridiagonal algebra
$T$ from
Theorem
     \ref{def:tdalg}, let $V$
denote 
a finite-dimensional
irreducible $T$-module on which the generators $A,A^*$
are diagonalizable, but do not form a tridiagonal pair.
Show that for $V$ the eigenspaces of $A$ and the
 eigenspaces of $A^*$
all have the same dimension.
\end{problem}

\begin{problem}
\label{prob:drin}
\rm
Referring to the tridiagonal pair $A,A^*$ in
Definition
   \ref{def:tdp},  
for $0 \leq i \leq d$ let $\theta_i$ (resp. $\theta^*_i$)
denote the eigenvalue of $A$ (resp. $A^*$) associated with
the eigenspace $V_i$ (resp. $V^*_i$).
Assume $V^*_0$ has dimension 1.
Observe that for $0 \leq i \leq d$ the space $V^*_0$ is an
eigenspace for 
\begin{eqnarray*}
(A^*-\theta^*_1I)
(A^*-\theta^*_2I)
\cdots 
(A^*-\theta^*_iI)
(A-\theta_{i-1}I)
\cdots 
(A-\theta_1I)
(A-\theta_0I);
\end{eqnarray*}
let $\zeta_i$ denote the corresponding eigenvalue.
Show that the tridiagonal pair $A,A^*$ is determined
up to isomorphism 
by the array
$(\lbrace \theta_i\rbrace_{i=0}^d, \lbrace \theta^*_i\rbrace_{i=0}^d, 
\lbrace \zeta_i\rbrace_{i=0}^d)$. We will call this array
a {\it parameter array} for $A,A^*$.
\end{problem}

\begin{problem}
\label{prob:cntrex}
\rm
Referring to Example
\ref{ex:cycliclp}, choose $b,c$ such that
$\varphi_i \not=0$ for $1 \leq i \leq d$ and
\begin{eqnarray*}
0 \not = \sum_{i=0}^d \varphi_1 \varphi_2 \cdots \varphi_i 
 p_{i+1}p_{i+2}\cdots p_d,
\end{eqnarray*}
where
 $p_j =(\theta_0-\theta_j)
(\theta^*_0-\theta^*_j)$ for $1 \leq j \leq d$.
Show that there exists a tridiagonal pair
with parameter array
$(\lbrace \theta_i\rbrace_{i=0}^d, \lbrace \theta^*_i\rbrace_{i=0}^d, 
\lbrace \varphi_1 \varphi_2 \cdots \varphi_i\rbrace_{i=0}^d)$.
Show that this tridiagonal pair is not a Leonard pair.
\end{problem}

\begin{problem}
\label{prob:tddual}
\rm
Referring to the tridiagonal pair $A,A^*$ in
Definition
\ref{def:tdp}, let $V'$ denote the dual space
of $V$; we recall $V'$ is the vector space over
$\K$ consisting of all $\K$-linear transformations
from $V$ to $\K$.
Define linear transformations
$B:V'\rightarrow V'$  and
$B^*:V'\rightarrow V'$ 
as follows: For $f \in V'$, $Bf$ (resp. $B^*f$) is that element of
$V'$ such that $(Bf)(v)= f(Av)$ (resp. $(B^*f)(v)=f(A^*v)$)
for all $v \in V$.
Show that the pair $B,B^*$ is a tridiagonal pair
on $V'$ which is isomorphic to $A,A^*$.
\end{problem}

\begin{problem}
\label{prob:bilform}
\rm
Referring to the tridiagonal pair $A,A^*$ in
Definition
\ref{def:tdp}, show that there exists
a nonzero symmetric bilinear form
$\langle \,,\,\rangle $ on $V$ such that
$\langle Au,v\rangle = 
\langle u,Av\rangle $
and 
$\langle A^*u,v\rangle = 
\langle u,A^*v\rangle $
for all $u,v \in V$.
Show that this form is nondegenerate 
and unique up to multiplication
by a nonzero scalar in $\K$.
For the case in which $A,A^*$ is a Leonard pair
the bilinear form is constructed in Section 15.
\end{problem}

\begin{problem}
\rm
Classify the tridiagonal pairs that have 
shape vector $(\rho,\rho, \ldots, \rho)$,
where $\rho$ is an integer at least 2.
See Section 34 for the definition of the shape vector.
\end{problem}

\begin{problem}
\label{prob:zeroint}
\rm
 Let 
$(A;A^*; \lbrace E_i\rbrace_{i=0}^d;  
\lbrace E^*_i\rbrace_{i=0}^d) $
denote a Leonard system in $\mathcal A$  and let
$V$ denote an irreducible $\mathcal A$-module.
For an integer $0 \leq r\leq d$ and for a subset
$L$ of $\lbrace 0,1,\ldots, d\rbrace$, show that
$\sum_{i=0}^r E^*_iV$ and
$\sum_{j \in L}E_jV$
 have zero intersection
if the cardinality of $L$ is at most $d-r$.
\end{problem}

\begin{problem}
\rm
Referring to the tridiagonal pair in Definition
   \ref{def:tdp}, find all the  linear transformations
   $S:V\to V$ such that $SV_i=V_i$
   and
  $S(V^*_0+V^*_1+\cdots + V^*_i)=
  V^*_d+V^*_{d-1}+\cdots + V^*_{d-i}$
   for $0 \leq i \leq d$.
Is $S$ diagonalizable? If so, find a basis for each
$V_i$ consisting of eigenvectors for $S$.
\end{problem}

%%%%%%%%%special case of Nomura ..height ..
%\begin{problem}
%\rm
%Let $A,A^*$ denote a tridiagonal pair
%with shape vector $(\rho_0, \rho_1, \ldots, \rho_d)$
%(see Section 34). Show that
%if each of $\rho_0, \rho_1$ is $1$ then
%$\rho_i=1$ for $0 \leq i \leq d$ or in other words
%$A,A^*$ is a Leonard pair.
%\end{problem}
%%%%%%%%%%%%%%%%%%%%%%above is special case of Nomure result%%%
%%%%%solved by Nomura and Terwilliger fall 06
%\begin{problem}
%\rm
%Referring to the tridiagonal pair in Definition
%   \ref{def:tdp}, 
%  let $\mathcal D$ denote the subalgebra
%  of $\mbox{\rm End}(V)$ generated by
%  $A$. Show that the following are equivalent.
%\begin{enumerate}
%\item There exists an element  $S\in \mathcal D$
%such that $SE^*_0V\subseteq E^*_dV$.
%\item Each of $V_i, V^*_i$ has dimension 1 for
%$0 \leq i \leq d$ or in other words the tridiagonal
%pair is a Leonard pair.
%\end{enumerate}
%\end{problem}
%%%%%%%%%%%%%%%%

\begin{problem}
\label{prob:tmat}
\rm
Let $\Phi=(A;A^*; \lbrace E_i\rbrace_{i=0}^d;  
\lbrace E^*_i\rbrace_{i=0}^d) $
denote a Leonard system in $\mathcal A$  and let
$V$ denote an irreducible 
$\mathcal A$-module.
Define
\begin{eqnarray*}
T=\nu^2\sum_{i=0}^d k^{-1}_i E_iE^*_0E_0E^*_i.
\end{eqnarray*}
Show that $T$ sends a $\Phi$-standard basis for $V$ to a
$\Phi^*$-standard basis for $V$. Show that $AT=TA^*$
if and only if the eigenvalue sequence of $\Phi$ and the
dual eigenvalue
sequence of $\Phi$ coincide.
\end{problem}

\begin{problem}
\label{prob:ntsplit}
\rm
Let $V$ denote a vector space over $\K$ with finite positive
dimension.
Consider a pair of linear transformations $A:V\to V$
and $A^*:V\to V$ that satisfy (i)--(iv) below:
    \begin{enumerate}
    \item Each of $A,A^*$ is diagonalizable.
    \item There exists an ordering $V_0, V_1,\ldots, V_d$ of the
    eigenspaces of $A$ such that
    \begin{eqnarray*}
    A^* V_i \subseteq V_0 + V_1+\cdots+ V_{i+1} \qquad \qquad (0 \leq i \leq d),
    \end{eqnarray*}
    where $V_{-1} = 0$, $V_{d+1}= 0$.
    \item There exists an ordering $V^*_0, V^*_1,\ldots, V^*_\delta$ of
    the
    eigenspaces of $A^*$ such that
    \begin{eqnarray*}
    A V^*_i \subseteq V^*_{0} + V^*_1+ \cdots+ V^*_{i+1} \qquad \qquad
    (0 \leq i \leq 
    \delta),
    \end{eqnarray*}
    where $V^*_{-1} = 0$, $V^*_{\delta+1}= 0$.
    \item There does not exist a 
    subspace $W$ of $V$ such  that $AW\subseteq W$,
    $A^*W\subseteq W$, $W\not=0$, $W\not=V$.
\end{enumerate}
\noindent Show that $d=\delta$. Now define
\begin{eqnarray*}
U_i = (V^*_0+V^*_1+\cdots +V^*_i)\cap (V_0+V_1+\cdots + V_{d-i})
\qquad \qquad (0 \leq i \leq d).
\end{eqnarray*}
\noindent Show that 
\begin{eqnarray*}
V = U_0+U_1+\cdots + U_d \qquad \qquad (\mbox{\rm direct sum}).
\end{eqnarray*}
\noindent Show that for $0 \leq i \leq d$ the dimensions of
$V_i$, $V^*_i$, $U_i$ coincide.
Let $\theta_i$ (resp. $\theta^*_i$) denote the eigenvalue
of $A$ (resp. $A^*$) associated with $V_i$ (resp. $V^*_i$).
Show that both
\begin{eqnarray*}
(A-\theta_{d-i}I)U_i &\subseteq & U_{i+1} \qquad \qquad (0 \leq i \leq d),
\\
(A^*-\theta^*_iI)U_i &\subseteq & U_{i-1} \qquad \qquad (0 \leq i \leq d), 
\end{eqnarray*}
where $U_{-1}=0$ and $U_{d+1}=0$.
\end{problem}

\begin{problem}
\label{prob:vs}
\rm Let $\lbrace \theta_i\rbrace_{i=0}^d$ denote a finite
sequence of mutually distinct scalars in $\K$ and assume
\begin{eqnarray}
\label{eq:3tr}
\frac{\theta_{i-2}-\theta_{i+1}}{\theta_{i-1}-\theta_i}
\end{eqnarray}
is 
independent of $i$ for $2\leq i \leq d-1$.
Let $\beta+1$ denote the common value of 
(\ref{eq:3tr}).
Consider the set of all sequences 
$(
 \lbrace \theta^*_i\rbrace_{i=0}^d;
 \lbrace \varphi_i\rbrace_{i=1}^d;
 \lbrace \phi_i\rbrace_{i=1}^d)$ that satisfy
(PA3), (PA4), and 
\begin{eqnarray*}
\mbox{\rm (PA5')}: \qquad \qquad  \theta^*_{i-2}-\theta^*_{i+1}
=(\beta+1)(\theta^*_{i-1}-\theta^*_i)
\qquad\qquad (2 \leq i \leq d-1).
\end{eqnarray*}
Observe that this set is a vector space over $\K$. Show that the
dimension of this vector space is $4$ provided $d\geq 2$.
Find an attractive basis for this vector space.
Find a geometric interpretation of this vector space.
How is this problem related to
Problem \ref{prob:psi}?
\end{problem}

\begin{problem} \rm
\label{prob:ae}
Let $A,A^*$ denote a Leonard pair on $V$. Find all
the 
linear transformations $A^{\varepsilon}:V\to V$
that satisfy (i)--(iv)  below:
\begin{enumerate}
\item
$A^{\varepsilon} \in \mbox{Span}\lbrace I,A,A^*,AA^*,A^*A\rbrace$; 
\item
$A \in \mbox{Span}\lbrace I,A^*,A^{\varepsilon},A^*A^{\varepsilon},
A^{\varepsilon}A^*\rbrace$;
\item
$A^* \in \mbox{Span}\lbrace I,A^{\varepsilon},A, A^{\varepsilon}A,
AA^{\varepsilon}\rbrace$;
\item Any two of $A,A^*,A^{\varepsilon}$ satisfy the
Askey-Wilson relations.
\end{enumerate}
\end{problem}

\begin{example}
\label{ex:condracah}
\rm
Referring to
Problem \ref{prob:ae},
assume $A,A^*$ has Racah type
as in Example
\ref{ex:racah}.
Show that conditions (i)--(iv) of the problem
are satisfied for 
 $A^{\varepsilon} = h^{-1}A + h^{*-1}A^*$,
where $h, h^*$ are from  
 Example
\ref{ex:racah}.
\end{example}

\noindent We now consider
the uniqueness of the element $A^{\varepsilon}$ 
in Problem
\ref{prob:ae}.

\begin{problem}
\rm
\label{prob:ae2}
Let $A,A^*$ denote a Leonard pair on $V$ and let
the linear transformation $A^{\varepsilon}:V\to V$
be as in Problem
\ref{prob:ae}. Show that 
 $A^{\varepsilon \dagger}$
 satisfies
the conditions (i)--(iv) in Problem
\ref{prob:ae},
 where 
$\dagger $ is the antiautomorphism from
Theorem \ref{thm:dagger}.
Now let  
 $A^{\varepsilon \prime}:V\to V$ denote any linear
transformation that satisfies
the conditions (i)--(iv) in Problem
\ref{prob:ae}.
Show that 
 $A^{\varepsilon \prime}$
is contained in 
$\mbox{Span}\lbrace
 A^{\varepsilon}, I\rbrace $
or 
$\mbox{Span}\lbrace 
 A^{\varepsilon \dagger}, I\rbrace $.
\end{problem}

\begin{problem}
\label{prob:aepseig}
\rm
Referring to Problem
\ref{prob:ae}, find the eigenvalues of $A^{\varepsilon}$.
Find necessary and sufficient conditions for 
$A^{\varepsilon}$ to be diagonalizable.
Show that 
$A^{\varepsilon}$ is diagonalizable if and only
if $A,A^*,
A^{\varepsilon}$ is a Leonard triple.
\end{problem}

\begin{problem}
\label{prob:flag}
\rm
Let $\Phi$ denote the Leonard system
from Definition \ref{eq:ourstartingpt} and
 let $V$ denote an irreducible
$\mathcal A$-module. 
By a {\it flag} on $V$ we mean a nested sequence
 $f_0\subseteq f_1 \subseteq \cdots \subseteq f_d$ 
of subspaces of $V$ such that $f_i$ has dimension $i+1$
for $0 \leq i \leq d$. 
First, find all the flags
 $f_0\subseteq f_1 \subseteq \cdots \subseteq f_d$ 
on $V$ such that each of
$A f_i$, $A^*f_i$ is contained in  $f_{i+1}$ 
for $0 \leq i \leq d-1$.
Show that the set of such flags has the structure of
an affine algebraic variety.
Secondly, find all the
 flags
 $f_0\subseteq f_1 \subseteq \cdots \subseteq f_d$ 
on $V$ such that each of
$A f_i$, $A^*f_i$ is contained in  $f_{i+1}$ 
for $0 \leq i \leq d-1$
 and $A^{\varepsilon}f_i\subseteq f_i$
for $0 \leq i \leq d$, where 
$A^{\varepsilon}$ is from
Problem
\ref{prob:ae}.
Show there exist at most two such flags for a given 
$A^{\varepsilon}$.
\end{problem}

\noindent To motivate the next two problems,
recall that the Lie algebra $\mathfrak{sl}_2$
has a basis $x,y,z$ such that
$\lbrack x,y\rbrack=2x+2y$,
$\lbrack y,z\rbrack=2y+2z$,
$\lbrack z,x\rbrack=2z+2x$. This basis is called
{\it equitable}.

\begin{problem}
\label{prob:equit}
\rm
Let $A,A^*$ denote a Leonard pair on $V$. Assume
$\beta=2$ where $\beta+1$ is the common value of
(\ref{eq:betaplusone}). Show that there exists
an irreducible 
$\mathfrak{sl}_2$-module structure on $V$
such that $A$ (resp. $A^*$) acts on $V$
as a linear combination of
$I,x,y,xy$ (resp. $I, y,z,yz$). 
Here $x,y,z$ is the equitable basis for 
$\mathfrak{sl}_2$.
\end{problem}

\begin{problem}
\label{prob:equit3}
\rm
Referring to Problem
\ref{prob:equit}, assume 
that $A,A^*$ extends to a Leonard triple $A,A^*,A^\varepsilon $.
Show that 
$A^\varepsilon$ acts on $V$ as
a linear combination of $I,z,x,zx$.
\end{problem}

\noindent In order to  motivate the next two problems
we recall the equitable presentation for $U_q(\mathfrak{sl}_2)$.
This presentation has generators $X,X^{-1}, Y,Z$  that
satisfy $XX^{-1} = X^{-1}X=1$  and
\begin{eqnarray*}
\frac{qXY-q^{-1}YX}{q-q^{-1}}=1,
\qquad \quad
\frac{qYZ-q^{-1}ZY}{q-q^{-1}}=1,
\qquad \quad
\frac{qZX-q^{-1}XZ}{q-q^{-1}}=1.
\end{eqnarray*}

\begin{problem}
\label{prob:equitq}
\rm
Let $A,A^*$ denote a Leonard pair on $V$.
Assume $q$ is not a root of unity, where 
$q^2+q^{-2}+1$ is the common value of
(\ref{eq:betaplusone}).
Show that there exists
an irreducible 
$U_q(\mathfrak{sl}_2)$-module structure on $V$
such that $A$ (resp. $A^*$) acts on $V$
as a linear combination of
$I,X,Y,XY$ (resp. $I, Y,Z,YZ$).
Here $X,Y,Z$ are the equitable generators for 
$U_q(\mathfrak{sl}_2)$.
\end{problem}

\begin{problem}
\label{prob:equitq3}
\rm
Referring to Problem
\ref{prob:equitq},
assume 
that $A,A^*$ extends to a Leonard triple $A,A^*,A^\varepsilon $.
Show that 
$A^\varepsilon$ or 
$A^{\varepsilon \dagger}$  
 acts on $V$ as
a linear combination of $I,Z,X,ZX$. Here
$\dagger$ is the antiautomorphism from
Theorem \ref{thm:dagger}.
\end{problem}

\begin{example}
\label{prob:uq}
\rm Let $A,A^*$ denote a Leonard pair of dual $q$-Krawtchouk type;
see Example \ref{ex:sbutnotss}. After a minor change of
variables and replacing $q$ by $q^2$, this Leonard pair
has a parameter
array of the form
\begin{eqnarray*}
\theta_i &=& \eta+uq^{d-2i}+vq^{2i-d},
\\
\theta^*_i &=& \eta^*+u^*q^{d-2i}
\end{eqnarray*}
for $0 \leq i \leq d$ and
\begin{eqnarray*}
\varphi_i &=& uu^*(q^i-q^{-i})(q^{i-d-1}-q^{d-i+1})q^{d-2i+1},
\\
\phi_i &=& 
vu^*(q^i-q^{-i})(q^{i-d-1}-q^{d-i+1})q^{d-2i+1}
\end{eqnarray*}
for $1 \leq i \leq d$.
Show that the underlying vector space
supports an irreducible $U_q(\mathfrak{sl}_2)$-module structure 
such that
\begin{eqnarray*}
&&A= \eta I+u X +v Y, \qquad \qquad 
A^*=\eta^*I+u^* Z,
\end{eqnarray*}
where $X,Y,Z$ are the equitable generators for
$U_q(\mathfrak{sl}_2)$.
\end{example}

\begin{problem}
\rm
\label{prob:sym}
 Let 
$A, A^*$
denote a Leonard pair of $q$-Racah type; see Example
\ref{ex:qracah}. After applying an affine transformation
to $A,A^*$, and after a minor change of variables
with $q$ replaced by $q^2$, this Leonard pair
has a parameter array of the form
\begin{eqnarray*}
\theta_i &=& aq^{2i-d}+a^{-1}q^{d-2i},
\\
\theta^*_i &=& bq^{2i-d}+b^{-1}q^{d-2i}
\end{eqnarray*}
for $0 \leq i \leq d$ and
\begin{eqnarray*}
\varphi_i&=&a^{-1}b^{-1}q^{d+1}(q^{-i}-q^i)(q^{d-i+1}-q^{i-d-1})
(q^{-i}-abcq^{i-d-1})
(q^{-i}-abc^{-1}q^{i-d-1}),
\\
\phi_i&=&ab^{-1}q^{d+1}(q^{-i}-q^i)(q^{d-i+1}-q^{i-d-1})
(q^{-i}-a^{-1}bcq^{i-d-1})
(q^{-i}-a^{-1}bc^{-1}q^{i-d-1})
\end{eqnarray*}
for $1 \leq i \leq d$.
In this notation the parameters
$\beta,\gamma,\gamma^*,\varrho,\varrho^*$, $\omega,\eta,\eta^*$
from Theorem \ref{lptheorem} satisfy 
\begin{eqnarray*}
\beta=q^2+q^{-2}, \quad & &  \quad \gamma=\gamma^*=0,
\qquad \qquad 
\varrho=\varrho^*=-(q^2-q^{-2})^2,
\\
-\frac{\omega}{(q-q^{-1})^2} 
&=&(a+a^{-1})(b+b^{-1})+(c+c^{-1})(q^{d+1}+q^{-d-1}),
\\
\frac{\eta}{(q-q^{-1})(q^2-q^{-2})}
&=&(a+a^{-1})(c+c^{-1})+(b+b^{-1})(q^{d+1}+q^{-d-1}),
\\
\frac{\eta^*}{(q-q^{-1})(q^2-q^{-2})}
&=&(b+b^{-1})(c+c^{-1})+(a+a^{-1})(q^{d+1}+q^{-d-1}).
\end{eqnarray*}
Upon defining $A^\varepsilon$ via the first equation
below, the Askey-Wilson relations can be expressed as 
\begin{eqnarray*}
\frac{q^{-1}AA^*-qA^*A}{q^2-q^{-2}} &=&
A^{\varepsilon} -
\frac{
(a+a^{-1})(b+b^{-1})+(c+c^{-1})(q^{d+1}+q^{-d-1})}{q+q^{-1}}I,
\\
\frac{q^{-1}A^*A^{\varepsilon}-qA^{\varepsilon}A^*}{q^2-q^{-2}} &=&
A -
\frac{
(b+b^{-1})(c+c^{-1})+(a+a^{-1})(q^{d+1}+q^{-d-1})}{q+q^{-1}}I,
\\
\frac{q^{-1}A^{\varepsilon}A-qAA^{\varepsilon}}{q^2-q^{-2}} &=&
A^* -
\frac{
(c+c^{-1})(a+a^{-1})+(b+b^{-1})(q^{d+1}+q^{-d-1})}{q+q^{-1}}I.
\end{eqnarray*}
The element $A^{\varepsilon}$ is not diagonalizable in
general.
The roots of its characteristic polynomial
are
\begin{eqnarray}
\theta^{\varepsilon}_i = cq^{2i-d}+c^{-1}q^{d-2i}
\qquad \qquad (0 \leq i \leq d).
\label{eq:ceig}
\end{eqnarray}
The element $A^{\varepsilon}$ is diagonalizable if and only
if the eigenvalues (\ref{eq:ceig}) are distinct 
if and only 
if $c^2$ is not among
$q^{2d-2}, q^{2d-4}, \ldots, q^{2-2d}$
if and only if  $A,A^*, A^{\varepsilon}$
is a Leonard triple.
So far we have been assuming $A,A^*$ is $q$-Racah type;
find similar results for the other types of Leonard pairs
described in Appendix 35.
\end{problem}

\begin{problem}
\label{prob:equitex}
\rm
Referring to 
Problem \ref{prob:sym},
show that the underlying vector space supports
an irreducible
 $U_q(\mathfrak{sl}_2)$-module structure
such that
\begin{eqnarray*}
A &=& a X + a^{-1} Y  + c b^{-1} \frac{XY-YX}{q-q^{-1}},
\\
A^* &=& b Z + b^{-1} X  + a c^{-1} \frac{ZX-XZ}{q-q^{-1}},
\\
A^\varepsilon &=& c Y + c^{-1} Z  + b a^{-1} \frac{YZ-ZY}{q-q^{-1}},
\end{eqnarray*}
where $X,Y,Z$ are the equitable generators for
 $U_q(\mathfrak{sl}_2)$.
\end{problem}

\begin{problem}
\label{prob:double}
\rm
Let 
$(\lbrace \theta_i\rbrace_{i=0}^d; \lbrace \theta^*_i\rbrace_{i=0}^d;
\lbrace \varphi_i\rbrace_{i=1}^d; \lbrace \phi_i\rbrace_{i=1}^d)$
denote a parameter array as in 
Definition
\ref{def:pa}.
Find all the parameter arrays
$(\lbrace \theta'_i\rbrace_{i=0}^d; \lbrace \theta^{*\prime}_i\rbrace_{i=0}^d;
\lbrace \varphi'_i\rbrace_{i=1}^d; \lbrace \phi'_i\rbrace_{i=1}^d)$
such that
$\theta^{*\prime}_i=\theta^*_i$ for $0 \leq i \leq d$
and 
$\varphi'_i \phi'_i=\varphi_i \phi_i$
for $1 \leq i \leq d$.
\end{problem}

\begin{problem}
\label{prob:eigstr}
\rm
Let $A,A^*$ denote a tridiagonal pair over $\K$. Assume
$\beta=2$ where $\beta+1$ is the common value of
(\ref{eq:betaplusone}). 
We are interested in the eigenvalues of 
$A+tA^*$  $(t \in \K)$.
Let $m_t \in \K\lbrack \lambda \rbrack$ denote the (monic) minimal polynomial
of $A+tA^*$.
Call $t$ {\it feasible} whenever 
there exist $a,b,c \in \K$ such that
\begin{eqnarray*}
m_t = \prod_{i=0}^d (\lambda - a - bi-ci^2).
\end{eqnarray*}
Note that $t=0$ is feasible. What other $t$ are feasible?
\end{problem}

\begin{problem}
\rm
\label{prob:bidi}
Let $V$ denote a vector space over $\K$ with finite positive
dimension.
Consider a pair of linear transformations $A:V\to V$
and $A^*:V\to V$ that satisfy (i)--(iii) below:
    \begin{enumerate}
    \item Each of $A,A^*$ is diagonalizable.
    \item There exists an ordering $\lbrace V_i\rbrace_{i=0}^d$
     (resp. $\lbrace V^*_i\rbrace_{i=0}^d$) of the 
    eigenspaces of $A$ (resp. $A^*$) such that
    \begin{eqnarray*}
    (A-\theta_iI)V^*_i &\subseteq & V^*_{i+1} \qquad \qquad (0 \leq i \leq d),
\\
    (A^*-\theta^*_iI)V_i &\subseteq & V_{i+1} \qquad \qquad (0 \leq i \leq d),
    \end{eqnarray*}
   where $\theta_i$ (resp. $\theta^*_i$) is the eigenvalue of
$A$ (resp. $A^*$) associated with $V_i$ (resp. $V^*_i$) and 
     $V_{d+1} = 0$, $V^*_{d+1}= 0$.
    \item For $0 \leq i \leq d/2$ the restrictions
\begin{eqnarray*}
(A-\theta_iI)
(A-\theta_{i+1}I)
\cdots
(A-\theta_{d-i-1}I)\vert_{V^*_i} : V^*_i& \to & V^*_{d-i},
\\
(A^*-\theta^*_iI)
(A^*-\theta^*_{i+1}I)
\cdots
(A^*-\theta^*_{d-i-1}I)\vert_{V_i} : V_i &\to & V_{d-i}
\end{eqnarray*}
are bijections.
\end{enumerate}
Call such a pair $A,A^*$ a {\it bidiagonal pair}. Classify
the bidiagonal pairs up to isomorphism. Show that
these objects are essentially in bijection with
the finite-dimensional modules for
$U_q(\mathfrak{sl}_2)$.
\end{problem}

\begin{note}
\rm
%(December 5, 2006)
(Jan 1, 2008)

\noindent 
George Brown is working on
Problem \ref{prob:diagzero}
and Problem
\ref{prob:bipdbip}.  

\noindent
John Caughman and his students are 
working on
Problem \ref{prob:abs},
Problem \ref{prob:raise}.
%
%\noindent
%Brian Curtin is working on
%Problem \ref{prob:lt},
%Problem \ref{prob:qc}.
%

\noindent Brian Curtin and Hassan Alnajjar are working on
Problem \ref{prob:lt},
Problem  \ref{prob:drin},
Problem \ref{prob:tddual},
Problem \ref{prob:bilform}.

\noindent
Ali Godjali is working on
%Problem \ref{prob:krawtriple},
Problem \ref{prob:4vand},
Problem \ref{prob:ntsplit}.

\noindent Darren Neubauer is working on
Problem 
\ref{prob:bidi}.

\noindent Kazumasa Nomura is working
on
Problem \ref{prob:bi-tri},
Problem \ref{prob:kki},
Problem \ref{prob:tmat},
Problem \ref{prob:ae}.
He has solved 
Problem \ref{prob:zeroint}.

\noindent
Raimundas Vidunas is working on
%related to
%Problem \ref{prob:psipre},
%Problem \ref{prob:psi},
Problem \ref{prob:raise},
Problem \ref{prob:vs}.
\end{note}

%\section{Acknowledgment} 

\noindent Paul Terwilliger \hfil\break
Department of Mathematics \hfil\break
University of Wisconsin \hfil\break
480 Lincoln Drive \hfil\break
Madison, Wisconsin, 53706 USA 
\hfil\break
{Email: \tt terwilli@math.wisc.edu} \hfil\break

%\noindent Paul Terwilliger,
%Department of Mathematics,
%University of Wisconsin,\hfil\break
%480 Lincoln Drive,
%Madison Wisconsin, 53706 USA.
%{Email: \tt terwilli@math.wisc.edu} \hfil\break


\begin{thebibliography}{10}
    
%\bibitem{AAR}
%    G.E.~Andrews, R.~Askey, R.~Roy.
%     \newblock{\em Special functions}.
%   \newblock{Encyclopedia of Mathematics and its Applications,
%   71}.
% \newblock {Cambridge University Press,
%   Cambridge,
%     1999}.
%%%%%%%%%%%%
%
%\bibitem{AWil}
%R.~Askey and J.A.~Wilson.
%\newblock A set of orthogonal polynomials that 
%generalize the {R}acah coefficients or $6-j$ symbols. 
%\newblock {\em SIAM J. Math. Anal.}, 10:1008--1016, 1979. 
%

%\bibitem{BanIto}
%E.~Bannai and T.~Ito.
%\newblock {\em Algebraic Combinatorics {I}: Association Schemes}.
%\newblock Benjamin/Cummings, London, 1984.

%\bibitem{GR}
%G.~Gasper and M.~Rahman.
%\newblock {\em Basic hypergeometric series}.
%\newblock Encyclopedia of Mathematics and its Applications, 35.
%\newblock Cambridge University Press, Cambridge, 1990.

%\bibitem{Zhed}
%Ya.~Granovskii, I.~Lutzenko, and A. Zhedanov.   
%\newblock  Mutual integrability, quadratic algebras, and dynamical symmetry.
%\newblock  {\em Ann. Physics}, 217(1):1--20, 1992.

%\bibitem{grove}
%L.~C.~Grove.
%\newblock {\em Classical groups and geometric algebra}.
%\newblock American Mathematical Society, Providence RI, 2002.

%\bibitem{Grun}
%F.~A.~Grunbaum and L.~Haine.
%\newblock A $q$-version of a theorem of Bochner.
%\newblock {\em J. Comput. Appl. Math.}, 68(1-2):103--114, 1996.


%\bibitem{TD00}
%T.~Ito, K.~Tanabe, and P.~Terwilliger.
%\newblock Some algebra related to ${P}$- and ${Q}$-polynomial association
%  schemes.
%\newblock {\em Codes and Association Schemes (Piscataway NJ, 1999)},
%\newblock {Amer.
%%Math. Soc., Providence RI}, 56:167--192, 2001.

%\bibitem{Kassel}
%C.~Kassel,
%\newblock {\em Quantum groups},
%\newblock Springer-Verlag,
%New York,
%    1995.

%\bibitem{KoeSwa}
%R.~Koekoek and R.~Swarttouw.
%\newblock {\em The Askey-scheme of hypergeometric orthogonal polyomials and its
%  $q$-analog}, volume 98-17 of {\em Reports of the faculty of Technical
%  Mathematics and Informatics}.
%\newblock Delft, The Netherlands, 1998.

%\bibitem{Koelink3}
%H.~T. Koelink.
%\newblock Askey-{W}ilson polynomials and the quantum ${\rm {s}{u}}(2)$ group:
%  survey and applications.
%\newblock {\em Acta Appl. Math.}, 44(3):295--352, 1996.

%\bibitem{Leon}
%D.~Leonard.
%\newblock Orthogonal polynomials, duality, and association
%schemes.
%\newblock {\em SIAM J. Math. Anal.}, 13(4):656--663, 1982.

%\bibitem{Hjal}
%H.~Rosengren.
%\newblock {\em Multivariable orthogonal polynomials as coupling
%coefficients for Lie and quantum algebra representations.}
%\newblock Centre for Mathematical Sciences, Lund University, Sweden,
%1999.

%\bibitem{rotman}
%J.~J.~Rotman.
%\newblock{\em Advanced modern algebra.}
%\newblock{Prentice Hall, Saddle River NJ 2002}.

%\bibitem{TersubI}
%P.~Terwilliger.
%\newblock The subconstituent algebra of an association scheme. 
%\newblock {\em J. Algebraic Combin.}, 1(4):363--388, 1992.

%\bibitem{LS99}
%P.~Terwilliger.
%\newblock Two linear transformations each tridiagonal with respect to an
%  eigenbasis of the other.
%\newblock {\em Linear Algebra Appl.},  330:149--203, 2001. 

%\bibitem{qSerre}
%P.~Terwilliger.
%\newblock Two relations that generalize the q-Serre relations and the 
%Dolan-Grady relations.
%\newblock {\em  Physics and
%combinatorics 1999 (Nagoya)}, World Scientific Publishing,
% River Edge, NJ, 2001 

%\bibitem{LS24}
%P.~Terwilliger.
%\newblock  Leonard pairs from 24 points of view.
%\newblock {\em Rocky mountain Journal of Mathematics.}, 32(2):827--888, 2002. 

%\bibitem{conform}
%P.~Terwilliger.
%\newblock Two linear transformations each tridiagonal with respect to an
%  eigenbasis of the other: the $TD$-$D$ and the $LB$-$UB$ canonical form.
%Preprint.

%\bibitem{lsint}
%P.~Terwilliger.
%\newblock Introduction to Leonard pairs.
%\newblock {OPSFA Rome 2001}.
%\newblock{\em J. Comput. Appl. Math.}, 153(2):463--475, 2003.

%\bibitem{Terint}
%P.~Terwilliger.
%\newblock Introduction to {L}eonard pairs and 
%	     {L}eonard systems.
%  \newblock {\em S\=urikaisekikenky\=usho K\=oky\=uroku},
%  \newblock 
%   (1109):67--79, 1999.
%    \newblock Algebraic combinatorics  (Kyoto, 1999).

%\bibitem{TLT:split}
%P.~Terwilliger.
%\newblock Two linear transformations each tridiagonal with respect to an
%  eigenbasis of the other: comments on the split decomposition.
%Preprint.

%\bibitem{TLT:array}
%P.~Terwilliger.
%\newblock Two linear transformations each tridiagonal with respect to an
%  eigenbasis of the other: comments on the parameter array.
%\newblock {\em Geometric and Algebraic Combinatorics 2,
%Oisterwijk, The Netherlands 2002}.
%\newblock Submitted.

%\bibitem{aw}
%P.~Terwilliger and R.~Vidunas.
%\newblock Leonard pairs and the Askey-Wilson relations.
%\newblock {\em J. Algebra Appl.} Submitted.
%Available at
%\newblock{ 
%\tt  http://front.math.ucdavis.edu/math.QA/0305356
%}


%\bibitem{Zhidd}
%A.~S. Zhedanov.
%\newblock ``{H}idden symmetry'' of {A}skey-{W}ilson polynomials.
%\newblock {\em Teoret. Mat. Fiz.}, 89(2):190--204, 1991.

%%%%%extrabib%%%%%
%{\small
%\begin{thebibliography}{10}



\bibitem{CKOnsn}
C.~Ahn and K.~Shigemoto.
\newblock Onsager algebra and integrable lattice models,
\newblock {\em Modern Phys. Lett. A} {\bf 6(38)} (1991) 3509--3515.

\bibitem{Albert}
\newblock G.~Albertini, B.~McCoy, J.~Perk.
  \newblock 
  Eigenvalue spectrum of the superintegrable chiral {P}otts
				 model,
\newblock in {\em Integrable systems in quantum field theory and statistical
						mechanics},
		\newblock {Adv. Stud. Pure Math.}
{\bf 19}
1--55.
\newblock Academic Press,
   Boston, MA,
	 1989.

\bibitem{agmpy}
G.~Albertini, B.~McCoy, J.~H.~H.~Perk, S.~Tang.
\newblock
Excitation spectrum and order parameter for
the integrable $N$-state chiral Potts model. 
\newblock {\em
Nuclear Phys. B}  {\bf 314}  (1989) 741--763. 

\bibitem{AAR}
G.~Andrews, R.~Askey, and R.~Roy.
\newblock  {\em Special functions},
\newblock Cambridge University Press, Cambridge, 1999. 

\bibitem{AWil}
R.~Askey and J.A.~Wilson.
\newblock A set of orthogonal polynomials that 
generalize the {R}acah coefficients or $6-j$ symbols. 
\newblock {\em SIAM J. Math. Anal.}, 10:1008--1016, 1979. 

\bibitem{atak}
N.~Atakishiyev and A.~Klimyk.
\newblock On $q$-orthogonal polynomials, dual to
little and big $q$-Jacobi polynomials.
\newblock {\em J. Math. Anal. Appl.}
{\bf 294} (2004) 246--254;
{\tt arXiv:math.CA/0307250}.

\bibitem{atak2}
 M.~N.~Atakishiyev and V.~Groza.
\newblock 
The quantum algebra {$U\sb q(\rm su\sb 2)$} and
     {$q$}-{K}rawtchouk families of polynomials.
\newblock {\em
J. Phys. A}.
{\bf 37}
     (2004)
      2625--2635.
							       
							    
\bibitem{atak3} 
M.~N.~Atakishiyev, N.~M.~Atakishiyev, and A.~Klimyk.
\newblock
Big {$q$}-{L}aguerre and {$q$}-{M}eixner polynomials and
     representations of the quantum algebra {$U\sb q({\rm su}\sb
  {1,1})$}.
\newblock{\em
     J. Phys. A}.
    {\bf 36}
   (2003)
   10335--10347.

\bibitem{auyangandperk}
H.~Au-Yang and J.~H.~H.~Perk.
\newblock Onsager's star-triangle
equation: master key to integrability. 
\newblock {\em
Integrable systems in quantum field theory and 
statistical mechanics
57--94},
\newblock
Adv. Stud. Pure Math. {\bf 19},
Academic Press, Boston, MA, 1989.

\bibitem{auyang2}
H.~Au-Yang, B.~McCoy, J.~H.~H.~Perk,
and S.~Tang.
\newblock
Solvable models in statistical mechanics and {R}iemann
  surfaces of genus greater than one,
\newblock in {\em Algebraic analysis, Vol.\ I},
29--39.
\newblock Academic Press,
Boston, MA,
1988.

\bibitem{BanIto}
E.~Bannai and T.~Ito.
\newblock {\em Algebraic Combinatorics {I}: Association Schemes},
\newblock Benjamin/Cummings, London, 1984.

\bibitem{Baz}
V.~V.~Bazhanov and Y.~G.~Stroganov.
   \newblock{Chiral {P}otts model as a descendant of the six-vertex model},
\newblock{\em J. Statist. Phys.}
	      {\bf 59}
		    (1990)
			     799--817.
				  
\bibitem{bcn}
A.~E. Brouwer, A.~M. Cohen, and A.~Neumaier.
\newblock {\em Distance-Regular Graphs},
\newblock Springer-Verlag, Berlin, 1989.

\bibitem{Cau}
J.~S. Caughman~{I}{V}.
\newblock The {T}erwilliger algebras of bipartite ${P}$- and ${Q}$-polynomial
  schemes,
\newblock {\em Discrete Math.} {\bf 196} (1999) 65--95.

\bibitem{charp}
V.~Chari and A.~Pressley.
\newblock {Q}uantum affine algebras.
\newblock {\em Commun. Math. Phys.}
{\bf 142} (1991), 261--283. 

\bibitem{ciccoli}
N.~Ciccoli and F.~Gavarini.
\newblock A quantum duality principle for
subgroups and homogeneous spaces, preprint
{\tt arXiv:math.QA/0312289}.



\bibitem{curt2hom}
B.~Curtin.
\newblock The Terwilliger algebra of a 2-homogeneous bipartite
distance-regular graph,
\newblock {\em J. Combin. Theory Ser. B}
 {\bf 81}
(2001) 125--141.


\bibitem{CurNom}
B.~Curtin and K.~Nomura.
\newblock Distance-regular graphs related to the quantum enveloping algebra of
  $sl(2)$,
\newblock {\em J. Algebraic Combin.} 
{\bf 12} (2000) 25--36.

\bibitem{Curspin}
B.~Curtin.
\newblock Distance-regular graphs which support a spin model are thin, in:
Proc.  16th British Combinatorial Conference
(London, 1997),
\newblock {\em Discrete Math.} {\bf 197/198} (1999) 205--216.

\bibitem{mlt}
B.~Curtin.
\newblock Modular Leonard triples. Preprint.

\bibitem{Dask}
C.~Daskaloyannis.
\newblock Quadratic {P}oisson algebras of two-dimensional classical
  superintegrable systems and quadratic associative algebras of
 quantum superintegrable systems,
\newblock{\em J. Math. Phys.}
	  {\bf 42}
(2001)
100--1119.
							   

\bibitem{DateRoan2}
E.~Date and S.S. Roan.
\newblock The structure of quotients of the {O}nsager algebra by closed ideals,
\newblock {\em J. Phys. A: Math. Gen.} {\bf 33} (2000) 3275--3296.

\bibitem{Davfirst}
B.~Davies.
\newblock
Onsager's algebra and superintegrability,
\newblock{\em 
J. Phys. A:Math. Gen.} {\bf 23} (1990) 2245--2261.


\bibitem{Dav}
B.~Davies.
\newblock Onsager's algebra and the {D}olan-{G}rady condition in the
  non-self-dual case,
\newblock {\em J. Math. Phys.} {\bf 32} (1991) 2945--2950.

\bibitem{McCoyd}
\newblock   T.~Deguchi, K.~Fabricius,  and B.~McCoy.
\newblock {The {${\rm sl}\sb 2$} loop algebra symmetry of the six-vertex
		       model at roots of unity},
\newblock in {\rm Proceedings of the Baxter Revolution in Mathematical Physics
				      (Canberra, 2000)},
		\newblock {J. Statist. Phys.}
   {\bf 102}
    (2001)
     701--736.


\bibitem{Dolgra}
L.~Dolan and M.~Grady.
\newblock Conserved charges from self-duality,
\newblock {\em Phys. Rev. D (3)} {\bf 25} (1982) 1587--1604.

\bibitem{fairlie}
D.~B.~Fairlie.
\newblock Quantum deformations of {SU}(2).
\newblock {\em J. Phys. A: Math. Gen.} {\bf 23} (1990) L183--L187.

\bibitem{gan}
T. Gannon.
\newblock  Modular data: the algebraic combinatorics
of conformal field theory.
\newblock{\em J. Algebraic Combin.} {\bf 22} (2005) 211--250.


\bibitem{GR}
G.~Gasper and M.~Rahman.
\newblock {\em Basic hypergeometric series}.
\newblock Encyclopedia of Mathematics and its Applications, 35,
\newblock Cambridge University Press, Cambridge, 1990.



\bibitem{vongehlen}
G.~von Gehlen and V.~Rittenberg.
\newblock
$Z_n$-symmetric quantum chains with infinite
set of conserved charges and $Z_n$ zero modes,
\newblock{\em Nucl. Phys. B}
{\bf 257}
(1985) 351--370.


\bibitem{go}
J.~T.~Go.
\newblock The {T}erwilliger algebra of the hypercube,
\newblock {\em European J. Combin.} {\bf 23} (2002) 399--429.

\bibitem{GYZnature}
Ya.~A. Granovski{\u\i} and A.~S. Zhedanov.
\newblock Nature of the symmetry group of the $6j$-symbol,
\newblock {\em Zh. \`Eksper. Teoret. Fiz.}  {\bf 94} (1988) 49--54.

\bibitem{GYLZmut}
Ya.~I. Granovski{\u\i}, I.~M. Lutzenko, and A.~S. Zhedanov.
\newblock Mutual integrability, quadratic algebras, and dynamical symmetry,
\newblock {\em Ann. Physics} {\bf 217} (1992) 1--20.

\bibitem{GYZTwisted}
Ya.~I. Granovski{\u\i} and A.~S. Zhedanov.
\newblock ``{T}wisted'' {C}lebsch-{G}ordan coefficients for ${\rm {s}{u}}\sb
  q(2)$,
\newblock {\em J. Phys. A} {\bf 25} (1992) L1029--L1032.

\bibitem{GYZlinear}
Ya.~I. Granovski{\u\i} and A.~S. Zhedanov.
\newblock Linear covariance algebra for ${\rm {s}{l}}\sb q(2)$,
\newblock {\em J. Phys. A} {\bf 26} (1993) L357--L359.

\bibitem{GYZspherical}
Ya.~I. Granovski{\u\i} and A.~S. Zhedanov.
\newblock Spherical $q$-functions,
\newblock {\em J. Phys. A} {\bf 26} (1993) 4331--4338.

\bibitem{grove}
L.~C.~Grove.
\newblock {\em Classical groups and geometric algebra}.
\newblock American Mathematical Society, Providence RI, 2002.


\bibitem{GH4}
F.~A.~ Gr{\"u}nbaum.
\newblock Some bispectral musings, in: 
\newblock  {\em The bispectral problem (Montreal, PQ, 1997)},
  Amer. Math. Soc., Providence, RI, 1998, 
 pp. 31--45.

\bibitem{GH5}
F.~A.~Gr{\"u}nbaum and L.~Haine.
\newblock Bispectral {D}arboux transformations: an extension of the {K}rall
  polynomials,
\newblock {\em Internat. Math. Res. Notices} {\bf 8} (1997) 359--392.

\bibitem{GH7}
F.~A.~Gr{\"u}nbaum and L.~Haine.
\newblock The $q$-version of a theorem of {B}ochner,
\newblock {\em J. Comput. Appl. Math.}  {\bf 68} (1996) 103--114.

\bibitem{GH6}
F.~A.~Gr{\"u}nbaum and L.~Haine.
\newblock Some functions that generalize the {A}skey-{W}ilson polynomials,
\newblock {\em Comm. Math. Phys.} {\bf 184} (1997) 173--202.

\bibitem{GH1}
F.~A.~Gr{\"u}nbaum and L.~Haine.
\newblock On a $q$-analogue of the string equation and a generalization of the
  classical orthogonal polynomials, in: 
\newblock  {\em Algebraic methods and $q$-special functions (Montr\'eal, QC,
  1996)}
   Amer. Math. Soc., Providence, RI, 1999, pp. 
  171--181.

\bibitem{GH3}
F.~A.~Gr{\"u}nbaum and L.~Haine.
\newblock The {W}ilson bispectral involution: some elementary examples, in:
\newblock  {\em Symmetries and integrability of difference equations
  (Canterbury, 1996)},
   Cambridge Univ. Press, Cambridge, 1999, 
  pp. 353--369. 

\bibitem{GH2}
F.~A.~Gr{\"u}nbaum, L.~Haine, and E.~Horozov.
\newblock Some functions that generalize the {K}rall-{L}aguerre polynomials,
\newblock {\em J. Comput. Appl. Math.} {\bf 106} (1999) 271--297.


\bibitem{havlicek}
M.~Havlí\v cek, A.~U.~Klimyk, and S.~Po\v sta.
\newblock
Representations of the cyclically
symmetric $q$-deformed algebra ${\rm so}\sb q(3)$.
\newblock{\em
J. Math. Phys. } {\bf 40}  (1999) 2135--2161.

\bibitem{TD00}
T.~Ito, K.~Tanabe, and P.~Terwilliger.
\newblock Some algebra related to ${P}$- and ${Q}$-polynomial association
  schemes,  in:
\newblock {\em Codes and Association Schemes (Piscataway NJ, 1999)}, Amer.
Math. Soc., Providence RI, 2001, pp.
     167--192; 
{\tt arXiv:math.CO/0406556}.

\bibitem{shape}
T.~Ito and P.~Terwilliger.
\newblock The shape of a tridiagonal pair.
\newblock {\em J. Pure Appl. Algebra},
	      {\bf 188}
		    (2004)
			     145--160;
{\tt arXiv:math.QA/0304244
}.

\bibitem{tdanduq}
T.~Ito and P.~Terwilliger.
\newblock {Tridiagonal pairs and the quantum affine 
algebra
$U_q({\widehat{sl}}_2)$.}
\newblock {\em Ramanujan J.}, submitted; 
{\tt arXiv:math.QA/0310042}.

\bibitem{qtet}
T.~Ito and P.~Terwilliger.
\newblock
{The $q$-tetrahedron algebra and its finite dimensional irreducible modules.}
\newblock {\em Comm. Algebra}, in press;
{\tt arXiv:math.QA/0602199}.


\bibitem{uglov2}
I.~T.~Ivanov and D.~B.~Uglov.
\newblock
{$R$}-matrices for the semicyclic 
representations of {$U\sb q   
\widehat{\rm sl}(2)$},
\newblock {\em Phys. Lett. A},
{\bf 167}
    (1992)
  459--464.


\bibitem{Kassel}
C.~Kassel.
\newblock {\em Quantum groups},
\newblock Springer-Verlag,
New York,
    1995.

\bibitem{Klish1}
S.~Klishevich and M.~Plyushchay.
\newblock {Dolan-{G}rady relations and noncommutative quasi-exactly
		   solvable systems},
\newblock{\em J. Phys. A},
			    {\bf36},
			(2003)
					   11299--11319.
\bibitem{Klish2}
S.~Klishevich and M.~Plyushchay.
\newblock{Nonlinear holomorphic supersymmetry on {R}iemann surfaces},
\newblock{\em Nuclear Phys. B},
	      {\bf 640}
		    (2002)
			     481--503.

\bibitem{Klish3}
S.~Klishevich and M.~Plyushchay.
 \newblock {Nonlinear holomorphic supersymmetry, {D}olan-{G}rady relations
		       and {O}nsager algebra},
\newblock{\em Nuclear Phys. B},
				{\bf 628},
				    (2002)
					       217--233.
						  
%%%
\bibitem{KoeSwa}
R.~Koekoek and R.~F.~Swarttouw.
\newblock {\em The Askey scheme of hypergeometric orthogonal
polyomials and its
  $q$-analog}, report 98-17, Delft University of Technology, The
  Netherlands, 1998.
  Available at
 \newblock{ \tt http://aw.twi.tudelft.nl/{\~{}}koekoek/research.html} 

\bibitem{Koelink3}
H.~T. Koelink.
\newblock Askey-{W}ilson polynomials and the quantum ${\rm {s}{u}}(2)$ group:
  survey and applications,
\newblock {\em Acta Appl. Math.} {\bf 44} (1996) 295--352.

\bibitem{cite37} 
H.~T. Koelink.
\newblock $q$-Krawtchouk polynomials as spherical functions on the Hecke
algebra of type $B$,
\newblock {\em Trans. Amer. Math. Soc.} {\bf 352} (2000) 4789--4813.

\bibitem{cite38}
H.~T. Koelink and J.~Van der Jeugt.
\newblock  Convolutions for orthogonal polynomials from Lie and 
quantum algebra representations,
\newblock {\em SIAM J. Math. Anal.} {\bf 29} (1998) 794--822.


\bibitem{cite39}
H.~T. Koelink and J.~Van der Jeugt.
\newblock  Bilinear generating functions for orthogonal polynomials,
\newblock {\em Constr. Approx.} {\bf 15} (1999) 481--497.

\bibitem{cite40}
T.~H.~Koornwinder.
\newblock Askey-Wilson polynomials as zonal spherical functions on
the $SU(2)$ quantum group,
\newblock {\em SIAM J. Math. Anal.} {\bf 24} (1993) 795--813.

\bibitem{koornonn}
T.~H.~Koornwinder and U.~Onn.
\newblock Lower-upper triangular decompositions, $q=0$ limits,
and $p$-adic interpretations of some $q$-hypergeometric orthogonal
polynomials,
\newblock 
preprint
{\tt arXiv:math.CA/0405309}.

\bibitem{kresh}
A.~Kresch and H.~Tamvakis.
\newblock Standard conjectures for the arithmetic
grassmannian $G(2,N)$ and Racah polynomials.
\newblock {\em Duke Math. J.} {\bf 110} (2001), 359--376.

\bibitem{Lee}
C.~W.~H.~Lee and S.~G.~Rajeev.
\newblock{A {L}ie algebra for closed strings, spin chains, and gauge
		   theories},
\newblock{\em J. Math. Phys.},
			    {\bf 39}
				  (1998)
					   5199--5230.
\bibitem{Leon}
D.~Leonard.
\newblock Orthogonal polynomials, duality, and association
schemes,
\newblock {\em SIAM J. Math. Anal.} {\bf 13} (1982) 656--663.

\bibitem{Marco}
J.~Marco and J. Parcet.
\newblock A new approach to the theory of classical hypergeometric
polynomials.
\newblock {\em Trans. Amer. Math. Soc.} To appear.


\bibitem{McCoy1}
\newblock B.~McCoy.
\newblock  Integrable models in statistical mechanics: the hidden field
		   with unsolved problems.
\newblock {\em Internat. J. Modern Phys. A},
					  {\bf 14}
						(1999)
					 3921--3933.

\bibitem{N:aw}
K.~Nomura,
\newblock
Tridiagonal pairs and the {A}skey-{W}ilson relations,
Linear Algebra Appl.\ 397 (2005) 99--106.

\bibitem{N:refine}
K.~Nomura,
A refinement of the split decomposition of a tridiagonal pair,
Linear Algebra Appl.\ 403 (2005) 1--23.

\bibitem{N:height1}
K.~Nomura,
Tridiagonal pairs of height one,
Linear Algebra Appl.\ 403 (2005) 118--142.

\bibitem{NT:balanced}
K.~Nomura, P.~Terwilliger,
Balanced Leonard pairs,
Linear Algebra Appl., in press;
{\tt arXiv:math.RA/0506219}.

\bibitem{NT:formula}
K.~Nomura, P.~Terwilliger,
Some trace formulae involving the split sequences of a Leonard pair,
Linear Algebra Appl.\ 413 (2006) 189--201;
{\tt arXiv:math.RA/0508407}.

\bibitem{NT:det}
K.~Nomura, P.~Terwilliger,
The determinant of $AA^*-A^*A$ for a Leonard pair $A,A^*$,
Linear Algebra Appl.\ 416 (2006) 880--889;
{\tt arXiv:math.RA/0511641}.

\bibitem{NT:mu}
K.~Nomura, P.~Terwilliger,
Matrix units associated with the split basis of a Leonard pair,
Linear Algebra Appl.\ 418 (2006) 775--787;
{\tt arXiv:math.RA/0602416}.

\bibitem{NT:span}
K.~Nomura, P.~Terwilliger,
Linear transformations that are tridiagonal with respect to
both eigenbases of a Leonard pair,
Linear Algebra Appl., in press;
{\tt arXiv:math.RA/0605316}.

\bibitem{NT:switch}
K.~Nomura, P.~Terwilliger,
The switching element for a Leonard pair,
Linear Algebra Appl., submitted for publication;
{\tt arXiv:math.RA/0608623}.

\bibitem{NT:aff}
K.~Nomura, P.~Terwilliger.
Affine transformations of a Leonard pair,
Linear Algebra Appl., submitted for publication;
{\tt arXiv:math.RA/0611783}. 


\bibitem{noumi1}   
M.~Noumi and K.~Mimachi.
\newblock
Askey-{W}ilson polynomials as spherical functions on
{${\rm
SU}\sb q(2)$},
\newblock {\em Quantum groups (Leningrad, 1990)},
\newblock Lecture Notes in Math.
{\bf 1510}, 
 98--103.
\newblock Springer,
		Berlin,
		1992.

\bibitem{noumi2}
M.~Noumi and K.~Mimachi.
\newblock
Spherical functions on a family of quantum {$3$}-spheres,
\newblock  {\em Compositio Math.},
{\bf 83}
(1992)
19--42.

\bibitem{noumi3}
M.~Noumi.
\newblock Quantum groups and {$q$}-orthogonal polynomials---towards a
realization of {A}skey-{W}ilson polynomials on {${\rm SU}\sb
q(2)$},
\newblock {\em  Special functions (Okayama, 1990)} 260--288.
 ICM-90 Satell. Conf. Proc.
\newblock  Springer,
	Tokyo,
	    1991.

\bibitem{noumi4}
M.~Noumi and K.~Mimachi.
Askey-{W}ilson polynomials and the quantum group {${\rm SU}\sb
 q(2)$}.
\newblock{
\em Proc. Japan Acad. Ser. A Math. Sci.},
{\bf 66}
(1990)
146--149.

\bibitem{odesskii}
M.~Odesskii.
\newblock
An analog of the Sklyanin algebra.
\newblock {\em  Funct. Anal. Appl.}, {\bf 20} (1986) 78--79.

\bibitem{Onsager}
L.~Onsager.
\newblock 
Crystal statistics. I. 
A two-dimensional model with an 
order-disorder transition.  
\newblock{\em Phys. Rev. (2)}  {\bf 65} (1944)
117--149.

\bibitem{perk}
J.~H.~H. Perk.
\newblock
Star-triangle relations, quantum Lax pairs,
and higher genus curves.
\newblock{\em
Proceedings of Symposia in Pure Mathematics}
{\bf 49} 341--354.
\newblock Amer. Math. Soc., Providence, RI, 1989.

\bibitem{roan}
S.~S.~Roan.
\newblock
Onsager's algebra, loop algebra and chiral Potts model,
\newblock Preprint MPI 91--70, Max Plank Institute for
Mathematics, Bonn, 1991.



\bibitem{Rosengren}
H.~Rosengren.
\newblock {\em Multivariable orthogonal polynomials as coupling
coefficients for Lie and quantum algebra representations.}
\newblock Ph.D. Thesis.
\newblock Centre for Mathematical Sciences, Lund University, Sweden,
1999.

\bibitem{Rosengren2}
H.~Rosengren.
\newblock 
An elementary approach to the $6j$-symbols (classical,
quantum, rational, trigonometric, and elliptic),
preprint
{\tt arXiv:math.CA/0312310}.

\bibitem{rotman}
J.~J.~Rotman.
\newblock{\em Advanced modern algebra.}
\newblock{Prentice Hall, Saddle River NJ 2002}.



\bibitem{uniform}
P.~Terwilliger.
\newblock
The incidence algebra of a uniform poset,
\newblock {\em Math and its applications} {\bf 20} (1990) 193--212.


\bibitem{TersubI}
P.~Terwilliger.
\newblock The subconstituent algebra of an association scheme I. 
\newblock {\em J. Algebraic Combin.} {\bf 1} (1992) 363--388.
   
\bibitem{tersub3}
P.~Terwilliger.
\newblock The subconstituent algebra of an association scheme III.
\newblock{\em
J. Algebraic Combin. }
{\bf 2}  (1993) 177--210.


\bibitem{LS99}
P.~Terwilliger.
\newblock Two linear transformations each tridiagonal with respect to an
  eigenbasis of the other.
  \newblock {\em Linear Algebra Appl.}  {\bf 330} (2001), 149--203;
{\tt arXiv:math.RA/0406555}.

  \bibitem{qSerre}
  P.~Terwilliger.
  \newblock Two relations that generalize the $q$-Serre relations and the
  Dolan-Grady relations. In
  \newblock {\em  Physics and
  Combinatorics 1999 (Nagoya)}, 377--398, World Scientific Publishing,
   River Edge, NJ, 2001; 
{\tt arXiv:math.QA/0307016}.

   \bibitem{LS24}
   P.~Terwilliger.
   \newblock  Leonard pairs from 24 points of view.
   \newblock {\em Rocky Mountain J. Math.} {\bf 32}(2) (2002), 827--888;
{\tt arXiv:math.RA/0406577}.

   \bibitem{conform}
   P.~Terwilliger.
   \newblock Two linear transformations each tridiagonal with respect to an
     eigenbasis of the other; the $TD$-$D$ and the $LB$-$UB$ canonical form.
\newblock {\em J. Algebra}, submitted;
{\tt arXiv:math.RA/0304077}.


    \bibitem{lsint}
    P.~Terwilliger.
    \newblock Introduction to Leonard pairs.
    \newblock {OPSFA Rome 2001}.
    \newblock{\em J. Comput. Appl. Math.} {\bf 153}(2) (2003),
    463--475.

\bibitem{Terint}
P.~Terwilliger.
\newblock Introduction to {L}eonard pairs and
  {L}eonard systems.
  \newblock {\em S\=urikaisekikenky\=usho K\=oky\=uroku},
 \newblock
 (1109):67--79, 1999.
\newblock Algebraic combinatorics  (Kyoto, 1999).

\bibitem{TLT:split}
P.~Terwilliger.
\newblock Two linear transformations each tridiagonal with respect to an
  eigenbasis of the other; comments on the split decomposition,
\newblock {\em OPSFA7} To appear;
{\tt arXiv:math.RA/0306290}.



 \bibitem{TLT:array}
 P.~Terwilliger.
 \newblock Two linear transformations each tridiagonal with respect to an
   eigenbasis of the other; comments on the parameter array.
 \newblock {\em Geometric and Algebraic Combinatorics 2,
 Oisterwijk, The Netherlands 2002}.
  \newblock  To appear, 
{\tt arXiv:math.RA/0306291}.

\bibitem{qrac}
P.~Terwilliger.
\newblock Leonard pairs and the $q$-Racah polynomials.
\newblock {\em Linear Algebra Appl.} To appear;
{\tt arXiv:math.QA/0306301}.

\bibitem{aw}
P.~Terwilliger and R.~Vidunas.
\newblock Leonard pairs and the Askey-Wilson relations.
\newblock {\em J. Algebra Appl.} To appear;
{\tt arXiv:math.QA/0305356}.

\bibitem{uglov1}   
D.~B.~Uglov and I.~T.~Ivanov.
\newblock {${\rm sl}(N)$} {O}nsager's algebra and integrability,
\newblock {\em J. Statist. Phys.},
 {\bf 82}
(1996)
87--113.

%\bibitem{uglov2}
%I.~T.~Ivanov and D.~B.~Uglov.
%\newblock
%{$R$}-matrices for the semicyclic 
%representations of {$U\sb q   
%\widehat{\rm sl}(2)$},
%\newblock  {\em Phys. Lett. A},
%{\bf 167}
%    (1992)
%  459--464.


\bibitem{Zhidd}
A.~S. Zhedanov.
\newblock ``{H}idden symmetry'' of {A}skey-{W}ilson polynomials,
\newblock {\em Teoret. Mat. Fiz.} {\bf 89} (1991) 190--204.

\bibitem{ZheCart}
A.~S. Zhedanov.
\newblock Quantum ${\rm {s}{u}}\sb q(2)$ algebra: ``{C}artesian'' version and
  overlaps,
\newblock {\em Modern Phys. Lett. A} {\bf 7} (1992) 1589--1593.

\bibitem{Zhidden}
A.~S. Zhedanov.
\newblock Hidden symmetry algebra and overlap coefficients for two ring-shaped
  potentials,
\newblock {\em J. Phys. A} {\bf 26} (1993) 4633--4641.

\bibitem{LPcm}
\newblock A.~S. Zhedanov, and A.~Korovnichenko.
	\newblock {``{L}eonard pairs'' in classical mechanics},
	   \newblock {\em J. Phys. A},
		  {\bf 5},
			(2002)
				 5767--5780.



\end{thebibliography}
\end{document}